\documentstyle[12pt]{article}

\large\normalsize
\newcommand{\eqref}[1]{(\ref{#1})}
\newtheorem{proposition}{Proposition}[section]
\newtheorem{definition}[proposition]{D\'efinition}
\newtheorem{remark and definition}[proposition]{Remarque et d\'efinition}
\newtheorem{theorem}[proposition]{Th\'eor\`eme}
\newtheorem{important example}[proposition]{Exemple important}
\newtheorem{principal mutations}[proposition]{Mutations principales}
\newtheorem{corollary}[proposition]{Corollaire}
\newtheorem{remark}[proposition]{Remarque}
\newtheorem{particular case}[proposition]{Cas particulier}
\newtheorem{power associative algebra}[proposition]{Alg\`ebre \`a puissances
associatives}
\newtheorem{flexible algebra}[proposition]{Alg\`ebre flexible}
\newtheorem{Jordan algebra}[proposition]{Alg\`ebre de Jordan}
\newtheorem{quadratic algebra}[proposition]{Alg\`ebre quadratique}
\newtheorem{alternative algebra}[proposition]{Alg\`ebre alternative}
\newtheorem{noncommutative Jordan algebra}[proposition]{Alg\`ebre de Jordan non commutative}
\newtheorem{Lie derivations algebra}[proposition]{Alg\`ebre de Lie des d\'erivations}
\newtheorem{note}[proposition]{Note}
\newtheorem{notes}[proposition]{Notes}
\newtheorem{Sckolem-Noether theorem}[proposition]{Th\'eor\`eme de Sckolem-Noether}
\newtheorem{Lie group and its Lie algebra}[proposition]{Groupe de Lie et son alg\`ebre de Lie}
\newtheorem{linear invertibility}[proposition]{Inversibilit\'e lin\'eaire}
\newtheorem{J-invertibility}[proposition]{J-Inversibilit\'e}
\newtheorem{Banach algebras}[proposition]{Alg\`ebres de Banach}
\newtheorem{Gelfand-Mazur theorem}[proposition]{Th\'eor\`eme de Gelfand-Mazur}
\newtheorem{DTLZ}[proposition]{Diviseurs topologiques lin\'eaires de z\'ero}
\newtheorem{GMK theorem}[proposition]{Th\'eor\`eme de Gelfand-Mazur-Kaplansky}
\newtheorem{cayley algebra}[proposition]{Alg\`ebre cayleyenne}
\newtheorem{CD process}[proposition]{Proc\'ed\'e de Cayley-Dickson}
\newtheorem{J-DTZ}[proposition]{J-diviseurs topologiques de z\'ero}
\newtheorem{isotopy}[proposition]{Isotopie}
\newtheorem{GCD process}[proposition]{Proc\'ed\'e de Cayley-Dickson g\'en\'eralis\'e}

\newtheorem{notation and remark}[proposition]{Notation et remarque}
\newtheorem{remarks}[proposition]{Remarques}
\newtheorem{note and remark}[proposition]{Note et remarque}
\newtheorem{lemma}[proposition]{Lemme}
\newtheorem{example}[proposition]{Exemple}
\newtheorem{examples}[proposition]{Exemples}
\newtheorem{definitions}[proposition]{D\'efinitions}
\newtheorem{definitions and notations}[proposition]{D\'efinitions et notations}
\newtheorem{definition and remarks}[proposition]{D\'efinition et remarques}
\newtheorem{principal isotopes}[proposition]{Isotopes Principales}

\newtheorem{remarks and open problem}[proposition]{Remarques et probl\`eme ouvert}

\def\cit{\hbox{\it l\hskip -5.5pt C\/}}
\def\oit{\hbox{\it l\hskip -5.5pt O\/}}
\def\rit{{ \mathbb{R}}}
\def\rit{\hbox{\it l\hskip -2pt R}}

\def\nit{\hbox{\it l\hskip -2pt N}}

\def\hit{\hbox{\it l\hskip -2pt H}}

\begin{document}

\sloppy
\renewcommand{\thesection}{\arabic{section}}
\renewcommand{\theequation}{\thesection.\arabic{equation}}

\setcounter{page}{1}
\begin{center} {\large\bf UNIVERSITE MOHAMMED V \\ FACULTE DES SCIENCES DE RABAT} \end{center}

\begin{center} {\large\bf Th\`ese} \end{center}

\begin{center} {\large\bf Pr\'esent\'ee pour obtenir la grade de \\
DOCTEUR ES SCIENCES MATHEMATIQUES}
\end{center}

\begin{center} {\large\bf PAR} \end{center}

\begin{center} {\large\bf ROCHDI ABDELLATIF} \end{center}

\begin{center} {\Large\bf Alg\`ebres non
associatives norm\'ees de division.}
\end{center}

\begin{center} {\Large\bf Classification des alg\`ebres r\'eelles de Jordan non commutatives de division lin\'eaire
de dimension $8.$}
\end{center}

\vspace{0.5cm} \hspace{2cm} {\large Soutenue le 05 Octobre 1994
devant le jury:}

\vspace{1cm} \hspace{0.1cm} {\small $M^r$ \ KERKOUR Ahmed
\hspace{7.8cm} Pr\'esident

\vspace{0.1cm} \hspace{0.9cm} {\scriptsize Professeur, Pr\'esident
de l'Universit\'e AL-AKHAWAYN}

\vspace{0.4cm} \hspace{0.1cm} {\small $M^r$ \ BENSLIMANE Mohamed
\hspace{6.5cm} Examinateur

\vspace{0.1cm} \hspace{0.9cm} {\scriptsize Professeur \`a la
Facult\'e des Sciences de T\'etouan}

\vspace{0.4cm} \hspace{0.1cm} {\small $M^r$ \ BOURASS Abdelhamid
\hspace{6.9cm} Examinateur

\vspace{0.1cm} \hspace{0.9cm} {\scriptsize Professeur \`a la
Facult\'e des Sciences de Rabat}

\vspace{0.4cm} \hspace{0.1cm} {\small $M^r$ \ CHIDAMI Mohammed
\hspace{7.1cm} Examinateur

\vspace{0.1cm} \hspace{0.9cm} {\scriptsize Professeur \`a la
Facult\'e des Sciences de Rabat}

\vspace{0.4cm} \hspace{0.1cm} {\small $M^r$ \ CUENCA MIRA Jos\'e
Antonio \hspace{5.9cm} Rapporteur

\vspace{0.1cm} \hspace{0.9cm} {\scriptsize Professeur \`a la
Facult\'e des Sciences de M\'alaga}

\vspace{0.4cm} \hspace{0.1cm} {\small $M^r$ \ ESSANNOUNI Hassane
\hspace{7cm} Rapporteur

\vspace{0.1cm} \hspace{0.9cm} {\scriptsize Professeur \`a la
Facult\'e des Sciences de Rabat}

\vspace{0.4cm} \hspace{0.1cm} {\small $M^r$ \ KAIDE El-Amin
\hspace{7.2cm} Directeur de recherche

\vspace{0.1cm} \hspace{0.9cm} {\scriptsize Professeur \`a la
Facult\'e des Sciences de Rabat}

\vspace{3cm} .

\vspace{6cm}
\begin{center} {\large\bf A la m\'emoire de mon p\`ere.

\vspace{0.1cm}
A ma m\`ere.

\vspace{0.1cm}
A ma femme et mes deux filles Fadwa et Salma.

\vspace{0.1cm} A mes fr\`eres et soeurs, beaux fr\`eres et belles
soeurs.

\vspace{0.1cm} A mes cousines.

\vspace{0.1cm} A ma belle famille.

\vspace{0.1cm} A tous les miens.}
\end{center}

\vspace{10cm} {\Large\bf Remerciements}

\vspace{1cm} \hspace{0.5cm} Je voudrais remercier le Professeur El
Amin KAIDI pour avoir dirig\'e mes recherches. Sa comp\'etence,
grande exp\'erience et disponibilit\'e \`a mon \'egard, sans
compter ses qualit\'es humaines, m'ont permis de mener \`a bien ce
travail. Qu'il trouve ici l'expression de ma profonde gratitude et
mon grand respect.

\vspace{0.6cm} \hspace{0.5cm} Ce travail a \'et\'e r\'ealis\'e, en
majorit\'e, \`a la Facult\'e des Sciences de M\'alaga, dans le
cadre de la coop\'eration qui existe entre le groupe
inter-universitaire (marocain) d'Alg\`ebre et Th\'eorie des
Nombres, et le groupe d'Alg\`ebre de la Facult\'e des Sciences de
M\'alaga. A ce propos, je remercie respectueusement le Professeur
Jos\'e Antonio CUENCA MIRA, mon co-directeur de recherche, pour
l'hospitalit\'e qu'il m'a offerte ainsi que les moyens de travail
qu'il a mis \`a ma disposition. Sa collaboration scientifique m'a
permis de r\'ealiser ce m\'emoire.

\vspace{0.6cm} \hspace{0.5cm} Je remercie \'egalement le
Professeur Ricardo DE LOS SANTOS VILLODRES pour sa collaboration
scientifique qui a \'et\'e tr\`es fructueuse, et pour les moments
amicaux que nous avons eu lors de ses visites \`a la Facult\'e des
Sciences de M\'alaga.

\vspace{0.6cm} \hspace{0.5cm} Je tiens \`a exprimer ma
reconnaissance au Professeur KERKOUR Ahmed, Pr\'esident de
l'Universit\'e Al-AKHAWAYN, d'avoir bien voulu pr\'esider le jury
de cette th\`ese.

\vspace{0.6cm} \hspace{0.5cm} Je remercie vivement les Professeurs
BENSLIMANE Mohamed, BOURASS Abdelhamid, CHIDAMI Mohammed et
ESSANNOUNI Hassane pour l'honneur qu'ils me font en acceptant de
participer au jury.

\vspace{0.6cm} \hspace{0.5cm} Je voudrais exprimer ma gratitude au
Professeur Antonio FERNANDEZ LOPEZ, Directeur du D\'epartement
d'Alg\`ebre, G\'eom\'etrie et Topologie de la Facult\'e des
Sciences de M\'alaga, pour sa collaboration, son soutien, ses
encouragements et les facilit\'es qu'il m'a offertes, et
\'egalement les Professeurs Mercedes SILES MOLINA, Eulalia GARCIA
RUS, Esperanza SANCHEZ CAMPOS, Alberto CASTELLON SERRANO,
C\'andido MARTIN GONZALEZ, Antonio SANCHEZ SANCHEZ et Aniceto
MUTILLO M\'AS, pour leur collaboration. L'attention et la
d\'elicatesse qu'ils m'ont prouv\'ees, m'ont \'enormement
touch\'e.

\vspace{0.6cm} \hspace{0.5cm} Je remercie chaleureusement le
Professeur Angel RODR\'IGUEZ PALACIOS pour l'int\'er\^et
particulier qu'il a toujours port\'e \`a mon travail, et
\'egalement les Professeurs Juan MARTINEZ MORENO et Miguel CABRERA
GARC\'IA.

\vspace{0.6cm} \hspace{0.5cm} Mes remerciements vont \'egalement
\`a Mademoiselle Elisa JAIME JIMENEZ, secr\'etaire du
d\'epartement cit\'e, pour son efficacit\'e et sa gentillesse.

\vspace{0.6cm} \hspace{0.5cm} J'adresse aussi un remerciement
sp\'ecial \`a notre cher biblioth\'ecaire Abdelkrim KOUBIDA.

\vspace{0.6cm} \hspace{0.5cm} Je ne saurais remercier assez ma
grande et, soutout, petite famille, qui a souffert de mes absences
r\'ep\'et\'ees, et a support\'e des moments difficiles sans me le
faire sentir. Elle a sans doute contribu\'e indirectement \`a ce
que je viens d'accomplir.

\vspace{18cm}
\begin{center}
{\Large\bf TABLE DES MATI\`ERES}
\end{center}

\vspace{0.6cm} \hspace{0.1cm} {\large\bf Introduction.}

\vspace{1.1cm} {\large\hspace{0.1cm}\bf Chapitre 1.
G\'en\'eralit\'es et pr\'erequis.}

\vspace{0.5cm} \hspace{0.8cm} {\bf\S 1. Alg\`ebres non
associatives} \hspace{6.8cm} ..... \hspace{0.2cm} 13

\vspace{0.3cm} \hspace{0.8cm} {\bf\S 2. Inversibilit\'e dans les
alg\`ebres non associatives} \hspace{2.8cm} ..... \hspace{0.2cm}
20

\vspace{0.3cm} \hspace{0.8cm} {\bf\S 3. Alg\`ebres quadratiques}
\hspace{7.3cm} ..... \hspace{0.2cm} 24

\vspace{1.1cm} {\large\hspace{0.1cm}\bf Chapitre 2. Alg\`ebres non
associatives norm\'ees de division.}

\vspace{0.5cm} \hspace{0.8cm} {\bf\S 1. Alg\`ebres norm\'ees non
associatives} \hspace{5cm} ..... \hspace{0.2cm} 30

\vspace{0.3cm} \hspace{0.8cm} {\bf\S 2. Diviseurs topologiques de
z\'ero dans une alg\`ebre norm\'ee} \hspace{0.7cm} .....
\hspace{0.2cm} 34

\vspace{0.3cm} \hspace{0.8cm} {\bf\S 3. Alg\`ebres complexes
norm\'ees de division lin\'eaire} \hspace{2.3cm} .....
\hspace{0.2cm} 39

\vspace{0.3cm} \hspace{0.8cm} {\bf\S 4. Alg\`ebres r\'eelles
norm\'ees de division lin\'eaire} \hspace{3cm} .....
\hspace{0.2cm} 44

\vspace{1.1cm} {\large\hspace{0.1cm}\bf Chapitre 3. Proc\'ed\'e de
Cayley-Dickson G\'en\'eralis\'e.}

\vspace{0.5cm} \hspace{0.8cm} {\bf\S 1. Les sous-alg\`ebres de
dimension $4$ dans une $\rit$-alg\`ebre de Jordan} \\
. \hspace{1.2cm} {\bf non commutative de division lin\'eaire de
dimension $8$} \hspace{1.4cm} ..... \hspace{0.2cm} 47

\vspace{0.3cm} \hspace{0.8cm} {\bf\S 2. Proc\'ed\'e de
Cayley-Dickson G\'en\'eralis\'e} \hspace{4.1cm} .....
\hspace{0.2cm} 53

\vspace{1.1cm} {\large\hspace{0.1cm}\bf Chapitre 4.
Classification des $\rit$-alg\`ebres de Jordan n.c. \\
.\hspace{2.7cm} de division lin\'eaire de dimension $8.$}

\vspace{0.5cm} \hspace{0.8cm} {\bf\S 1. Isotopie vectorielle}
\hspace{8cm} ..... \hspace{0.2cm} 61

\vspace{0.3cm} \hspace{0.8cm} {\bf\S 2. Probl\`emes
d'isomorphisme} \hspace{6.4cm} ..... \hspace{0.2cm} 64

\vspace{0.3cm} \hspace{0.8cm} {\bf\S 3. Th\'eor\`eme de
classification} \hspace{6.5cm} ..... \hspace{0.2cm} 68

\vspace{2cm} {\large\hspace{0.1cm}\bf Chapitre 5.
$\rit$-alg\`ebres de Jordan n.c. de division lin\'eaire \\ .
\hspace{1.7cm} de dimension $8$ ayant un automorphisme non
trivial.}

\vspace{0.5cm} \hspace{0.8cm} {\bf\S 1. Etude des
$\rit$-alg\`ebres de Jordan n.c. de division lin\'eaire de} \\
. \hspace{1.2cm} {\bf dimension $8$ qui poss\`edent une
d\'erivation non triviale} \hspace{1.3cm} ..... \hspace{0.2cm} 84

\vspace{0.3cm} \hspace{0.8cm} {\bf\S 2. Caract\'erisation des
$\rit$-alg\`ebres de Jordan n.c. de division}
\\ . \hspace{1.2cm} {\bf lin\'eaire de dimension $8$ ayant un
automorphisme non trivial} \hspace{0.2cm} ..... \hspace{0.2cm} 99

\vspace{0.8cm} \hspace{0.1cm} {\large\bf Index} \hspace{11.1cm}
......... \hspace{0.2cm} 111

\vspace{0.7cm} \hspace{0.1cm} {\large\bf R\'ef\'erences}
\hspace{9.8cm} ......... \hspace{0.2cm} 113

\vspace{18cm} \normalsize
{\Large\bf Introduction}

\vspace{0.5cm} \hspace{0.3cm} Dans une alg\`ebre $A,$ non
n\'ecessairement associative et non n\'ecessairement unitaire sur
un corps commutatif $K,$ un \'el\'ement $x\neq 0$ est
lin\'eairement inversible si les op\'erateurs lin\'eaires de
multiplication

\[  L_x:y\mapsto xy \hspace{0.2cm} \mbox{ et } \hspace{0.2cm} R_x:y\mapsto yx \]

sont inversibles dans l'alg\`ebre \ $End_K(A),$ \ des op\'erateurs
lin\'eaires de $A.$ L'alg\`ebre $A$ est dite de division
lin\'eaire si tout \'el\'ement non nul de $A$ est lin\'eairement
inversible. Elle est dite de division lin\'eaire \`a gauche si
pour tout \'el\'ement non nul $x$ de $A,$ $L_x$ est inversible
dans $End_K(A).$ Si $K=\rit$ ou $\cit,$ l'alg\`ebre $A$ est dite
norm\'ee si l'espace vectoriel $A$ est muni d'une norme \ $||.||$
\ sous-multiplicative i.e. \ $||xy||\leq ||x||$ $||y||$ \ pour
tous $x,y\in A.$

\vspace{0.3cm} \hspace{0.3cm} L'\'etude des alg\`ebres norm\'ees
de division lin\'eaire a connu son premier succ\'es en 1941
gr\^ace \`a Mazur et Gelfand qui prouv\`erent que le corps \
$\cit$ des nombres complexes est l'unique, \`a isomorphisme
pr\`es, $\cit$-alg\`ebre associative norm\'ee de division [BD 73],
[Ri 60]. Ce r\'esultat a \'et\'e \'etendu aux alg\`ebres non
n\'ecessairement associatives par Kaidi ([Kai 77] Teorema {\bf
1.6}) qui prouva en 1977 que les \ $\cit$-alg\`ebres norm\'ees
compl\`etes de division lin\'eaire \`a gauche sont \`a
isomorphisme pr\`es $\cit.$ Cependant, le probl\`eme de la
d\'etermination des \ $\cit$-alg\`ebres norm\'ees (non
n\'ecessairement compl\`etes) de division lin\'eaire, est encore
ouvert. Nous avons apport\'e une contribution modeste, \`a ce
probl\`eme, en assurant la validit\'e de ce dernier r\'esultat
dans des situations apparemment plus g\'en\'erales que celles du
cas complet. Nous avons montr\'e que les \ $\cit$-alg\`ebres
norm\'ees sans diviseurs topologiques lin\'eaires de z\'ero \`a
gauche (d.t.l.z.g.) et qui contiennent des \'el\'ements
lin\'eairement inversibles \`a gauche, sont isomorphes \`a \
$\cit$ (Th\'eor\`eme {\bf 2.51}). Nous avons montr\'e \'egalement
que les \ $\cit$-alg\`ebres norm\'ees de division lin\'eaire \`a
gauche dans lesquelles l'ensemble des d.t.l.z.g. est une partie
compl\`ete, sont isomorphes \`a \ $\cit$ (Th\'eor\`eme {\bf
2.54}).

\vspace{0.3cm} \hspace{0.3cm} Le probl\`eme de la d\'etermination
des $\rit$-alg\`ebres norm\'ees (compl\`etes) de division
lin\'eaire est r\'esolu au cas des alg\`ebres alternatives et de
Jordan [Kai 77]. Frobenius prouva en 1877 que les \
$\rit$-alg\`ebres associatives alg\'ebriques de division sont
isomorphes \`a $\rit,$ $\cit$ ou $\hit$ (le corps r\'eel des
quaternions de Hamilton) [E-R 91]. Kaplansky [Kap 49] acheva ce
premier travail et donna \'egalement une extension, pour les
alg\`ebres norm\'ees sans diviseurs topologiques de z\'ero
(d.t.z.), en montrant en 1949 que les \ $\rit$-alg\`ebres
associatives norm\'ees sans d.t.z. non nuls sont isomorphes \`a \
$\rit,$ $\cit$ ou $\hit.$ Albert [A 49] prouva en 1949 que les \
$\rit$-alg\`ebres alternatives alg\'ebriques de division sont
isomorphes \`a \ $\rit,$ $\cit,$ $\hit$ ou \ $\oit$ (l'alg\`ebre
r\'eelle des octonions de Cayley-Dickson). Il s'est int\'eress\'e
\'egalement \`a l'\'etude des \ $\rit$-alg\`ebres absolument
valu\'ees [A 47, 49]. Cette derni\`ere \'etude fut achev\'ee par
Urbanik et Wright qui prouv\`erent, en 1960, que les \
$\rit$-alg\`ebres absolument valu\'ees unitaires sont isomorphes
\`a \ $\rit,$ $\cit,$ $\hit$ ou \ $\oit$ [UW 60]. R\'ecemment
Cabrera et Rodriguez [CR] ont donn\'e une nouvelle
d\'emonstration, simple, du Th\'eor\`eme de Kaplansky, et \`a
l'aide de celui d'Albert, ils prouv\`erent que les \
$\rit$-alg\`ebres alternatives norm\'ees sans d.t.z. non nuls sont
isomorphes \`a \ $\rit,$ $\cit,$ $\hit$ ou $\oit.$

\vspace{0.3cm} \hspace{0.3cm} Wright avait conjectur\'e [Wr 53]
que les \ $\rit$-alg\`ebres norm\'ees de division lin\'eaire sont
de dimension finie. Cette conjecture s'est aver\'ee extr\^emement
difficile dans sa g\'en\'eralit\'e et on est actuellement loin
d'une r\'eponse affirmative, cependant on a des r\'esultats
partiels. R\'ecemment, Cuenca [Cu 92] a donn\'e des exemples d'une
classe d'alg\`ebres r\'eelles norm\'ees compl\`etes de dimension
infinie de division lin\'eaire \`a gauche, et Rodriguez [Rod
92$_1$] les a compl\`etement d\'ecrites. Ant\'erieurement, on a
confirm\'e la validit\'e de la conjecture de Wright pour les
alg\`ebres alternatives et de Jordan, et dans le cas Jordan non
commutatif, on a montr\'e que l'alg\`ebre est quadratique
(cayleyenne) [Kai 77]. L'existence des alg\`ebres de Jordan non
commutatives quadratiques de division lin\'eaire de dimension
infinie est actuellement un probl\`eme ouvert. D'apr\`es le
Th\'eor\`eme de Hopf [H 40], Kervaire [Ke 58] et Milnor-Bott [BM
58], qui affirme que $1,2,4,8$ sont les seules possibilit\'es pour
la dimension d'une \ $\rit$-alg\`ebre de division lin\'eaire de
dimension finie, on est amen\'e d'une mani\`ere naturelle \`a
\'etudier ces alg\`ebres en dimension finie. Osborn [Os 62] a
initi\'e la th\'eorie des alg\`ebres quadratiques et a
d\'etermin\'e toutes les \ $\rit$-alg\`ebres quadratiques de
division lin\'eaire de dimension $4$ et une classe particuli\`ere
de \ $\rit$-alg\`ebres quadratiques (non alternatives) de division
lin\'eaire de dimension $8.$ Kaidi [Kai 77] prouva que les \
$\rit$-alg\`ebres norm\'ees de Jordan non commutatives de division
lin\'eaire qui satisfont \`a l'identit\'e \ $(x,x,[x,y])=0$ \ sont
quadratiques, alternatives et isomorphes \`a \ $\rit,$ $\cit,$
$\hit$ ou $\oit.$ En particulier, les \ $\rit$-alg\`ebres de
Jordan norm\'ees de division lin\'eaire sont isomorphes \`a \
$\rit$ ou $\cit.$ Les \ $\rit$-alg\`ebres norm\'ees de Jordan non
commutatives de division lin\'eaire dans lesquelles deux
\'el\'ements qui n'appartiennent pas \`a la m\^eme sous-alg\`ebre
de dimension $2$ engendrent une sous-alg\`ebre de dimension $4,$
sont de dimension finie et isomorphes \`a \ $\rit,$ $\cit,$
$\hit^{(\lambda)}$ ou $\oit^{(\lambda)},$ $\lambda\neq\frac{1}{2}$
(les mutations de \ $\rit,$ $\cit,$ $\hit$ et $\oit$) [Kai 77],
[Roc 87]. Benkart, Britten et Osborn [BBO 82] r\'eduis\`erent, en
1982, la d\'etermination des \ $\rit$-alg\`ebres flexibles de
division lin\'eaire de dimension finie \`a celle des alg\`ebres de
Jordan non commutatives

\vspace{1cm} \hspace{0.3cm} Motiv\'es par les r\'esultats
pr\'ec\'edents, nous nous sommes int\'eress\'e au probl\`eme de la
d\'etermination des alg\`ebres r\'eelles de Jordan non
commutatives de division lin\'eaire de dimension $8.$ Nous avons
suivi deux approches. La premi\`ere consiste \`a construire ces
derni\`eres alg\`ebres, \`a partir d'une alg\`ebre de dimension
$4,$ par "duplication", et la seconde, \`a faire une d\'eformation
appropri\'ee du produit de l'alg\`ebre des octonions de
Cayley-Dickson.

\vspace{0.3cm} \hspace{0.3cm} Dans la premi\`ere approche, un
r\'esultat important a \'et\'e \'etabli dans [CDKR] assurant
l'existence d'une sous-alg\`ebre de dimension $4$ dans une
$\rit$-alg\`ebre de Jordan non commutatives de division lin\'eaire
de dimension $8$ (Th\'eor\`eme {\bf 3.7}). Pour le doublage, nous
avons g\'en\'eralis\'e le proc\'ed\'e de Cayley-Dickson [KR 92].
Soient $(B,\overline{.})$ une $K$-alg\`ebre cayleyenne (car(K)=0)
et $\gamma, \alpha, \delta\in K$ avec $\gamma\neq 0,$ alors le
produit

\[ (x,y)(x',y')=(x.^{\alpha}x'+\gamma\overline{y'}y,
y\overline{x'}+y'x+\frac{\delta}{2}[y',y]) \]

\vspace{0.1cm} munit l'espace vectoriel \ $B\times B$ \ d'une
structure de $K$-alg\`ebre cayleyenne, qu'on appelle extension
cayleyenne g\'en\'eralis\'ee de $(B,\overline{.})$ d'indice
$(\gamma,\alpha,\delta)$ et qu'on note \
$E_{\gamma,\alpha,\delta}(B).$ Nous avons \'etudi\'e ces
alg\`ebres et donn\'e une condition n\'ecessaire et suffisante
pour qu'elles soient de Jordan non commutatives (Proposition {\bf
3.17}), puis des conditions n\'ecessaires et suffisantes, lorsque
$K=\rit$ et $B$ de dimension $4,$ pour qu'elles soient de division
lin\'eaire (Corollaire {\bf 4.14}). Ceci nous a permi l'obtention
d'une nouvelle famille d'alg\`ebres r\'eelles de Jordan non
commutatives de division lin\'eaire de dimension $8:$

\[ \Big( E_{-1,\alpha,\delta}(\hit^{(\lambda)})\Big)^{(\mu)} \hspace{0.2cm} \mbox{
o\`u} \hspace{0.2cm} \lambda, \mu\neq\frac{1}{2}, \
\alpha>\frac{1}{2} \hspace{0.2cm} \mbox{ et } \hspace{0.2cm}
(2\alpha-1)\delta^2<4. \]

\vspace{0.1cm} Malheureusement (Note {\bf 5.26}), ce premier
proc\'ed\'e s'est av\'er\'e insuffisant pour la d\'etermination de
toutes les \ $\rit$-alg\`ebres de Jordan non commutatives de
division lin\'eaire de dimension $8$ [Roc$_1$], [Roc$_2$]
cependant, il nous a permis, comme on le verra ult\'erieurement,
l'\'etude des alg\`ebres qui poss\`edent une d\'erivation non
triviale et de r\'epondre affirmativement \`a une question pos\'ee
par Benkart et Osborn en 1981 dans [BO 81$_1$].

\vspace{5cm} \hspace{0.3cm} Pour la deuxi\`eme approche, on a
donn\'e une m\'ethode de construction de $\rit$-alg\`ebres de
Jordan non commutatives de division lin\'eaire de dimension $8,$
\`a partir de la donn\'ee de l'une d'elles [CDKR]. Soient
$(V,\wedge)$ une $\rit$-alg\`ebre anti-commutative de dimension
finie $\geq 1,$ $(.|.)$ une forme bilin\'eaire sym\'etrique
d\'efinie n\'egative sur $V$ et $\varphi$ un automorphisme de
l'espace vectoriel $V.$ On pose: \ $x\Delta y=\varphi^*\Big(
\varphi(x)\wedge\varphi(y)\Big),$ \hspace{0.2cm} o\`u $x,y\in V$ \
($\varphi^*$ \'etant l'automorphisme adjoint de $\varphi$) et on
d\'esigne par \ $(V,(.|.),\wedge)$ et $(V,(.|.),\Delta),$
respectivement, les alg\`ebres cayleyennes construites \`a partir
des alg\`ebres anti-commutatives $(V,\wedge)$ et $(V,\Delta),$ et
de la forme bilin\'eaire sym\'etrique $(.|.).$ On a \'etabli le
r\'esultat important suivant (Proposition {\bf 4.1}):

\vspace{0.3cm} {\bf Proposition} [CDKR]. {\em $(V,(.|.),\wedge)$
est flexible de division lin\'eaire si et seulement si
$(V,(.|.),\Delta)$ est flexible de division lin\'eaire.}

\vspace{0.4cm} \hspace{0.3cm} L'alg\`ebre $(V,(.|.),\Delta)$ est
dite obtenue \`a partir de $A=(V,(.|.),\wedge)$ et $\varphi,$ par
isotopie vectorielle, et est not\'ee $A(\varphi)$ ({\bf 4.2 1)}).
A l'aide de ce dernier proc\'ed\'e et des r\'esultats
pr\'ec\'edents, nous avons pu d\'eterminer toutes les \
$\rit$-alg\`ebres de Jordan non commutatives de division
lin\'eaire de dimension $8,$ puis nous avons r\'esolu le
probl\`eme d'isomorphisme (Corollaire {\bf 4.10}). Nous avons
\'etabli le Th\'eor\`eme de classification suivant (Th\'eor\`eme
{\bf 4.17}):

\vspace{0.3cm} {\bf Th\'eor\`eme} [CDKR]. {\em Les alg\`ebres
r\'eelles de Jordan non commutatives de division lin\'eaire de
dimension $8$ s'obtiennent \`a partir de l'alg\`ebre r\'eelle \
$\oit=\Big( W,(.|.),\wedge\Big)$ de Cayley-Dickson par isotopie
vectorielle et sont \`a isomorphisme pr\`es \ $\oit(s)$ o\`u $s$
est un automorphisme sym\'etrique de l'espace euclidien $\Big(
W,-(.|.)\Big),$ d\'efini positif. De plus \ $\oit(s)\simeq$
$\oit(s')$ \ ($s,s'$ \'etant deux automorphismes sym\'etriques de
l'espace euclidien $\Big( W,-(.|.)\Big),$ d\'efinis positifs) si
et seulement si il existe $f\in G_2$ (le groupe des automorphismes
de l'alg\`ebre $\oit$) tel que \ $\tilde{s'}=f\circ\tilde{s}\circ
f^{-1}$ \ ($\tilde{s}$ \'etant le prolongement naturel de $s$ \`a
$\oit:$ \ $\alpha+u\mapsto \alpha+s(u)$).}

\vspace{0.3cm} \hspace{0.3cm} Benkart et Osborn [BO 81$_2$] ont
donn\'e toutes les possibilit\'es pour l'alg\`ebre de Lie $Der(A)$
des d\'erivations d'une $\rit$-alg\`ebre, $A,$ de division
lin\'eaire de dimension $8:$ \begin{enumerate} \item $G_2$
compacte, \item $su(3),$ \item $su(2)\oplus su(2),$ \item
$su(2)\oplus N$ o\`u $N$ est une alg\`ebre ab\'elienne de
dimension $\leq 1,$ \item $N$ une alg\`ebre ab\'elienne de
dimension $\leq 2.$ \end{enumerate}

\vspace{0.1cm} Ils ont obtenu ensuite dans [BO 81$_1$] une
classification compl\`ete pour les alg\`ebres r\'eelles de
division lin\'eaire de dimension $8$ dont l'alg\`ebre de Lie des
d\'erivations est $G_2$ compacte, $su(3)$ ou $su(2)\oplus su(2).$
Ils ont donn\'e \'egalement des exemples d'alg\`ebres r\'eelles de
division lin\'eaire pour chacun des autres cas de l'alg\`ebre de
Lie des d\'erivations puis ils ont pos\'e, parmi d'autres, le
probl\`eme de l'existence d'une alg\`ebre r\'eelle de division
lin\'eaire dont l'alg\`ebre de Lie des d\'erivations est $su(2)$
et dont la d\'ecomposition en $su(2)$-modules irr\'eductibles est
de la forme: $1+1+3+3.$

\vspace{0.2cm} \hspace{0.3cm} Ces travaux nous ont permi
d'appliquer avec succ\`es les r\'esultats obtenus sur la
duplication et d'apporter des \'eclaircissements sur l'\'etude des
\ $\rit$-alg\`ebres de Jordan non commutatives de division
lin\'eaire de dimension $8$ ayant une d\'erivation non triviale.
Nous avons montr\'e que ces alg\`ebres ne peuvent pas avoir
$su(3)$ comme alg\`ebre de Lie des d\'erivations (Remarque {\bf
5.10}), et prouv\'e, pour les alg\`ebres de Lie $G_2$ compacte et
$su(2)\oplus su(2),$ les deux r\'esultats suivants (Th\'eor\`emes
{\bf 5.12, 5.13}):

\vspace{0.3cm} {\bf Th\'eor\`eme} [Roc$_1$]. {\em Soit $A$ une
$\rit$-alg\`ebre de Jordan non commutative de division lin\'eaire
de dimension $8.$ Alors $Der(A)=G_2$ compacte si et seulement si
$A\simeq$ $\oit^{(\lambda)},$ $\lambda\neq\frac{1}{2}.$}

\vspace{0.3cm} {\bf Th\'eor\`eme} [Roc$_1$]. {\em Soit $A$ une \
$\rit$-alg\`ebre de Jordan non commutative de division lin\'eaire
de dimension $8.$ Alors $Der(A)=su(2)\oplus su(2)$ compacte si et
seulement si $A\simeq$ $\Big(
E_{-1,\alpha,0}(\hit)\Big)^{(\lambda)}$ avec $1\neq\alpha
>\frac{1}{2}$ et $\lambda\neq\frac{1}{2}.$}

\vspace{0.4cm} \hspace{0.3cm} Nous avons donn\'e ensuite (Remarque
{\bf 5.14}) une r\'eponse affirmative \`a la question
pr\'ec\'edente \`a l'aide d'une alg\`ebre de Jordan non
commutative, obtenue par le proc\'ed\'e de Cayley-Dickson
g\'en\'eralis\'e: \ $\Big( E_{-1,\alpha,\delta}(\hit)\Big)$ o\`u
$\alpha >\frac{1}{2},$ $\delta\neq 0$ et $(2\alpha-1)\delta^2<4.$

\vspace{0.3cm} \hspace{0.3cm} Nous nous sommes int\'eress\'es
ensuite au groupe des automorphismes d'une alg\`ebre r\'eelle de
Jordan non commutatives de division lin\'eaire de dimension $8.$
Nous avons donn\'e un exemple d'une telle alg\`ebre dont le groupe
des automorphismes est trivial (Note {\bf 5.26}), puis nous avons
caract\'eris\'e les alg\`ebres ayant un automorphisme non trivial
(Th\'eor\`eme {\bf 5.25}):

\vspace{0.5cm} {\bf Th\'eor\`eme} [Roc$_1$]. {\em Soient
$\oit=\Big( W,(.|.),\wedge\Big)$ l'alg\`ebre r\'eelle de
Cayley-Dickson et $s$ un automorphisme sym\'etrique de l'espace
euclidien $\Big( W,-(.|.)\Big),$ d\'efini positif. Alors les deux
propri\'et\'es suivantes sont \'equivalentes: \begin{enumerate}
\item $Aut\Big( \oit(s)\Big) $ n'est pas trivial. \item
$\tilde{s}$ laisse stable une sous-alg\`ebre de $\oit$ de
dimension $4.$ \end{enumerate}}

\vspace{0.4cm} \hspace{0.3cm} Ceci met en \'evidence l'immensit\'e
de la classe des $\rit$-alg\`ebres de Jordan non commutatives de
division lin\'eaire de dimension $8$ [Roc$_1$], [Roc$_2$].

\vspace{0.5cm} \hspace{0.3cm} Nous avons divis\'e ce m\'emoire en
cinq chapitres. Le premier est consacr\'e aux g\'en\'eralit\'es
sur les alg\`ebres non associatives. Il contient \'egalement
quelques r\'esultats faisant partie de notre travail (Lemmes {\bf
1.3, 1.29}, Corollaire {\bf 1.30}). Dans le deuxi\`eme chapitre,
nous exposons d'abord quelques r\'esultats utiles de la th\'eorie
de base des alg\`ebres norm\'ees compl\`etes non associatives,
puis \'etudions les \ $\cit$-alg\`ebres norm\'ees de division
lin\'eaire (\`a gauche). Nous exposons enfin, les r\'esultats
connus ou parus r\'ecemment sur les \ $\rit$-alg\`ebres norm\'ees
(compl\`etes) de division lin\'eaire. Dans le troisi\`eme chapitre
nous d\'emontrons, en premier lieu, le r\'esultat concernant
l'existence d'une sous-alg\`ebre de dimension $4$ dans une \
$\rit$-alg\`ebre de Jordan non commutative de division lin\'eaire
de dimension $8,$ puis nous \'etudions des propri\'et\'es
particuli\`eres satisfaites par ces sous-alg\`ebres (Proposition
{\bf 3.14}). Nous introduisons ensuite le proc\'ed\'e de
Cayley-Dickson "G\'en\'eralis\'e" suivi par l'exposition de ses
propri\'et\'es fondamentales. Dans le quatri\`eme chapitre nous
introduisons le proc\'ed\'e "d'isotopie vectorielle", puis nous
traitons le probl\`eme d'isomorphisme. Nous d\'emontrons ensuite
le Th\'eor\`eme de classification des \ $\rit$-alg\`ebres de
Jordan non commutatives de division lin\'eaire de dimension $8.$
Enfin, dans le dernier chapitre, nous \'etudions les \
$\rit$-alg\`ebres de Jordan non commutatives de division
lin\'eaire de dimension $8,$ dont l'alg\`ebre de Lie des
d\'erivations n'est pas triviale, et caract\'erisons le cas, plus
g\'en\'eral, o\`u le groupe des automorphismes n'est pas trivial.

\vspace{1cm} \hspace{10cm} M\'alaga, fin Mars 1994.

\vspace{5cm}
\section{G\'en\'eralit\'es et pr\'erequis}

\vspace{0.5cm} \hspace{0.3cm} Dans ce chapitre $K$ d\'esignera un
corps commutatif de caract\'eristique nulle.

\vspace{0.7cm} \subsection{Alg\`ebres non associatives}

\vspace{0.2cm}
\begin{definitions and notations} . \begin{enumerate} \item On appelle $K$-alg\`ebre tout $K$-espace vectoriel $A$
muni d'une application $K$-bilin\'eaire \ $A\times A
\longrightarrow A\quad (x,y)\mapsto xy$ \ appel\`ee produit de
$A.$ L'alg\`ebre $A$ est dite associative (resp. commutative,
anti-commutative) si son produit est associatif (resp. commutatif,
anti-commutatif). Les d\'efinitions d'\'el\'ement unit\'e (\`a
gauche), sous-alg\`ebre, id\'eal, isomorphisme, automorphisme,
sont les m\^emes que dans le cas associatif. L'alg\`ebre $A$ est
dite simple si elle ne poss\`ede pas d'id\'eaux bilat\`eres non
nuls propres. On note \ $\rit,$ $\cit,$ $\hit$ et $\oit,$
respectivement, les \ $\rit$-alg\`ebres des nombres r\'eels, des
nombres complexes, des quaternions de Hamilton et des octonions de
Cayley-Dickson. \item Soit $A$ une $K$-alg\`ebre et soient
$x,y,z\in A,$ on note \ $[x,y]$ le commutateur $xy-yx,$ de $x$ et
$y$ \ et \ $(x,y,z)$ l'associateur $(xy)z-x(yz),$ de $x,y$ et $z.$
Les sous-ensembles de $A$ suivants: \[ N(A)=\{x\in A:
(x,A,A)=(A,x,A)=(A,A,x)=0\} \ \mbox{ et } \] \[ Z(A)=\{x\in
N(A):[X,A]=0
\] sont des sous-alg\`ebres associatives de $A,$ de plus $Z(A)$
est commutative. On les appelle, respectivement, le noyau et le
centre de $A.$ Si $A$ est unitaire, son \'el\'ement unit\'e est
not\'e $1.$ Elle est dite centrale si \ $Z(A)=K.1.$ Si $x\in A,$
on note \ $L_x,$ $R_x$ les op\'erateurs lin\'eaires de
multiplication par $x:$ \[ \mbox{ \`a gauche } \ L_x:y\mapsto xy \
\mbox{ et \`a droite } \ R_x:y\mapsto yx,
\] qu'on appelle op\'erateurs de multiplication par $x,$ \`a
gauche et \`a droite. On note \ $End_K(A)$ \ la $K$-alg\`ebre
associative et unitaire des op\'erateurs lin\'eaires de $A.$ Si
$\lambda\in K,$ on appelle mutation $\lambda$ de $A,$ et on note
$A^{(\lambda)},$ l'alg\`ebre ayant pour espace vectoriel
sous-jacent $A$ et pour produit \ $x.^{\lambda}y=\lambda
xy+(1-\lambda)yx.$ Si $\lambda, \mu\in K,$ on a \ $\Big(
A^{(\lambda)}\Big)^{(\mu)}=A^{(\alpha)}$ o\`u
$\alpha=2\lambda\mu-\lambda-\mu+1.$ L'alg\`ebre \
$A^{(\frac{1}{2})}$ est commutative, not\'ee simplement $A^+,$ on
l'appelle sym\'etrisation de $A.\Box$ \end{enumerate}
\end{definitions and notations}

\vspace{0.1cm} \begin{definition} Soit $A$ une $K$-alg\`ebre, on
dit qu'une partie $S,$ non vide de $A,$ engendre lin\'eairement
$A,$ si $S$ est une partie g\'en\'eratrice de l'espace vectoriel
$A.$ On note $[S]_A$ la sous-alg\`ebre de $A$ engendr\'ee par $S.$
La sous-alg\`ebre engendr\'ee par un \'el\'ement $a\in A$ est
not\'ee $[a]_A.$ Si $A$ est unitaire et si $x_1,\dots,x_n$ sont
des \'el\'ements de $A-\{1\},$ alors la sous-alg\`ebre de $A$
engendr\'ee par la partie \ $\{1,x_1,\dots,x_n\}$ est not\'ee
$K_A[x_1,\dots,x_n].$ La dimension de l'alg\`ebre $A$ est la
dimension de l'espace vectoriel $A.$ Si ${\cal B}=\{u_i:i\in I\}$
est une base de l'alg\`ebre $A,$ c'est \`a dire base de l'espace
vectoriel $A,$ alors pour tous $i,j\in I$ on a

\vspace{0.1cm}
\begin{eqnarray} u_iu_j=\sum_{k\in I}\lambda_{ijk}u_k \end{eqnarray}

\vspace{0.1cm} o\`u les $\lambda_{ijk}$ sont des \'el\'ements de
$K,$ nuls, sauf pour un nombre fini d'indices $k\in I.$ Les
relations {\bf (1.1)} s'appelent la table de multiplication de $A$
relativement \`a la base ${\cal B}.$ R\'eciproquement, si ${\cal
B}=\{u_i:i\in I\}$ est une base d'un $K$-espace vectoriel $A,$
\'etant donn\'ee une famille \ $\{\lambda_{ijk}\in K: i,j,k\in
K\}$ telle que, pour tous $i,j\in K$ fix\'es, les $\lambda_{ijk}$
o\`u $k\in I,$ sont nuls, sauf pour un nombre fini, il existe
alors sur $A$ une seule structure de $K$-alg\`ebre pour laquelle
les relations {\bf (1.1)} sont satisfaites {\em ([Bou 70] {\bf A
III} p. 10)}.$\Box$
\end{definition}

\vspace{0.2cm}
\begin{lemma} {\em ([Roc$_1$] p. 2).} Soit $A$ une $K$-alg\`ebre et soit $S$ une partie non vide de $A.$
Alors pour tout $\lambda\in K-\{\frac{1}{2}\},$ la sous-alg\`ebre
de $A^{(\lambda)}$ engendr\'ee par $S$ coincide avec la mutation
$\lambda$ de la sous-alg\`ebre de $A$ engendr\'ee par $S:$ \
$[S]_{A^{(\lambda)}}=([S]_A)^{(\lambda)}.\Box$
\end{lemma}

\vspace{0.3cm}
\begin{power associative algebra} Les puissances \`a gauche d'un \'el\'ement $a$ d'une
$K$-alg\`ebre $A$ sont d\'efinis par \ $a^1=a$ \ et \
$a^{n+1}=aa^n$ pour tout $n\in\nit^*.$ L'alg\`ebre $A$ est dite
\`a puissances associatives si la sous-alg\`ebre $[a]_A$ est
associative pour tout $a\in A,$ ce qui est \'equivalent \`a dire
que $a^na^m=a^{n+1}$ pour tous $n,m\in\nit^*$ et pour tout $a\in
A.$ L'alg\`ebre r\'eelle \ $\stackrel{*}{\cit}$ de Mc Clay {\em
[A]}, ayant pour espace vectoriel sous-jacent $\cit$ et pour
produit \ $z\odot z'=\overline{z}$ $\overline{z'},$ $\overline{z}$
\'etant le conjugu\'e de $z,$ est commutative et n'est pas \`a
puissances associatives.$\Box$
\end{power associative algebra}

\vspace{0.2cm}
\begin{note} Dans {\em ([A 48] p. 554)} Albert d\'ecouvrit un ensemble minimal d'identit\'es
assurant l'associativit\'e des puissances et prouva qu'une
$K$-alg\`ebre $A$ est \`a puissances associatives si et seulement
si elle satisfait aux deux identit\'es \ $(x,x,x)=(x,x,x^2)=0$
pour tout $x\in A.$ Notons que ces deux identit\'es entra\^inent \
$(x^2,x,x)=0$ pour tout $x\in A.\Box$
\end{note}

\vspace{0.2cm} \hspace{0.3cm} Soit $A$ une $K$-alg\`ebre et soient
$x,y\in A.$ Pour motif de simplification, nous aurons l'occasion
d'utiliser, les notations suivantes:
\begin{enumerate} \item $xy+yx:=x\bullet y.$
\item $L_x+R_x:=V_x.$ \item $L_xL_y+L_yL_x:=L_{x,y}.$ \item
$R_xR_y+R_yR_x:=R_{x,y}.$
\end{enumerate}

\vspace{0.2cm}
\begin{proposition} Soit $A$ une $K$-alg\`ebre qui satisfait \`a l'identit\'e $(x,x,x)=0$ pour tout $x\in A.$ Alors
les trois propri\'et\'es suivantes sont \'equivalentes:
\begin{enumerate} \item $A$ est \`a puissances associatives.
\item $A$ satisfait \`a l'identit\'e \[ (x,y,x\bullet
z)+(x,z,x\bullet y)+(y,x,x\bullet z)+(x,x,y\bullet
z)+(z,x,x\bullet y)+(y,z,x^2)+(z,y,x^2)=0
\] pour tous $x,y,z\in A.$ \item $A$ satisfait \`a l'identit\'e \[ (x\bullet
z,y,x)+(x\bullet y,z,x)+(x\bullet z,x,y)+(y\bullet
z,x,x)+(x\bullet y,x,z)+(x^2,z,y)+(x^2,y,z)=0
\] pour tous $x,y,z\in A.$
\end{enumerate}
\end{proposition}

\vspace{0.1cm} {\bf Preuve.} 1) $\Rightarrow$ 2) est \'etablie
dans ([Sc 66] p. 129). Pour 2) $\Rightarrow$ 1) on fait $y=z=x$ et
on utilise la Note {\bf 1.5}. Enfin 1) $\Leftrightarrow$ 3)
s'obtient de la m\^eme fa\c{c}on.$\Box$

\vspace{0.2cm}
\begin{corollary} Soit $A$ une $K$-alg\`ebre \`a puissances associatives, alors pour tout entier $m\geq 2$ et pour tout
$x\in A,$ \ $L_{x^m}$ et $R_{x^m}$ appartiennent \`a la
sous-alg\`ebre de $End_K(A)$ engendr\'ee par $L_x,$ $R_x,$
$L_{x^2}$ et $R_{x^2}$ et on a: \begin{enumerate} \item
$3R_{x^{m+1}}=(2V_{x^m}-L_{x,x^{m-1}})V_x+(2R_{x^2}-R_x^2)V_{x^{m-1}}-2L_xR_{x^m}-L_{x^{m-1}}R_{x^2}.$
\item
$3L_{x^{m+1}}=(2V_{x^m}-R_{x,x^{m-1}})V_x+(2L_{x^2}-L_x^2)V_{x^{m-1}}-2R_xL_{x^m}-R_{x^{m-1}}L_{x^2}.$
\end{enumerate}
\end{corollary}

\vspace{0.1cm} {\bf Preuve.} Les deux identit\'es 1) et 2)
s'obtiennent des identit\'es 2) et 3) de la Proposition {\bf 1.6}
pr\'ec\'edente en faisant $z=x^{m-1}.$ La premi\`ere proposition
s'obtient alors des identit\'es 1) et 2) par r\'ecurrence sur
$m.\Box$

\vspace{0.3cm} \hspace{0.3cm} On montre dans ([Sc 66] Lemma {\bf
5.3}) le r\'esultat int\'eressant suivant:

\vspace{0.2cm}
\begin{lemma} Soit $A$ une $K$-alg\`ebre \`a puissances associatives et sans diviseurs de z\'ero. Si $A$ contient un
idempotent non nul $e,$ alors $A$ est unitaire d'unit\'e $e.\Box$
\end{lemma}

\vspace{0.2cm}
\begin{definition} {\em (Alg\`ebre alg\'ebrique).} Soit $A$ une $K$-alg\`ebre \`a puissances associatives unitaire
et soit $K[X]$ l'alg\`ebre des polyn\^omes \`a une
ind\'etermin\'ee \`a coefficients dans $K.$ Un \'el\'ement $a\in
A$ est dit alg\'ebrique s'il existe un polyn\^ome $P\in
K[X]-\{0\}$ tel que \ $P(a)=0,$ ce qui est \'equivalent \`a dire
que la dimension de la $K$-alg\`ebre $K_A[a]$ est finie.
L'alg\`ebre $A$ est dite alg\'ebrique si tous ces \'el\'ements
sont alg\'ebriques.$\Box$
\end{definition}

\vspace{0.3cm}
\begin{flexible algebra} Une $K$-alg\`ebre $A$ est dite flexible si elle satisfait \`a l'une des
deux identit\'es \'equivalentes suivantes: \begin{enumerate} \item
$(x,y,x)=0$ pour tous $x,y\in A.$ \item $[L_x,R_x]=0$ pour tout
$x\in A.$
\end{enumerate}

Il est clair qu'une $K$-alg\`ebre associative ou commutative est
flexible.$\Box$
\end{flexible algebra}

\vspace{0.3cm}
\begin{proposition} Dans une $K$-alg\`ebre $A,$ les deux op\'erateurs lin\'eaires \[ L_x(L_x+R_x)-L_{x^2} \ \mbox{ et }
\ R_x(L_x+R_x)-R_{x^2} \] coincident. On les note $U_x.$
\end{proposition}

\vspace{0.1cm} {\bf Preuve.} La proposition s'obtient par une
lin\'earisation de l'identit\'e $(x,x,x)=0$ et en tenant compte de
la flexibilit\'e de $A.\Box$

\vspace{0.3cm}
\begin{alternative algebra} Une $K$-alg\`ebre $A$ est dite alternative si elle satisfait aux deux identit\'es
suivantes \ $(y,x,x)=(x,x,y)=0$ pour tous $x,y\in A.\Box$
\end{alternative algebra}

\vspace{0.3cm} \hspace{0.3cm} Il est bien connu, selon un
Th\'eor\`eme d'Artin, qu'une $K$-alg\`ebre $A$ est alternative si
et seulement si deux \'el\'ement quelconques de $A$ engendrent une
sous-alg\`ebre associative ([Sc 66] p. 29). De plus $A$ satisfait
aux trois identit\'es de Moufang suivantes ([Sc 66] p. 28):

\vspace{3cm}
\begin{enumerate} \item $x(y(xz))=((xy)x))z$
\hspace{0.3cm} (\`a gauche). \vspace{0.1cm} \item
$((zx)y)x=z((xy)x)$ \hspace{0.3cm} (\`a droite). \vspace{0.1cm}
\item $x(yz)x=(xy)(zx)$   \hspace{0.3cm} (moyenne).
\end{enumerate}

\vspace{0.3cm} \hspace{0.3cm} Les alg\`ebres alternatives sont
proches des associatives, comme le montre le Th\'eor\`eme de Artin
suivant ([Sc 66] p 29):

\vspace{0.3cm}
\begin{Jordan algebra} Une $K$-alg\`ebre commutative $A$ est dite de Jordan si elle satisfait \`a l'identit\'e de Jordan
suivante \ {\bf (J)} $(x^2,y,x)=0$ pour tout $x\in A.\Box$
\end{Jordan algebra}

\vspace{0.3cm}
\begin{examples} . \begin{enumerate} \item La sym\'etrisation $A^+$ d'une $K$-alg\`ebre associative est de Jordan.
\item Soit $V$ un $K$-espace vectoriel muni d'une forme
bilin\'eaire sym\'etrique $f,$ alors l'espace vectoriel \ $K\times
V$ \ muni du produit

\[ (\alpha,x)(\beta,y)=\Big( \alpha\beta+f(x,y),\alpha y+\beta x\Big) \]

est une $K$-alg\`ebre de Jordan unitaire, appel\'ee l'alg\`ebre de
Jordan associ\'ee \`a la forme bilin\'eaire sym\'etrique $f,$
not\'ee \ $J(V,f).\Box$ \end{enumerate}
\end{examples}

\vspace{0.3cm}
\begin{noncommutative Jordan algebra} Une alg\`ebre flexible $A$ est dite de Jordan non commutative si elle
satisfait \`a l'identit\'e {\bf (J)} de Jordan. Toute alg\`ebre
alternative est, d'apr\`es le Th\'eor\`eme d'Artin, de Jordan non
commutative.$\Box$
\end{noncommutative Jordan algebra}

\vspace{0.3cm} \hspace{0.3cm} On trouve dans [BKo 66] les deux
r\'esultats importants suivants:

\vspace{0.3cm}
\begin{proposition} Soit $A$ une $K$-alg\`ebre flexible, alors $A$ est de Jordan non commutative si et seulement si sa
sym\'etrisation $A^+$ est de Jordan.$\Box$
\end{proposition}

\vspace{0.3cm}
\begin{proposition} Soit $A$ une $K$-alg\`ebre de Jordan non commutative, alors $A$ est \`a puissances
associatives. De plus, pour tout entier $m\geq 2$ et pour tout
$x\in A,$ $L_{x^m}$ et $R_{x^m}$ appartiennent \`a la
sous-alg\`ebre de $End_K(A)$ engendr\'ee par $L_x,$ $R_x,$
$L_{x^2},$ $R_{x^2}$ et on a: \begin{enumerate} \item
$R_{x^{m+1}}=R_{x^m}(L_x+R_x)-U_xR_{x^{m-1}},$ \item
$L_{x^{m+1}}=L_{x^m}(L_x+R_x)-U_xL_{x^{m-1}}.\Box$ \end{enumerate}
\end{proposition}

\vspace{0.3cm}
\begin{Lie derivations algebra} Une $K$-alg\`ebre anti-commutative $L$ est dite de Lie si elle satisfait \`a l'identit\'e
de Jacobi suivante: \ $(xy)z+(yz)x+(zx)y=0$ pour tous $x,y,z\in
L.$ Soit maintenant $A$ une $K$-alg\`ebre, on appelle d\'erivation
de $A,$ tout op\'erateur lin\'eaire $\partial$ de $A$ satisfaisant
\`a l'une des trois conditions \'equicalentes suivantes:
\begin{enumerate}
\item $\partial(xy)=(\partial x)y+x\partial y$ pour tous $x,y\in
A.$ \item $[\partial, L_x]=L_{\partial x}$ pour tout $x\in A.$
\item $[\partial, R_x]=R_{\partial x}$ pour tout $x\in A.\Box$
\end{enumerate}
\end{Lie derivations algebra}

\vspace{0.3cm}
\begin{notes} Soit $A$ une $K$-alg\`ebre. \begin{enumerate} \item Si $A$ est
unitaire, alors $\partial 1=0.$ \item L'ensemble de toutes les
d\'erivations de $A,$ muni de sa structure naturelle de $K$-espace
vactoriel et du produit \ $[f,g]=fg-gf$ \ est une alg\`ebre de Lie
not\'ee $Der(A)$ et appel\'ee alg\`ebre de Lie des d\'erivations
de $A$ {\em [J 62], [Sc 66]}. \item On note $Aut(A)$ le groupe des
automorphismes de $A,$ celui de l'alg\`ebre r\'eelle $\oit$ de
Cayley-Dickson est not\'e $G_2$ {\em [GG 73]}. \item Si $A$ est
associative, les op\'erateurs lin\'eaires \ $L_x-R_x,$ $x\in A,$
sont des d\'erivations de $A$ dites int\'erieures. Si de plus, $A$
est unitaire, les op\'erateurs lin\'eaires \ $L_{x^{-1}}R_x,$ o\`u
$x$ est un \'el\'ement inversible de $A,$ sont des automorphismes
de $A,$ dits int\'erieurs.$\Box$
\end{enumerate}
\end{notes}

\vspace{0.3cm} \hspace{0.3cm} Nous avons les deux r\'esultats bien
connus suivants:

\vspace{0.3cm}
\begin{Sckolem-Noether theorem} Tout automorphisme d'une $K$-alg\`ebre associative simple centrale de dimension finie, est
int\'erieur {\em ([Pi 82] p. 230)}.$\Box$
\end{Sckolem-Noether theorem}

\vspace{0.3cm}
\begin{theorem} Toute d\'erivation d'une $K$-alg\`ebre associative simple centrale de dimension finie, est
int\'erieure {\em [J 37]}.$\Box$
\end{theorem}

\vspace{0.3cm} \hspace{0.3cm} Nous avons \'etabli le r\'esultat
pr\'eliminaire utile suivant:

\vspace{0.3cm}
\begin{lemma} Soit $A$ une $K$-alg\`ebre et soit $\lambda\in K,$ alors: \ $Der(A)\subseteq Der(A^{(\lambda)}).$ Si
de plus, $\lambda\neq\frac{1}{2},$ alors \
$Der(A)=Der(A^{(\lambda)}).$
\end{lemma}

\vspace{0.1cm} {\bf Preuve.} L'op\'erateur de multiplication \`a
gauche par un \'el\'ement $a\in A,$ dans $A^{(\lambda)}$ s'\'ecrit
\ $L_a^{(\lambda)}=\lambda L_a+(1-\lambda)R_a.$ Si $\partial\in
Der(A),$ on a

\begin{eqnarray*} [\partial, L_a^{(\lambda)}] &=& \lambda[\partial,
L_a]+(1-\lambda)[\partial, R_a] \\
&=& \lambda L_{\partial a}+(1-\lambda)R_{\partial a} \\
&=& L_{\partial a}^{(\lambda)}
\end{eqnarray*}

i.e. $\partial\in Der(A^{(\lambda)}).$ Si de plus,
$\lambda\neq\frac{1}{2},$ l'\'egalit\'e
$Der(A)=Der(A^{(\lambda)})$ est \'etablie du fait que \ $A=\Big(
A^{(\lambda)}\Big)^{(\mu)}$ o\`u
$\mu=\lambda(2\lambda-1)^{-1}.\Box$

\vspace{0.3cm}
\begin{Lie group and its Lie algebra} Soit $M$ une vari\'et\'e de dimension $n,$ alors les espaces tangents \`a $M:$ \
$T_a(M)$ o\`u $a$ parcourt $M,$ forment une vari\'et\'e
diff\'erentiable $T(M)$ de dimension $2n$ qui se projette
canoniquement sur $M.$ La projection \ $\pi:T(M)\rightarrow M$ \
associe \`a tout vecteur $B,$ son "point d'application" i.e. un
point $a\in M$ tel que $B\in T(M),$ de sorte que \
$T_a(M)=\pi^{-1}(a).$ Les applications diff\'erentiables \
$X:M\rightarrow T(M)$ \ $a\mapsto x_a$ \ telles que \ $\pi\circ
X=Id_M$ \ i.e. \ $X_a\in T_a(M)$ \ s'appellent champs de vecteurs
sur $M$ {\em [Po 85]}. Si $X$ et $Y$ sont deux champs de vecteurs
sur $M,$ le "crochet de Lie" de $X$ et $Y$ est d\'efini par \
$[X,Y]_a(f)=X_a(Yf)-Y_a(Xf)$ {\em [War 83]}. Un groupe de Lie est
un groupe $G$ muni d'une structure de vari\'et\'e diff\'erentiable
et tel que les applications \ $G\times G\rightarrow G$ \
$(x,y)\mapsto xy$ \ et \ $G\rightarrow G$ \ $x\mapsto x^{-1}$ \
sont diff\'erentiables. Soit $f:G\rightarrow G$ une fonction
diff\'erentiable et soit $a\in G,$ la diff\'erentielle de $f$ au
point $a$ est l'application lin\'eaire \ $df:T_a(G)\rightarrow
T_{f(a)}(G)$ \ d\'efinie de la mani\`ere suivante. Si $B\in
T_a(G),$ alors \ $df(B)\in T_{f(a)}(G)$ et si $g$ est une fonction
diff\'erentiable sur un voisinage de $f(a),$ on a \
$df(B)(g)=B(g\circ f).$ Un champs de vecteurs $X$ sur $G$ est dit
invariant \`a gauche si pour tout $a\in G,$ on a \ $(dL_a)\circ
X=X\circ L_a$ ($L_a:G\rightarrow G$ \'etant la translation $a$ \`a
gauche) {\em [War 83]}. L'espace vectoriel des champs invariants
\`a gauche, not\'e ${\cal I}(G),$ muni du crochet de Lie, est une
alg\`ebre de Lie appel\'ee l'alg\`ebre de Lie du groupe de Lie $G$
{\em ([Po 85], [War 83])}.$\Box$
\end{Lie group and its Lie algebra}

\vspace{0.3cm}
\begin{proposition} Soit $A$ une $\rit$-alg\`ebre de dimension finie, alors son groupe d'automorphismes $Aut(A)$
est un groupe de Lie dont l'alg\`ebre de Lie et l'alg\`ebre
$Der(A)$ coincident: \ ${\cal I}\Big( Aut(A)\Big)=Der(A).$
\end{proposition}

\vspace{0.1cm} {\bf Preuve.} ([Po 85] p. 202).$\Box$

\vspace{1cm}
\subsection{Inversibilit\'e dans les alg\`ebres non associatives}

\vspace{0.5cm}
\begin{linear invertibility} Soit $A$ une $K$-alg\`ebre, un \'el\'ement $x\in A$ est dit lin\'eairement
inversible (l.i.) si les op\'erateurs lin\'eaires $L_x,$ $R_x,$ de
multiplication par $x,$ sont inversibles dans $End_K(A).$
L'ensemble des \'el\'ements l.i. de $A$ est not\'e $L-inv(A).$ Un
\'el\'ement $x\in A$ tel que $L_x$ (resp. $R_x$) est inversible
dans $End_K(A)$ est dit lin\'eairement inversible \`a gauche
(l.i.g.) (resp. \`a droite). L'ensemble des \'el\'ements l.i.g.
(resp. \`a droite) est not\'e $L-inv_g(A)$ (resp. $L-inv_d(A)$).
Si $A$ n'est pas r\'eduite \`a $\{0\},$ elle est dite de division
lin\'eaire \`a gauche (d.l.g.) (resp. \`a droite) si
$L-inv_g(A)=A-\{0\}$ (resp. $L-inv_d(A)=A-\{0\}$). $A$ est dite de
division lin\'eaire (d.l.) si elle l'est \`a gauche et \`a
droite.$\Box$
\end{linear invertibility}

\vspace{0.2cm}
\begin{proposition} Soit $A$ une $K$-alg\`ebre de dimension finie, alors les quatre propri\'et\'es suivantes
sont \'equivalentes: \begin{enumerate} \item $A$ est de division
lin\'eaire. \item $A$ est de division lin\'eaire \`a gauche. \item
$A$ est de division lin\'eaire \`a droite. \item $A$ est sans
diviseurs de z\'ero.$\Box$
\end{enumerate}
\end{proposition}

\vspace{0.3cm} \hspace{0.3cm} Nous avons le r\'esultat, bien
connu, suivant:

\vspace{0.3cm}
\begin{theorem} Soit $K_0$ un corps alg\'ebriquement clos et soit $A\neq\{0\}$ une $K_0$-alg\`ebre unitaire sans
diviseurs de z\'ero de dimension finie. Alors $A\simeq K_0.\Box$
\end{theorem}

\vspace{0.3cm}
\begin{note} Dans {\em ([Sc 66] p. 134)} on montre qu'une
$K$-alg\`ebre \`a puissances associatives sans diviseurs de z\'ero
de dimension finie, contient un \'el\'ement unit\'e.$\Box$
\end{note}

\vspace{0.3cm} \hspace{0.3cm} Il est facile de voir qu'une
$K$-alg\`ebre associative de division lin\'eaire est unitaire.
Nous allons montrer que ce r\'esultat reste valable notament pour
les alg\`ebres de Jordan non commutatives.

\vspace{0.3cm}
\begin{lemma} Soit $A$ une $K$-alg\`ebre de Jordan non commutative telle que $L-inv_g(A)\neq\emptyset.$ Alors $A$
contient un idempotent non nul.
\end{lemma}

\vspace{0.1cm} {\bf Preuve.} Soit $a\in L-inv_g(A),$ il existe
$e\in A-\{0\}$ tel que $ae=a$ et on a \ $a^2=(ae)a=a(ea).$ Donc
$ea=a$ car $L_a$ est injectif. D'apr\`es la relation 1) de la
Proposition {\bf 1.17}, on a \begin{eqnarray*} ae^3 &=& R_{e^3}(a) \\
&=& \Big( R_{e^2}(L_e+R_e)-U_eR_e\Big)(a) \\
&=& (ea+ae)e^2-U_e(ae) \\
&=& 2ae^2-\Big( R_e(L_e+R_e)-R_{e^2}\Big)(a) \\
&=& 2ae^2-(2ae-ae^2) \\
&=& a(3e^3-2e) \end{eqnarray*}

Donc $e^3=3e^2-2e$ car $L_a$ est injectif. Si $e^2=e,$ c'est fini.
Sinon on pose $e_0=\frac{1}{2}(e^2-e)$ et on a \begin{eqnarray*} e_0^2 &=& \frac{1}{4}(e^4-2e^3+e^2) \\
&=& \frac{1}{4}(e^3-e^2) \\
&=& \frac{1}{2}(e^2-e) \\
&=& e_0. \end{eqnarray*}

Donc $e_0$ est un idempotent non nul de $A.\Box$

\vspace{0.3cm}
\begin{corollary} Soit $A$ une $K$-alg\`ebre de Jordan non commutative sans diviseurs de z\'ero et telle que
$L-inv_g(A)\neq\emptyset.$ Alors $A$ contient un \'el\'ement
unit\'e.
\end{corollary}

\vspace{0.1cm} {\bf Preuve.} $A$ contient un idempotent non nul.
La Proposition {\bf 1.17} et le Lemme {\bf 1.8} ach\`event la
d\'emonstration.$\Box$

\vspace{0.3cm}
\begin{remark} Ce dernier r\'esultat ne persiste pas en g\'en\'eral pour les alg\`ebres non associatives. En effet,
l'alg\`ebre $\stackrel{*}{\cit}$ de Mc Clay est de division
lin\'eaire (commutative ayant un idempotent non nul) non
unitaire.$\Box$
\end{remark}

\vspace{0.3cm}
\begin{J-invertibility} Soit $A$ une $K$-alg\`ebre de Jordan non commutative unitaire. Un \'el\'ement $x\in A$ est
dit inversible au sens de Jacobson, ou $J$-inversible, s'il existe
$y\in A,$ dit inverse de $x,$ tel que \begin{enumerate} \item
$xy=yx=1$ et \item $x^2y=yx^2=x$ {\em [Mc 65]}.
\end{enumerate}
\end{J-invertibility}

\vspace{0.3cm} \hspace{0.3cm} Si $A$ est alternative, les
relations 1) et 2) de {\bf 1.32} se r\'eduisent \`a 1). Les
propri\'et\'es fondamentales sont donn\'ees par les deux
r\'esultats suivants:

\vspace{0.3cm}
\begin{theorem} {\em (Jacobson)} Soit $A$ une $K$-alg\`ebre de Jordan unitaire et soit $x$ un \'el\'ement de $A.$ Alors
les propri\'et\'es suivantes sont \'equivalentes:
\begin{enumerate} \item $x$ est $J$-inversible dans $A.$ \item
$1\in U_x(A).$ \item L'op\'erateur $U_x=2L_x^2-L_{x^2}$ est
inversible dans $End_K(A).$
\end{enumerate}

Dans ces conditions, l'inverse $y$ est unique et on a \
$y=U_x^{-1}(x)$ {\em [J 68]}.$\Box$
\end{theorem}

\vspace{0.3cm}
\begin{theorem} {\em (Mc Crimmon)} Soit $A$ une $K$-alg\`ebre de Jordan non commutative unitaire. Alors un \'el\'ement
$x\in A$ est $J$-inversible et poss\`ede pour inverse $y\in A$ si
et seulement si $y$ est l'inverse de $x$ dans l'alg\`ebre de
Jordan $A^+$ {\em [Mc 65]}.$\Box$
\end{theorem}

\vspace{0.3cm}
\begin{remark} Le Th\'eor\`eme {\bf 1.34} montre que
l'inversibilit\'e au sens de Jacobson dans une alg\`ebre $A$ de
Jordan non commutative unitaire est la m\^eme que dans l'alg\`ebre
de Jordan $A^+.$ On notera l'inverse de $x$ par $x^{-1}$ et
l'ensemble des \'el\'ements $J$-inversibles de $A,$ par
$inv(A).\Box$
\end{remark}

\vspace{0.3cm}
\begin{definition} Soit $A$ une $K$-alg\`ebre de
Jordan non commutative unitaire. L'alg\`ebre $A$ est dite de
$J$-division si $inv(A)=A-\{0\}.$ Une telle alg\`ebre est simple
et il en est de m\^eme pour l'alg\`ebre de Jordan $A^+.\Box$
\end{definition}

\vspace{0.3cm}
\begin{theorem} {\em (Th\'eor\`eme de Artin avec inverses)} Dans une $K$-alg\`ebre alternative unitaire, deux \'el\'ements
quelconques et leur inverse, s'il existe, engendrent une
sous-alg\`ebre associative {\em ([Kai 77] p. 39)}.$\Box$
\end{theorem}

\vspace{0.3cm}
\begin{proposition} Soit $A$ une $K$-alg\`ebre de Jordan non commutative unitaire, alors tout \'el\'ement $x\in A$
lin\'eairement inversible \`a gauche (resp. \`a droite) est
$J$-inversible, d'inverse $x^{-1},$ lin\'eairement inversible \`a
droite (resp. \`a gauche) et on a \ $R_{x^{-1}}=U_x^{-1}l_x$ \
(resp. $L_{x^{-1}}=U_x^{-1}R_x$) {\em ([Kai 77] p. 31)}.$\Box$
\end{proposition}

\vspace{0.3cm}
\begin{corollary} Soit $A$ une $K$-alg\`ebre alternative unitaire et soit $x\in A.$ Alors les deux propri\'et\'es
suivantes sont \'equivalentes: \begin{enumerate} \item $x$ est
$J$-inversible. \item $x$ est lin\'eairement inversible.
\end{enumerate}
\end{corollary}

\vspace{0.1cm} {\bf Preuve.} L'implication 2) $\Rightarrow$ 1)
d\'ecoule de la Proposition {\bf 1.38}. Soit maintenant $x\in
inv(A)$ et soit $z\in A.$ D'apr\`es le Th\'eor\`eme {\bf 1.37},
$x(x^{-1}z)=(xx^{-1})z=z$ i.e. $x\in L-inv_g(A).$ De m\^eme $x\in
L-inv_d(A).$ Ceci montre l'implication 1) $\Rightarrow$ 2).$\Box$

\vspace{0.3cm} \hspace{0.3cm} L'inversibilit\'e lin\'eaire \`a
gauche coincide avec l'inversibilit\'e au sens de Jacobson pour
les alg\`ebres alternatives (unitaires) {\bf 1.39} et entra\^ine
cette derni\`ere pour les alg\`ebres de Jordan non commutatives
(unitaires) {\bf 1.38}. Cependant ces deux notions
d'inversibilit\'e ne sont pas \'equivalen tes pour les alg\`ebres
de Jordan non commutatives, et les contre-exemples existent en
abondance ([Kai] p. 19), [Pete 81].$\Box$

\vspace{0.3cm}
\begin{definition} Soit $A$ une $K$-alg\`ebre de Jordan n.c. unitaire. Une sous-alg\`ebre $B$ de $A,$ qui
contient l'\'el\'ement unit\'e de $A,$ est dite pleine si \
$inv(B)=B\cap inv(A).\Box$
\end{definition}

\vspace{0.3cm}
\begin{proposition} Tout \'el\'ement $a$ d'une $K$-alg\`ebre $A$ de Jordan n.c. unitaire, est contenu
dans une sous-alg\`ebre associative, commutative et pleine {\em
([Kai 77] p. 33)}.$\Box$
\end{proposition}

\vspace{0.3cm}
\begin{note} La plus petite sous-alg\`ebre pleine de $A$ qui contient $a,$ not\'ee $K(a),$ est appel\'ee la
sous-alg\`ebre pleine engendr\'ee par $a.$ Evidemment $K(a)$
contient $K_A[a].\Box$
\end{note}

\vspace{1cm}
\subsection{Alg\`ebres quadratiques}

\vspace{0.3cm}
\begin{definitions and notations} .\begin{enumerate} \item Une $K$-alg\`ebre unitaire $A$ est dite quadratique si
$1,x,x^2$ sont lin\'erement d\'ependants pour tout $x\in A.$ \item
Il est bien connu qu'une $K$-alg\`ebre quadratique $A$ s'obtient
\`a partir d'une $K$-alg\`ebre anti-commutative $(V,\wedge)$ et
d'une forme bilin\'eaire $(.,.)$ sur $V,$ en munissant \
$A=K\times V$ \ de sa structure naturelle de $K$-espace vectoriel
et du produit

\[ (\alpha+x)(\beta+y)=\Big( \alpha\beta+(x,y)\Big)+(\alpha y+\beta x+x\wedge y). \]

On note, de la m\^eme mani\`ere, la forme bilin\'eaire \[ A\times
A\rightarrow K \hspace{0.2cm}
(\alpha+x,\beta+y)\mapsto\alpha\beta+(x,y), \] qu'on appelle la
forme bilin\'eaire associ\'ee \`a $A.$ $(V,\wedge)$ est appel\'ee
l'alg\`ebre anti-commutative associ\'ee \`a $A,$ ses \'el\'ements
sont appel\'es vecteurs, ceux de $K$ scalaires. On note $A$ par
$\Big( V,(.,.),\wedge\Big).$ \item L'alg\`ebre $J(V,f)$ {\em
(Exemple {\bf 1.14 2)}} est quadratique.$\Box$
\end{enumerate}
\end{definitions and notations}

\vspace{3cm}
\begin{remarks} .\begin{enumerate} \item Une $K$-alg\`ebre quadratique et flexible, est de Jordan non commutative.
\item Si $A=\Big( V,(.,.),\wedge\Big)$ est une $K$-alg\`ebre
quadratique, sa sym\'etrisation $A^+$ est de Jordan et quadratique
associ\'ee \`a la forme bilin\'eaire sym\'etrique \[
(x|y)=\frac{1}{2}\Big( (x,y)+(y,x)\Big), \]

dite sym\'etrisation de $(.,.).\Box$ \end{enumerate}
\end{remarks}

\vspace{0.3cm}
\begin{definition} Une $K$-alg\`ebre quadratique $A$ est dite $Q$-simple, si l'alg\`ebre de Jordan $A^+$ est simple.
On dira \'egalement que $A$ est de division, si l'alg\`ebre de
Jordan $A^+$ est de division.$\Box$
\end{definition}

\vspace{0.3cm}
\begin{proposition} {\em (Kaidi).} Soit $A=\Big( V,(.,.),\wedge\Big)$ une \ $\rit$-alg\`ebre quadratique et soit
$q:V\rightarrow\rit$ \ $x\mapsto q(x)=(x|x)$ \ la forme
quadratique associ\'ee \`a la sym\'etrisation $(.|.)$ de $(.,.).$
Alors $A$ est de division si et seulement si $q$ est d\'efinie
n\'egative.
\end{proposition}

\vspace{0.1cm} {\bf Preuve.} ([Kai 77] p.
98).$\Box$

\vspace{0.3cm}
\begin{proposition} Soit $A$ une \ $\rit$-alg\`ebre quadratique flexible $Q$-simple de dimension finie. Alors
$A$ est de division lin\'eaire si et seulement si $A^{(\lambda)}$
est de division lin\'eaire pour tout
$\lambda\in\rit-\{\frac{1}{2}\}.$
\end{proposition}

\vspace{0.1cm} {\bf Preuve.} [A 48].$\Box$

\vspace{0.3cm}
\begin{remark} L'alg\`ebre r\'eelle $\hit$ des quaternions de Hamilton est quadratique flexible $Q$-simple et de
division. cepandant, sa sym\'erisation $\hit^+$ n'est pas de
division lin\'eaire. Price {\em ([Pr 51] p. 294)} a montr\'e que
si $A$ est une $K$-alg\`ebre associative centrale de dimension
finie impaire, alors l'alg\`ebre de Jordan $A^+$ est de division
lin\'eaire.$\Box$
\end{remark}

\vspace{3cm}
\begin{theorem} {\em (Osborn).} Soit $A=\Big( V,(.,.),\wedge\Big)$ une \ $K$-alg\`ebre quadratique, avec
la propri\'et\'e que la sous-alg\`ebre engendr\'ee par un
\'el\'ement quelconque de $A,$ soit un corps. Alors les
propri\'et\'es suivantes sont \'equivalentes:
\begin{enumerate}
\item $A$ est sans diviseurs de z\'ero. \item $A$ ne poss\`ede pas
de sous-alg\`ebres de dimension $3.$ \item Pour deux vecteurs
$x,y$ de $A,$ lin\'eairement ind\'ependants, les vecteurs
$x,y,x\wedge y$ sont lin\'eairement ind\'ependants. \item Il
n'existe pas de vecteurs $x,y$ de $A,$ lin\'eairement
ind\'ependants tels que \ $x\wedge y=x$ ou $x\wedge y=0.$
\end{enumerate}
\end{theorem}

\vspace{0.1cm} {\bf Preuve.} [Os 62].$\Box$

\vspace{0.5cm} \hspace{0.3cm} Osborn [Os 62] a donn\'e une
classification pour toutes les alg\`ebres quadratiques de division
lin\'eaire de dimension $4,$ sur un corps commutatif $K_0$ de
caract\'eristique $\neq 2,$ en montrant le r\'esultat suivant:

\vspace{0.3cm}
\begin{theorem} {\em (Osborn).} Dans une $K$-alg\`ebre anti-commutative $(V,\wedge)$ admissible pour la division de
dimension $3,$ il existe une base \ $\{x,y,z\}$ \ telle que
\begin{eqnarray} y\wedge z=x, \hspace{0.1cm} z\wedge x=\alpha
y+\beta z \mbox{ et } x\wedge y=\gamma z \mbox{ o\`u }
\beta\in\{0,1\} \mbox{ et } \alpha, \gamma\in K-\{0\}.
\end{eqnarray} Inversement, si les \'el\'ements de la base d'une alg\`ebre
anti-commutative $(V,\wedge)$ satisfont \`a {\bf (1.2)}, alors
$(V,\wedge)$ est admissible pour la division si et seulement si la
forme quadratique \
$\lambda_1^2+\beta\lambda_1\lambda_2+\alpha\gamma\lambda_2^2+\alpha\lambda_3^2+\beta\lambda_3\lambda_4+\gamma\lambda_4^2$
\ est non d\'eg\'en\'er\'ee.
\end{theorem}

\vspace{0.1cm} {\bf Preuve.} [Os 62].$\Box$

\vspace{0.5cm} \hspace{0.3cm} Osborn [Os 62] a distingu\'e une
classe importante d'alg\`ebres. Nous donnons, pour cela, la
d\'efinition suivante:

\vspace{0.3cm}
\begin{definition} On dit qu'une $K$-alg\`ebre $A$ satisfait \`a la propri\'et\'e d'Osborn si deux \'el\'ements
quelconques de $A$ qui n'appartiennent pas \`a la m\^eme
sous-alg\`ebre de dimension $2,$ engendrent une sous-alg\`ebre de
dimension $4$ {\em [Os 62]}.$\Box$ \end{definition}

\vspace{0.3cm}
\begin{cayley algebra} Une $K$-alg\`ebre unitaire $A$ est dite cayleyenne
si elle est munie d'une involution multiplicative (ou
anti-automorphisme involutif) \ $s:x\mapsto\overline{x}$ \ telle
que \ $x+\bar{x},$ $x\bar{x}\in K.1$ {\em ([Bou 70] A III p. 15)}.
L'involution $s$ est unique, appel\'ee la conjugaison cayleyenne
de $A.$ On v\'erifie facilement que $A$ est quadratique i.e.
$A=(V,(.,.),\wedge)$ et que $s$ est d\'efinie par \ $A=K\oplus
V\rightarrow A$ \ $\alpha+u\mapsto s(\alpha+u)=\alpha-u.\Box$
\end{cayley algebra}

\vspace{0.3cm}
\begin{definitions} Soit $A$ une $K$-alg\`ebre quadratique munie d'une forme bilin\'eaire sym\'etrique $(.|.)$ et
soit $f\in End_K(A).$ On dit que \begin{enumerate} \item $f$ est
isom\'etrique, par rapport \`a $(.|.),$ si $(f(x)|f(y))=(x|y)$
pour tous $x,y\in A.$ \item $f$ est anti-sym\'etrique, par rapport
\`a $(.|.),$ si $(f(x)|y))=-(x|f(y))$ pour tous $x,y\in A.$ \item
$(.|.)$ est une forme trace sur $A$ si $(xy|z)=(x|yz)$ pour tous
$x,y,z\in A.\Box$ \end{enumerate}
\end{definitions}

\vspace{0.3cm}
\begin{lemma} {\em (Osborn).} Soit $B=(V,[.,.),\wedge)$ une $K$-alg\`ebre quadratique. Alors
\begin{enumerate}
\item $B$ est cayleyenne si et seulement si $(.,.)$ est
sym\'etrique. \item $B$ est flexible si et seulement si $(.,.)$
est sym\'etrique et l'une des trois propri\'et\'es \'equivalentes
suivantes a lieu: \begin{enumerate} \item $(.,.)$ est une forme
trace sur $B.$ \item $(.,.)$ est une forme trace sur $(V,\wedge).$
\item $(x\wedge y,x)=0$ pour tout $x\in V.\Box$
\end{enumerate} \end{enumerate}
\end{lemma}

\vspace{0.3cm}
\begin{remark} L'alg\`ebre $J(V,f)$ est cayleyenne.$\Box$
\end{remark}

\vspace{0.3cm}
\begin{note} Si $A=(V,(.|.),\wedge)$ est une $K$-alg\`ebre cayleyenne, on dira que deux \'el\'ements $x,y\in A$ sont
orthogonaux si \ $(x|y)=0$ et on notera \ $S^{\perp}$ \ le
sous-espace vectoriel de $A$ orthogonal \`a une partie $S$ de $A.$
On a \'evidemment \ $V=(K 1)^{\perp}.\Box$
\end{note}

\vspace{3cm}
\begin{remarks} Si $(B,\overline{.})=(V,(.|.),\wedge)$ est une $K$-alg\`ebre cayleyenne, alors
\begin{enumerate} \item La conjugaison cayleyenne \ $\overline{.}$ \ est isom\'etrique par rapport \`a la forme
bilin\'eaire sym\'etrique $(.|.).$ \item Pour tous $x,y\in B,$ on
a \ $\overline{[x,y]}=[x,\overline{y}]=-[x,y].$ Si, de plus, $B$
est flexible, les identit\'e suivantes ont lieu pour tous
$x,y,z\in B$ et pour tout $\alpha\in K:$ \begin{enumerate} \item
$(x.^{\alpha}y|z)=(x|y.^{\alpha}z)=(y|z.^{\alpha}x),$ \item
$([x,y]|z)=(x|[y,z])=(y|[z,x]).$ \end{enumerate} \item Pour tout
$\lambda\in K,$ l'application \ $B^{(\lambda)}\rightarrow
B^{(1-\lambda)}$ \ $\alpha+u\mapsto\alpha-u$ \ est un isomorphisme
d'alg\`ebres.$\Box$
\end{enumerate}
\end{remarks}

\vspace{0.2cm}
\begin{CD process} Soit $(B,\overline{.})$ une $K$-alg\`ebre cayleyenne et soit $\gamma\in K,$ avec $\gamma\neq 0.$
On appelle extension cayleyenne de $(B,\overline{.}),$ d'indice
$\gamma,$ l'alg\`ebre cayleyenne not\'ee \ $E_\gamma(B)$ \ ayant
pour espace vectoriel sous-jacent \ $B\times B,$ \ muni de produit
\[ (x,y)(x',y')=(xx'+\gamma\bar{y'}y,y\bar{x'}+y'x) \] et de la conjugaison cayleyenne \ $s(x,y)=(\bar{x},-y)$
{\em ([Bou 70] A III p. 16)}. $B\times\{0\}$ est une
sous-alg\`ebre de $E_\gamma(B),$ isomorphe \`a $B,$ qu'on
identifie \`a $B.$ Si $f=(0,1),$ alors tout \'el\'ement $(x,y)$ de
$E_\gamma(B)$ s'\'ecrit d'une mani\`ere unique sous la forme
$x+yf$ i.e. $E_\gamma(B)=B\oplus Bf.$ On v\'erifie dans {\em [Bou
70]} que $E_\gamma(B)$ est associative si et seulement si $B$ est
associative et commutative, et que $E_\gamma(B)$ est alternative
si et seulement si $B$ est associative. Ce proc\'ed\'e de
construction d'alg\`ebres cayleyennes par extension est appel\'e
la proc\'ed\'e de Cayley-Dickson. L'alg\`ebre \ $A_0=K$ \ munie de
l'application identique est une $K$-alg\`ebre cayleyenne, ainsi en
consid\'erant des scalaires $\gamma_i,$ $i\geq 1$ non nuls, et en
posant \ $A_i=E_{\gamma_i}(A_{i-1}),$ on obtient des
$K$-alg\`ebres cayleyennes. $A_1$ est associative et commutative.
$A_2$ est associative, simple et centrale {\em [Sc 66]} appel\'ee
une alg\`ebre des quaternions sur $K.$ $A_3$ est alternative,
simple et telle que \ $N(A_3)=K.1$ {\em [Sc 66]}. Elle est
appel\'ee une alg\`ebre des octonions sur $K.$ Si $K=\rit$ et
$\gamma_i=-1,$ on retrouve: \begin{enumerate} \item Le corps des
nombres complexes \ $\cit=A_1.$ \item L'alg\`ebre r\'eelles de
division des quaternions de Hamilton \[ \hit=A_2=\cit\oplus\cit j.
\] Elle poss\`ede une base \ ${\cal B}_H=\{1,i,j,k\},$ appel\'ee
la base canonique de $\hit,$ pour laquelle la table de
multiplication est donn\'ee par: \ $1$ est l'\'el\'ement unit\'e,
\[ i^2=j^2=k^2=-1, \ ij=-ji=k, \ jk=-kj=i, \ ki=-ik=j. \]
\item L'alg\`ebre r\'eelle de division des octonions de
Cayley-Dickson \[ \oit=A_3=\hit\oplus\hit f. \] Elle poss\`ede une
base \ ${\cal B}_O=\{1,i,j,k,f,if,jf,kf\},$ appel\'ee la base
canonique de $\oit,$ pour laquelle la table de multiplication est
donn\'ee par

\vspace{0.3cm} \[ \begin{tabular}{ccccccccc} \\
\multicolumn{1}{c}{} & \multicolumn{1}{c}{$1$} &
\multicolumn{1}{c}{$i$} & \multicolumn{1}{c}{$j$} &
\multicolumn{1}{c}{$k$} & \multicolumn{1}{c}{$f$} &
\multicolumn{1}{c}{$if$} & \multicolumn{1}{c}{$jf$} &
\multicolumn{1}{c}{$kf$} \\ \cline{2-9} \multicolumn{1}{c|}{$1$} &
\multicolumn{1}{|c}{$1$} & \multicolumn{1}{c}{$i$} &
\multicolumn{1}{c}{$j$} & \multicolumn{1}{c|}{$k$} &
\multicolumn{1}{|c}{$f$} & \multicolumn{1}{c}{$if$} &
\multicolumn{1}{c}{$jf$} & \multicolumn{1}{c|}{$kf$} \\
\multicolumn{1}{c|}{$i$} & \multicolumn{1}{|c}{$i$} &
\multicolumn{1}{c}{$-1$} & \multicolumn{1}{c}{$k$} &
\multicolumn{1}{c|}{$-j$} & \multicolumn{1}{|c}{$if$} &
\multicolumn{1}{c}{$-f$} & \multicolumn{1}{c}{$-kf$} &
\multicolumn{1}{c|}{$jf$} \\
\multicolumn{1}{c|}{$j$} & \multicolumn{1}{|c}{$j$} &
\multicolumn{1}{c}{$-k$} & \multicolumn{1}{c}{$-1$} &
\multicolumn{1}{c|}{$i$} & \multicolumn{1}{|c}{$jf$} &
\multicolumn{1}{c}{$kf$} & \multicolumn{1}{c}{$-f$} &
\multicolumn{1}{c|}{$-if$} \\
\multicolumn{1}{c|}{$k$} & \multicolumn{1}{|c}{$k$} &
\multicolumn{1}{c}{$j$} & \multicolumn{1}{c}{$-i$} &
\multicolumn{1}{c|}{$-1$} & \multicolumn{1}{|c}{$kf$} &
\multicolumn{1}{c}{$-jf$} & \multicolumn{1}{c}{$if$} &
\multicolumn{1}{c|}{$-f$} \\ \cline{2-9} \multicolumn{1}{c|}{$f$}
& \multicolumn{1}{|c}{$f$} & \multicolumn{1}{c}{$-if$} &
\multicolumn{1}{c}{$-jf$} & \multicolumn{1}{c|}{$-kf$} &
\multicolumn{1}{|c}{$-1$} & \multicolumn{1}{c}{$i$} &
\multicolumn{1}{c}{$j$} & \multicolumn{1}{c|}{$k$} \\
\multicolumn{1}{c|}{$if$} & \multicolumn{1}{|c}{$if$} &
\multicolumn{1}{c}{$f$} & \multicolumn{1}{c}{$-kf$} &
\multicolumn{1}{c|}{$jf$} & \multicolumn{1}{|c}{$-i$} &
\multicolumn{1}{c}{$-1$} & \multicolumn{1}{c}{$-k$} &
\multicolumn{1}{c|}{$j$} \\
\multicolumn{1}{c|}{$jf$} & \multicolumn{1}{|c}{$jf$} &
\multicolumn{1}{c}{$kf$} & \multicolumn{1}{c}{$f$} &
\multicolumn{1}{c|}{$-if$} & \multicolumn{1}{|c}{$-j$} &
\multicolumn{1}{c}{$k$} & \multicolumn{1}{c}{$-1$} &
\multicolumn{1}{c|}{$-i$} \\
\multicolumn{1}{c|}{$kf$} & \multicolumn{1}{|c}{$kf$} &
\multicolumn{1}{c}{$-jf$} & \multicolumn{1}{c}{$if$} &
\multicolumn{1}{c|}{$f$} & \multicolumn{1}{|c}{$-k$} &
\multicolumn{1}{c}{$-j$} & \multicolumn{1}{c}{$i$} &
\multicolumn{1}{c|}{$-1$} \\ \cline{2-9}
\end{tabular} \]

\vspace{0.3cm} $\{i,j,k,f,if,jf,kf\}$ est appel\'ee la base
canonique de l'alg\`ebre anti-commutative des vecteurs aasoci\'ee
\`a \ $\oit.$ On la notera $\{e_1,\dots,e_7\}.\Box$
\end{enumerate}
\end{CD process}

\vspace{0.3cm} \hspace{0.3cm} Nous avons les r\'esultats, bien
connus, suivants:

\vspace{0.3cm}
\begin{theorem} {\em (Frobenius)} Les uniques alg\`ebres r\'eelles
associatives alg\'ebriques de division sont de dimension finie et
isomorphes \`a $\rit,$ $\cit$ ou $\hit$ {\em [E-R 91]}.$\Box$
\end{theorem}

\vspace{0.3cm} \begin{theorem} {\em (Zorn)} Les uniques alg\`ebres
r\'eelles alternatives de division de dimension finie sont
isomorphes \`a \ $\rit,$ $\cit,$ $\hit$ ou $\oit$ {\em ([E-R 91]
p. 262)}.$\Box$
\end{theorem}

\vspace{0.3cm} \begin{theorem} {\em (Albert)} Les uniques
alg\`ebres r\'eelles alternatives alg\'ebriques de division sont
de dimension finie et isomorphes \`a \ $\rit,$ $\cit,$ $\hit$ ou
$\oit$ {\em ([A 49] p. 767)}.$\Box$
\end{theorem}

\vspace{4cm}
\section{Alg\`ebres non associatives norm\'ees de division}

\vspace{0.5cm} \hspace{0.3cm} Dans ce chapitre $K$ d\'esignera
$\rit$ ou $\cit.$

\vspace{0.7cm} \subsection{Alg\`ebres norm\'ees non associatives}

\vspace{0.3cm} \begin{Banach algebras} On dit qu'une $K$-alg\`ebre
$A$ est norm\'ee, si l'espace vectoriel $A$ est muni d'une norme \
$||.||$ \ sous multiplicative i.e. \ $||xy||\leq ||x||$ $||y||$ \
pour tous $x,y\in A,$ ce qui est \'equivalent \`a dire que la
norme $||.||$ rend continu le produit de $A$ et, plus
particuli\`erement, les op\'erateurs de multiplication de $A.$
\begin{enumerate}
\item Si l'espace vectoriel norm\'e $(A,||.||)$ est de Banach,
l'alg\`ebre $A$ est dite norm\'ee compl\`ete. \item Si, de plus,
$A$ est associative, elle est dite de Banach. \item L'alg\`ebre
$A$ est dite absolument valu\'ee si l'espace vectoriel $A$ est
muni d'une norme \ $||.||$ \ multiplicatiove i.e. \ $||xy||=||x||$
$||y||$ \ pour tous $x,y\in A.$ \item On note $\hat{B}$ la
compl\'etion d'une alg\`ebre norm\'ee $B$ {\em [BD 73]}.$\Box$
\end{enumerate}
\end{Banach algebras}

\vspace{0.2cm}
\begin{proposition} Dans une alg\`ebre de Banach unitaire $A,$ l'ensemble $inv(A)$ des \'el\'ements inversibles est
un ouvert de $A$ {\em [BD 73]}.$\Box$
\end{proposition}

\vspace{0.2cm}
\begin{definition} Soient $E, F$ deux espaces topologiques et soit, pour tout $x\in E,$ $\varphi(x)$ une partie de $F.$
L'application \ $\varphi:E\rightarrow{\cal P}(F)$ \ est dite
semi-continue sup\'erieurement si pour tout $x\in E$ et pour tout
voisinage $V$ de $\varphi(x),$ il existe un voisinage $U$ de $x$
tel que \ $\varphi(U)\subset V.\Box$
\end{definition}

\vspace{0.2cm}
\begin{note} Soit $A$ une $\cit$-alg\`ebre associative norm\'ee unitaire. Il est bien connu que
\begin{enumerate} \item Pour tout $x\in A,$ le spectre \
$sp_A(x)=\{\lambda\in\cit: x-\lambda.1\notin inv(A)\}$ \ de $x$
dans $A$ est non vide. \item Si, de plus, $A$ est de Banach, alors
$sp_A(x)$ est un compact de $\cit$ pour tout $x\in A$ et
l'application \ $A\rightarrow\{\mbox{compacts de } \ $\cit$\}$ \
$x\mapsto sp_A(x)$ \ est semi-continue sup\'erieurement {\em [BD
73]}.$\Box$
\end{enumerate}
\end{note}

\vspace{0.1cm} \hspace{0.3cm} Un des r\'esultats classiques de la
th\'eorie des $\cit$-alg\`ebres de Banach, est le suivant:

\vspace{0.2cm} \begin{Gelfand-Mazur theorem} Soit $A$ une
$\cit$-alg\`ebre associative norm\'ee de division. Alors $A\simeq$
$\cit$ {\em [BD 73]}.
\end{Gelfand-Mazur theorem}

\vspace{0.1cm} {\bf Preuve.} Soit $x\in A$ et soit $\lambda\in
sp_A(x)$ {\bf 2.4}, on a \ $x-\lambda 1\notin inv(A)=A-\{0\}$ i.e.
$x=\lambda 1.$ L'application \ $\cit\rightarrow A$ \
$\lambda\mapsto\lambda 1$ \ est alors un isomorphisme
d'alg\`ebres.$\Box$

\vspace{0.2cm}
\begin{note} Soient $(E,||.||)$ un $K$-espace vectoriel norm\'e, $S(E)=\{x\in E:||x||=1\}$ la sph\`ere unit\'e de $E,$ et
$End_K(E)$ la $K$-alg\`ebre associative et unitaire des
op\'erateurs lin\'eaires de $E.$ On d\'esigne par $BL(E),$ la
sous-alg\`ebre de $End_K(E),$ des op\'erateurs lin\'eaires
continus de $E.$ Muni de la norme op\'erationnelle \[
|||f|||=\sup_{x\in S(E)}||f(x)||, \] $BL(E)$ est une $K$-alg\`ebre
norm\'ee {\em ([Ber 73] p. 167-168)}. De plus, $(BL(E),|||.|||)$
est de Banach, si $(E,||.||)$ est un espace de Banach.$\Box$
\end{note}

\vspace{0.1cm} \hspace{0.3cm} Les r\'esultats pr\'eliminaires
suivants, dans le cadre des op\'erateurs lin\'eaires dans un
espace de Banach, nous seront utiles:

\vspace{0.2cm}
\begin{definition} Soit $(E,||.||)$ un $K$-espace vectoriel norm\'e et soit $T\in BL(E).$ On dit que $T$ est born\'e
inf\'erieurement s'il existe $m>0$ tel que \ $||T(x)||\geq m||x||$
pour tout $x\in E.$ On dit alors que $T$ est born\'e
inf\'erieurement par $m.\Box$
\end{definition}

\vspace{0.2cm}
\begin{lemma} Soit $(E,||.||)$ un espace vectoriel norm\'e et soit $T\in BL(E),$ $T$ bijectif. Alors $T$ est born\'e
inf\'erieurement si et seulement si $T^{-1}\in BL(E).\Box$
\end{lemma}

\vspace{0.2cm}
\begin{lemma} {\em [Rod 92$_1$].} Soit $(E,||.||)$ un espace de Banach et soit $T\in BL(E).$ Alors les deux propri\'et\'es
suivantes sont \'equivalentes: \begin{enumerate} \item $T$ est
born\'e inf\'erieurement. \item $T$ est injectif d'image ferm\'ee.
\end{enumerate}
\end{lemma}

\vspace{0.1cm} {\bf Preuve.} 1) $\Rightarrow$ 2). $T$ est
\'evidemment injectif, et si $(y_n)_{n}=(T(x_n))_n$ est une suite
d'\'el\'ements de $T(E)$ convergente vers un \'el\'ement $y,$
alors $(x_n)_n$ est une suite de Cauchy de $E,$ car $T$ est
born\'e inf\'erieurement, convergente vers un \'el\'ement $x\in
E,$ et on a \ $y=T(x)\in T(E).$

\vspace{0.2cm} 2) $\Rightarrow$ 1). Si $m=\inf_{x\in
S(E)}||T(x)||$ est nul, alors il existe une suite $(x_n)_n$
d'\'el\'ements de $S(E)$ telle que $(T(x_n))_n$ soit convergente
vers $0.$ Comme $T\Big( S(E)\Big)$ est un ferm\'e, on a \
$0=\lim_nT(x_n)\in T\Big( S(E)\Big),$ ce qui est absurde car $T$
est injectif. Donc $m>0$ i.e. $T$ est born\'e
inf\'erieurement.$\Box$

\vspace{0.2cm}
\begin{lemma} Soit $(E,||.||)$ un espace de Banach et soit $T\in BL(E),$ born\'e inf\'erieurement par $m$ et non
surjectif. Alors la boule ouverte de $BL(E)$ centr\'ee en $T$ et
de rayon $m$ est constitu\'ee uniquement d'op\'erateurs born\'es
inf\'erieurement et non surjectifs. En particulier, \[ \{f\in
BL(E): f \mbox{ born\'e inf\'erieurement et non surjectif}\} \]
est un ouvert de $BL(E)$ {\em ([Rod 92$_1$] Lemma {\bf 1}).}$\Box$
\end{lemma}

\vspace{0.2cm}
\begin{proposition} Soient $E$ un espace vectoriel norm\'e, $\hat{E}$ son compl\'et\'e et \ $\varphi:BL(E)\rightarrow
BL(\hat{E})$ \ $f\mapsto\overline{f}$ \ (le prolongement par
continuit\'e de $f$). Alors: \begin{enumerate} \item $\varphi$ est
un homomorphisme isom\'etrique d'alg\`ebres. \item Si $T\in
BL(E),$ alors $\overline{T}$ est born\'e inf\'erieurement si et
seulement si $T$ est born\'e inf\'erieurement. Dans ces
conditions, $\overline{T}$ est inversible dans $BL(\hat{E})$ si et
seulement si $T$ a une image dense dans $E.$ \end{enumerate}
\end{proposition}

\vspace{0.1cm} {\bf Preuve.} {\bf 2.} La premi\`ere partie est
cons\'equence de la continuit\'e de $T$ et de la densit\'e de $E$
dans $\hat{E}.$ Si $\overline{T}$ est bijectif, pour tout $x\in E$
il existe $a\in\hat{E},$ limite d'une suite d'\'el\'ements $a_n\in
E,$ tel que \ $\overline{T}(a)=\lim_nT(a_n)=x.$ Donc $T(E)$ est
dense dans $E.$ R\'eciproquement, si $T$ est born\'e
inf\'erieurement et $T$ a une image dense dans $E,$ alors
$\overline{T}$ est born\'e inf\'erieurement et a une image dense
dans $\hat{E},$ i.e. $\overline{T}$ est injectif et surjectif
(Lemme {\bf 2.10}).$\Box$

\vspace{0.2cm}
\begin{remark} Soit $A$ une alg\`ebre norm\'ee et soit $a\in A,$ alors le prolongement par continuit\'e
$\overline{L}_a,$ de $L_a$ \`a $\hat{A},$ n'est autre que
l'op\'erteur de multiplication \`a gauche par $a$ dans
$\hat{A}.\Box$ \end{remark}

\vspace{0.2cm}
\begin{note} {\em (Alg\`ebres norm\'ees de Jordan non commutatives).} Dans {\em [Kai 77]} Kaidi a d\'efini le spectre
d'un \'el\'ement $x,$ d'une $\cit$-alg\`ebre $A$ de Jordan non
commutative unitaire, comme dans le cas associatif: \
$sp_A(x)=\{\lambda\in\cit: x-\lambda 1\notin inv(A)\}.$ Il est
clair que \ $sp_A(x)=sp_{C(x)}(x)$ \ o\`u \ $\cit(x)$ est la
sous-alg\`ebre (associative et commutative) pleine engendr\'ee par
$x$ {\em (Note {\bf 1.42})}. Ainsi, le spectre d'un \'el\'ement
d'une $\cit$-alg\`ebre norm\'ee unitaire de Jordan non
commutative, est non vide.$\Box$
\end{note}

\vspace{0.2cm} \hspace{0.3cm} On en d\'eduit imm\'ediatement la
g\'en\'eralisation suivante du Th\'eor\`eme de Gelfand-Mazur:

\vspace{0.2cm}
\begin{theorem} {\em (Kaidi).} Soit $A$ une \ $\cit$-alg\`ebre de Jordan non commutative norm\'ee unitaire de division.
Alors $A\simeq$ $\cit$ {\em ([Kai 77] p. 92)}.$\Box$
\end{theorem}

\vspace{0.2cm} \hspace{0.3cm} Nous obtenons, \`a l'aide de ce
r\'esultat, le suivant:

\vspace{0.2cm}
\begin{corollary} {\em (Rochdi).} Soit $A$ une $\cit$-alg\`ebre de Jordan non commutative norm\'ee de division
lin\'eaire \`a gauche. Alors $A\simeq$ $\cit.$
\end{corollary}

\vspace{0.1cm} {\bf Preuve.} Le Corollaire {\bf 1.30} et la
Proposition {\bf 1.38} montrent que $A$ est unitaire et de
division.$\Box$

\vspace{0.2cm}
\begin{remark} Il est int\'eressant de savoir si une \ $\cit$-alg\`ebre commutative, \`a puissances
associatives norm\'ee de division lin\'eaire, est isomorphe \`a \
$\cit.\Box$
\end{remark}

\vspace{0.2cm} \hspace{0.3cm} Nous nous limitons maintenant \`a
exposer deux r\'esultats utiles, dans le cadre des alg\`ebres
norm\'ees non associatives:

\vspace{0.2cm}
\begin{proposition} {\em [Kai 77].} Soit $A$ une $K$-alg\`ebre norm\'ee compl\`ete, alors
l'ensemble \ $L-inv_g(A)$ est un ouvert de $A.$
\end{proposition}

\vspace{0.1cm} {\bf Preuve.} L'application \ $L:A\rightarrow
BL(A)$ \ $x\mapsto L_x$ \ est continue, car $|||L_x|||\leq ||x||,$
et pour tout $x\in A,$ on a

\begin{eqnarray*} x\in L-inv_g(A) &\Leftrightarrow& L_x\in inv\Big( End_K(A)\Big) \\
&\Leftrightarrow& L_x\in inv\Big( BL(A)\Big) \ \mbox{ d'apr\`es le
Th\'eor\`eme d'isomorphisme de Banach } \end{eqnarray*}

Donc \ $L-inv_g(A)=L^{-1}\Big( inv(BL(A))\Big)$ est un ouvert de
$A$ (Proposition {\bf 2.2}).$\Box$

\vspace{0.2cm}
\begin{proposition} {\em (Rochdi).} Soit $A$ une $\cit$-alg\`ebre norm\'ee, alors
l'application \ $\varphi:A\rightarrow\{\mbox{compacts de} \
$\cit$\}$ \ $x\mapsto sp_{BL(\hat{A})}(\overline{L}_x)$ \ est
semi-continue sup\'erieurement.
\end{proposition}

\vspace{0.1cm} {\bf Preuve.} Les deux applications \ $A\rightarrow
BL(A)$ \ $x\mapsto L_x$ \ et \ $BL(A)\rightarrow BL(\hat{A})$ \
$f\mapsto\overline{f}$ \ sont continues. De plus, l'application \
$BL(\hat{A})\rightarrow\{\mbox{compacts de} \ $\cit$\}$ \
$u\mapsto sp_{BL(\hat{A})}(u)$ \ est semi-continue
sup\'erieurement (Note {\bf 2.4}). L'application $\varphi,$
compos\'ee de ces trois derni\`eres, est semi-continue
sup\'erieurement.$\Box$

\vspace{1cm}
\subsection{Diviseurs topologiques de z\'ero dans une alg\`ebre norm\'ee}

\vspace{0.3cm} \begin{DTLZ} Soit $A$ une $K$-alg\`ebre norm\'ee.
\begin{enumerate}
\item Un \'el\'ement $a\in A$ est dit diviseur topologique
lin\'eaire de z\'ero \`a gauche (d.t.l.z.g.) (resp. diviseur
topologique lin\'eaire de z\'ero \`a droite (d.t.l.z.d.)) s'il
existe une suite $(x_n)_n,$ d'\'el\'ements de la sph\`ere unit\'e
de $A,$ telle que \ $ax_n\rightarrow 0$ \ (resp. $x_na\rightarrow
0$). L'\'el\'ement $a$ est dit diviseur topologique lin\'eaire de
z\'ero (d.t.l.z.) s'il est d.t.l.z.g. ou d.t.l.z.d. Il est clair
qu'un diviseur de z\'ero est un d.t.l.z. \item Si $A$ est
associative (m\^eme alternative), le mot "lin\'eaire" est
supprim\'e dans les d\'efinitions. \item Une alg\`ebre absolument
valu\'ee est sans diviseurs de z\'ero et ne contient aucun
d.t.l.z. non nul. \item On note $D_g(A)$ l'ensemble des d.t.l.z.g.
de $A.$
\end{enumerate}
\end{DTLZ}

\vspace{0.3cm} \begin{proposition} Soit $A$ une $\rit$-alg\`ebre
norm\'ee compl\`ete, alors les deux propri\'et\'es suivantes sont
\'equivalentes: \begin{enumerate} \item $A$ est de division
lin\'eaire \`a gauche. \item $A$ est sans d.t.l.z.g. non nuls et
$L-inv_g(A)\neq\emptyset.$ \end{enumerate} \end{proposition}

\vspace{0.1cm} {\bf Preuve.} 1) $\Rightarrow$ 2) est cons\'equence
du Th\'eor\`eme d'isomorphisme de Banach et 2) $\Rightarrow$ 1)
est \'etablie dans [Kai 77].$\Box$

\vspace{0.4cm} \hspace{0.3cm} Le r\'esultat suivant nous sera
utile:

\vspace{0.3cm} \begin{theorem} Soit $E$ un espace vectoriel
norm\'e et soit $T\in BL(E),$ alors propri\'et\'es suivantes sont
\'equivalentes:
\begin{enumerate} \item $T$ est un d.t.z.g. de $BL(E).$
\item Il existe une suite $(x_n)_n$ d'\'el\'elements de $S(E)$
telle que \ $T(x_n)\rightarrow 0.$ \item $T$ n'est pas born\'e
inf\'erieurement.
\end{enumerate}
\end{theorem}

\vspace{0.1cm} {\bf Preuve.} ([Ber 73] p. 241).$\Box$

\vspace{0.3cm} \begin{remark} Il est clair, en vertu du
Th\'eor\`eme {\bf 2.22}, qu'un \'el\'ement $a$ d'une alg\`ebre
norm\'ee $A$ est un d.t.l.z.g. si et seulement si $L_a$ est un
d.t.z.g. de $BL(A).\Box$
\end{remark}

\vspace{0.4cm} \hspace{0.3cm} Un d.t.z. d'une alg\`ebre
associative norm\'ee unitaire n'est pas inversible [BD 73], mais
il peut l'\^etre dans une extension de l'alg\`ebre:

\vspace{0.3cm} \begin{example} Soit $H$ un $K$-espace de Hilbert
de dimension hilbertienne d\'enombrable et soit ${\cal
B}=(e_n)_{n\geq 1}$ une base orthonorm\'ee totale de $H.$ On
consid\`ere le sous-espace vectoriel $E$ de $H$ engendr\'e
(alg\'ebriquement) par ${\cal B}$ et on d\'esigne par $f$
l'op\'erateur lin\'eaire de $E$ d\'efini par \
$f(e_n)=\frac{1}{n}e_n$ pour tout $n\geq 1.$ Alors $f\in BL(E)$ et
$f$ est inversible dans $End_K(E).$ De plus, $f$ est un d.t.z.g.
de $BL(E)$ car $f(e_n)\rightarrow 0$ i.e. $f^{-1}$ est non
continu.$\Box$
\end{example}

\vspace{0.3cm} \hspace{0.3cm} Un des r\'esultats classiques de la
th\'eorie des alg\`ebre r\'eelles de Banach, est le suivant:

\vspace{0.2cm} \begin{GMK theorem} Les uniques \ $\rit$-alg\`ebres
associatives norm\'ees de division sont isomorphes \`a \ $\rit,$
$\cit$ ou $\hit$ {\em ([BD 73], [Ri 60])}.$\Box$
\end{GMK theorem}

\vspace{0.3cm} \hspace{0.3cm} Ce th\'eor\`eme est une
cons\'equence du r\'esultat suivant d\^u \`a Kaplansky:

\vspace{0.2cm} \begin{theorem} {\em (Kaplansky).} Les uniques \
$\rit$-alg\`ebres associatives norm\'ees sans d.t.z. non nuls sont
\`a isomorphisme pr\`es \ $\rit,$ $\cit$ ou $\hit$ {\em [Kap
49]}.$\Box$
\end{theorem}

\vspace{0.4cm} \hspace{0.3cm} Dans une alg\`ebre $A$ (non
n\'ecessairement associative) un centralisateur d\'efini
partiellement sur $A$ d\'esigne un op\'erateur lin\'eaire $f,$
d\'efini dans un id\'eal non nul de $A,$ not\'e $dom(f),$ \`a
valeur dans $A,$ tel que \ $f(xy)=f(x)y$ \ et \ $f(yx)=yf(x)$ \
pour tous $x\in dom(f)$ et $y\in A.$ Si $A$ est premi\`ere (i.e.
le produit de deux id\'eaux non nuls de $A,$ est non nul), la
relation \ $\sim$ \ d\'efinie sur l'ensemble de tous les
centralisateurs d\'efinis partiellement sur $A$ par

\[ f\sim g \ \mbox{ si et seulement si } f \mbox{ et } g \ \mbox{ coincident sur } \ dom(f)\cap dom(g), \]

est d'\'equivalence. La somme et la compos\'ee de deux
centralisateurs d\'efinis partiellement sur $A,$ en tant
qu'op\'erateurs d\'efinis partiellement, sont \'egalement des
centralisateurs d\'efinis partiellement sur $A.$ Ces op\'erations
sont compatibles avec \ $\sim$ \ et l'anneau quotient \ ${\cal
C}(A)$ ainsi obtenu est appel\'e "le centro{\" i}de \'etendu
(extended centroid) de $A$". R\'ecemment, Cabrera et Rodriguez
[CR] ont donn\'e une d\'emonstration simple du Th\'eor\`eme de
Kaplansky. Vu son originalit\'e, nous en exposons ici les grandes
lignes:

\vspace{0.2cm} \begin{theorem} {\em (Cabrera-Rodriguez).} Les
seules $\rit$-alg\`ebres alternatives norm\'ees sans d.t.z. non
nuls sont \`a isomorphisme pr\`es \ $\rit,$ $\cit,$ $\hit$ ou
$\oit.$
\end{theorem}

\vspace{0.1cm} {\bf Preuve.} La premi\`ere \'etape de la
d\'emonstration consiste \`a r\'eduire, moyennant le Th\'eor\`eme
d'Albert-Frobenius {\bf 1.61}, le probl\`eme au cas associatif et
commutatif. Ensuite on montre qu'une \ $\rit$-alg\`ebre $A$
associative, commutative et sans diviseurs de z\'ero, est
identifi\'ee \`a une sous-alg\`ebre de son centro{\" i}de \'etendu
${\cal C}(A).$ Ce centro{\" i}de \'etendu est une extension du
corps de base \ $\rit$ ([ErMO 75] Theorem {\bf 2.1}). enfin, \`a
l'aide de l'absence des d.t.z. non nuls, on montre que le
centro{\" i}de \'etendu de $A$ est muni d'une norme d'alg\`ebre,
et le Th\'eor\`eme de Gelfand-Mazur-Kaplansky ach\`eve la
d\'emonstration.$\Box$

\vspace{0.2cm} \begin{corollary} Soit $A$ une \ $\rit$-alg\`ebre
\`a puissances associatives norm\'ee sans d.t.l.z.g. non nuls.
Alors $A$ est quadratique et de division.
\end{corollary}

\vspace{0.1cm} {\bf Preuve.} La sous-alg\`ebre engendr\'ee par un
\'el\'ement quelconque de $A$ est isomorphe \`a \ $\rit$ ou
$\cit,$ i.e. $A$ poss\`ede un idempotent non nul. Comme $A$ est
sans diviseurs de z\'ero, elle est unitaire (Lemme {\bf 1.8}).
donc $A$ est quadratique et de division.$\Box$

\vspace{0.4cm} \hspace{0.3cm} Nous introduisons maintenant une
notions appropri\'ee de d.t.z. pour les alg\`ebres norm\'ees de
Jordan non commutatives:

\vspace{0.2cm} \begin{proposition} Soit $A$ une alg\`ebre
associative norm\'ee et soit $a\in A.$ Alors les affirmations
suivantes sont \'equivalentes: \begin{enumerate} \item $a$ est un
d.t.z. de $A.$ \item Il existe une suite $(x_n)_{n\geq 0}$
d'\'el\'ements de $S(A)$ telle que \ $ax_na\rightarrow 0.$
\end{enumerate}
\end{proposition}

\vspace{0.1cm} {\bf Preuve.} [Kai 77].$\Box$

\vspace{0.2cm} \begin{remark} $a$ est un d.t.z. de $A$ si et
seulement si il existe une suite $(x_n)_{n\geq 0}$ d'\'el\'ements
de $S(A)$ telle que \ $U_a(x_n)\rightarrow 0.\Box$
\end{remark}

\vspace{0.3cm} \hspace{0.3cm} Ceci conduit \`a la d\'efinition
suivante:

\vspace{0.2cm} \begin{J-DTZ} {\em [Kai 77].} Soit $A$ une
$K$-alg\`ebre norm\'ee de Jordan non commutative. Un \'el\'ement
$a\in A$ est dit $J$-diviseur topologique de z\'ero ($J$-d.t.z.)
s'il existe une suite $(x_n)_{n\geq 0}$ d'\'el\'ements de $S(A)$
telle que \[ U_a(x_n)=a(ax_n+x_na)-a^2x_n\rightarrow 0 \] i.e.
$U_a$ est un d.t.z.g. de $BL(A).$ Si $A$ est une alg\`ebre
quadratique, un \'el\'ement de $A$ est dit $J$-d.t.z. s'il est un
$J$-d.t.z. de l'alg\`ebre de Jordan $A^+.\Box$
\end{J-DTZ}

\vspace{0.5cm} \hspace{0.3cm} Les r\'esultats classiques
subsistent pour cette nouvelle notion de d.t.z. ([Kai 77], [KS]).
Par exemple, le r\'esultat suivant:

\vspace{0.2cm} \begin{proposition} {\em (Kaidi).} Soit $A$ une
alg\`ebre de Jordan non commutative norm\'ee unitaire. Alors un
\'el\'ement inversible de $A$ ne peut \^etre un $J$-d.t.z.
\end{proposition}

\vspace{0.1cm} {\bf Preuve.} ([Kai 77] p. 86).$\Box$

\vspace{0.2cm} \begin{remark} Si $A$ est une alg\`ebre norm\'ee de
Jordan non commutative, alors un d.t.l.z. de $A$ n'est pas
n\'ecessairement un $J$-d.t.z. de $A.$ En effet, dans l'alg\`ebre
r\'eelle de Jordan \ $\hit^+,$ l'\'el\'ement $i$ est inversible
mais non l.i. inversible. On en d\'eduit que $i$ est un d.t.l.z.
mais pas un $J$-d.t.z. {\em (Proposition {\bf 2.32})}.$\Box$
\end{remark}

\vspace{0.3cm} \hspace{0.3cm} Cette \'etude a abouti \`a une
extension du Th\'eor\`eme de Kaplansky aux alg\`ebres de Jordan
non commutatives:

\vspace{0.2cm} \begin{theorem} {\em ([Kai 77], [KS]).} Soit $A$
une $K$-alg\`ebre de Jordan non commutative norm\'ee sans
$J$-d.t.z. non nuls. Alors $A\simeq$ $\cit$ si $K=\cit$ ou $A$ est
quadratique de division si $K=\rit.\Box$
\end{theorem}

\vspace{0.3cm} \hspace{0.3cm} Nous avons \'egalement les
th\'eor\`emes de structure suivants:

\vspace{0.2cm} \begin{theorem} {\em [Kai 77].} Soit $A$ une
$\rit$-alg\`ebre quadratique, munie d'une norme d'espace
vectoriel. Alors les affirmations suivantes sont \'equivalentes:
\begin{enumerate} \item $A$ est norm\'ee de division. \item $A$ est norm\'ee sans $J$-d.t.z. non nuls.
\item $A$ est associ\'ee \`a une forme quadratique continue
d\'efinie n\'egative, et \`a une alg\`ebre anti-commutative dont
le produit est continu.$\Box$
\end{enumerate}
\end{theorem}

\vspace{0.2cm} \begin{theorem} {\em [Kai 77].} Soit $A$ une
$\rit$-alg\`ebre de Jordan non commutative, munie d'une norme
d'espace vectoriel. Alors les affirmations suivantes sont
\'equivalentes:
\begin{enumerate} \item $A$ est norm\'ee de division. \item $A$ est norm\'ee sans $J$-d.t.z. non nuls.
\item $A$ est quadratique associ\'ee \`a une forme trace continue
d\'efinie n\'egative, et \`a une alg\`ebre anti-commutative dont
le produit est continu.$\Box$
\end{enumerate}
\end{theorem}

\vspace{0.2cm} \begin{theorem} {\em [Kai 77].} Soit $A$ une
$\rit$-alg\`ebre de Jordan. Alors les affirmations suivantes sont
\'equivalentes:
\begin{enumerate} \item $A$ est norm\'ee de division. \item $A$ est norm\'ee sans $J$-d.t.z. non nuls.
\item $A$ est l'alg\`ebre de Jordan associ\'ee \`a une forme
quadratique d\'efinie n\'egative.$\Box$
\end{enumerate}
\end{theorem}

\vspace{1cm}
\subsection{Alg\`ebre complexes norm\'ees de division lin\'eaire}

\vspace{0.5cm} \hspace{0.3cm} Le Th\'eor\`eme de Gelfand-Mazur a
connu une extension (non associative) en 1977:

\vspace{0.2cm} \begin{theorem} {\em (Kaidi).} Soit $A$ une
$\cit$-alg\`ebre norm\'ee compl\`ete de division lin\'eaire \`a
gauche. Alors $A\simeq$ $\cit.$
\end{theorem}

\vspace{0.1cm} {\bf Preuve.} ([Kai 77] p. 80).$\Box$

\vspace{0.4cm} \hspace{0.3cm} Nous n'avons pu \'eliminer de ce
r\'esultat, l'hypoth\`ese de la compl\'etude. Contrairement au cas
associatif, il nous semble que le probl\`eme de savoir si une \
$\cit$-alg\`ebre norm\'ee de d.l. est isomorphe \`a \ $\cit$ est
tr\`es difficile. Nous exposons ici les r\'esultats que nous avons
pu \'etablir, assurant la validit\'e du Th\'eor\`eme {\bf 2.38}
dans des situations appremment plus g\'en\'erales que celles du
cas complet. Pour cela, nous avons besoin des r\'esultats
pr\'eliminaires suivants:

\vspace{0.2cm} \begin{lemma} {\em [Ri 60].} Soit $(A,||.||)$ une
alg\`ebre norm\'ee, alors l'application

\[ \Phi:A\rightarrow\rit\hspace{0.4cm} x\mapsto\inf_{y\in S(A)}||xy||=\inf_{y\neq 0}\frac{||xy||}{||y||} \]

est continue.$\Box$
\end{lemma}

\vspace{0.3cm} \hspace{0.3cm} Nous avons obtenu le r\'esultat
suivant, introuvable apparemment dans la litt\'erature:

\vspace{0.2cm} \begin{proposition} Soit $(A,||.||)$ une
$K$-alg\`ebre norm\'ee, alors: \begin{enumerate} \item $D_g(A)$
est un ferm\'e de $A,$ et on a \ $D_g(A)=A\cap D_g(\hat{A}).$
\item $L-inv_g(\hat{A})$ contient \ $(A\setminus D_g(A))\cap
L-inv_g(A).$ Si, de plus, $A$ est de division lin\'eaire \`a
gauche, on a \ $A\setminus D_g(A)=A\cap L-inv_g(\hat{A}).$
\end{enumerate}
\end{proposition}

\vspace{0.1cm} {\bf Preuve.} \begin{enumerate} \item l'application
$\Phi$ dans le Lemme {\bf 2.38} est continue. Donc \[
D_g(A)=\{x\in A:\Phi(x)=0\}=\Phi^{-1}\{0\} \] est un ferm\'e de
$A.$ De plus, la densit\'e de $S(A)$ dans $S(\hat{A})$ montre que
\ $D_g(A)=A\cap D_g(\hat{A}).$ \item Soit $a\in(A\setminus
D_g(A))\cap L-inv_g(A),$ alors $L_a\in BL(A)$ est born\'e in
f\'erieurement et a une image dense dans $A.$ Son prolongement
$\overline{L}_a$ est inversible dans $BL(\hat{A})$ (Proposition
{\bf 2.11}) i.e. $a\in L-inv_g(\hat{A}).$ Si $a\in A\cap
L-inv_g(\hat{A}),$ alors $\overline{L}_a$ est inversible dans
$End_K(\hat{A})$ et $\overline{L}_a^{-1}\in BL(\hat{A}),$ en vertu
du le Th\'eor\`eme d'isomorphisme de Banach. Donc $\overline{L}_a$
est born\'e inf\'erieurement (Lemme {\bf 2.8}). Ainsi $a\in
A\setminus D_g(A)=A\cap L-inv_g(\hat{A}).\Box$
\end{enumerate}

\vspace{0.2cm} \begin{isotopy} Deux alg\`ebres $A$ et $B,$ sur un
corps commutatif $K_0,$ sont dites isotopes s'il existe trois
bijections lin\'eaires $u,v,w:A\rightarrow B$ telles que \
$w(xy)=u(x)v(y)$ pour tous $x,y\in A.$ Les deux alg\`ebres
r\'eelles $\stackrel{*}{\cit}$ et $\cit$ sont isotopes.$\Box$
\end{isotopy}

\vspace{0.2cm} \begin{lemma} {\em [Kai 77].} Toute \
$\cit$-alg\`ebre isotope \`a \ $\cit$ est isomorphe \`a \
$\cit.\Box$
\end{lemma}

\vspace{0.2cm} \begin{note} {\em [Kai 77].} Soit $A$ une alg\`ebre
norm\'ee dans laquelle il existe un \'el\'ement $a$ l.i.g. pour
lequel $L_a^{-1}$ est continu. Alors l'alg\`ebre $A^a,$ isotope
\`a $A,$ ayant pour espace vectoriel sous-jacent $A$ et pour
produit \ $x\odot y=L_a^{-1}(x)L_a^{-1}(y),$ est norm\'ee {\em
(}$||x\odot y||\leq |||L_a|||^2||x||$ $||y||${\em )} unitaire \`a
gauche {\em (}d'unit\'e \`a gauche $a^2${\em )}.$\Box$
\end{note}

\vspace{0.3cm} \hspace{0.3cm} Nous pouvons exposer, dans le reste
de ce paragraphe, les r\'esultats obtenus:

\vspace{0.2cm} \begin{lemma} Soit $A$ une $K$-alg\`ebre norm\'ee
telle que \ $(A\setminus D_g(A))\cap L-inv_g(A)\neq\emptyset.$
Alors pour tout $a\in (A\setminus D_g(A))\cap L-inv_g(A),$
$L_a^{-1}$ est continue.
\end{lemma}

\vspace{0.1cm} {\bf Preuve.} Cons\'equence du Lemme {\bf 2.8} et
du Th\'eor\`eme {\bf 2.21}.$\Box$

\vspace{0.2cm} \begin{remarks} Soit $A^a$ l'alg\`ebre introduite
dans la {\em Note {\bf 2.42}}. \begin{enumerate} \item La
bicontinuit\'e de $L_a$ montre, facilement, l'\'egalit\'e \
$D_g(A^a)=L_a(D_g(A)).$ \item Si $A$ est de d.l.g., il en est de
m\^eme pour $A^a.$ En effet, l'op\'erateur de multiplication \`a
gauche par $x\in A$ dans $A^a,$ s'\'ecrit \
$L_x^{\odot}=L_{L_a^{-1}(x)}\circ L_a^{-1}.\Box$ \end{enumerate}
\end{remarks}

\vspace{0.2cm} \begin{lemma} Soit $A$ une alg\`ebre norm\'ee de
division lin\'eaire \`a gauche. Alors \
$D_g(\hat{A})=\hat{A}\setminus L-inv_g(\hat{A}).$
\end{lemma}

\vspace{0.1cm} {\bf Preuve.} On note $\Omega=\{f\in BL(\hat{A}): f
\mbox{ born\'e inf\'erieurement et non surjectif}\}$ et $L$
l'application \ $\hat{A}\rightarrow BL(\hat{A})$ \ $x\mapsto L_x.$
Soit maintenant $a\in\hat{A},$ on a

\begin{eqnarray*} a\in\hat{A}\setminus D_g(\hat{A}) &\Leftrightarrow& L_a \mbox{ born\'e inf\'erieurement, et }
a\in\hat{A}\setminus L-inv_g(\hat{A}) \\
&\Leftrightarrow& L_a \mbox{ n'est pas bijectif. }
\end{eqnarray*} Ainsi

\[ a\in\Big( \hat{A}\setminus D_g(\hat{A})\Big)\cap\Big(
\hat{A}\setminus L-inv_g(\hat{A})\Big)=\hat{A}\setminus\Big(
D_g(\hat{A})\cup L-inv_g(\hat{A})\Big)\Leftrightarrow
L_a\in\Omega. \]

Donc \ $\hat{A}\setminus\Big( D_g(\hat{A})\cup
L-inv_g(\hat{A})\Big)=L^{-1}(\Omega)$ \ est un ouvert de $\hat{A}$
et, par cons\'equent, \ $D_g(\hat{A})\cup L-inv_g(\hat{A})$ \ est
un ferm\'e de $\hat{A}.$ D'autre part

\begin{eqnarray*} A &=& D_g(A)\cup(A\setminus D_g(A)) \\
&=& \Big( A\cap D_g(A)\Big)\cup\Big( A\cap L-inv_g(\hat{A})\Big)
\\
&\subseteq& D_g(\hat{A})\cup L-inv_g(\hat{A}). \end{eqnarray*}

Donc \ $D_g(\hat{A})\cup L-inv_g(\hat{A})=\hat{A}.$ De plus \
$D_g(\hat{A})\cap L-inv_g(\hat{A})=\emptyset.\Box$

\vspace{0.3cm} \begin{lemma} Soit $A$ une \ $\cit$-alg\`ebre
norm\'ee unitaire \`a gauche, de division lin\'eaire \`a gauche.
Alors \begin{enumerate} \item $A=\cit 1+D_g(A)$ et \item
$\hat{A}=\cit 1+\overline{D_g(A)}$ \ {\em (}la fermeture dans
$\hat{A}${\em )}.
\end{enumerate}
\end{lemma}

\vspace{0.1cm} {\bf Preuve.} . \begin{enumerate} \item Soit $x\in
A,$ alors
\begin{eqnarray*} \lambda\in sp_{BL(\hat{A})}(\overline{L}_x)
&\Leftrightarrow&
\overline{L}_x-\lambda 1\notin inv\Big( BL(\hat{A})\Big) \\
&\Leftrightarrow& \overline{L}_{x-\lambda 1}\notin
inv\Big( End_C(\hat{A})\Big) \ \mbox{ d'apr\`es le TIB } \\
&\Leftrightarrow& x-\lambda 1\notin
L-inv_g(\hat{A}) \\
&\Leftrightarrow& x-\lambda 1\in D_g(\hat{A}) \\
&\Leftrightarrow& x-\lambda 1\in D_g(A) \end{eqnarray*}

i.e. $A=\cit 1+D_g(A).$ \item Soit maintenant $x\in\hat{A},$ $x$
est limite d'une suite d'\'el\'ements $x_n$ de $A,$ pour tout
$n\in\nit,$ soit $\lambda_n\in
sp_{BL(\hat{A})}(\overline{L}_{x_n}).$ Il existe $u_n\in D_g(A)$
tel que $x_n=\lambda_n1+u_n.$ Si $V$ est un voisinage (born\'e) de
$sp_{BL(\hat{A})}(\overline{L}_{x_n}),$ il existe un rang $n_0$
tel que \ $\lambda_n\in
sp_{BL(\hat{A})}(\overline{L}_{x_n})\subset V$ pour tout $n\geq
n_0.$ La suite $(\lambda_n)_n$ est donc born\'ee et contient une
sous-suite $(\lambda_{\varphi(n)})_n$ convergente vers un nombre
complexe $\lambda.$ Donc $(u_{\varphi(n)})_n$ est convergente vers
un \'el\'ement $u\in\overline{D_g(A)}.$ Comme $x$ est limite des
\'el\'ements $x_{\varphi(n)}$ on a $x=\lambda 1+u\in\cit
1+\overline{D_g(A)}$ i.e. $\hat{A}=\cit 1+\overline{D_g(A)}.\Box$
\end{enumerate}

\vspace{0.3cm} \begin{lemma} Soit $A$ une alg\`ebre norm\'ee de
division lin\'eaire \`a gauche et telle que $D_g(A)$ soit une
partie compl\`ete. Alors $A$ est isotope \`a une alg\`ebre
norm\'ee unitaire \`a gauche de d.l.g. $A^a.$ De plus, l'ensemble
des d.t.l.z.g. de $A^a$ est une partie compl\`ete.
\end{lemma}

\vspace{0.1cm} {\bf Preuve.} $D_g(A)$ est une partie ferm\'ee de
$\hat{A}$ distincte de $\hat{A},$ donc $A\neq D_g(A).$ Soit alors
$a\in\Big( A\setminus D_g(A)\Big)\cap L-inv_g(A)=A\setminus
D_g(A)\neq\emptyset,$ le Lemme {\bf 2.43} montre que $L_a^{-1}$
est continu. En tenant compte des Remarques {\bf 2.44}, l'isotope
$A^a$ de $A$ introduite dans la Note {\bf 2.42} r\'epond \`a la
question.$\Box$

\vspace{0.3cm} \begin{remark} Soit $A$ une alg\`ebre norm\'ee de
division lin\'eaire. On ne sait, toujours, pas si la compl\'etion
$\hat{A}$ de $A$ contient un \'el\'ement lin\'eairement
inversible.$\Box$
\end{remark}

\vspace{0.3cm} \hspace{0.3cm} Nous avons obtenu auparavant [Roc
87], les deux r\'esultats suivants:

\vspace{0.3cm} \begin{theorem} Soit $A$ une $\cit$-alg\`ebre
norm\'ee de division lin\'eaire et sans d.t.l.z.g. non nuls. Alors
$A\simeq$ $\cit.$
\end{theorem}

\vspace{0.1cm} {\bf Preuve.} ([Roc 87] p. 34).$\Box$

\vspace{0.3cm} \begin{theorem} Soit $A$ une $\cit$-alg\`ebre
norm\'ee unitaire \`a gauche et sans d.t.l.z.g. non nuls. Alors
$A\simeq$ $\cit.$
\end{theorem}

\vspace{0.1cm} {\bf Preuve.} ([Roc 87] p. 34-35).$\Box$

\vspace{0.3cm} \hspace{0.3cm} Nous avons les r\'esultats nouveaux
suivants:

\vspace{0.3cm} \begin{theorem} {\em (Rochdi).} Soit $A$ une \
$\cit$-alg\`ebre norm\'ee sans d.t.l.z.g. non nuls et telle que \
$L-inv_g(A)\neq\emptyset.$ Alors $A\simeq$ $\cit.$
\end{theorem}

\vspace{0.1cm} {\bf Preuve.} Soit $a\in L-inv_g(A),$ alors le
Lemme {\bf 2.43} montre que $L_a^{-1}$ est continu et l'isotope
$A^a$ de $A$ introduite dans la Note {\bf 2.42} est norm\'ee,
unitaire \`a gauche. En outre, $A^a$ est sans d.t.l.z.g. (Remarque
{\bf 2.44 1)}). Le r\'esultat est alors cons\'equence du
Th\'eor\`eme {\bf 2.50} et du Lemme {\bf 2.41}.$\Box$

\vspace{0.3cm} \hspace{0.3cm} Le Th\'eor\`eme {\bf 2.51}
pr\'ec\'edent g\'en\'eralise apparemment les Th\'eor\`emes {\bf
2.49, 2.50}, et \'egalement le Th\'eor\`eme {\bf 2.37}. D'autre
part, dans le Th\'eor\`eme {\bf 2.51}, l'hypoth\`ese de
l'existence d'un \'el\'ement l.i.g. est n\'ecessaire car il existe
des exemples de $\cit$-alg\`ebres norm\'ees sans d.t.l.z.g. non
nuls (absolument valu\'ees compl\`etes) de dimension infinie ([Roc
87] p. 35), ([Rod 92$_1$] Remarks {\bf 3. i)}).

\vspace{0.3cm} \begin{corollary} Soit $A$ une $\cit$-alg\`ebre
absolument valu\'ee telle que \ $L-inv_g(A)\neq\emptyset.$ Alors
$A\simeq$ $\cit.\Box$
\end{corollary}

\vspace{0.3cm} \begin{corollary} Soit $A$ une $\cit$-alg\`ebre
norm\'ee \`a puissances associatives et sans d.t.l.z.g. Alors
$A\simeq$ $\cit.\Box$
\end{corollary}

\vspace{0.3cm} \hspace{0.3cm} Nous avons \'egalement le r\'esultat
suivant:

\vspace{0.3cm} \begin{theorem} {\em (Rochdi).} Soit $A$ une
$\cit$-alg\`ebre norm\'ee de d.l.g. et telle que \ $D_g(A)$ soit
une partie compl\`ete. Alors $A\simeq$ $\cit.$
\end{theorem}

\vspace{0.1cm} {\bf Preuve.} L'isotope $A^a$ de $A$ construite \`a
partir d'un \'el\'ement $a\in A\setminus D_g(A),$ est norm\'ee,
unitaire \`a gauche de d.l.g., et on a:

\begin{eqnarray*} A^a &=& \cit 1+D_g(A^a) \hspace{0.3cm} \mbox{ Lemme {\bf 2.46}} \\
&=& \cit 1+\overline{D_g(A^a)} \hspace{0.2cm} \mbox{ car }
D_g(A^a) \mbox{ est
compl\`ete } \\
&=& \hat{A} \hspace{0.3cm} \mbox{ Lemme {\bf 2.46}}
\end{eqnarray*}

i.e. $A^a$ est compl\`ete. Le r\'esultat est alors cons\'equence
du Lemme {\bf 2.41} et du Th\'eor\`eme {\bf 2.37}.$\Box$

\vspace{0.3cm} \begin{remark} Soit $A$ une $\cit$-alg\`ebre
norm\'ee de d.l.g. On aimerait savoir si $D_g(A)$ est une partie
compl\`ete.$\Box$
\end{remark}

\vspace{4cm}
\subsection{Alg\`ebre r\'eelles norm\'ees de division lin\'eaire}

\vspace{0.6cm} \hspace{0.3cm} Contrairement au cas complexe
(complet), la d\'etermination des alg\`ebres r\'eelles norm\'ees
compl\`etes de division lin\'eaire est un probl\`eme apparemment
plus compliqu\'e. En 1953, Wright [Wr 53] a conjectur\'e que les
$\rit$-alg\`ebres norm\'ees de division lin\'eaire, sont de
dimension finie. Cette conjecture est \'etablie pour les
alg\`ebres faiblement alternatives (...) et pour les alg\`ebres
quadratiques qui satisfont \`a la propri\'et\'e d'Osborn (..),
mais elle ne l'est pas encore pour les alg\`ebres de Jordan non
commutatives. Dans ce paragraphe, nous exposons les r\'esultats
\'etablis dans le cadre des $\rit$-alg\`ebres norm\'ees de
division lin\'eaire.

\vspace{0.3cm} \begin{proposition} Soit $A=\Big(
V,(.|.),\times\Big)$ une $\rit$-alg\`ebre cayleyenne de division.
Si $u_1,\dots,u_n$ sont des vecteurs non nuls tels que
$(u_i|u_j)=0$ si $i\neq j,$ alors

\[ V=\Big( \bigoplus_{1\leq i\leq n}\rit u_i\Big)\oplus\Big( \bigcap_{1\leq i\leq n}W(u_i)\Big) \]

avec $W(u)=\{x\in A:(x|u)=0\}.$
\end{proposition}

\vspace{0.1cm} {\bf Preuve.} ([Kai 77] p. 117).$\Box$

\vspace{0.3cm} \begin{definition} {\em [Kai 77].} Une
$K$-alg\`ebre $A,$ de Jordan non commutatove, est dite faiblement
alternative si elle satisfait \`a l'identit\'e \ $(x,x,[x,y])=0$
pour tous $x,y\in A.$ Une alg\`ebre alternative ou de Jordan est
\'evidemment faiblement alternative.$\Box$
\end{definition}

\vspace{0.3cm} \begin{proposition} {\em [Kai 77].} Soit $A=\Big(
V,(.|.),\times\Big)$ une $\rit$-alg\`ebre quadratique, faiblement
alternative de division et sans diviseurs de z\'ero. Alors $A$ est
alternative.
\end{proposition}

\vspace{0.1cm} {\bf Preuve.} Si $A$ est de dimension $>2,$ on a
$A=\rit\oplus V$ et $\dim_R(V)>1.$ Soient $e_1$ un vecteur norm\'e
($e_1^2=-1$) de $V$ et $e_2$ un vecteur norm\'e de $W(e_1),$ alors
$e_1, e_2, e_1\times e_2:=e_3$ sont lin\'eairement ind\'ependnats
et on a \ $A=\rit 1\oplus\rit e_1\oplus \rit e_2\oplus\rit
e_3\oplus L$ o\`u $L=W(e_1)\cap W(e_2)\cap W(e_3).$ Comme $A$ est
faiblement alternative et $e_1e_2=e_1\times
e_2=\frac{1}{2}[e_1,e_2],$ on a

\begin{eqnarray*} 0 &=& (e_1,e_1,[e_1,e_2]) \\
&=& -e_1e_2-e_1\Big( e_1(e_1e_2)\Big) \\
&=& -e_1\Big( e_2+e_1(e_1e_2)\Big) \end{eqnarray*}

i.e. $e_2=-e_1(e_1e_2).$ De m\^eme $e_2(e_2e_1)=-e_1.$ Le
caract\`ere quadratique de $A$ et la flexibilit\'e assurent que
$\rit 1\oplus\rit e_1\oplus \rit e_2\oplus\rit e_3$ est une
sous-alg\`ebre de $A,$ isomorphe \`a $\hit.$ Ainsi, deux
\'el\'ements quelconques de $A$ engendrent une sous-alg\`ebre
associative, et le Th\'eor\`eme d'Artin ach\`eve la
d\'emonstration.$\Box$

\vspace{0.3cm} \begin{theorem} {\em [Kai 77].} Soit $A$ une
$\rit$-alg\`ebre faiblement alternative. Alors les trois
affirmations suivantes sont \'equvalentes: \begin{enumerate} \item
$A$ est norm\'ee de d.l.g. \item $A$ est norm\'ee sans d.t.l.z.g.
non nuls. \item $A$ est \`a isomorphisme pr\`es \ $\rit,$ $\cit,$
$\hit$ ou $\oit.$
\end{enumerate}
\end{theorem}

\vspace{0.1cm} {\bf Preuve.} Chacune des affirmations {\bf 1.} et
{\bf 2.} entra\^ine que $A$ est quadratique (de division) et sans
diviseurs de z\'ero. Donc $A$ est alternative (Proposition {\bf
2.58}) et le Th\'eor\`eme {\bf 1.61} de Albert montre que $A$ est
isomorphe \`a \ $\rit,$ $\cit,$ $\hit$ ou $\oit.\Box$

\vspace{0.3cm} \begin{corollary} Les uniques $\rit$-alg\`ebres de
Jordan, de division lin\'eaire, sont \`a isomorphisme pr\`es \
$\rit$ ou $\cit.\Box$
\end{corollary}

\vspace{0.3cm} \hspace{0.3cm} Des \'etudes sur les alg\`ebres
r\'eelles absolument valu\'ees, d\^ues \`a Albert [A 47, 49],
Wright [Wr 53] et Urbanik-Wright [UW 60], ont abouti en 1960 aux
r\'esultats suivants:

\vspace{0.3cm} \begin{theorem} Les uniques \ $\rit$-alg\`ebres
absolument valu\'ee unitaires, sont \`a isomorphisme pr\`es \
$\rit,$ $\cit,$ $\hit$ ou $\oit.$
\end{theorem}

\vspace{0.1cm} {\bf Preuve.} ([UW 60] p. 863).$\Box$

\vspace{0.3cm} \begin{theorem} Les uniques \ $\rit$-alg\`ebres
absolument valu\'ee commutatives, sont \`a isomorphisme pr\`es \
$\rit,$ $\cit$ ou $\stackrel{*}{\cit}.$
\end{theorem}

\vspace{0.1cm} {\bf Preuve.} ([UW 60] p. 865).$\Box$

\vspace{0.3cm} \hspace{0.3cm} Le Th\'eor\`eme {\bf 2.61} s'\'etend
facilement aux alg\`ebres absolument valu\'ees contenant un
\'el\'ement lin\'eairement inversible ([ElM 81] p. 245-246):

\vspace{0.3cm} \begin{theorem} Les uniques \ $\rit$-alg\`ebres
absolument valu\'ee qui contiennent un \'el\'ement l.i., sont de
dimension finie $1,2,4$ ou $8$ et isotopes \`a \ $\rit,$ $\cit,$
$\hit$ ou $\oit.$
\end{theorem}

\vspace{0.1cm} {\bf Preuve.} Soit $A$ une telle alg\`ebre et soit
$a$ un \'el\'ement l.i., qu'on peut supposer norm\'e. Alors
l'alg\`ebre $A^{\odot},$ isotope de $A,$ ayant pour produit \
$x\odot y=R_a^{-1}(x)L_a^{-1}(y)$ est unitaire, d'unit\'e $e=a^2,$
et on v\'erifie facilemtn qu'elle est absolument valu\'ee. Le
Th\'eor\`eme {\bf 2.61} montre que $A^{\odot}$ est de dimension
finie $1,2,4$ ou $8$ et isomorphe \`a \ $\rit,$ $\cit,$ $\hit$ ou
$\oit.\Box$

\vspace{0.3cm} \hspace{0.3cm} Ce dernier r\'esultat persiste si
l'on \'echange l'hypoth\`ese de l'existence d'un \'el\'ement l.i.
par celle de la flexibilit\'e ([ElM 81]. 244).

\vspace{0.3cm} \hspace{0.3cm} Ult\'erieurement, El-Mallah a
donn\'e une classification pour les \ $\rit$-alg\`ebres absolument
valu\'ees de dimension finie, qui satisfont \`a l'identit\'e
$(x,x,x)=0$ et montre que de telles alg\`ebres sont flexibles [Elm
87].

\vspace{0.3cm} \hspace{0.3cm} Wright a conjectur\'e, auparavant,
que les $\rit$-alg\`ebres norm\'ees de d.l. sont de dimension
finie [Wr 53]. R\'ecemment, Cuenca [Cu 92] a donn\'e des exemples
d'alg\`ebres norm\'ees de d.l.g. (absolument valu\'ees compl\`etes
et unitaires \`a gauche) de dimension infinie. De son c\^ot\'e,
Rodriguez ([Rod 92] Theorem {\bf 2} p. 941) les a compl\`etement
d\'ecrites.

\vspace{0.3cm} \hspace{0.3cm} Un des r\'esultats fondamentaux pour
les alg\`ebres r\'eelles de division lin\'eaire, de dimension
finie, \'etabli en 1958, est le suivant:

\vspace{0.3cm} \begin{theorem} {\em (Hopf-Kervaire-Milnor-Bott).}
Si l'espace vectoriel r\'eel $\rit^n$ poss\`ede un produit
bilin\'eaire sans diviseurs de z\'ero, alors $n=1,2,4$ ou $8.\Box$
\end{theorem}

\vspace{0.3cm} \hspace{0.3cm} Il est facile de voir que $\rit$ est
l'unique, \`a isomorphisme pr\`es, $\rit$-alg\`ebre de division
lin\'eaire de dimension $1.$ D'autre part, les $\rit$-alg\`ebres
de divison lin\'eaire de dimension $2$ ont \'et\'e "classifi\'ees"
[AK 83]. Dans le m\^eme th\`eme, des \'etudes ont \'et\'e faites
en dimension $4$ et $8$ [Wr 53], [Os 62], [Cz 76], [Kai 77], [Ok
80], [OO 81$_{1,2}$], [BBO 82], [AHK 86, 87, 92], [Roc 87], [KR
92], et ont donn\'e des r\'esultats partiels. Une $\rit$-alg\`ebre
unitaire de division lin\'eaire de dimension finie $\geq 2$
contient une sous-alg\`ebre isomorphe \`a \ $\cit$ [Y 81], [Petr
87], cependant la d\'etermination des $\rit$-alg\`ebres unitaires
de division lin\'eaire de dimension $4$ est encore un probl\`eme
ouvert.

\vspace{0.3cm} \begin{theorem} {\em [Kai 77].} Les
$\rit$-alg\`ebres quadratiques sans diviseurs de z\'ero, qui
satisfont \`a la propri\'et\'e d'Osborn, sont de dimension finie
$1,2,4$ ou $8.\Box$
\end{theorem}

\vspace{20cm}
\section{Proc\'ed\'e de Cayley-Dickson g\'en\'eralis\'e}

\vspace{0.7cm} \subsection{Les sous-alg\`ebres de dimension $4$
dans une $\rit$-alg\`ebre de Jordan n.c. de division lin\'eaire de
dimension $8$}

\vspace{0.5cm} \hspace{0.3cm} Soit maintenant $\Big(
W,(.|.),\wedge\Big)$ une alg\`ebre r\'eelle de Jordan non
commutative de division lin\'eaire de dimension $8.$ L'application

\[ \langle .,.\rangle:W\times W\rightarrow\rit\hspace{0.3cm} (x,y)\mapsto \langle x,y\rangle=-(x|y) \]

est une forme trace d\'efinie positive qui munit $W$ d'une
structure d'espace euclidien \ $(W,-(.|.)).$ Si $u\in W-\{0\},$ on
note $W(u)$ l'orthogonal \`a $u$ dans $W$ et $||u||=\sqrt{\langle
u,u\rangle}.$

\vspace{0.3cm} \begin{lemma} Sit $u\in W-\{0\},$ alors
l'application lin\'eaire \ $L_u^*:W\rightarrow W$ \ $y\mapsto
u\wedge y$ \ est anti-sym\'etrique par rapport \`a
$\langle.,.\rangle$ et induit une bijection \ $f_u:W(u)\rightarrow
W(u).$ De plus, l'application lin\'eaire sym\'etrique $f_u^2$ est
d\'efinie n\'egative par rapport \`a $\langle.,.\rangle.$
\end{lemma}

\vspace{0.1cm} {\bf Preuve.} Le fait que $L_u^*$ soit
anti-sym\'etrique est cons\'equence de la propri\'et\'e trace de
$(.|.).$ Pour la m\^eme raison et d'apr\`es la propri\'et\'e {\bf
4.} du Th\'eor\`eme {\bf 1.49}, $L_u^*$ induit une bijection $f_u$
de $W(u)$ dans $W(u).$ Si $y_1\in W(u)$ est un vecteur propre
norm\'e ($||y_1||=1$) de l'application lin\'eaire sym\'etrique
$f_u^2$ associ\'e \`a la valeur propre $\lambda_1,$
on a: \begin{eqnarray*} \lambda_1 &=& \lambda_1||y_1||^2 \\
&=& -(\lambda_1y_1|y_1) \\
&=& -(u\wedge(u\wedge y_1)|y_1) \\
&=& -(u\wedge y_1|u\wedge y_1) \\
&=& -||u\wedge y_1||^2<0. \end{eqnarray*}

Donc $f_u^2$ est d\'efinie n\'egative par rapport \`a
$\langle.,.\rangle.$ L'\'el\'ement $u\wedge y_1$ est \'egalement
un vecteur propre de $f_u^2$ associ\'e \`a la valeur propre
$\lambda_1.\Box$

\vspace{0.3cm} \begin{lemma} $f_u^2$ poss\`ede au plus trois
valeurs propres distinctes.
\end{lemma}

\vspace{0.1cm} {\bf Preuve.} Si $y_2$ est un vecteur propre
(norm\'e) de $f_u^2$ orthogonal aux vecteurs $y_1$ et $u\wedge
y_1$ associ\'e \`a la valeur propre $\lambda_2,$ alors le vecteur
propre $u\wedge y_2$ est orthogonal aux vecteurs $y_1$ et $u\wedge
y_1.$ On construit alors une base orthogonale \[ y_1,
\hspace{0.1cm} u\wedge y_1, \hspace{0.2cm} y_2, \hspace{0.1cm}
u\wedge y_2, \hspace{0.2cm} y_3, \hspace{0.1cm} u\wedge y_3
\] de $W(u)$ form\'ee de vecteurs propres de $f_u^2$ associ\'es
respectivement aux valeurs propres \ $\lambda_1, \lambda_1,
\lambda_2, \lambda_2, \lambda_3, \lambda_3.\Box$

\vspace{0.3cm} \begin{note} Les vecteurs propres orthonormaux \
$y_1,$ $y_2,$ $y_3$ \ peuvent \^etre compl\'et\'es par les
vecteurs orthonormaux \[ z_1=(-\lambda_1)^{-\frac{1}{2}}u\wedge
y_1, \ z_2=(-\lambda_2)^{-\frac{1}{2}}u\wedge y_2, \
z_3=(-\lambda_3)^{-\frac{1}{2}}u\wedge y_3 \] en une base
orthonorm\'ee de $W(u).$ On d\'esignera par \
$vect\{x_1,\dots,x_n\}$ le sous-espace vectoriel de $A$ engendr\'e
par des \'el\'ements $x_1,\dots,x_n$ de $A.\Box$
\end{note}

\vspace{0.3cm} \begin{remark} Si $\lambda$ est une valeur propre
de $f_u^2,$ alors le sous-espace propre correspondant $E_\lambda$
est de dimension $2,$ $4$ ou $6$ et se d\'ecompose en somme
directe orthogonale de sous-espaces de dimension $2$ stables par
$f_u.\Box$
\end{remark}

\vspace{0.3cm} \begin{lemma} Il existe des vecteurs orthonormaux
$u_0, v_0\in W$ tels que $v_0$ soit un vecteur propre de
$f_{u_0}^2$ et $u_0$ soit un vecteur propre de $f_{v_0}^2$
associ\'e \`a la m\^eme valeur propre.
\end{lemma}

\vspace{0.1cm} {\bf Preuve.} Soit $S=\{x\in W: ||x||=1\}$ la
sph\`ere unit\'e de $W$ et soit $H$ l'application \[ S\times
S\rightarrow\rit^+ \ (x,y)\mapsto ||x\wedge y||. \] L'application
$H$ est continue sur le compact $S\times S$ et il existe $u_0,
v_0\in S$ tels que \ $||u_0\wedge v_0||=\sup H(S\times S).$ Selon
la Note {\bf 3.3}, il existe une base orthonorm\'ee \ $y_1, z_1,
y_2, z_2, y_3, z_3$ \ de $W(u_0)$ pour laquelle la matrice de
$f_{u_0}$ s'\'ecrit:

\vspace{0.3cm}
\[ \left(
\begin{array}{llllll}
\hspace{0.5cm} 0 & -(-\lambda_1)^{\frac{1}{2}} &                                                  & & & \\
(-\lambda_1)^{\frac{1}{2}} & \hspace{0.8cm} 0  &                                                  & & & \\
                 &                             & \hspace{0.5cm} 0           & -(-\lambda_2)^{\frac{1}{2}} & & \\
                 &                             & (-\lambda_2)^{\frac{1}{2}} & \hspace{0.8cm} 0            & & \\
                 &                             &                            &                             &
\hspace{0.5cm} 0 & -(-\lambda_3)^{\frac{1}{2}} \\
                 &                             &                            &                             &
(-\lambda_3)^{\frac{1}{2}} & \hspace{0.8cm}  0 \\
\end{array}
\right)  \]

\vspace{0.3cm} o\`u les $\lambda_i<0,$ $i=1,2,3.$ Si
$v=\sum_{1\leq i\leq 3}(\alpha_i y_i+\beta_i z_i)$ est un
\'el\'ement arbitraire de $S,$ on a: \[ f_{u_0}(v)=\sum_{1\leq
i\leq 3}\Big(
\alpha_i(-\lambda_i)^{\frac{1}{2}}z_i-\beta_i(-\lambda_i)^{\frac{1}{2}}
y_i\Big) \ \mbox{ et } \] \begin{eqnarray*} ||f_{u_0}(v)||^2 &=&
-\sum_{1\leq
i\leq 3}\lambda_i(\alpha_i^2+\beta_i^2) \\
&\leq& \sup\{-\lambda_i:1\leq i\leq 3\}\sum_{1\leq i\leq
3}(\alpha_i^2+\beta_i^2) \\
&=& \sup\{-\lambda_i:1\leq i\leq 3\}:=-\lambda. \end{eqnarray*}

\vspace{0.2cm}
Ce qui montre que la borne sup\'erieure de $H$ est
atteinte sur les vecteurs norm\'es du sous-espace propre
$E_\lambda$ de $f_{u_0}$ associ\'e \`a plus grande, en valeur
absolue, valeur propre. En effet, $vect\{y_{i_0},z_{i_0}\}:=E$ est
stable par $f_{u_0}$ et si
$v_0=\alpha_{i_0}y_{i_0}+\beta_{i_0}z_{i_0}$ est un vecteur
norm\'e arbitraire de $E,$ alors \
$f_{u_0}(v_0)=\alpha_{i_0}(-\lambda)^{\frac{1}{2}}z_{i_0}-\beta_{i_0}(-\lambda)^{\frac{1}{2}}
y_{i_0}$ i.e.
$||f_{u_0}(v_0)||^2=-\lambda(\alpha_{i_0}^2+\beta_{i_0}^2)=-\lambda.$
Par cons\'equent $v_0$ est un vecteur propre de $f_{u_0}^2$ et les
m\^emes consid\'erations sont valables pour $f_{v_0}.$ Ceci
\'etablit la premi\`ere partie du Lemme. Comme $||u_0\wedge
v_0||=-\lambda,$ la deuxi\`eme partie est \'egalement
\'etablie.$\Box$

\vspace{0.3cm} \begin{remarks} . \begin{enumerate} \item La borne
sup\'erieure de $H$ est atteinte en fait sur le compact \[
\{(u,v)\in S\times S:\langle u,v\rangle=0\}=S\times S\cap \langle
.,.\rangle^{-1}\{0\}:=U \] et on a: \ $\sup H(S\times S)=\sup
H(U)=||u_0\wedge v_0||:=M.$ \ D'autre part, il existe
$(u_1,v_1)\in U$ tel que \ $\inf H(U)=||u_1\wedge v_1||:=m>0.$

\vspace{0.1cm} \hspace{0.3cm} La sous-alg\`ebre $\rit[u_0,v_0]$ de
$A$ engendr\'ee par $u_0$ et $v_0$ a pour base \
$\{1,u_0,v_0,w_0\},$ o\`u $w_0=(-\lambda)^{-\frac{1}{2}}u_0\wedge
v_0,$ et est isomorphe \`a \
$\hit^{(\frac{1+\sqrt{-\lambda}}{2})}=\hit^{(\frac{1+M}{2})}:$ \[ \begin{tabular}{ccccc} \\
\multicolumn{1}{c}{} & \multicolumn{1}{c}{$1$} &
\multicolumn{1}{c}{$u_0$} & \multicolumn{1}{c}{$v_0$} &
\multicolumn{1}{c}{$w_0$} \\
\cline{2-5} \multicolumn{1}{c|}{$1$} & \multicolumn{1}{|c}{$1$} &
\multicolumn{1}{c}{$u_0$} & \multicolumn{1}{c}{$v_0$} &
\multicolumn{1}{c|}{$w_0$} \\
\multicolumn{1}{c|}{$u_0$} & \multicolumn{1}{|c}{$u_0$} &
\multicolumn{1}{c}{$-1$} & \multicolumn{1}{c}{$\sqrt{-\lambda} \
w_0$} &
\multicolumn{1}{c|}{$-\sqrt{-\lambda} \ v_0$} \\
\multicolumn{1}{c|}{$v_0$} & \multicolumn{1}{|c}{$v_0$} &
\multicolumn{1}{c}{$-\sqrt{-\lambda} \ w_0$} &
\multicolumn{1}{c}{$-1$} &
\multicolumn{1}{c|}{$\sqrt{-\lambda} \ u_0$} \\
\multicolumn{1}{c|}{$w_0$} & \multicolumn{1}{|c}{$w_0$} &
\multicolumn{1}{c}{$\sqrt{-\lambda} \ v_0$} &
\multicolumn{1}{c}{$-\sqrt{-\lambda} \ u_0$} &
\multicolumn{1}{c|}{$-1$} \\ \cline{2-5}
\end{tabular} \]

Par analogie, on a \ $\rit[u_1,v_1]\simeq\hit^{(\frac{1+m}{2})}.$
\item Si $(u,v)\in U$ o\`u $u,v\in\rit[u,v],$ on a \
$\rit[u,v]=\rit[u_0,v_0]\simeq\hit^{(\frac{1+M}{2})}$ \ i.e.
$||u\wedge v||=M.\Box$
\end{enumerate}
\end{remarks}

\vspace{0.2cm} \hspace{0.3cm} On obtient le r\'esultat suivant:

\vspace{0.3cm} \begin{theorem} Soit $A$ une $\rit$-alg\`ebre de
Jordan non commutative de division lin\'eaire de dimension $8.$
Alors $A$ contient une sous-alg\`ebre de dimension $4,$ isomorphe
\`a $\hit^{(\lambda)}$ o\`u
$\lambda\in\rit-\{{\frac{1}{2}}\}.\Box$
\end{theorem}

\vspace{0.3cm} \begin{remark} Les alg\`ebres r\'eelles de d.l. de
dimension $8$ que nous connaissons contiennent une sous-alg\`ebre
de dimension $4,$ mais nous ne savons pas encore si toutes les
$\rit$-alg\`ebres quadratiques de d.l. de dimension $8$
contiennent une sous-alg\`ebre de dimension $4.$ Cette question
est apparemment un probl\`eme ouvert.$\Box$
\end{remark}

\vspace{0.3cm} \begin{theorem} Soit $A$ une $\rit$-alg\`ebre de
Jordan non commutative de division lin\'eaire de dimension $8.$
Alors les deux propri\'et\'es suivantes sont \'equivalentes:
\begin{enumerate}
\item $A$ satisfait \`a la propri\'et\'e d'Osborn. \item Deux
sous-alg\`ebres de $A,$ de dimension $4,$ sont isomorphes.
\end{enumerate}
\end{theorem}

\vspace{0.1cm} {\bf Preuve.} Nous utilisons, pour la
d\'emonstration, les notations du Lemme {\bf 3.5}.

\vspace{0.2cm} 2) $\Rightarrow$ 1). Les sous-alg\`ebres
$\hit^{(\frac{1+M}{2})}$ et $\hit^{(\frac{1+m}{2})}$ (Remarque \bf
3.6 1)}) sont isomorphes si et seulement si $m=M$ ([Roc 87] p.
42). Ainsi, pour tout $(u,v)\in U,$ on a \ $m=||u\wedge v||=M$ \
i.e. la sous-alg\`ebre $\rit[u,v]$ est de dimension $4.$ Donc $A$
satisfait \`a la propri\'et\'e d'Osborn.

\vspace{0.2cm} 1) $\Rightarrow$ 2). Soit $v=\sum_{1\leq i\leq
3}(\alpha_iy_i+\beta_iz_i)$ un vecteur norm\'e de $A$ tel que
$(u_0,v)\in U,$ alors $v$ et $f_{u_0}^2(v)=\sum_{1\leq i\leq
3}\lambda_i(\alpha_iy_i+\beta_iz_i)$ sont colin\'eaires, car $u_0$
et $v$ satisfont \`a la propri\'et\'e d'Osborn:
$u_0\wedge(v\wedge(u_0\wedge v))=v\wedge(u_0\wedge(v\wedge
u_0))=0.$ Donc $f_{u_0}^2$ admet une unique valeur propre, et on
a: $M=||u_0\wedge v_0||=||u_0\wedge v||.$ Soit maintenant
$(u,v)\in U,$ il existe $v_0'\in W(u)\cap\rit[u_0,v_0]-\{0\}$ et
$u_0'\in W(v_0')\cap\rit[u_0,v_0]-\{0\},$ qu'on peut choisir
norm\'es, et on a \begin{eqnarray*} M &=& ||u_0'\wedge v_0'|| \
\mbox{ (Remarque {\bf 3.6 2)} } \\
&=& ||u\wedge v_0'|| \ \mbox{ car } (u,v_0')\in U \\
&=& ||u\wedge v||. \end{eqnarray*}

Donc $\rit[u,v]\simeq \hit^{(\frac{1+M}{2})}$ l'unique, \`a
isomorphisme pr\`es, sous-alg\`ebre de $A$ de dimension $4.\Box$

\vspace{0.5cm} \hspace{0.3cm} Nous avons besoin, pour la suite,
des deux r\'esultats pr\'eliminaires suivants:

\vspace{1cm} \begin{lemma} Soit $\lambda\in K$ et soit $D$ une
$K$-alg\`ebre pour laquelle il existe un \'el\'ement $u$ tel que
$L_u-R_u\neq 0.$ Alors \ $D^{(\lambda)}$ est commutative si et
seulement si $\lambda=\frac{1}{2}.$
\end{lemma}

\vspace{0.1cm} {\bf Preuve.} Pour tout $x\in D,$ les op\'erateurs
de multiplication \`a gauche et \`a droite par $x$ dans
$D^{(\lambda)}$ sont: $L_x^{(\lambda)}=\lambda L_x+(1-\lambda)R_x$
\ et \ $R_x^{(\lambda)}=\lambda R_x+(1-\lambda)L_x.$ On a
\begin{eqnarray*} L_x^{(\lambda)}-R_x^{(\lambda)} &=& (2\lambda-1)(L_x-R_x). \end{eqnarray*}
Donc $D^{(\lambda)}$ est commutative si et seulement si
$\lambda=\frac{1}{2}.\Box$

\vspace{0.3cm} \begin{lemma} Soient $\lambda,\mu\in K$ et $D$ une
$K$-alg\`ebre alternative pour laquelle il existe un \'el\'ement
$u$ tel que $(L_u-R_u)^2\neq 0.$ Alors \ $D^{(\lambda)}$ est
alternative si et seulement si $\lambda=0$ ou $1.$ Si, de plus,
$D$ est cayleyenne, on a $D^{(\lambda)}\simeq D^{(\mu)}$ si et
seulement si $\lambda=\mu$ ou $\lambda=1-\mu.$
\end{lemma}

\vspace{0.1cm} {\bf Preuve.} Nous utilisons, pour la
d\'emonstration, les notations du Lemme {\bf 3.10}. On a
\begin{eqnarray*} L_{x^2}^{(\lambda)} &=& \lambda L_{x^2}+(1-\lambda)R_{x^2} \\
&=& \lambda L_x^2+(1-\lambda)R_x^2 \end{eqnarray*} et
\begin{eqnarray*} (L_x^{(\lambda)})^2 &=&
\lambda^2L_x^2+2\lambda(1-\lambda)L_xR_x+(1-\lambda)^2R_x^2.
\end{eqnarray*} Donc
\begin{eqnarray*} L_{x^2}^{(\lambda)}-(L_x^{(\lambda)})^2 &=& (\lambda-\lambda^2)L_x^2-2\lambda(1-\lambda)L_xR_x+
\lambda(1-\lambda)R_x^2 \\
&=& \lambda(1-\lambda)(L_x-R_x)^2.
\end{eqnarray*}
De m\^eme \
$R_{x^2}^{(\lambda)}-(R_x^{(\lambda)})^2=\lambda(1-\lambda)(R_x-L_x)^2.$
Ainsi, $D^{(\lambda)}$ est alternative si et seulement si
$\lambda=0$ ou $\lambda=1.$

\vspace{0.2cm} \hspace{0.3cm} Soient maintenant
$\lambda,\mu\in\rit$ que l'on peut supposer, en vertu du Lemme
{\bf 3.10}, distincts de $\frac{1}{2},$ alors \
$D^{(\lambda)}\simeq D^{(\mu)}\Leftrightarrow\Big(
D^{(\lambda)}\Big)^{(\frac{\mu}{2\mu-1})}=D^{(\frac{\lambda+\mu-1}{2\mu-1})}\simeq
\Big( D^{(\mu)}\Big)^{(\frac{\mu}{2\mu-1})}=D \ \mbox{ alternative
}.$ \ Donc $\frac{\lambda+\mu-1}{2\mu-1}=0$ ou $1$ i.e.
$\lambda=\mu$ ou $\lambda=1-\mu.$ La r\'eciproque, de cette
derni\`ere proposition, est cons\'equence de la Remarque {\bf 1.57
3)}.$\Box$

\vspace{0.3cm} \begin{theorem} Les uniques $\rit$-alg\`ebres de
Jordan non commutatives de division lin\'eaire, de dimension finie
$\leq 4,$ sont \`a isomorphisme pr\`es \ $\rit,$ $\cit$ ou
$\hit^{(\lambda)},$ $\lambda\neq\frac{1}{2}.$ De plus \
$\hit^{(\lambda)}\simeq \hit^{(\mu)},$
$\lambda,\mu\neq\frac{1}{2}$ si et seulement si $\lambda=\mu$ ou
$\lambda=1-\mu.$
\end{theorem}

\vspace{0.1cm} {\bf Preuve.} ([Roc 87] p. 42).$\Box$

\vspace{0.3cm} \begin{theorem} Les uniques $\rit$-alg\`ebres de
Jordan non commutatives sans diviseurs de z\'ero, qui satisfont
\`a la propri\'et\'e d'Osborn, sont de dimension finie $1,2,4$ ou
$8$ et \`a isomorphisme pr\`es \ $\rit,$ $\cit,$
$\hit^{(\lambda)}$ ou \ $\oit^{(\lambda)},$
$\lambda\neq\frac{1}{2}.$ De plus \ $D^{(\lambda)}\simeq
D^{(\mu)},$ $D=\hit$ ou $\oit$ et $\lambda,\mu\neq\frac{1}{2}$ si
et seulement si $\lambda=\mu$ ou $\lambda=1-\mu.$
\end{theorem}

\vspace{0.1cm} {\bf Preuve.} Une telle alg\`ebre $A,$ est unitaire
quadratique de division et est de dimension finie $1,2,4$ ou $8$
(Th\'eor\`eme {\bf 2.65}). Le Th\'eor\`eme {\bf 3.12} permet de
supposer que $A$ est de dimension $8,$ et l'on d\'esigne par
$\hit^{(\lambda)},$ $\lambda\neq \frac{1}{2}$ la classe
d'isomorphisme des sous-alg\`ebres de $A$ de dimension $4.$
L'alg\`ebre $A^{(\frac{\lambda}{2\lambda-1})}$ est de division
lin\'eaire (Proposition {\bf 1.47}), et si $x,y$ sont deux
vecteurs non nuls orthogonaux de $A,$ on a
\begin{eqnarray*} \rit_{A^{(\frac{\lambda}{2\lambda-1})}}[x,y] &=&
\Big( \rit_A[x,y]\Big)^{(\frac{\lambda}{2\lambda-1})}
\ \mbox{ Lemme {\bf 1.3}} \\
&=& \Big( \hit^{(\lambda)}\Big)^{(\frac{\lambda}{2\lambda-1})}
\\
&=& \hit \ \mbox{ associative }. \end{eqnarray*}

D'apr\`es le Th\'eor\`eme d'Artin,
$A^{(\frac{\lambda}{2\lambda-1})}$ est alternative, et d'apr\`es
le Th\'eor\`eme d'Albert, $A^{(\frac{\lambda}{2\lambda-1})}\simeq$
$\oit$ i.e. $A=\Big(
A^{(\frac{\lambda}{2\lambda-1})}\Big)^{(\lambda)}\simeq$
$\oit^{(\lambda)}.\Box$

\vspace{0.3cm} \begin{proposition} Soit $A=\Big(
W,(.|.),\wedge\Big)$ une $\rit$-alg\`ebre de Jordan non
commutative de division lin\'eaire, de dimension $8.$ Si $B$ est
une sous-alg\`ebre de $A,$ de dimension $4$ telle que \
$B^{\perp}B^{\perp}\subseteq B.$ Alors il existe $(u_0,x_0)\in
B\times B^{\perp}$ tel que $u_0$ et $x_0$ engendrent une
sous-alg\`ebre de $A$ de dimension $4.$
\end{proposition}

\vspace{0.1cm} {\bf Preuve.} Si $x\in S(B^{\perp}),$ alors
l'application $f_x^2$ induit une bijection \ $g_x:B\cap
W\rightarrow B\cap W$ \ sym\'etrique. Si $x_i,$ $i=1,2,3$
d\'esigne une base orthonorm\'ee de $B\cap W$ form\'ee de vecteurs
propres de $g_x$ associ\'es respectivement aux valeurs propres
$\lambda_i,$ $i=1,2,3$ alors les vecteurs \ $x_1,$ $x'_1,$ $x_2,$
$x'_2,$ $x_3,$ $x'_3$ o\`u
$x'_i=(-\lambda_i)^{-\frac{1}{2}}x\wedge x_i$ constituent une base
orthonorm\'ee de $W(x)$ form\'ee de vecteurs propres de $f_x^2$
associ\'es respectievement aux valeurs propres \[ \lambda_1, \
\lambda_1, \hspace{0.2cm} \lambda_2, \ \lambda_2, \hspace{0.2cm}
\lambda_3, \ \lambda_3. \]

\vspace{0.2cm} \hspace{0.3cm} Soit maintenant $u$ un \'el\'ement
de $S(B\cap W),$ alors l'application $f_u$ induit une bijection \
$h_u:B^{\perp}\rightarrow B^{\perp}$ \ anti-sym\'etrique. Si $y_i,
z_i$ o\`u $i\in\{2,3\}$ d\'esigne une base orthonorm\'ee de
$B^{\perp}$ form\'ee de vecteurs propres de $h_u^2,$ avec
$z_i=(-\mu_1)^{-\frac{1}{2}}u\wedge y_i,$ associ\'es
respectivement aux valeurs propres \ $\mu_1,$ $\mu_1,$ \ $\mu_2,$
$\mu_2,$ \ on peut la compl\'eter en une base de $W(u)$ form\'ee
de vecteurs propres de $f_u^2,$ par des vecteurs orthonormaux \
$y_1,$ $z_1$ \ de $B$ orthogonaux \`a $u,$ avec
$z_1=(-\mu_1)^{-\frac{1}{2}}u\wedge y_1,$ associ\'es \`a la valeur
propre $\mu_1.$

\vspace{0.2cm} \hspace{0.3cm} Soit maintenant $H$ l'application
continue \ $S(B\cap W)\times S(B^{\perp})\rightarrow\rit^+$ \
$(u,x)\mapsto ||u\wedge x||,$ \ il existe $(u_0,x_0)\in S(B\cap
W)\times S(B^{\perp})$ tel que \ $\sup H=||u_0\wedge x_0||.$ On
utilise alors les notations pr\'ec\'edentes pour les \'el\'ements
propres de $h_{u_0}^2$ et $g_{x_0}.$ Si $x\in S(B^{\perp}),$ il
existe des scalaires $a_i, b_i$ o\`u $i\in\{1,2\},$ avec
$\sum_{1\leq i\leq 2}(a_i^2+b_i^2)=1,$ tels que $x=\sum_{1\leq
i\leq 2}(a_iy_i+b_iz_i)$ et on a \ $u_0\wedge x=\sum_{1\leq i\leq
2}\Big(
a_i(-\mu_i)^{\frac{1}{2}}z_i-b_i(-\mu_i)^{\frac{1}{2}}y_i\Big).$ \
Ainsi \begin{eqnarray*} ||u_0\wedge x||^2 &=& \sum_{1\leq i\leq
2}(-\mu_i)(a_i^2+b_i^2) \\
&\leq& \sup_{1\leq i\leq 2}(-\mu_i). \end{eqnarray*} Donc $x_0$
est un vecteur propre de $h_{u_0}^2$ associ\'e \`a la plus grande,
en valeur absolue, valeur propre.

\vspace{0.2cm} \hspace{0.3cm} Soit maintenant $u\in S(B\cap W),$
il existe des scalaires $\gamma_i,$ $i\in\{1,2,3\}$ avec
$\sum_{1\leq i\leq 3}\gamma_i^2=1,$ tels que $u=\sum_{1\leq i\leq
3}\gamma_ix_i,$ et on a \ $x_0\wedge u=\sum_{1\leq i\leq
3}\gamma_i(-\lambda_i)^{\frac{1}{2}}x'_i.$ Ainsi \begin{eqnarray*}
||x_0\wedge u||^2 &=& \sum_{1\leq i\leq
3}(-\lambda_i)\gamma_i^2 \\
&\leq& \sup_{1\leq i\leq 2}(-\lambda_i). \end{eqnarray*} Donc
$u_0$ est un vecteur propre de $g_{x_0}$ associ\'e \`a la plus
grande, en valeur absolue, valeur propre. Les vecteurs $u_0$ et
$x_0$ engendrent alors une sous-alg\`ebre de $A$ de dimension
$4.\Box$

\vspace{0.7cm} \subsection{Proc\'ed\'e de Cayley-Dickson
g\'en\'eralis\'e}

\vspace{0.5cm} \begin{GCD process} Soit $(B,\overline{.})=\Big(
W,(.|.),\wedge\Big)$ une $K$-alg\`ebre cayleyenne et soient
$\gamma, \alpha, \beta, \delta, \theta\in K$ avec $\gamma\neq 0.$
Alors le produit \[ (x,y)(x',y')=\Big(
x.^{\alpha}x'+\frac{\beta}{2}([x,y']+[y,x'])+\gamma\overline{y'}y,
y\overline{x'}+y'x+\frac{\delta}{2}[y',y]+\frac{\theta}{2}[x',x]\Big)
\]

munit l'espace vectoriel $B\times B$ d'une structure de
$K$-alg\`ebre cayleyenne associ\'ee \`a la m\^eme forme
bilin\'eaire sym\'etrique \ $\Big( (x,y),
(x',y')\Big)=(x|x')+\gamma(y|\overline{y'}),$ \ que celle de
l'extension cayleyenne $E_\gamma(B)$ de $(B,\overline{.}),$
d'indice $\gamma.$ On l'appelle extension cayleyenne
"g\'en\'eralis\'ee" de $(B,\overline{.}),$ d'indice
$(\gamma,\alpha,\beta,\delta,\theta)$ et on note \
$E_{\gamma,\alpha,\beta,\delta,\theta}(B).$ Les alg\`ebres \
$E_{\gamma,\alpha,0,\delta,0}(B)$ et $E_{\gamma,1,0,\delta,0}(B)$
not\'ees respectivement $E_{\gamma,\alpha,\delta}(B)$ et
$E_{\gamma,\delta}(B).\Box$
\end{GCD process}

\vspace{4cm} \begin{remarks} .
\begin{enumerate} \item L'application \ $B^{(\alpha)}\rightarrow E_{\gamma,\alpha,\beta,\delta,\theta}(B)$ \
$x\mapsto (x,0)$ \ est un monomorphisme d'alg\`ebres. \item
$E_{\gamma,\alpha,\beta,\delta,\theta}(B)=E_{\gamma}(B)$ \ si et
seulement si \ $\alpha-1=\beta=\delta=\theta=0$ ou $B$ est
commutative. \item Si le polyn\^ome \ $X^2+\gamma$ \ poss\`ede une
racine $\omega,$ dans $K,$ alors l'application \[
E_{\gamma,\alpha,\beta,\delta,\theta}(B)\rightarrow
E_{-1,\alpha,\beta\omega^{-1},\delta\omega^{-1},\theta\omega}(B)
\hspace{0.2cm} (x,y)\mapsto (x,\omega y) \] est un isomorphisme
d'alg\`ebres. \item L'op\'erateur \ $L_{(0,1)}^2$ de
$E_{\gamma,\alpha,\beta,\delta,\theta}(B):=A$ est une
homoth\'etie, \'egale \`a $\gamma I_A.\Box$
\end{enumerate}
\end{remarks}

\vspace{0.3cm} \begin{proposition} Soient $(B,\overline{.})$ une
$K$-alg\`ebre cayleyenne, \ $\gamma,\alpha,\beta,\delta,\theta\in
K,$ avec $\gamma\neq 0,$ \ et $A$ l'extension cayleyenne
g\'en\'eralis\'ee de $(B,\overline{.})$ d'indice
$(\gamma,\alpha,\beta,\delta,\theta).$ Alors
\begin{enumerate} \item $A$ est associative si et seulement si $B$ est associative et commutative.
\item Si $\beta=\gamma\theta$ et $\alpha\neq\frac{1}{2},$ alors
$A$ est flexible si et seulement si $B$ est flexible.
\end{enumerate}
\end{proposition}

\vspace{0.1cm} {\bf Preuve.} . \begin{enumerate} \item Si $A$ est
aassociative, alors pour tous $x,y\in B,$ on a \begin{eqnarray*}
(0,0) &=& \Big(
(x,0),(0,1),(y,0)\Big) \\
&=& (\frac{\beta}{2}[x,y-\overline{y}],[x,\overline{y}])
\end{eqnarray*} i.e. $B$ est commutative, et on a $A=E_\gamma(B).$ Par
cons\'equent $B$ est associative. R\'eciproquement, si $B$ est
associative et commutative, alors $A=E_\gamma(B)$ est associative.
\item On suppose que $\alpha\neq\frac{1}{2}.$ Si $A$ est flexible,
alors pour tous $x,y\in B,$ on a
\begin{eqnarray*}
(0,0) &=& \Big(
(x,0),(y,0),(x,0)\Big) \\
&=& \Big(
(x,y,x)^{(\alpha)},\frac{\theta}{2}([y,x]\overline{x}+[x,x.^{\alpha}y]-[x,y]x-[y.^{\alpha}x,x])\Big)
\end{eqnarray*}
o\`u $(x,y,x)^{(\alpha)}$ d\'esigne l'associateur de $x,y$ et $x$
dans $B^{(\alpha)}.$ Donc $B^{(\alpha)}$ est flexible et, comme
$\alpha\neq\frac{1}{2},$ $B$ est flexible.

\vspace{0.2cm} \hspace{0.3cm} On suppose maintenant que
$\beta=\gamma\theta$ et $(B,\overline{.})$ flexible. Soient alors
$(x,y),$ $(x',y'),$ $(x",y")$ trois \'el\'ements de $A$ et
$(.|.),$ $(.,.)$ respectivement la forme trace associ\'ee \`a $B$
et la forme bilin\'eaire sym\'etrique associ\'ee \`a $A.$ En
utilisant les trois identit\'es de la Remarque {\bf 1.57 2)}, on a
\[ \Big( (x,y)(x',y'),(x",y")\Big) = \] \[ \Big(
x.^{\alpha}x'+\frac{\beta}{2}([x,y']+[y,x'])+\gamma\overline{y}'y\mid
x"\Big)+\gamma\Big(
y\overline{x}'+y'x+\frac{\delta}{2}[y',y]+\frac{\theta}{2}[x',x]\mid\overline{y}"\Big)=
\] \[ \Big(
x\mid
x'.^{\alpha}x"+\frac{\beta}{2}[y',x"]+\frac{\gamma\theta}{2}[\overline{y}",x']+\gamma\overline{y}"y'\Big)+
\gamma\Big(
y\mid\frac{\delta}{2}[\overline{y}",y']+\frac{\theta}{2}[x',x"]+x"\overline{y}'+\overline{x}'\overline{y}"\Big)=
\] \[ \Big(
x\mid
x'.^{\alpha}x"+\frac{\beta}{2}([x',y"]+[y',x"])+\gamma\overline{y}"y'\Big)+
\gamma\Big( y\mid
\overline{y'\overline{x}"+y"x'+\frac{\delta}{2}[y",y']+\frac{\theta}{2}[x",x']}\Big)=
\] \[ \Big( (x,y),(x',y')(x",y")\Big) \] Donc $A$ est
flexible.$\Box$
\end{enumerate}

\vspace{0.3cm} \begin{remarks} . \begin{enumerate} \item Si $B$
est une $K$-alg\`ebre des quaternions, de division, alors
l'alg\`ebre $E_{-1,\alpha,\beta,\delta,\theta}(B)$ est alternative
si et seulement si $\alpha=1$ et $\beta=\delta=\theta=0.$ \item Si
$\beta\neq\gamma\theta,$ alors
$E_{\gamma,\alpha,\beta,\delta,\theta}(B)$ n'est pas
n\'ecessairement flexible. En effet, l'alg\`ebre r\'eelle
cayleyenne \ $E_{-1,1,0,0,1}(\hit)$ n'est pas flexible.$\Box$
\end{enumerate}
\end{remarks}

\vspace{10cm} \subsection{Automorphismes et d\'erivations des
alg\`ebres obtenues par le proc\'ed\'e de Cayley-Dickson
g\'en\'eralis\'e}

\vspace{0.5cm} \begin{proposition} Soient $(B,\overline{.})$ une
$K$-alg\`ebre cayleyenne, $\gamma,\alpha,\beta,\delta,\theta\in
K,$ avec $\gamma\neq 0,$ et $A$ l'extension cayleyenne
g\'en\'eralis\'ee de $(B,\overline{.})$ d'indice
$(\gamma,\alpha,\beta,\delta,\theta).$ Si \ $\partial$ est une
d\'erivation de $B$ qui commute avec l'involution \
$\overline{.},$ \ alors l'application \ $D_\partial:A\rightarrow
A$ \hspace{0.2cm} $(x,y)\mapsto(\partial x,\partial y)$ \ est une
d\'erivation de $A.$ Si, de plus, les d\'erivations de $B$
commutent avec l'involution \ $\overline{.},$ alors l'application
\ $Der(B)\rightarrow Der(A)$ \hspace{0.2cm} $\partial\mapsto
D_\partial$ \ est un monomorphisme d'alg\`ebres de Lie.
\end{proposition}

\vspace{0.1cm} {\bf Preuve.} Soient $(x,y), (x',y')$ deux
\'el\'ements de $A,$ on a

\vspace{0.2cm}
\[ D_\partial\Big( (x,y)(x',y')\Big)=\Big(
\partial(x.^{\alpha}x'+\frac{\beta}{2}([x,y']+[y,x'])+\gamma\overline{y'}y),
\partial(y\overline{x'}+y'x+\frac{\delta}{2}[y',y]+\frac{\theta}{2}[x',x])\Big)
\] \[
=\Big( (\partial x).^{\alpha}x'+x.^{\alpha}\partial
x'+\frac{\beta}{2}([\partial x,y']+[x,\partial y']+[\partial
y,x']+[y,\partial x'])+\gamma\overline{\partial
y'}y+\gamma\overline{y}'\partial y, \hspace{1.5cm} . \] \[
.\hspace{3cm} (\partial y)\overline{x}'+y\overline{\partial
x'}+(\partial y')x+y'\partial x+\frac{\delta}{2}([\partial
y',y]+[y',\partial y])+\frac{\theta}{2}([\partial
x',x]+[x',\partial x])\Big)
\]
\[ =(\partial x,\partial y)(x',y')+(x,y)(\partial x',\partial
y')=\Big( D_\partial(x,y)\Big)(x',y')+(x,y)D_\partial(x',y'). \]
On dira que $D_\partial$ est une extension naturelle de $\partial$
\`a $A.\Box$

\vspace{0.5cm} \begin{proposition} Soit $(B,\overline{.})$ une
$K$-alg\`ebre cayleyenne des quaternions, de division. Alors les
d\'erivations de $B$ commutent avec l'involution \ $\overline{.}$
. Si, de plus, $\gamma, \alpha, \delta\in K,$ avec
$\gamma\delta\neq 0$ et $\alpha\neq\frac{1}{2},$ alors toute
d\'erivation de $E_{\gamma,\alpha,\delta}(B)$ est une extension
naturelle d'une d\'erivation de $B,$ qui commute avec l'involution
\ $\overline{.},$ i.e. \ $Der(B)\rightarrow Der\Big(
E_{\gamma,\alpha,\delta}(B)\Big)$ \hspace{0.2cm} $\partial\mapsto
D_\partial$ \ est un isomorphisme d'alg\`ebres de Lie.
\end{proposition}

\vspace{0.1cm} {\bf Preuve.} Les d\'erivations de $B$ sont
int\'erieures (Th\'eor\`eme {\bf 1.21}) de la forme \ $L_x-R_x,$
$x\in B,$ donc commutent avec $\overline{.}$ . On suppose
maintenant $\delta\neq 0.$ Si $D$ est une d\'erivation de
$E_{\gamma,\alpha,\delta}(B),$ il existe $\partial, f\in End_K(B)$
et $x_0,y_0\in B$ tels que \ $(a,0)=(\partial a, f(a),$ pour tout
$a\in B,$ et $D(0,1)=(x_0,y_0).$ Ainsi, pour tous $a,b\in B,$ on a
\begin{eqnarray*} D(a,b) &=& D\Big( (a,0)+(b,0)(0,1)\Big) \\
&=& D(a,0)+\Big( D(b,0)\Big)(0,1)+(b,0)D(0,1) \\
&=& (\partial a,f(a))+(\partial b,f(b))(0,1)+(b,0)(x_0,y_0) \\
&=& \Big( \partial a+\gamma f(b)+b.^{\alpha}x_0, f(a)+\partial
b+y_0b\Big). \end{eqnarray*} L'\'egalit\'e $D\Big(
(a,b)(c,d)\Big)=D\Big( (a,b)\Big)(c,d)+(a,b)D(c,d)$ donne:

\vspace{0.1cm}
\begin{enumerate} \item $\partial\in
Der(B^{(\alpha)}),$ \ en faisant $b=d=0.$ \item
$f(a.^{\alpha}c)=f(a)\overline{c}+f(c)a,$ \ en faisant $b=d=0.$
\item $\gamma f(b\overline{c})+(b\overline{c}).^{\alpha}x_0=\gamma
f(b).^{\alpha}c+\gamma\overline{f(c)}b+(b.^{\alpha}x_0).^{\alpha}c,$
\ en faisant $a=d=0.$ \item $\partial(b\overline{c})=(\partial
b)\overline{c}+b\overline{\partial c}+\frac{\delta}{2}[f(c),b],$ \
en faisant $a=d=0.$ \item
$\gamma\partial(\overline{d}b)+\frac{\gamma\delta}{2}f([d,b])+\frac{\delta}{2}[d,b].^{\alpha}x_0=\gamma\Big(
\overline{d}\partial b+\overline{\partial
d}.b+(y_0+\overline{y}_0)\overline{d}b\Big),$ \ en faisant
$a=c=0.$ \item $\gamma f(\overline{d}b)+\frac{\delta}{2}\partial
[d,b]+\frac{\delta}{2}y_0[d,b]=$ \[ d\Big( \gamma
f(b)+b.^{\alpha}x_0\Big)+\frac{\delta}{2}[d,\partial
b+y_0b]+b\Big(
\gamma\overline{f(d)}+\overline{x}_0.^{\alpha}\overline{d}\Big)+\frac{\delta}{2}[\partial
d+y_0d.b], \] en faisant $a=c=0.$

\vspace{0.2cm} De plus $D(1,0)=(0,0)$ donne \item $\partial
1=f(1)=0.$
\end{enumerate}

Le fait que $\alpha$ soit distincte de $\frac{1}{2}$ et {\bf 1.}
montrent que $\partial\in Der(B).$ Ainsi, en faisant $b=1$ dans
{\bf 4.}, il r\'esulte

\vspace{0.2cm} {\bf 1'.} $\partial\overline{c}=\overline{\partial
c}$ \ et \ $\frac{\delta}{2}[f(c),b]=0.$ Comme $B$ est centrale et
$\delta\neq 0,$ on a

\vspace{0.2cm} {\bf 2'.} $f(B)\subseteq K.$

\vspace{0.2cm} \hspace{0.3cm} Soit maintenant $a$ un vecteur de
$B,$ il existe un vecteur $c$ de $B$ lin\'eairement ind\'ependant
\`a $a$. Les \'egalit\'es {\bf 2.} et {\bf 2'.} entrainent alors
que $f(a)=0,$ et on d\'eduit de la seconde \'egalit\'e de {\bf 7.}
que $f\equiv 0.$ De {\bf 5.} et {\bf 6.} il r\'esulte que $x_0$ et
$y_0$ sont des vecteurs et on obtient, en faisant $(b,d)=(1,1)$
dans {\bf 5.} et {\bf 6.} puis $d=1$ dans {\bf 7.}, l'\'egalit\'e
\ $[\frac{\delta}{2}y_0+(1-\alpha)x_0,b]=0.$ Donc \
$\frac{\delta}{2}y_0+(1-\alpha)x_0=0,$ et en faisant
$(b,c)=(1,x_0)$ dans {\bf 3.}, on obtient $2x_0^2=0$ i.e.
$x_0=y_0=0.$ Ainsi $D=D_\partial.\Box$

\vspace{0.3cm}
\begin{remark} On a \[ \bigcap_{D\in Der(E_{\gamma,\alpha,\delta}(B))}\ker D=
\bigcap_{\partial\in Der(B)}\ker D_\partial=K\times K.\Box \]
\end{remark}

\vspace{0.3cm}
\begin{corollary} L'inclusion $Der(\hit)\subset Der\Big( E_{-1,\alpha,\delta}(\hit)\Big)$ a lieu pour tous r\'eels
$\alpha, \delta.$ Si, de plus, $\alpha\neq\frac{1}{2}$ et
$\delta\neq 0,$ on a l'\'egalit\'e: $Der(\hit)=Der\Big(
E_{-1,\alpha,\delta}(\hit)\Big).\Box$
\end{corollary}

\vspace{0.3cm}
\begin{proposition} Soient $B=\Big( V,(.|.),\wedge\Big)$ une $K$-alg\`ebre des quaternions, de division,
$\gamma, \alpha, \delta\in K,$ avec $\gamma\neq 0,$ et soit $A$
l'extension cayleyenne de $B$ d'indice $(\gamma,\alpha,\delta).$
On suppose que $\delta\neq 0$ et $\alpha\neq\frac{1}{2}.$ Alors \
$K\times\{0\},$ $\{0\}\times K,$ $V\times\{0\}$ et $\{0\}\times V$
\ sont des sous-$Der(A)$-modules irr\'eductibles de $A$ (leur
somme directe).
\end{proposition}

\vspace{0.1cm} {\bf Preuve.} D'apr\`es la Proposition {\bf 3.20},
$Der(A)=\{D_\partial:\partial\in Der(B)\},$ ainsi les sous-modules
$K\times\{0\}$ et $\{0\}\times K$ sont annul\'es par les
\'el\'emenst de $Der(A).$ D'autre part, comme les d\'erivations de
$B$ sont int\'erieures, et vu que $V$ est dimension $3$ sur $K,$
ces d\'erivations ne laissent invariant aucun sous-$Der(B)$-module
propre de $V.$ Donc les d\'erivations $D_\partial$ ne laissent
invariant aucun sous $Der(A)$-module propre de $V\times\{0\}$ ou
de $\{0\}\times V.$ Ces deux derniers sont \`a leur tour des
sous-$Der(A)$-modules irr\'eductibles de $A.\Box$

\vspace{0.3cm}
\begin{proposition} Soit $B=\Big( V,(.|.),\wedge\Big)$ une $K$-alg\`ebre cayleyenne et soient
$\gamma, \alpha, \beta, \delta, \theta\in K,$ avec $\gamma\neq 0.$
Si $f$ est un automorphisme de $B$ qui commute avec l'involution
$\overline{.},$ alors l'application \
$\Phi_f:E_{\gamma,\alpha,\beta,\delta,\theta}(B)\rightarrow
E_{\gamma,\alpha,\beta,\delta,\theta}(B)$ \
$(x,y)\mapsto(f(x),f(y))$ \ est un automorphisme. Si, de plus, les
automorphismes de $B$ commutent avec l'involution $\overline{.},$
alors l'application \ $Aut(B)\rightarrow Aut\Big(
E_{\gamma,\alpha,\beta,\delta,\theta}(B)\Big)$ \ $f\mapsto\Phi_f$
\ est un monomorphisme de groupes.
\end{proposition}

\vspace{0.1cm} {\bf Preuve.} Soient $(x,y),$ $(x',y')$ deux
\'el\'ements de $E_{\gamma,\alpha,\beta,\delta,\theta}(B),$ on a
\begin{eqnarray*} \Phi_f\Big( (x,y)(x',y')\Big) &=&
\Big( f(x.^{\alpha}x'+\frac{\beta}{2}([x,y']+[y,x'])+\gamma\overline{y'}y),
f(y\overline{x'}+y'x+\frac{\delta}{2}[y',y]+\frac{\theta}{2}[x',x])\Big) \\
&=& \Big(
f(x).^{\alpha}f(x')+\frac{\beta}{2}([f(x),f(y')]+[f(y),f(x')])+\gamma\overline{f(y')}f(y)),
\\ && \hspace{2.5cm} f(y)\overline{f(x')}+f(y')f(x)+\frac{\delta}{2}[f(y'),f(y)]+\frac{\theta}{2}[f(x'),f(x)])\Big) \\
&=& (f(x),f(y))(f(x'),f(y')) \\
&=& \Phi_f(x,y)\Phi_f(x',y'). \end{eqnarray*}

On dira que $\Phi_f$ est une extension naturelle de $f$ \`a
$E_{\gamma,\alpha,\beta,\delta,\theta}(B).\Box$

\vspace{0.3cm}
\begin{lemma} Soient $A$ et $B$ deux $K$-alg\`ebres et
$\Phi:A\rightarrow B$ un isomorphisme d'alg\`ebres. Alors pour
toute d\'erivation $\partial$ de $A,$ \ $\Phi\partial\Phi^{-1}$
est une d\'erivation de $B$ et l'application \ $Der(A)\rightarrow
Der(B)$ \ $\partial\mapsto\Phi\partial\Phi^{-1}$ \ est un
isomorphisme d'alg\`ebres de Lie. Si, de plus, $A=\Big(
V,(.|.),\times\Big)$ et $B=\Big( W,\langle.,.\rangle,\wedge\Big)$
sont quadratiques, alore $\Phi(1)=1$ et $\Phi(V)=W.$
\end{lemma}

\vspace{0.1cm} {\bf Preuve.} Pour tous $x,y\in A,$ on a
$\Phi(xy)=\Phi(x)\Phi(y)$ i.e. \ $L_{\Phi(x)}=\Phi L_x\Phi^{-1}.$

\vspace{0.2cm} \hspace{0.3cm} Soient maintenant $\partial\in
Der(A)$ et $x\in B,$ on a
\begin{eqnarray*} [\Phi\partial\Phi^{-1},L_x] &=& [\Phi\partial\Phi^{-1},\Phi L_{\Phi^{-1}(x)}\Phi^{-1}] \\
&=& \Phi[\partial, L_{\Phi^{-1}(x)}]\Phi^{-1} \\
&=& \Phi L_{\partial \Phi^{-1}(x)}\Phi^{-1} \\
&=& L_{\Phi\partial\Phi^{-1}(x)} \end{eqnarray*} i.e.
$\Phi\partial\Phi^{-1}\in Der(B).$ Donc l'application \
$H:Der(A)\rightarrow Der(B)$ \ $\partial\mapsto
\Phi\partial\Phi^{-1}$ \ est bien d\'efinie et on v\'erifie
facilement que c'est un isomorphisme d'alg\`ebres de Lie. Si $A$
et $B$ sont quadratiques, alors $\Phi(1)=1$ et si $x\in V,$ on a \
$\Phi(x)^2=\Phi(x^2)=x^2\in K 1.$ Si, de plus, $\Phi(x)\in K 1,$
alors $x$ est colin\'eaire \`a $1,$ car $\Phi$ est injective. Donc
$\Phi(V)\subseteq W,$ et on a $\Phi^{-1}(W)\subseteq V$ i.e.
$\Phi(V)=W.\Box$

\vspace{0.3cm}
\begin{proposition} Soit $(B,\overline{.})=\Big( V,(.|.),\times\Big)$ une $K$-alg\`ebre cayleyenne des
quaternions, de division, $\gamma, \delta\in K,$ avec $\gamma\neq
0$ et soit $A$ l'extension cayleyenne de $B$ d'indice $(\gamma,
\delta).$ Alors les automorphismes de $B$ commute avec
l'involution $\overline{.}.$ Si, de plus, $\delta\neq 0,$ alors
tout automorphisme de $A$ est une extension naturelle d'un
automorphisme de $B$ qui commute l'involution $\overline{.},$ i.e.
\ $Aut(B)\rightarrow Aut(A)$ \ $f\mapsto\Phi_f$ \ est un
isomorphisme de groupes.
\end{proposition}

\vspace{0.1cm} {\bf Preuve.} Les automorphismes de $B$ sont
int\'erieurs d'apr\`es le Th\'eor\`eme {\bf 1.20} de
Sckolem-Noether, de la forme \
$L_{x^{-1}}R_x=L_{\overline{x}}R_x,$ $x\in B$ avec
$x\overline{x}=1.$ Donc commutent avec $\overline{.}$ . On suppose
maintenant $\delta\neq 0.$ Si $\Phi$ est un automorphisme de $A,$
alors $\Phi^{-1}D\Phi\in Der(A)$ pour toute d\'erivation $D$ de
$A,$ et on a \ $\Phi^{-1}D\Phi(0,1)=(0,0)$ i.e.
$D\Phi(0,1)=(0,0).$ Ainsi \ \[ \Phi(0,1)\in\bigcap_{D\in
Der(A)}\ker D=K\times K. \] Comme $\Phi(0,1)$ est un vecteur, il
existe $\alpha\in K$ tel que $\Phi(0,1)=(0,\alpha),$ et on a:
\begin{eqnarray*} (\gamma,0) &=& \Phi(\gamma,0) \\
&=& \Phi\Big( (0,1)^2\Big) \\
&=& \Big( \Phi(0,1)\Big)^2 \\
&=& (0,\alpha)^2 \\
&=& (\gamma\alpha^2,0). \end{eqnarray*} Donc $\Phi(0,1)=\pm(0,1).$

\vspace{0.2cm} \hspace{0.3cm} D'autre part, il existe $f,g\in
End_K(B)$ tels que $\Phi(a,0)=(f(a),g(a))$ pour tout $a\in B.$
Ainsi, pour tous $a,b\in B,$ on a:

\vspace{7cm} \begin{enumerate} \item \begin{eqnarray*} \Phi(a,b) &=& \Phi\Big( (a,0)+(b,0)(0,1)\Big) \\
&=& \Phi(a,0)+\Phi(b,0)\phi(0,1) \\
&=& \Big( f(a)+\varepsilon\gamma g(b), g(a)+\varepsilon f(b)\Big).
\end{eqnarray*}

De plus \item $(f(1),g(1))=\Phi(1,0)=(1,0).$ \item $f(V)\subseteq
V,$ \ $g(B)\subseteq V$ \ (d'apr\`es {\bf 1.}).
\end{enumerate}

\vspace{0.2cm} \hspace{0.3cm} On en d\'eduit que \
$f(\overline{a})=\overline{f(a)}$ \ et \
$g(\overline{a})=\overline{g(a)}=-g(a)$ pour tout $a\in B,$ \`a
l'aide de quoi l'\'egalit\'e \ $\Phi\Big(
-\gamma(\overline{a},0)(b,0)+(0,b)(0,a)\Big)=-\gamma\Phi(\overline{a},0)\Phi(b,0)+\Phi(0,b)\Phi(0,a),$
\ pour $a,b\in B,$ donne \
$\frac{\varepsilon\gamma\delta}{2}g([a,b])=0,$ \ i.e. $g\equiv 0.$
On en d\'eduit de {\bf 1.}, et en tenant compte de du fait que $B$
est de dimension finie, que $f$ est bijective. L'\'egalit\'e \
$\Phi\Big( (0,b)(0,c)\Big)=\Phi(0,b)\Phi(0,c)$ o\`u $b,c\in B,$
donne alors \ $\Big( \gamma
f(\overline{c}b),\frac{\varepsilon\delta}{2}f([c,b])\Big)=\Big(
\gamma\overline{f(c)}f(b),\frac{\delta}{2}[f(c),f(b)] \Big).$ Ce
qui montre que $f$ est un homomorphisme d'alg\`ebres et
$\varepsilon=1.\Box$

\vspace{0.3cm}
\begin{corollary} Soit $\delta$ un nombre r\'eel arbitraire. Alors
\ $Aut(\hit)\subset Aut\Big( E_{-1,\delta}(\hit)\Big).$ Si, de
plus, $\delta\neq 0,$ on a l'\'egalit\'e: $Aut(\hit)=Aut\Big(
E_{-1,\delta}(\hit)\Big).\Box$
\end{corollary}

\vspace{0.3cm} \begin{theorem} Soit $(B,\overline{.})$ une
$K$-alg\`ebre cayleyenne des quaternions, de division, \ $\gamma,
\gamma', \delta, \delta'\in K,$ avec $\gamma\gamma'\neq 0$ tels
que les deux polyn\^omes \ $X^2+\gamma$ et $X^2+\gamma'$
poss\`edent respectivement une racine $\omega$ et $\omega'$ dans
$K.$ Alors \ $E_{\gamma',\delta'}(B)\simeq E_{\gamma,\delta}(B)$
et seulement si $\delta'\omega=\pm\delta\omega'.$
\end{theorem}

\vspace{0.1cm} {\bf Preuve.} On peut supposer $\delta\delta'\neq
0,$ en vertu de la Remarque {\bf 3.18 1)}, et on a \
$E_{\gamma',\delta'}(B)\simeq E_{\gamma,\delta}(B)$ si et
seulement si $E_{-1,\omega'^{-1}\delta'}(B):={\cal A}'\simeq
E_{-1,\omega^{-1}\delta}(B):={\cal A}$ d'apr\`es la Remarque {\bf
III 3.16 3)}. Soit alors $\Phi:{\cal A}\rightarrow{\cal A}'$ un
isomorphisme, alors $\Phi^{-1}D\Phi\in Der({\cal A}),$ pour toute
d\'erivation $D$ de ${\cal A}'.$ Donc $D\Phi(0,1)=(0,0)$ pour
toute d\'erivation $D$ de ${\cal A}'$ i.e. $\Phi(0,1)\in K_times
K.$ Comme $\Phi(0,1)$ est un vecteur unitaire de ${\cal A}',$ on
a: \ $\Phi(0,1)=\pm(0,1).$ un m\^eme raisonnement que celui dans
la Proposition {\bf 3.26} montre qu'il existe $f\in End_(B),$ qui
commute avec $\overline{.},$ tel que \ $\Phi(a,b)=(f(a),\pm f(b))$
pour tout $a,b\in B.$ L'\'egalit\'e $\Phi\Big(
(0,b)(0,c)\Big)=\Phi(0,b)\Phi(0,c),$ o\`u $b,c\in B,$ donne alors
$\delta'\omega=\pm\delta\omega'.\Box$

\vspace{0.3cm} \begin{corollary} $E_{-1,\delta'}(\hit)\simeq$
$E_{-1,\delta}(\hit)$ si et seulement si \
$\delta'=\pm\delta.\Box$
\end{corollary}

\vspace{1cm}
\section{Classification des $\rit$-alg\`ebres de Jordan n.c. de division lin\'eaire de dimension $8$}

\vspace{0.6cm} \subsection{Isotopie vectorielle}

\vspace{0.4cm} \hspace{0.3cm} Soient $(V,\wedge)$ une
$\rit$-alg\`ebre anti-commutative de dimension $\geq 1,$ $(.|.)$
une forme bilin\'eaire sym\'etrique d\'efinie n\'egative sur $V$
et $\varphi$ un automorphisme de l'espace vectoriel $V.$ On pose \
$x\Delta y=\varphi^*\Big( \varphi(x)\wedge\varphi(y)\Big),$
$x,y\in V,$ \ $\varphi^*$ \'etant l'automorphisme adjoint de
$\varphi,$ et on d\'esigne par \ $\Big( V,(.|.),\wedge\Big),$ et
$\Big( V,(.|.),\Delta\Big)$ \ respectivement, les alg\`ebres
cayleyennes construites \`a partir des alg\`ebres
anti-commutatives $(V,\wedge)$ et $(V,\Delta),$ et de la forme
bilin\'eaire sym\'etrique $(.|.).$ Nous avons le r\'esultat cl\'e
suivant:

\vspace{0.3cm} \begin{proposition} $\Big( V,(.|.),\wedge\Big)$ est
flexible de division lin\'eaire si et seulement si \ $\Big(
V,(.|.),\Delta\Big)$ est flexible de division lin\'eaire.
\end{proposition}

\vspace{0.1cm} {\bf Preuve.} Si $\Big( V,(.|.),\Delta\Big)$ est
flexible, alors $(.|.)$ est une forme trace sur $V,$ i.e. $\Big(
V,(.|.),\wedge\Big)$ est flexible. Si, de plus, $\Big(
V,(.|.),\Delta\Big)$ est de division lin\'eaire, on v\'erifie
facilement que la propri\'et\'e {\bf 4.} du Th\'eor\`eme {\bf
1.49} a lieu dans $\Big( V,(.|.),\Delta\Big)$ aussi bien que dans
$\Big( V,(.|.),\wedge\Big).$ L'implication 2) $\Rightarrow$ 1)
s'\'etablit de la m\^eme fa\c{c}on.$\Box$

\vspace{0.3cm} \begin{definition and remarks} .
\begin{enumerate} \item On dira que l'alg\`ebre $\Big( V,(.|.),\Delta\Big)$ est obtenue, \`a partir de l'alg\`ebre
$A=\Big( V,(.|.),\wedge\Big)$ et de l'automorphisme $\varphi,$ par
isotopie vectorielle. On la note $A(\varphi).$ Cette notion
d'isotopie vectorielle est une relation d'\'equivalence dans la
classe ${\cal C}_n,$ $n\geq 2$ des alg\`ebres r\'eelles
cayleyennes de dimension $n,$ dont l'espace r\'eel $V=\rit^{n-1}$
des vecteurs associ\'e, est muni d'une m\^eme forme bilin\'eaire
sym\'etrique $(.|.)$ d\'efinie n\'egative. En outre, on v\'erifie
facilement que pour toute alg\`ebre $B,$ de cette classe, et pour
deux automorphismes $\varphi$ et $\psi$ de l'espace des vecteurs
assoc\'e \`a $B,$ on a
\begin{eqnarray} \Big( B(\varphi)\Big)(\psi) &=&
B(\varphi\psi). \end{eqnarray}

\item Les alg\`ebres r\'eelles de Jordan non commutatives de
division lin\'eaire de dimension finie $n=2$ et $4$ s'obtiennent,
respectivement, \`a partir de $\cit$ et $\hit,$ par isotopie
vectorielle. Plus pr\'ecis\'ement, de telles alg\`ebres coincident
avec $\cit$ pour $n=2$ et sont les mutations $\hit^{(\lambda)}$ de
$\hit$ o\`u $\lambda\in\rit-\{\frac{1}{2}\}.$ L'alg\`ebre
$\hit^{(\lambda)}$ n'est autre que l'isotope vectorielle
$\hit(\varphi)$ o\`u $\varphi$ est l'homoth\'etie
$(2\lambda-1)^{\frac{1}{3}}I_V$ de l'espace $V$ des vecteurs
associ\'e \`a l'alg\`ebre r\'eelle $\hit.\Box$
\end{enumerate}
\end{definition and remarks}

\vspace{0.2cm} \begin{example} Soient $\oit=\Big(
V,(.|.),\times\Big)$ l'alg\`ebre r\'eelle de Cayley-Dickson de
division, $\lambda, \mu\in\rit-\{\frac{1}{2}\}$ et $\varphi$
l'automorphisme de $V$ dont la matrice, par rapport \`a la base
canonique $\{e_1,\dots,e_7\}$ de $V$ est donn\'ee par:

\vspace{0.5cm}
\[ \left(
\begin{array}{lllllll}
\lambda' &          &          &      &          &          &          \\
         & \lambda' &          &      &          &          &          \\
         &          & \lambda' &      &          &          &          \\
         &          &          & \mu' &          &          &          \\
         &          &          &      & \lambda' &          &          \\
         &          &          &      &          & \lambda' &          \\
         &          &          &      &          &          & \lambda'
\end{array}
\right)  \]

\vspace{0.5cm} o\`u
$\lambda'=(2\mu-1)^{\frac{1}{3}}(2\lambda-1)^{\frac{1}{3}}$ \ et \
$\mu'=(2\mu-1)^{\frac{1}{3}}(2\lambda-1)^{-\frac{2}{3}}.$ Alors \
$\oit(\varphi)=\Big( E_{-1}(\hit^{(\lambda)})\Big)^{(\mu)}.\Box$
\end{example}

\vspace{0.3cm} \hspace{0.3cm} Nous donnons maintenant un
r\'esultat important, qui constitue \'egalement un exemple:

\vspace{0.3cm} \begin{proposition} Les alg\`ebres r\'eelles \
$E_{-1,\delta}(\hit)$ o\`u $0\leq\delta<2$ s'obtiennent, \`a
partir de l'alg\`ebre r\'eelle \ $\oit=\Big(V,(.|.),\times\Big)$
de Cayley-Dickson, par isotopie vectorielle. \end{proposition}

\vspace{0.1cm} {\bf Preuve.} Soit $\alpha\in ]\frac{1}{2},1],$
alors l'endomorphisme $\varphi$ de $V$ dont la matrice, par
rapport \`a la base canonique ${\cal B}=\{e_1,\dots,e_7\}$ de $V$
est donn\'ee par:

\vspace{0.5cm}
\[ \left(
\begin{array}{lllllll}
(2\alpha-1)^{-1} &                  &                  & 0 & (1-\alpha^2)^{\frac{1}{2}} &               &          \\
                 & (2\alpha-1)^{-1} &                  & 0 &                            &
(1-\alpha^2)^{\frac{1}{2}} &          \\
                 &                  & (2\alpha-1)^{-1} & 0 &                            &               &
(1-\alpha^2)^{\frac{1}{2}} \\
\hspace{0.6cm} 0 & \hspace{0.6cm} 0 & \hspace{0.6cm} 0 & 1 & \hspace{0.6cm} 0 & \hspace{0.6cm} 0 & \hspace{0.6cm} 0 \\
\hspace{0.6cm} \alpha'          &                  &                  & 0 & \hspace{0.6cm} \alpha &     &          \\
                 & \hspace{0.6cm} \alpha' &         & 0 &                            & \hspace{0.6cm} \alpha &      \\
                 &                  & \hspace{0.6cm} \alpha'           & 0 &    &    & \hspace{0.6cm} \alpha
\end{array}
\right)  \]

\vspace{0.5cm} o\`u
$\alpha=(1-\alpha^2)^{\frac{1}{2}}(\alpha+1)^{-1}(2\alpha-1)^{-1},$
est un automorphisme. On pose alors $A=\oit(\varphi)$ et on
obtient, par rapport \`a la base ${\cal B},$ la table:

\vspace{0.3cm}
\[ \begin{tabular}{ccccccccc} \\
\multicolumn{1}{c}{} & \multicolumn{1}{c}{$1$} &
\multicolumn{1}{c}{$e_1$} & \multicolumn{1}{c}{$e_2$} &
\multicolumn{1}{c}{$e_3$} & \multicolumn{1}{c}{$e_4$} &
\multicolumn{1}{c}{$e_5$} & \multicolumn{1}{c}{$e_6$} &
\multicolumn{1}{c}{$e_7$}
\\ \cline{2-9}
\multicolumn{1}{c|}{$1$} & \multicolumn{1}{|c}{$1$} &
\multicolumn{1}{c|}{$e_1$} & \multicolumn{1}{|c}{$e_2$} &
\multicolumn{1}{c|}{$e_3$} & \multicolumn{1}{|c}{$e_4$} &
\multicolumn{1}{c|}{$e_5$} & \multicolumn{1}{|c}{$e_6$} &
\multicolumn{1}{c|}{$e_7$}
\\ \multicolumn{1}{c|}{$e_1$} & \multicolumn{1}{|c}{} &
\multicolumn{1}{c|}{$-1$} & \multicolumn{1}{|c}{$\beta e_3$} &
\multicolumn{1}{c|}{$-\beta e_2$} & \multicolumn{1}{|c}{$e_5$} &
\multicolumn{1}{c|}{$-e_4$} & \multicolumn{1}{|c}{$-e_7$} &
\multicolumn{1}{c|}{$e_6$} \\ \cline{2-9}
\multicolumn{1}{c|}{$e_2$} & \multicolumn{1}{|c}{} &
\multicolumn{1}{c|}{} & \multicolumn{1}{|c}{$-1$} &
\multicolumn{1}{c|}{$\beta e_1$} & \multicolumn{1}{|c}{$e_6$} &
\multicolumn{1}{c|}{$e_7$} & \multicolumn{1}{|c}{$-e_4$} &
\multicolumn{1}{c|}{$-e_5$} \\
\multicolumn{1}{c|}{$e_3$} & \multicolumn{1}{|c}{} &
\multicolumn{1}{c|}{} & \multicolumn{1}{|c}{} &
\multicolumn{1}{c|}{$-1$} & \multicolumn{1}{|c}{$e_7$} &
\multicolumn{1}{c|}{$-e_6$} & \multicolumn{1}{|c}{$e_5$} &
\multicolumn{1}{c|}{$-e_4$} \\ \cline{4-9}
\multicolumn{1}{c|}{$e_4$} & \multicolumn{1}{|c}{} &
\multicolumn{1}{c}{} & \multicolumn{1}{c}{} &
\multicolumn{1}{c|}{} & \multicolumn{1}{|c}{$-1$} &
\multicolumn{1}{c|}{$e_1$} & \multicolumn{1}{|c}{$e_2$} &
\multicolumn{1}{c|}{$e_3$} \\
\multicolumn{1}{c|}{$e_5$} & \multicolumn{1}{|c}{} &
\multicolumn{1}{c}{} & \multicolumn{1}{c}{} & \multicolumn{1}{c}{}
& \multicolumn{1}{|c}{} & \multicolumn{1}{c|}{$-1$} &
\multicolumn{1}{|c}{$-e_3-\lambda e_7$} &
\multicolumn{1}{c|}{$e_2+\lambda e_6$} \\ \cline{6-9}
\multicolumn{1}{c|}{$e_6$} & \multicolumn{1}{|c}{} &
\multicolumn{1}{c}{} & \multicolumn{1}{c}{} & \multicolumn{1}{c}{}
& \multicolumn{1}{c}{} & \multicolumn{1}{c}{} &
\multicolumn{1}{|c}{$-1$} & \multicolumn{1}{c|}{$-e_1-\lambda
e_5$} \\
\multicolumn{1}{c|}{$e_7$} & \multicolumn{1}{|c}{} &
\multicolumn{1}{c}{} & \multicolumn{1}{c}{} & \multicolumn{1}{c}{}
& \multicolumn{1}{c}{} & \multicolumn{1}{c}{} &
\multicolumn{1}{|c}{} & \multicolumn{1}{c|}{$-1$} \\ \cline{2-9}
\end{tabular} \]

\vspace{0.4cm} o\`u $\beta=2(\alpha+1)^{-1}(2\alpha-1)^{-2}>0$ et
$\lambda=(4\alpha^2-1)(1-\alpha^2)^{\frac{1}{2}}\geq 0.$ On
consid\`ere enfin l'automorphisme $\psi$ de $V$ dont la matrice,
par rapport \`a la base ${\cal B}$ est donn\'ee par:

\vspace{0.5cm}
\[ \left(
\begin{array}{lllllll}
\beta^{-\frac{1}{3}} &                      &                  &     &     &               &          \\
         & \beta^{-\frac{1}{3}} &                      &     &     &               &          \\
         &                      & \beta^{-\frac{1}{3}} &     &     &               &          \\
 &       &                      &                      \beta^{\frac{1}{6}} & &       &    \\
 &       &                      &                      & \beta^{\frac{1}{6}} &       &    \\
 &       &                      &                      &                      & \beta^{\frac{1}{6}} &       \\
 &       &                      &  &                   &   & \beta^{\frac{1}{6}}
\end{array}
\right)  \]

\vspace{0.5cm} La table de multiplication de l'alg\`ebre
$A(\psi)=\oit(\varphi\psi),$ par rapport \`a la base ${\cal B},$
est la m\^eme que celle de l'alg\`ebre \
$E_{-1,\delta_{\alpha}}(\hit)$ o\`u
$\delta_\alpha=\lambda\beta^{\frac{1}{2}}\geq 0.$ De plus
$\delta_\alpha^2=2(1-\alpha)(2\alpha+1)^2$ est une fonction, en
$\alpha,$ continue et strictement d\'ecroissante sur
$]\frac{1}{2},1].$ Elle atteint, en vertu du Th\'eor\`eme des
valeurs interm\'ediaires, toutes les valeurs de l'intervalle \[
[\delta_1^2,\lim_{\alpha\rightarrow\frac{1}{2}^+}\delta_\alpha^2[=[0,4[.\Box
\]

\vspace{5cm} \subsection{Probl\`emes d'isomorphisme}

\vspace{0.5cm}
\begin{note} Soit $A=\Big( V,(.|.),\wedge\Big)$ est une $\rit$-alg\`ebre quadratique et soit $f$ un automorphisme
de l'espace vectoriel r\'eel $V,$ on d\'esigne par \ $\tilde{f}$
l'application \ $\alpha+u\mapsto\alpha+f(u)$ \ $A\rightarrow A$ \
qui est un automorphisme de l'espace vectoriel r\'eel $A.$ On
l'appelera prolongement naturel de $f$ \`a $A.\Box$
\end{note}

\vspace{0.3cm}
\begin{proposition} Soient $A=\Big( V,(.|.),\Delta\Big)$ et $B=\Big( V,(.|.),\wedge\Big)$ deux
alg\`ebres r\'eelles cayleyennes de la classe ${\cal C}_n$ et soit
$f$ un automorphisme de l'espace vectoriel r\'eel $V.$ Alors
\begin{enumerate}
\item Les alg\`ebres $A$ et $B$ sont isomorphes si et seulement si
elle sont isom\'etriquement vectoriellement isotopes. \item
$\tilde{f}\in Aut(A)$ si et seulement si $A(f)=A$ et $f$ est une
isom\'etrie de l'espace euclidien $(V,-(.|.)).$ \end{enumerate}
\end{proposition}

\vspace{0.1cm} {\bf Preuve.} . \begin{enumerate} \item Soit
$g:A\rightarrow B$ un isomorphisme d'alg\`ebres et soient $x,y\in
A.$ On a: \begin{eqnarray*} (x|y)+g(x\Delta y) &=& g\Big(
(x|y)+x\Delta y\Big) \\
&=& g(xy) \\
&=& g(x)g(y) \\
&=& (g(x)|g(y))+g(x)\wedge g(y). \end{eqnarray*}

Donc $((g(x)|g(y))=(x|y)$ et $g(x\Delta y)=g(x)\wedge g(y)$ i.e.
$g$ est une isom\'etrie de l'espace vectoriel $A$ qui est
prolongement naturel d'une isom\'etrie $g_0$ de l'espace euclidien
$(V,-(.|.)).$ De plus, pour tous $x,y\in V:$ \begin{eqnarray*}
x\Delta y &=& g_0^{-1}(g_0(x)\wedge g_0(y)) \\
&=& g_0^*(g_0(x)\wedge g_0(y)) \end{eqnarray*} i.e. $A=B(g_0).$

\vspace{0.2cm} \hspace{0.3cm} R\'eciproquement, si $A=B(h)$ o\`u
$h$ est une isom\'etrie de l'espace euclidien $(V,-(.|.)),$ alors
le prolongement naturel $\tilde{h},$ de $h$ \`a $\rit\oplus V,$
est un isomorphisme de l'alg\`ebre $A$ dans l'alg\`ebre $B.$
\item La proposition {\bf 2.} est cons\'equence de {\bf 1.}.$\Box$
\end{enumerate}

\vspace{0.3cm}
\begin{lemma} L'application \ $f:E_{-1,-\delta}(\hit)\rightarrow E_{-1,\delta}(\hit)$ \ $(x,y)\mapsto(x,-y)$ \ o\`u
$\delta\in\rit$ est un isomorphisme d'alg\`ebres. En particulier,
$E_{-1,-\delta}(\hit)$ et $E_{-1,\delta}(\hit)$ sont
vectoriellement isotopes.
\end{lemma}

\vspace{0.1cm} {\bf Preuve.} Soient $x,y,x',y'\in \hit,$ on a
\begin{eqnarray*} f\Big( (x,y)(x',y')\Big) &=& f\Big( xx'-\overline{y'}y,
y\overline{x'}+y'x-\frac{\delta}{2}[y',y]\Big) \\
&=& \Big( xx'-(\overline{-y'})(-y),
(-y)\overline{x'}+(-y')x+\frac{\delta}{2}[(-y'),(-y)]\Big) \\
&=& (x,-y)(x',-y') \\
&=& f(x,y)f(x',y').\Box \end{eqnarray*}

\vspace{0.3cm}
\begin{corollary} Les alg\`ebres r\'eelles \ $E_{-1,\delta}(\hit)$ \ o\`u
$|\delta|<2$ s'obtiennent, \`a partir de l'alg\`ebre r\'eelle
$\oit$ de Cayley-Dickson, par isotopie vectorielle.
\end{corollary}

\vspace{0.1cm} {\bf Preuve.} Le Lemme {\bf 4.7} et l'\'egalit\'e
{\bf (4.3)} dans ({\bf 4.2 1)}) permettent de supposer $\delta\geq
0.$ La Proposition {\bf 4.4} ach\`eve alors la
d\'emonstration.$\Box$

\vspace{0.3cm} \begin{proposition} Soit \ $\oit=\Big(
V,(.|.),\times\Big)$ l'alg\`ebre r\'eelle de Cayley-Dickson et
soit $\varphi$ un automorphisme de l'espace vectoriel $V.$ Alors
les deux propri\'et\'es suivantes sont \'equivalentes:
\begin{enumerate} \item $\oit(\varphi)=\oit.$ \item $\tilde{\varphi}\in G_2.$
\end{enumerate}
\end{proposition}

\vspace{0.1cm} {\bf Preuve.} L'implication 2) $\Rightarrow$ 1) est
cons\'equence de la Proposition {\bf 4.6 2)}.

\vspace{0.2cm} 1) $\Rightarrow$ 2). On note \ $\oit(\varphi)=\Big(
V,(.|.),\wedge\Big)$ et, pour tous $x,y\in V,$ on a \[
\varphi^*(\varphi(x)\times \varphi(y))=x\wedge y=x\times y. \] Si
$x_1,\dots,x_7$ est une base orthonorm\'ee de $V$ form\'ee de
vecteurs propres de $\varphi\varphi^*$ associ\'es, respectivement,
aux valeurs propres $\lambda_1,\dots,\lambda_7,$ alors:

\begin{enumerate} \item $\lambda_i>0.$ \item $(\varphi^*(x_i)|\varphi^*(x_j))=-\lambda_i\delta_{ij}$ o\`u $\delta_{ij}$
est le symbole de Kronecker. \item
$\varphi^*(x_i)\times\varphi^*(x_j)=\lambda_i\lambda_j\varphi^*(x_i\times
x_j).$ \end{enumerate}

On suppose $i\neq j$ et on a \begin{eqnarray*} \varphi^*\Big(
\lambda_ix_i\times\varphi\varphi^*(x_i\times x_j)\Big) &=&
\varphi^*(x_i)\wedge\varphi^*(x_i\times x_j) \\
&=& \varphi^*(x_i)\times\varphi^*(x_i\times x_j) \\
&=& (\lambda_i\lambda_j)^{-1}\varphi^*(x_i)\times\Big(
\varphi^*(x_i)\times\varphi^*(x_j)\Big) \ \mbox{ d'apr\`es
l'\'egalit\'e {\bf 3.} ci-dessus }
\\ &=& (\lambda_i\lambda_j)^{-1}\Big(
\varphi^*(x_i)\Big)^2\varphi^*(x_j) \\
&=& -\lambda_j^{-1}\varphi^*(x_j) \ \mbox{ d'apr\`es l'\'egalit\'e
{\bf 2.} ci-dessus } \\
&=& \varphi^*(-\lambda_j^{-1}x_j) \\ &=& \varphi^*\Big(
\lambda_ix_i\times(\lambda_i\lambda_j)^{-1}(x_i\times x_j)\Big).
\end{eqnarray*}

\vspace{0.2cm} Donc \
$\lambda_ix_i\times\varphi\varphi^*(x_i\times
x_j)=\lambda_ix_i\times(\lambda_i\lambda_j)^{-1}(x_i\times x_j)$
et on a: \begin{eqnarray} \varphi\varphi^*(x_i\times x_j) &=&
(\lambda_i\lambda_j)^{-1}x_i\times x_j. \end{eqnarray} Par
cons\'equent \ $x_i, x_i\times x_1,\dots,x_1\times x_7$ \ est une
base orthonorm\'ee de $V$ form\'ee de vecteurs propres de
$\varphi\varphi^*$ associ\'es, respectivement, aux valeurs propres
\ $\lambda_i,$
$(\lambda_i\lambda_1)^{-1},\dots,(\lambda_i\lambda_7)^{-1}.$ Si
$k\neq i,j$ l'\'egalit\'e {\bf (4.4)} montre que
\begin{eqnarray}
\varphi\varphi^*\Big( (x_i\times x_k)\times(x_i\times x_j)\Big)
&=& \Big(
(\lambda_i\lambda_k)^{-1}(\lambda_i\lambda_j)^{-1}\Big)^{-1}(x_i\times
x_k)\times(x_i\times x_j) \end{eqnarray}

et l'on distingue les deux cas suivants:

\vspace{0.4cm} {\bf Premier cas.} Si les vecteurs propres $x_i$ et
$x_j\times x_k$ ne sont pas orthogonaux, donc associ\'es \`a la
m\^eme valeur propre, on a $\lambda_i=(\lambda_j\lambda_k)^{-1}.$

\vspace{0.2cm} {\bf Deuxi\`eme cas.} Si $x_i$ et $x_j\times x_k$
sont orthogonaux, on a:
\begin{eqnarray*} \Big( (x_i\times x_k)\times(x_j\times x_j)\Big)\mid x_k\times x_j\Big) &=&
-\Big( (x_ix_k)(x_jx_i)\mid x_kx_j \Big) \ (x_i\times x_k, x_i\times x_j \mbox{ orthogonaux }) \\
&=& -\Big( (x_i(x_kx_j)x_i\mid x_kx_j \Big) \ \mbox{ Identit\'e moyenne de Moufang) } \\
&=& -\Big( (x_i(x_kx_j)\mid x_i(x_kx_j) \Big) \\
&=& -\Big( x_i\times(x_k\times x_j)\Big)^2 \ (x_i, x_k\times x_j
\mbox{ orthogonaux }). \end{eqnarray*}

\vspace{0.2cm} Donc les deux vecteurs propres \ $(x_i\times
x_k)\times(x_j\times x_j)$ {\bf (4.5)} et $x_k\times x_j$ ne sont
pas orthogonaux, donc associ\'es \`a la m\^eme valeur propre, et
on a \ $(\lambda_i\lambda_k)^{-1}(\lambda_i\lambda_j)^{-1}=
(\lambda_j\lambda_k)^{-1}.$

\vspace{0.2cm} \hspace{0.3cm} Dans les deux cas
$(\lambda_i\lambda_j\lambda_k)^2=1$ et comme $\lambda_i>0,$
d'apr\`es l'\'egalit\'e {\bf 1.}, on a $\varphi\varphi^*=I_V.$
Ainsi $\tilde{\varphi}\in G_2.\Box$

\vspace{3cm} \begin{corollary} Soit \ $\oit$ l'alg\`ebre r\'eelle
de Cayley-Dickson et soient $\varphi,\psi$ deux automorphismes de
l'espace vectoriel r\'eel $\rit^7.$ Alors:
\begin{enumerate} \item $\oit(\varphi)\simeq$ $\oit$ si et seulement si \ $\varphi\in O_7(\rit)$
(le groupe des isom\'etries de l'espace euclidien $\rit^7$). \item
$\oit(\varphi)=\oit(\psi)$ si et seulement si \
$\overline{\varphi\psi^{-1}}\in G_2.$ \item $\oit(\varphi)\simeq$
$\oit(\psi)$ si et seulement si il existe $f\in O_7(\rit)$ tel que
\ $\overline{\varphi f\psi^{-1}}\in G_2.$
\end{enumerate}
\end{corollary}

\vspace{0.1cm} {\bf Preuve.} \begin{enumerate} \item $\Leftarrow/$
\ $\tilde{\varphi}:\oit(\varphi)\rightarrow\oit$ \ est un
isomorphisme.

\vspace{0.2cm} $\Rightarrow/$ \ Il existe $f\in O_7(\rit)$ tel que
$\oit(\varphi f)=\oit$ et on a $\varphi f\in O_7(\rit).$ Donc
$\varphi=(\varphi f)f^{-1}\in O_7(\rit).$ \item On a
\begin{eqnarray*} \overline{\varphi\psi^{-1}}\in G_2 &\Leftrightarrow& \oit(\varphi\psi^{-1})=\oit \
\mbox{ Proposition {\bf 4.9} } \\
&\Leftrightarrow& \oit(\varphi)=\Big(
\oit(\varphi\psi^{-1})\Big)(\psi)=\oit(\psi). \end{eqnarray*}
\item On a \begin{eqnarray*} \oit(\varphi)\simeq\oit(\psi)
&\Leftrightarrow& \mbox{ Il existe } f\in O_7(\rit): \
\oit(\varphi)=\Big( \oit(\psi)\Big)(f) \ \mbox{ Proposition {\bf 4.6 1)} } \\
&\Leftrightarrow& \mbox{ Il existe } f\in O_7(\rit): \
\oit(\varphi)=\oit(\psi f) \\
&\Leftrightarrow& \mbox{ Il existe } f\in O_7(\rit): \
\oit(\varphi)=\oit(\psi f^{-1}) \\
&\Leftrightarrow& \mbox{ Il existe } f\in O_7(\rit): \
\overline{\varphi f\psi^{-1}}\in G_2 \ \mbox{ (proposition {\bf
2.} pr\'ec\'edente) }.\Box
\end{eqnarray*}
\end{enumerate}

\vspace{7cm} \subsection{Th\'eor\`eme de classification}

\vspace{0.5cm} \hspace{0.3cm} Soient $A=\Big( W,(.|.),\wedge\Big)$
une $\rit$-alg\`ebre de Jordan n.c. de division lin\'eaire de
dimension $8$ et soit $B$ une sous-alg\`ebre de $A$ de dimension
$4.$ Il existe une base orthonorm\'ee \ $\{1,u,y_1,z_1\}$ \ de $B$
que l'on peut compl\'eter, selon la note {\bf 3.3}, en une base
orthonorm\'ee \ ${\cal B}=\{1,u,y_1,z_1,y_2,z_2,y_3,z_3\}$ \ de
$A$ telle que \[ u\wedge y_i=a_iz_i, \hspace{0.2cm} u\wedge
z_i=-a_iy_i, i=1,2,3 \ \mbox{ et } \ y_1\wedge z_1=a_1u \] o\`u
les $a_i$ sont des param\`etres $>0.$ Il existe \'egalement des
param\`etres \[ \theta_{ij}, \ \omega_{ij}, \ \pi_{ij}, \
\alpha_{ijk}, \ \beta_{ijk}, \ \gamma_{ijk}, \ \eta_{ijk}, \
\lambda_{ijk}, \ \mu_{ijk} \] o\`u $i,j,k\in\{1,2,3\}$ tels que

\vspace{0.2cm} \begin{eqnarray*} y_i\wedge y_j &=&
\theta_{ij}u+\sum_{1\leq k\leq 3}(\alpha_{ijk}y_k+\beta_{ijk}z_k) \\
y_i\wedge z_j &=& \omega_{ij}u+\sum_{1\leq k\leq
3}(\gamma_{ijk}y_k+\eta_{ijk}z_k) \\
z_i\wedge z_j &=& \pi_{ij}u+\sum_{1\leq k\leq
3}(\lambda_{ijk}y_k+\mu_{ijk}z_k) \end{eqnarray*}

\vspace{0.2cm} La propri\'et\'e trace de $(.|.)$ et
l'anti-commutativit\'e de "$\wedge$" montrent que $\alpha_{ijk}$
et $\mu_{ijk}$ sont altern\'es en $i,j$ et $k,$
$\beta_{ijk}=-\beta_{jik}=\gamma_{jki},$
$\eta_{ijk}=-\eta_{ikj}=\lambda_{jki},$
$\theta_{ij}=\pi_{ij}=\beta_{i11}=\gamma_{11i}=\eta_{11i}=0,$
$\omega_{ii}=a_i$ et $\omega_{ij}=0$ si $i\neq j.$ Ainsi les
param\`etres pr\'ec\'edents se r\'eduisent aux suivants: $a_1,$
$a_2,$ $a_3,$ $\alpha_{123},$ $\beta_{123},$ $-\beta_{132},$
$\eta_{123},$ $\beta_{231},$ $-\eta_{213},$ $\eta_{312},$
$\mu_{123},$ $\beta_{133},$ $-\eta_{212},$ $-\eta_{313},$
$\eta_{223},$ $-\eta_{323},$ $\beta_{122},$ $-\beta_{232},$
$\beta_{233},$ \ qu'on note respectivement \[ a, \ b, \ c, \
\alpha, \ \beta, \ \gamma, \ \mu, \ \lambda, \ \eta, \ \sigma, \
\delta, \ \nu, \ \pi, \ \rho, \ \theta, \ \omega, \ \pi_*, \
\theta_*, \ \omega_*.
\] On obtient, par rapport \`a la base ${\cal B},$ une premi\`ere
table de multiplication de $A,$ et on se limite \`a la partie
triangulaire sup\'erieure, vu l'anti-commutativit\'e de
"$\wedge$":

\vspace{0.2cm}
\[ \begin{tabular}{ccccccccc} \\
\multicolumn{1}{c}{} & \multicolumn{1}{c}{$1$} &
\multicolumn{1}{c}{$u$} & \multicolumn{1}{c}{$y_1$} &
\multicolumn{1}{c}{$z_1$} & \multicolumn{1}{c}{$y_2$} &
\multicolumn{1}{c}{$z_2$} & \multicolumn{1}{c}{$y_3$} &
\multicolumn{1}{c}{$z_3$}
\\ \cline{2-9}
\multicolumn{1}{c|}{$1$} & \multicolumn{1}{|c}{$1$} &
\multicolumn{1}{c|}{$u$} & \multicolumn{1}{|c}{$y_1$} &
\multicolumn{1}{c|}{$z_1$} & \multicolumn{1}{|c}{$y_2$} &
\multicolumn{1}{c|}{$z_2$} & \multicolumn{1}{|c}{$y_3$} &
\multicolumn{1}{c|}{$z_3$}
\\
\multicolumn{1}{c|}{$u$} & \multicolumn{1}{|c}{} &
\multicolumn{1}{c|}{$-1$} & \multicolumn{1}{|c}{$az_1$} &
\multicolumn{1}{c|}{$-ay_1$} & \multicolumn{1}{|c}{$bz_2$} &
\multicolumn{1}{c|}{$-by_2$} & \multicolumn{1}{|c}{$cz_3$} &
\multicolumn{1}{c|}{$-cy_3$}
\\ \cline{2-9}
\multicolumn{1}{c|}{$y_1$} & \multicolumn{1}{|c}{} &
\multicolumn{1}{c|}{} & \multicolumn{1}{|c}{$-1$} &
\multicolumn{1}{c|}{$au$} & \multicolumn{1}{|c}{$$} &
\multicolumn{1}{c|}{$$} & \multicolumn{1}{|c}{$-\alpha y_2-\gamma
z_2+\nu z_3$} & \multicolumn{1}{c|}{$-\beta y_2-\mu z_2-\nu y_3$}
\\
\multicolumn{1}{c|}{$z_1$} & \multicolumn{1}{|c}{} &
\multicolumn{1}{c|}{} & \multicolumn{1}{|c}{} &
\multicolumn{1}{c|}{$-1$} & \multicolumn{1}{|c}{$$} &
\multicolumn{1}{c|}{$$} & \multicolumn{1}{|c}{$-\lambda y_2-\sigma
z_2+\rho z_3$} & \multicolumn{1}{c|}{$-\eta y_2-\delta z_2-\rho
z_3$}
\\ \cline{4-9}
\multicolumn{1}{c|}{$y_2$} & \multicolumn{1}{|c}{} &
\multicolumn{1}{c}{} & \multicolumn{1}{c}{} &
\multicolumn{1}{c|}{} & \multicolumn{1}{|c}{$-1$} &
\multicolumn{1}{c|}{$$} & \multicolumn{1}{|c}{$$} &
\multicolumn{1}{c|}{$$}
\\
\multicolumn{1}{c|}{$z_2$} & \multicolumn{1}{|c}{} &
\multicolumn{1}{c}{} & \multicolumn{1}{c}{} & \multicolumn{1}{c}{}
& \multicolumn{1}{|c}{} & \multicolumn{1}{c|}{$-1$} &
\multicolumn{1}{|c}{$$} & \multicolumn{1}{c|}{$$}
\\ \cline{6-9}
\multicolumn{1}{c|}{$y_3$} & \multicolumn{1}{|c}{} &
\multicolumn{1}{c}{} & \multicolumn{1}{c}{} & \multicolumn{1}{c}{}
& \multicolumn{1}{c}{} & \multicolumn{1}{c}{} &
\multicolumn{1}{|c}{$-1$} & \multicolumn{1}{c|}{$$}
\\
\multicolumn{1}{c|}{$z_3$} & \multicolumn{1}{|c}{} &
\multicolumn{1}{c}{} & \multicolumn{1}{c}{} & \multicolumn{1}{c}{}
& \multicolumn{1}{c}{} & \multicolumn{1}{c}{} &
\multicolumn{1}{|c}{} & \multicolumn{1}{c|}{$-1$}
\\ \cline{2-9}
\end{tabular} \]

\vspace{0.1cm}
{\small \[ \begin{tabular}{ccccc} \\
\multicolumn{1}{c}{} & \multicolumn{1}{c}{$y_2$} &
\multicolumn{1}{c}{$z_2$} & \multicolumn{1}{c}{$y_3$} & \multicolumn{1}{c}{$z_3$} \\
\cline{2-5} \multicolumn{1}{c|}{$y_1$} &
\multicolumn{1}{|c}{$\pi_*z_2+\alpha y_3+\beta z_3$} & \multicolumn{1}{c|}{$-\pi_*y_2+\gamma y_3+\mu z_3$}
& \multicolumn{1}{|c}{} & \multicolumn{1}{c|}{} \\
\multicolumn{1}{c|}{$z_1$} &
\multicolumn{1}{|c}{$\pi z_2+\lambda y_3+\eta z_3$} & \multicolumn{1}{c|}{$-\pi y_2+\sigma y_3+\delta z_3$}
& \multicolumn{1}{|c}{} & \multicolumn{1}{c|}{} \\
\cline{2-5} \multicolumn{1}{c|}{$y_2$} & \multicolumn{1}{|c}{} &
\multicolumn{1}{c|}{} & \multicolumn{1}{|c}{$\alpha y_1+\lambda
z_1-\theta_*z_2+\omega_*z_3$}
& \multicolumn{1}{c|}{$\beta y_1+\eta z_1-\theta z_2-\omega_*y_3$} \\
\multicolumn{1}{c|}{$z_2$} & \multicolumn{1}{|c}{} &
\multicolumn{1}{c|}{} & \multicolumn{1}{|c}{$\gamma y_1+\sigma
z_1+\theta_*y_2+\omega z_3$} & \multicolumn{1}{c|}{$\mu y_1+\delta z_1+\theta y_2-\omega z_2$} \\
\cline{2-5}
\end{tabular} \]}

\[ \begin{tabular}{ccc} \\
\multicolumn{1}{c}{} & \multicolumn{1}{c}{$z_2$} &
\multicolumn{1}{c}{$z_3$} \\ \cline{2-3}
\multicolumn{1}{c|}{$y_2$} & \multicolumn{1}{|c}{$bu+\pi_*y_1+\pi
z_1+\theta_*y_3+\theta z_3$} &
\multicolumn{1}{c|}{} \\
\multicolumn{1}{c|}{$y_3$} & \multicolumn{1}{|c}{} &
\multicolumn{1}{c|}{$cu+\nu y_1+\rho z_1+\omega_*y_2+\omega z_2$} \\
\cline{2-3}
\end{tabular} \]

\vspace{0.3cm} On note $H_i$ le sous-espace vectoriel \
$vect\{y_i,z_i\}$ \ pour $i\in\{2,3\}.$

\vspace{0.5cm} Si $(\pi_*^2+\pi^2)^{\frac{1}{2}}:=\pi'\neq 0,$ on
pose

\begin{eqnarray*} y_1' &=& \pi'^{-1}(\pi y_1-\pi_*z_1), \\
z_1' &=& \pi'^{-1}(\pi_* y_1+\pi z_1). \end{eqnarray*}

Les vecteurs $y_1'$ et $z_1'$ sont orthonormaux et on a \[ u\wedge
y_1'=az_1', \hspace{0.2cm} u\wedge y_1'=-ay_1', \hspace{0.2cm}
y_1'\wedge z_1'=y_1\wedge z_1 \ \mbox{ et } \ L_{y_1'}(H_2)\subset
H_3. \] Ce qui permet de supposer $\pi_*=0.$

\vspace{0.5cm} Si $(\theta_*^2+\theta^2){\frac{1}{2}}:=\theta'\neq
0,$ on pose

\begin{eqnarray*} y_3' &=& \theta'^{-1}(\theta y_3-\theta_*z_3), \\
z_3' &=& \theta'^{-1}(\theta_* y_3+\theta z_3).
\end{eqnarray*}

Les vecteurs $y_3'$ et $z_3'$ sont orthonormaux et on a
\[ u\wedge y_3'=cz_3', \hspace{0.2cm} u\wedge y_3'=-cy_3',
\hspace{0.2cm} y_3'\wedge z_3'=y_3\wedge z_3 \ \mbox{ et } \
L_{y_3'}(H_2)\subset vect\{y_1,z_1,z_3'\}.
\] Ce qui permet de supposer $\theta_*=0$ et, pour des raisons analogues, $\omega_*=0.$
On r\'eduit ainsi la multiplication de $A$ \`a $16$ param\`etres:

\vspace{0.2cm}
{\scriptsize \[ \begin{tabular}{ccccccccc} \\
\multicolumn{1}{c}{} & \multicolumn{1}{c}{$1$} &
\multicolumn{1}{c}{$u$} & \multicolumn{1}{c}{$y_1$} &
\multicolumn{1}{c}{$z_1$} & \multicolumn{1}{c}{$y_2$} &
\multicolumn{1}{c}{$z_2$} & \multicolumn{1}{c}{$y_3$} &
\multicolumn{1}{c}{$z_3$}
\\ \cline{2-9}
\multicolumn{1}{c|}{$1$} & \multicolumn{1}{|c}{$1$} &
\multicolumn{1}{c|}{$u$} & \multicolumn{1}{|c}{$y_1$} &
\multicolumn{1}{c|}{$z_1$} & \multicolumn{1}{|c}{$y_2$} &
\multicolumn{1}{c|}{$z_2$} & \multicolumn{1}{|c}{$y_3$} &
\multicolumn{1}{c|}{$z_3$}
\\
\multicolumn{1}{c|}{$u$} & \multicolumn{1}{|c}{} &
\multicolumn{1}{c|}{$-1$} & \multicolumn{1}{|c}{$az_1$} &
\multicolumn{1}{c|}{$-ay_1$} & \multicolumn{1}{|c}{$bz_2$} &
\multicolumn{1}{c|}{$-by_2$} & \multicolumn{1}{|c}{$cz_3$} &
\multicolumn{1}{c|}{$-cy_3$}
\\ \cline{2-9}
\multicolumn{1}{c|}{$y_1$} & \multicolumn{1}{|c}{} &
\multicolumn{1}{c|}{} & \multicolumn{1}{|c}{$-1$} &
\multicolumn{1}{c|}{$au$} & \multicolumn{1}{|c}{$\alpha y_3+\beta
z_3$} & \multicolumn{1}{c|}{$\gamma y_3+\mu z_3$} &
\multicolumn{1}{|c}{$-\alpha y_2-\gamma z_2+\nu z_3$} &
\multicolumn{1}{c|}{$-\beta y_2-\mu z_2-\nu y_3$}
\\
\multicolumn{1}{c|}{$z_1$} & \multicolumn{1}{|c}{} &
\multicolumn{1}{c|}{} & \multicolumn{1}{|c}{} &
\multicolumn{1}{c|}{$-1$} & \multicolumn{1}{|c}{$\pi z_2+\lambda
y_3+\eta z_3$} & \multicolumn{1}{c|}{$-\pi y_2+\sigma y_3+\delta
z_3$} & \multicolumn{1}{|c}{$-\lambda y_2-\sigma z_2+\rho z_3$} &
\multicolumn{1}{c|}{$-\eta y_2-\delta z_2-\rho z_3$}
\\ \cline{4-9}
\multicolumn{1}{c|}{$y_2$} & \multicolumn{1}{|c}{} &
\multicolumn{1}{c}{} & \multicolumn{1}{c}{} &
\multicolumn{1}{c|}{} & \multicolumn{1}{|c}{$-1$} &
\multicolumn{1}{c|}{$bu+\pi z_1+\theta z_3$} &
\multicolumn{1}{|c}{$\alpha y_1+\lambda z_1$} &
\multicolumn{1}{c|}{$\beta y_1+\eta z_1-\theta z_2$}
\\
\multicolumn{1}{c|}{$z_2$} & \multicolumn{1}{|c}{} &
\multicolumn{1}{c}{} & \multicolumn{1}{c}{} & \multicolumn{1}{c}{}
& \multicolumn{1}{|c}{} & \multicolumn{1}{c|}{$-1$} &
\multicolumn{1}{|c}{$\gamma y_1+\sigma z_1+\omega z_3$} &
\multicolumn{1}{c|}{$\mu y_1+\delta z_1+\theta y_2-\omega y_3$}
\\ \cline{6-9}
\multicolumn{1}{c|}{$y_3$} & \multicolumn{1}{|c}{} &
\multicolumn{1}{c}{} & \multicolumn{1}{c}{} & \multicolumn{1}{c}{}
& \multicolumn{1}{c}{} & \multicolumn{1}{c}{} &
\multicolumn{1}{|c}{$-1$} & \multicolumn{1}{c|}{$cu+\nu y_1+\rho
z_1+\omega z_2$}
\\
\multicolumn{1}{c|}{$z_3$} & \multicolumn{1}{|c}{} &
\multicolumn{1}{c}{} & \multicolumn{1}{c}{} & \multicolumn{1}{c}{}
& \multicolumn{1}{c}{} & \multicolumn{1}{c}{} &
\multicolumn{1}{|c}{} & \multicolumn{1}{c|}{$-1$}
\\ \cline{2-9}
\end{tabular} \]}
\begin{center} {\bf Table 1} \end{center}

\vspace{0.5cm} \hspace{0.3cm} R\'eciproquement, un calcul direct
montre qu'une alg\`ebre r\'eelle dont la multiplication est
donn\'ee par la Table {\bf 1} est de Jordan non commutative.

\vspace{0.3cm} \hspace{0.3cm} Nous avons besoin, pour la suite,
des r\'esultats pr\'eliminaires suivants:

\vspace{0.3cm} \begin{lemma} Si \ $\nu=0,$ alors les conditions: \
$\beta\gamma-\alpha\mu,$ \ $\beta\lambda-\alpha\eta,$ \
$\gamma\lambda-\alpha\sigma>0$ \ sont n\'ecessaires pour que
l'alg\`ebre $A$ soit de division lin\'eaire. \end{lemma}

\vspace{0.1cm} {\bf Preuve.} On suppose que $A$ est de division
lin\'eaire et on distingue les deux cas suivants:

\vspace{0.3cm} {\bf 1) Si $\alpha=0,$} alors $\beta\gamma\neq 0$
et on a: \[ (\beta xu-cy_1)(xy_2+y_3)=(\beta bx^2+c\gamma)z_2 \
\mbox{ et } \ (\gamma xu-by_3)(xy_1+y_2)=(\gamma ax^2+b\lambda)z_1
\] pour tout $x\in\rit.$ Comme les termes de gauche dans les deux
derni\`eres \'egalit\'es sont non nuls, les deux trin\^omes en
$x:$ \ $\beta bx^2+c\gamma$ \ et \ $\gamma ax^2+b\lambda$ \ ont
des discriminants n\'egatifs, i.e. $\beta\gamma,$ $\gamma\lambda,$
$\beta\lambda=\gamma^{-2}(\beta\gamma)(\lambda\gamma))>0.$

\vspace{0.3cm} {\bf 2) Si $\alpha\neq 0,$} on consid\`ere
l'automorphisme $\varphi$ de $W$ dont la matrice, par rapport \`a
la base ${\cal B}=\{u,y_1,z_1,y_2,z_2,y_3,z_3\},$ est donn\'ee
par:

\vspace{0.5cm}
\[ \left(
\begin{array}{lllllll}
1 &   &                     &   &                    &     &                   \\
  & 1 & -\lambda\alpha^{-1} &   &                    &     &                   \\
  &   & \hspace{0.5cm}  1   &   &                    &     &                   \\
  &   &                     & 1 & -\gamma\alpha^{-1} &     &                   \\
  &   &                     &   & \hspace{0.5cm} 1   &     &                   \\
  &   &                     &   &                    &  1  & -\beta\alpha^{-1} \\
  &   &                     &   &                    &     & \hspace{0.5cm} 1
\end{array}
\right)  \]

\vspace{0.5cm} La multiplication de l'alg\`ebre $A(\varphi),$ par
rapport \`a la base ${\cal B}$ est donn\'ee par la table:

\vspace{0.3cm}
{\scriptsize \[ \begin{tabular}{ccccccccc} \\
\multicolumn{1}{c}{} & \multicolumn{1}{c}{$1$} &
\multicolumn{1}{c}{$u$} & \multicolumn{1}{c}{$y_1$} &
\multicolumn{1}{c}{$z_1$} & \multicolumn{1}{c}{$y_2$} &
\multicolumn{1}{c}{$z_2$} & \multicolumn{1}{c}{$y_3$} &
\multicolumn{1}{c}{$z_3$} \\ \cline{2-9} \multicolumn{1}{c|}{$1$}
& \multicolumn{1}{|c}{$1$} & \multicolumn{1}{c|}{$u$} &
\multicolumn{1}{|c}{$y_1$} & \multicolumn{1}{c|}{$z_1$} &
\multicolumn{1}{|c}{$y_2$} & \multicolumn{1}{c|}{$z_2$} &
\multicolumn{1}{|c}{$y_3$} & \multicolumn{1}{c|}{$z_3$} \\
\multicolumn{1}{c|}{$u$} & \multicolumn{1}{|c}{} &
\multicolumn{1}{c|}{$-1$} & \multicolumn{1}{|c}{$az_1$} &
\multicolumn{1}{c|}{$-ay_1$} & \multicolumn{1}{|c}{$bz_2$} &
\multicolumn{1}{c|}{$-by_2$} & \multicolumn{1}{|c}{$cz_3$} &
\multicolumn{1}{c|}{$-cy_3$} \\ \cline{2-9}
\multicolumn{1}{c|}{$y_1$} & \multicolumn{1}{|c}{} &
\multicolumn{1}{c|}{} & \multicolumn{1}{|c}{$-1$} &
\multicolumn{1}{c|}{$au$} & \multicolumn{1}{|c}{$\alpha y_3$} &
\multicolumn{1}{c|}{$\mu' z_3$} & \multicolumn{1}{|c}{$-\alpha
y_2$} & \multicolumn{1}{c|}{$-\mu' z_2$} \\
\multicolumn{1}{c|}{$z_1$} & \multicolumn{1}{|c}{} &
\multicolumn{1}{c|}{} & \multicolumn{1}{|c}{} &
\multicolumn{1}{c|}{$-1$} & \multicolumn{1}{|c}{$\pi z_2+\eta'
z_3$} & \multicolumn{1}{c|}{$-\pi y_2+\sigma' y_3+\delta' z_3$} &
\multicolumn{1}{|c}{$-\sigma' z_2+\rho z_3$} &
\multicolumn{1}{c|}{$-\eta' y_2-\delta' z_2-\rho z_3$} \\
\cline{4-9} \multicolumn{1}{c|}{$y_2$} & \multicolumn{1}{|c}{} &
\multicolumn{1}{c}{} & \multicolumn{1}{c}{} &
\multicolumn{1}{c|}{} & \multicolumn{1}{|c}{$-1$} &
\multicolumn{1}{c|}{$bu+\pi z_1+\theta z_3$} &
\multicolumn{1}{|c}{$\alpha y_1$} & \multicolumn{1}{c|}{$\eta'
z_1-\theta z_2$} \\
\multicolumn{1}{c|}{$z_2$} & \multicolumn{1}{|c}{} &
\multicolumn{1}{c}{} & \multicolumn{1}{c}{} & \multicolumn{1}{c}{}
& \multicolumn{1}{|c}{} & \multicolumn{1}{c|}{$-1$} &
\multicolumn{1}{|c}{$\sigma' z_1+\omega z_3$} &
\multicolumn{1}{c|}{$\mu' y_1+\delta' z_1+\theta y_2-\omega y_3$}
\\ \cline{6-9} \multicolumn{1}{c|}{$y_3$} & \multicolumn{1}{|c}{}
& \multicolumn{1}{c}{} & \multicolumn{1}{c}{} &
\multicolumn{1}{c}{} & \multicolumn{1}{c}{} & \multicolumn{1}{c}{}
& \multicolumn{1}{|c}{$-1$} & \multicolumn{1}{c|}{$cu+\rho
z_1+\omega z_2$} \\
\multicolumn{1}{c|}{$z_3$} & \multicolumn{1}{|c}{} &
\multicolumn{1}{c}{} & \multicolumn{1}{c}{} & \multicolumn{1}{c}{}
& \multicolumn{1}{c}{} & \multicolumn{1}{c}{} &
\multicolumn{1}{|c}{} & \multicolumn{1}{c|}{$-1$}
\\ \cline{2-9}
\end{tabular} \]}
\begin{center} {\bf Table 2} \end{center}

\vspace{0.3cm} o\`u \ $\mu'=\mu-\beta\gamma\alpha^{-1},$
$\sigma'=\sigma-\gamma\lambda\alpha^{-1},$
$\eta'=\eta-\beta\lambda\alpha^{-1},$
$\delta'=\delta-(\gamma\eta+\lambda\mu+\beta\sigma)\alpha^{-1}+2\beta\gamma\lambda
\alpha^{-2}.$ \ $A(\varphi)$ \'etant de division lin\'eaire, on a
$\sigma'\neq 0,$ car sinon on aurait \
$L_{y_3}(vect\{z_1,z_2\})\subseteq vect\{z_3\}.$ De plus, les
\'egalit\'es

\[ (xu+z_1)\Big( \sigma' y_2+(bx+\pi)y_3\Big)=\Big(
bcx^2+(b\rho+\pi c)x+\sigma'\eta'+\pi\rho)z_3, \]
\[ (xu+y_1)(-az_2+bxy_3)=(bcx^2-\alpha\mu')z_3, \ \mbox{ et } \]
\[ (xu+z_2)(-\sigma'y_1+axy_3)=(acx^2+\alpha\omega x+\sigma'\mu')z_3 \]

ont lieu pour tout $x\in\rit.$ Comme les termes de gauche dans les
trois derni\`eres \'egalit\'es sont non nuls, les trois trin\^omes
en $x,$ \[ bcx^2+(b\rho+\pi c)x+\sigma'\eta', \hspace{0.3cm}
bcx^2-\alpha\mu', \hspace{0.3cm} acx^2+\alpha\omega x+\sigma'\mu'
\] ont des discriminants n\'egatifs \ i.e. \
$\beta\gamma-\alpha\mu,$ \ $\beta\lambda-\alpha\eta,$ \
$\gamma\lambda-\alpha\sigma>0.\Box$

\vspace{0.3cm} \begin{note} Suivant les notations de la Table {\bf
2}, il existe $\varepsilon\in\{1,-1\}$ tel que \
$|\alpha|=\varepsilon\alpha,$ \ $|\mu'|=-\varepsilon \mu',$ \
$|\sigma'|=-\varepsilon\sigma'$ \ et \
$|\eta'|=-\varepsilon\eta'.\Box$
\end{note}

\vspace{0.3cm} \begin{lemma} Si \ $\nu=\pi=\theta=0,$ alors $A$
est division lin\'eaire si et seulement si \[
\beta\gamma-\alpha\mu, \ \beta\lambda-\alpha\eta, \
\gamma\lambda-\alpha\sigma>0 \ \mbox{ et } \]
\[ c(\alpha\delta-\beta\sigma-\lambda\mu+\gamma\eta)^2+b(\beta\gamma-\alpha\mu)\rho^2+a(\beta\lambda-\alpha\eta)\omega^2
< 4c(\beta\lambda-\alpha\eta)(\mu\sigma-\gamma\delta). \]
\end{lemma}

\vspace{0.1cm} {\bf Preuve.} On suppose que $A$ est division
lin\'eaire et on distingue les deux cas suivants:

\vspace{0.5cm} {\bf 1) Si \ $a=b=c=1,$} on distingue les deux
sous-cas suivants:

\vspace{0.3cm} {\bf i)} \ $\alpha\neq 0.$ On consid\`ere, en
premier lieu, la m\^eme automorphisme $\varphi$ que celui du Lemme
{\bf 4.11} pr\'ec\'edent et, en second lieu, l'automorphisme
$\psi$ de $W$ dont la matrice, par rapport \`a la base ${\cal B},$
est donn\'ee par:

\vspace{0.5cm}
\[ \left(
\begin{array}{lllllll}
1 &   &         &   &        &   &       \\
  & q &         &   &        &   &       \\
  &   & q^{-1}  &   &        &   &       \\
  &   &         & r &        &   &       \\
  &   &         &   & r^{-1} &   &       \\
  &   &         &   &        & s &       \\
  &   &         &   &        &   & s^{-1}
\end{array}
\right)  \]

\vspace{0.5cm} o\`u \begin{eqnarray*} q &=& -\alpha^{-1}\Big(
(\beta\lambda-\alpha\eta)(\gamma\lambda-\alpha\sigma)\Big)^{\frac{1}{2}}, \\
r &=& -\alpha^{-1}\Big(
(\beta\gamma-\alpha\mu)(\gamma\lambda-\alpha\sigma)\Big)^{\frac{1}{2}},
\\ s &=& -\alpha^{-1}\Big(
(\beta\lambda-\alpha\eta)(\beta\gamma-\alpha\mu)\Big)^{\frac{1}{2}}.
\end{eqnarray*}

\vspace{0.5cm} On obtient la table de multiplication de
l'alg\`ebre $A_1=\Big( A(\varphi)\Big)(\psi):$

\vspace{0.3cm}
\[ \begin{tabular}{ccccccccc} \\
\multicolumn{1}{c}{} & \multicolumn{1}{c}{$1$} &
\multicolumn{1}{c}{$u$} & \multicolumn{1}{c}{$y_1$} &
\multicolumn{1}{c}{$z_1$} & \multicolumn{1}{c}{$y_2$} &
\multicolumn{1}{c}{$z_2$} & \multicolumn{1}{c}{$y_3$} &
\multicolumn{1}{c}{$z_3$} \\ \cline{2-9} \multicolumn{1}{c|}{$1$}
& \multicolumn{1}{|c}{$1$} & \multicolumn{1}{c|}{$u$} &
\multicolumn{1}{|c}{$y_1$} & \multicolumn{1}{c|}{$z_1$} &
\multicolumn{1}{|c}{$y_2$} & \multicolumn{1}{c|}{$z_2$} &
\multicolumn{1}{|c}{$y_3$} & \multicolumn{1}{c|}{$z_3$} \\
\multicolumn{1}{c|}{$u$} & \multicolumn{1}{|c}{} &
\multicolumn{1}{c|}{$-1$} & \multicolumn{1}{|c}{$z_1$} &
\multicolumn{1}{c|}{$-y_1$} & \multicolumn{1}{|c}{$z_2$} &
\multicolumn{1}{c|}{$-y_2$} & \multicolumn{1}{|c}{$z_3$} &
\multicolumn{1}{c|}{$-y_3$} \\ \cline{2-9}
\multicolumn{1}{c|}{$y_1$} & \multicolumn{1}{|c}{} &
\multicolumn{1}{c|}{} & \multicolumn{1}{|c}{$-1$} &
\multicolumn{1}{c|}{$u$} & \multicolumn{1}{|c}{$\alpha_* y_3$} &
\multicolumn{1}{c|}{$z_3$} & \multicolumn{1}{|c}{$-\alpha_*
y_2$} & \multicolumn{1}{c|}{$-z_2$} \\
\multicolumn{1}{c|}{$z_1$} & \multicolumn{1}{|c}{} &
\multicolumn{1}{c|}{} & \multicolumn{1}{|c}{} &
\multicolumn{1}{c|}{$-1$} & \multicolumn{1}{|c}{$z_3$} &
\multicolumn{1}{c|}{$y_3+\delta_* z_3$} &
\multicolumn{1}{|c}{$-z_2+\rho' z_3$} &
\multicolumn{1}{c|}{$-y_2-\delta_* z_2-\rho' y_3$} \\
\cline{4-9} \multicolumn{1}{c|}{$y_2$} & \multicolumn{1}{|c}{} &
\multicolumn{1}{c}{} & \multicolumn{1}{c}{} &
\multicolumn{1}{c|}{} & \multicolumn{1}{|c}{$-1$} &
\multicolumn{1}{c|}{$u$} & \multicolumn{1}{|c}{$\alpha_* y_1$} &
\multicolumn{1}{c|}{$z_1$} \\
\multicolumn{1}{c|}{$z_2$} & \multicolumn{1}{|c}{} &
\multicolumn{1}{c}{} & \multicolumn{1}{c}{} & \multicolumn{1}{c}{}
& \multicolumn{1}{|c}{} & \multicolumn{1}{c|}{$-1$} &
\multicolumn{1}{|c}{$z_1+\omega' z_3$} &
\multicolumn{1}{c|}{$y_1+\delta_* z_1-\omega' y_3$}
\\ \cline{6-9} \multicolumn{1}{c|}{$y_3$} & \multicolumn{1}{|c}{}
& \multicolumn{1}{c}{} & \multicolumn{1}{c}{} &
\multicolumn{1}{c}{} & \multicolumn{1}{c}{} & \multicolumn{1}{c}{}
& \multicolumn{1}{|c}{$-1$} & \multicolumn{1}{c|}{$u+\rho'
z_1+\omega' z_2$} \\
\multicolumn{1}{c|}{$z_3$} & \multicolumn{1}{|c}{} &
\multicolumn{1}{c}{} & \multicolumn{1}{c}{} & \multicolumn{1}{c}{}
& \multicolumn{1}{c}{} & \multicolumn{1}{c}{} &
\multicolumn{1}{|c}{} & \multicolumn{1}{c|}{$-1$}
\\ \cline{2-9}
\end{tabular} \]

\vspace{0.5cm} o\`u \begin{eqnarray*} \alpha_* &=& -\alpha^{-2}(\beta\gamma-\alpha\mu)(\beta\lambda-\alpha\eta)
(\gamma\lambda-\alpha\sigma)<0, \\
\delta_* &=&
-\alpha^{-3}\delta'(\beta\gamma-\alpha\mu)^{-1}(\beta\lambda-\alpha\eta)^{-1}
(\gamma\lambda-\alpha\sigma)^{-1}, \\
\rho' &=& -\alpha\rho\Big(
(\beta\lambda-\alpha\eta)(\gamma\lambda-\alpha\sigma)\Big)^{-\frac{1}{2}},
\\ \omega' &=& -\alpha\omega\Big(
(\beta\gamma-\alpha\mu)(\gamma\lambda-\alpha\sigma)\Big)^{-\frac{1}{2}}.
\end{eqnarray*}

\vspace{0.2cm} \hspace{0.3cm} On consid\`ere, finalement,
l'automorphisme $f$ de $W$ dont la matrice, par rapport \`a la
base ${\cal B},$ est:

\vspace{0.5cm}
\[ \left(
\begin{array}{lllllll}
-(-\alpha_*)^{\frac{1}{6}} &                         &                           &                         &  &  & \\
                           & \alpha_*^{-\frac{1}{3}} &                           &                         &  &  & \\
                           &                         & (-\alpha_*)^{\frac{1}{6}} &                         &  &  & \\
                           &                         &                           & \alpha_*^{-\frac{1}{3}} &  &  & \\
                           &                         &             &             & (-\alpha_*)^{\frac{1}{6}}  &  & \\
                           &                         &             &      &      & \alpha_*^{-\frac{1}{3}}    &    \\
                           &                         &             &      &   &  & (-\alpha_*)^{\frac{1}{6}}
\end{array}
\right)  \]

\vspace{0.5cm} On obtient la multiplication de l'alg\`ebre
$A_2=A_1(f):$

\vspace{0.5cm}
\[ \begin{tabular}{ccccccccc} \\
\multicolumn{1}{c}{} & \multicolumn{1}{c}{$1$} &
\multicolumn{1}{c}{$u$} & \multicolumn{1}{c}{$y_1$} &
\multicolumn{1}{c}{$z_1$} & \multicolumn{1}{c}{$y_2$} &
\multicolumn{1}{c}{$z_2$} & \multicolumn{1}{c}{$y_3$} &
\multicolumn{1}{c}{$z_3$} \\ \cline{2-9} \multicolumn{1}{c|}{$1$}
& \multicolumn{1}{|c}{$1$} & \multicolumn{1}{c|}{$u$} &
\multicolumn{1}{|c}{$y_1$} & \multicolumn{1}{c|}{$z_1$} &
\multicolumn{1}{|c}{$y_2$} & \multicolumn{1}{c|}{$z_2$} &
\multicolumn{1}{|c}{$y_3$} & \multicolumn{1}{c|}{$z_3$} \\
\multicolumn{1}{c|}{$u$} & \multicolumn{1}{|c}{} &
\multicolumn{1}{c|}{$-1$} & \multicolumn{1}{|c}{$z_1$} &
\multicolumn{1}{c|}{$-y_1$} & \multicolumn{1}{|c}{$z_2$} &
\multicolumn{1}{c|}{$-y_2$} & \multicolumn{1}{|c}{$z_3$} &
\multicolumn{1}{c|}{$-y_3$} \\ \cline{2-9}
\multicolumn{1}{c|}{$y_1$} & \multicolumn{1}{|c}{} &
\multicolumn{1}{c|}{} & \multicolumn{1}{|c}{$-1$} &
\multicolumn{1}{c|}{$u$} & \multicolumn{1}{|c}{$y_3$} &
\multicolumn{1}{c|}{$-z_3$} & \multicolumn{1}{|c}{$-y_2$} & \multicolumn{1}{c|}{$z_2$} \\
\multicolumn{1}{c|}{$z_1$} & \multicolumn{1}{|c}{} &
\multicolumn{1}{c|}{} & \multicolumn{1}{|c}{} &
\multicolumn{1}{c|}{$-1$} & \multicolumn{1}{|c}{$-z_3$} &
\multicolumn{1}{c|}{$-y_3+\delta^* z_3$} &
\multicolumn{1}{|c}{$z_2-\rho' z_3$} &
\multicolumn{1}{c|}{$y_2-\delta^* z_2+\rho' y_3$} \\
\cline{4-9} \multicolumn{1}{c|}{$y_2$} & \multicolumn{1}{|c}{} &
\multicolumn{1}{c}{} & \multicolumn{1}{c}{} &
\multicolumn{1}{c|}{} & \multicolumn{1}{|c}{$-1$} &
\multicolumn{1}{c|}{$u$} & \multicolumn{1}{|c}{$y_1$} &
\multicolumn{1}{c|}{$-z_1$} \\
\multicolumn{1}{c|}{$z_2$} & \multicolumn{1}{|c}{} &
\multicolumn{1}{c}{} & \multicolumn{1}{c}{} & \multicolumn{1}{c}{}
& \multicolumn{1}{|c}{} & \multicolumn{1}{c|}{$-1$} &
\multicolumn{1}{|c}{$-z_1-\omega' z_3$} &
\multicolumn{1}{c|}{$-y_1+\delta^* z_1+\omega' y_3$}
\\ \cline{6-9} \multicolumn{1}{c|}{$y_3$} & \multicolumn{1}{|c}{}
& \multicolumn{1}{c}{} & \multicolumn{1}{c}{} &
\multicolumn{1}{c}{} & \multicolumn{1}{c}{} & \multicolumn{1}{c}{}
& \multicolumn{1}{|c}{$-1$} & \multicolumn{1}{c|}{$u-\rho'
z_1-\omega' z_2$} \\
\multicolumn{1}{c|}{$z_3$} & \multicolumn{1}{|c}{} &
\multicolumn{1}{c}{} & \multicolumn{1}{c}{} & \multicolumn{1}{c}{}
& \multicolumn{1}{c}{} & \multicolumn{1}{c}{} &
\multicolumn{1}{|c}{} & \multicolumn{1}{c|}{$-1$}
\\ \cline{2-9}
\end{tabular} \]

\vspace{0.5cm} o\`u \ $\delta^*=(-\alpha_*)^{\frac{1}{2}}\delta_*.$ On pose alors \ $e_4=u,$ $e_i=y_i$ \ et \
$e_{i+4}=-z_i,$ $i\in\{1,2,3\}$ et on obtient, par rapport \`a la base \ $1, e_1,\dots, e_7$ \ la table:

\vspace{0.5cm}
\[ \begin{tabular}{ccccccccc} \\
\multicolumn{1}{c}{} & \multicolumn{1}{c}{$1$} &
\multicolumn{1}{c}{$e_1$} & \multicolumn{1}{c}{$e_2$} &
\multicolumn{1}{c}{$e_3$} & \multicolumn{1}{c}{$e_4$} &
\multicolumn{1}{c}{$e_5$} & \multicolumn{1}{c}{$e_6$} &
\multicolumn{1}{c}{$e_7$} \\ \cline{2-9} \multicolumn{1}{c|}{$1$}
& \multicolumn{1}{|c}{$1$} & \multicolumn{1}{c|}{$e_1$} &
\multicolumn{1}{|c}{$e_2$} & \multicolumn{1}{c|}{$e_3$} &
\multicolumn{1}{|c}{$e_4$} & \multicolumn{1}{c|}{$e_5$} &
\multicolumn{1}{|c}{$e_6$} & \multicolumn{1}{c|}{$e_7$} \\
\multicolumn{1}{c|}{$e_1$} & \multicolumn{1}{|c}{} &
\multicolumn{1}{c|}{$-1$} & \multicolumn{1}{|c}{$e_3$} &
\multicolumn{1}{c|}{$-e_2$} & \multicolumn{1}{|c}{$e_5$} &
\multicolumn{1}{c|}{$-e_4$} & \multicolumn{1}{|c}{$-e_7$} &
\multicolumn{1}{c|}{$e_6$} \\ \cline{2-9}
\multicolumn{1}{c|}{$e_2$} & \multicolumn{1}{|c}{} &
\multicolumn{1}{c|}{} & \multicolumn{1}{|c}{$-1$} &
\multicolumn{1}{c|}{$e_1$} & \multicolumn{1}{|c}{$e_6$} &
\multicolumn{1}{c|}{$e_7$} & \multicolumn{1}{|c}{$-e_4$} & \multicolumn{1}{c|}{$-e_5$} \\
\multicolumn{1}{c|}{$e_3$} & \multicolumn{1}{|c}{} &
\multicolumn{1}{c|}{} & \multicolumn{1}{|c}{} &
\multicolumn{1}{c|}{$-1$} & \multicolumn{1}{|c}{$e_7$} &
\multicolumn{1}{c|}{$-e_6+\rho' e_7$} &
\multicolumn{1}{|c}{$e_5+\omega' e_7$} &
\multicolumn{1}{c|}{$-e_4-\rho' e_5-\omega' e_6$} \\
\cline{4-9} \multicolumn{1}{c|}{$e_4$} & \multicolumn{1}{|c}{} &
\multicolumn{1}{c}{} & \multicolumn{1}{c}{} &
\multicolumn{1}{c|}{} & \multicolumn{1}{|c}{$-1$} &
\multicolumn{1}{c|}{$e_1$} & \multicolumn{1}{|c}{$e_2$} &
\multicolumn{1}{c|}{$e_3$} \\
\multicolumn{1}{c|}{$e_5$} & \multicolumn{1}{|c}{} &
\multicolumn{1}{c}{} & \multicolumn{1}{c}{} & \multicolumn{1}{c}{}
& \multicolumn{1}{|c}{} & \multicolumn{1}{c|}{$-1$} &
\multicolumn{1}{|c}{$-e_3-\delta^* e_7$} &
\multicolumn{1}{c|}{$e_2+\rho' e_3+\delta^* e_6$}
\\ \cline{6-9} \multicolumn{1}{c|}{$e_6$} & \multicolumn{1}{|c}{}
& \multicolumn{1}{c}{} & \multicolumn{1}{c}{} &
\multicolumn{1}{c}{} & \multicolumn{1}{c}{} & \multicolumn{1}{c}{}
& \multicolumn{1}{|c}{$-1$} & \multicolumn{1}{c|}{$-e_1+\omega'
e_3-\delta^* e_5$} \\
\multicolumn{1}{c|}{$e_7$} & \multicolumn{1}{|c}{} &
\multicolumn{1}{c}{} & \multicolumn{1}{c}{} & \multicolumn{1}{c}{}
& \multicolumn{1}{c}{} & \multicolumn{1}{c}{} &
\multicolumn{1}{|c}{} & \multicolumn{1}{c|}{$-1$}
\\ \cline{2-9}
\end{tabular} \]

\vspace{0.5cm} Si $\rho=\omega=0,$ alors $A_2\simeq$
$E_{-1,\delta^*}(\hit)$ et est de division lin\'eaire si et
seulement si \ $|\delta^*|<2$ ([KR 92] Proposition {\bf 3.9}). On
a

\begin{eqnarray*} \delta^{*2}<4 &\Leftrightarrow&
\alpha^4\delta'^2(\beta\gamma-\alpha\mu)^{-1}(\beta\lambda-\alpha\eta)^{-1}(\gamma\lambda-\alpha\sigma)^{-1}<4
\\ &\Leftrightarrow&
\Big(
\alpha^2\delta-(\gamma\eta+\lambda\mu+\beta\sigma)\alpha+2\beta\lambda\gamma)\Big)^2<
4(\beta\gamma-\alpha\mu)(\beta\lambda-\alpha\eta)(\gamma\lambda-\alpha\sigma)
\\ &\Leftrightarrow&
\Big(
\alpha(\alpha\delta-\beta\sigma-\lambda\mu+\gamma\eta)+2\gamma(\beta\lambda-\alpha\eta)\Big)^2<
\\ && \hspace{1.5cm} 4(\beta\lambda-\alpha\eta)\Big(
(\alpha^2(\sigma\mu-\gamma\delta)+\gamma\alpha(\alpha\delta-\beta\sigma-\lambda\mu+\gamma\eta)+
\gamma^2(\beta\lambda-\alpha\eta)\Big) \\
&\Leftrightarrow&
(\alpha\delta-\beta\sigma-\lambda\mu+\gamma\eta)^2<4(\beta\lambda-\alpha\eta)(\sigma\mu-\gamma\delta).
\end{eqnarray*}

\vspace{0.5cm} Si \ $(\rho,\omega)\neq(0,0),$ on pose \
$\delta_1=(\rho'^2+\delta^{*2})^{\frac{1}{2}},$ \
$\delta_2=(\delta_1^2+\omega'^2)^{\frac{1}{2}}$ \ et

\begin{eqnarray*} e_1' &=& \delta_1^{-1}(\delta^*e_1+\rho'e_4) \ \mbox{ si } \rho\neq 0 \ (=e_1 \ \mbox{ si } \rho=0) \\
e_2' &=& \delta_2^{-1}(-\omega'\rho'\delta_1^{-1}e_1+\delta_1e_2+\omega'\delta^*\delta_1^{-1}e_4) \ \mbox{ si }
\rho\neq 0 \ \Big( =\delta_2^{-1}(\delta^*e_2+\omega'e_4) \ \mbox{ si } \rho=0\Big) \\
e_3' &=& \delta_2^{-1}(\delta^*e_3+\omega'e_5-\rho'e_6) \\
e_4' &=& \delta_2^{-1}(-\rho' e_1-\omega'e_2+\delta^*e_4) \\
e_5' &=& \delta_2^{-1}(-\omega'\delta^*\delta_1^{-1}e_3+\delta_1e_5+\omega'\rho'\delta_1^{-1}e_6) \ \mbox{ si }
\rho\neq 0 \ \Big( =\delta_2^{-1}(-\omega'e_3+\delta^*e_5) \ \mbox{ si } \rho=0\Big) \\
e_6' &=& \delta_1^{-1}(\rho' e_3+\delta^*e_6) \ \mbox{ si } \rho\neq 0 \ (=e_6 \ \mbox{ si } \rho=0) \\
e_7' &=& e_7.
\end{eqnarray*}

\vspace{0.5cm} On obtient la table:

\vspace{0.5cm}
\[ \begin{tabular}{ccccccccc} \\
\multicolumn{1}{c}{} & \multicolumn{1}{c}{$1$} &
\multicolumn{1}{c}{$e_1'$} & \multicolumn{1}{c}{$e_2'$} &
\multicolumn{1}{c}{$e_3'$} & \multicolumn{1}{c}{$e_4'$} &
\multicolumn{1}{c}{$e_5'$} & \multicolumn{1}{c}{$e_6'$} &
\multicolumn{1}{c}{$e_7'$} \\ \cline{2-9} \multicolumn{1}{c|}{$1$}
& \multicolumn{1}{|c}{$1$} & \multicolumn{1}{c|}{$e_1'$} &
\multicolumn{1}{|c}{$e_2'$} & \multicolumn{1}{c|}{$e_3'$} &
\multicolumn{1}{|c}{$e_4'$} & \multicolumn{1}{c|}{$e_5'$} &
\multicolumn{1}{|c}{$e_6'$} & \multicolumn{1}{c|}{$e_7'$} \\
\multicolumn{1}{c|}{$e_1'$} & \multicolumn{1}{|c}{} &
\multicolumn{1}{c|}{$-1$} & \multicolumn{1}{|c}{$e_3'$} &
\multicolumn{1}{c|}{$-e_2'$} & \multicolumn{1}{|c}{$e_5'$} &
\multicolumn{1}{c|}{$-e_4'$} & \multicolumn{1}{|c}{$-e_7'$} &
\multicolumn{1}{c|}{$e_6'$} \\ \cline{2-9}
\multicolumn{1}{c|}{$e_2'$} & \multicolumn{1}{|c}{} &
\multicolumn{1}{c|}{} & \multicolumn{1}{|c}{$-1$} &
\multicolumn{1}{c|}{$e_1'$} & \multicolumn{1}{|c}{$e_6'$} &
\multicolumn{1}{c|}{$e_7'$} & \multicolumn{1}{|c}{$-e_4'$} & \multicolumn{1}{c|}{$-e_5'$} \\
\multicolumn{1}{c|}{$e_3'$} & \multicolumn{1}{|c}{} &
\multicolumn{1}{c|}{} & \multicolumn{1}{|c}{} &
\multicolumn{1}{c|}{$-1$} & \multicolumn{1}{|c}{$e_7'$} &
\multicolumn{1}{c|}{$-e_6'$} & \multicolumn{1}{|c}{$e_5'$} &
\multicolumn{1}{c|}{$-e_4'$} \\
\cline{4-9} \multicolumn{1}{c|}{$e_4'$} & \multicolumn{1}{|c}{} &
\multicolumn{1}{c}{} & \multicolumn{1}{c}{} &
\multicolumn{1}{c|}{} & \multicolumn{1}{|c}{$-1$} &
\multicolumn{1}{c|}{$e_1'$} & \multicolumn{1}{|c}{$e_2'$} &
\multicolumn{1}{c|}{$e_3'$} \\
\multicolumn{1}{c|}{$e_5'$} & \multicolumn{1}{|c}{} &
\multicolumn{1}{c}{} & \multicolumn{1}{c}{} & \multicolumn{1}{c}{}
& \multicolumn{1}{|c}{} & \multicolumn{1}{c|}{$-1$} &
\multicolumn{1}{|c}{$-e_3'-\delta_2e_7'$} &
\multicolumn{1}{c|}{$e_2'+\delta_2e_6'$}
\\ \cline{6-9} \multicolumn{1}{c|}{$e_6'$} & \multicolumn{1}{|c}{}
& \multicolumn{1}{c}{} & \multicolumn{1}{c}{} &
\multicolumn{1}{c}{} & \multicolumn{1}{c}{} & \multicolumn{1}{c}{}
& \multicolumn{1}{|c}{$-1$} & \multicolumn{1}{c|}{$-e_1'-\delta_2e_5'$} \\
\multicolumn{1}{c|}{$e_7'$} & \multicolumn{1}{|c}{} &
\multicolumn{1}{c}{} & \multicolumn{1}{c}{} & \multicolumn{1}{c}{}
& \multicolumn{1}{c}{} & \multicolumn{1}{c}{} &
\multicolumn{1}{|c}{} & \multicolumn{1}{c|}{$-1$}
\\ \cline{2-9}
\end{tabular} \]

\vspace{0.5cm} Donc $A_2\simeq$ $E_{-1,\delta_2}(\hit)$ et est de
division lin\'eaire si et seulement si \ $|\delta_2|<2.$ On a

\begin{eqnarray*} \delta_2^2 &=& \omega'^2+\rho'^2+\delta^{*2} \\
&=&
\alpha^2\omega^2(\beta\gamma-\alpha\mu)^{-1}(\gamma\lambda-\alpha\sigma)^{-1}+\alpha^2\rho^2(\beta\lambda-\alpha\eta)^{-1}
(\gamma\lambda-\alpha\sigma)^{-1}
\\ && \hspace{5cm} +\alpha^4\delta'^2(\beta\gamma-\alpha\mu)^{-1}(\beta\lambda-\alpha\eta)^{-1}
(\gamma\lambda-\alpha\sigma)^{-1}.
\end{eqnarray*} Ainsi \begin{eqnarray*} \delta_2^2<4 &\Leftrightarrow& \alpha^2(\beta\lambda-\alpha\eta)\omega^2+
\alpha^2(\beta\gamma-\alpha\mu)\rho^2+\alpha^4\delta'^2<4(\beta\gamma-\alpha\mu)(\beta\lambda-\alpha\eta)
(\gamma\lambda-\alpha\sigma) \\
&\Leftrightarrow&
(\alpha\delta-\beta\sigma-\lambda\mu+\gamma\eta)^2+(\beta\lambda-\alpha\eta)\omega^2+(\beta\gamma-\alpha\mu)\rho^2
<4(\beta\lambda-\alpha\eta)(\sigma\mu-\gamma\delta). \\
\end{eqnarray*}

En tenant compte de la d\'efinition de l'automorphisme $\psi,$ du
d\'ebut de la d\'emonstration, $A$ est de division lin\'eaire si
et seulement si \ $\beta\lambda-\alpha\eta,$
$\beta\gamma-\alpha\mu,$ $\gamma\lambda-\alpha\sigma>0$ \ et \
$A_2$ est de division lin\'eaire. Ce qui \'etablit le r\'esultat
dans cette premi\`ere situation.

\vspace{0.3cm} {\bf ii)} \ $\alpha=0.$ On pose \ $y_3'=z_3,$ \
$z_3'=-y_3$ et on obtient la table

\vspace{0.5cm}
\[ \begin{tabular}{ccccccccc} \\
\multicolumn{1}{c}{} & \multicolumn{1}{c}{$1$} &
\multicolumn{1}{c}{$u$} & \multicolumn{1}{c}{$y_1$} &
\multicolumn{1}{c}{$z_1$} & \multicolumn{1}{c}{$y_2$} &
\multicolumn{1}{c}{$z_2$} & \multicolumn{1}{c}{$y_3'$} &
\multicolumn{1}{c}{$z_3'$} \\ \cline{2-9} \multicolumn{1}{c|}{$1$}
& \multicolumn{1}{|c}{$1$} & \multicolumn{1}{c|}{$u$} &
\multicolumn{1}{|c}{$y_1$} & \multicolumn{1}{c|}{$z_1$} &
\multicolumn{1}{|c}{$y_2$} & \multicolumn{1}{c|}{$z_2$} &
\multicolumn{1}{|c}{$y_3'$} & \multicolumn{1}{c|}{$z_3'$} \\
\multicolumn{1}{c|}{$u$} & \multicolumn{1}{|c}{} &
\multicolumn{1}{c|}{$-1$} & \multicolumn{1}{|c}{$z_1$} &
\multicolumn{1}{c|}{$-y_1$} & \multicolumn{1}{|c}{$z_2$} &
\multicolumn{1}{c|}{$-y_2$} & \multicolumn{1}{|c}{$z_3'$} &
\multicolumn{1}{c|}{$-y_3'$} \\ \cline{2-9}
\multicolumn{1}{c|}{$y_1$} & \multicolumn{1}{|c}{} &
\multicolumn{1}{c|}{} & \multicolumn{1}{|c}{$-1$} &
\multicolumn{1}{c|}{$u$} & \multicolumn{1}{|c}{$\beta y_3'$} &
\multicolumn{1}{c|}{$\mu y_3'-\gamma z_3'$} & \multicolumn{1}{|c}{$-\beta y_2-\mu z_2$} &
\multicolumn{1}{c|}{$\gamma z_2$} \\
\multicolumn{1}{c|}{$z_1$} & \multicolumn{1}{|c}{} &
\multicolumn{1}{c|}{} & \multicolumn{1}{|c}{} &
\multicolumn{1}{c|}{$-1$} & \multicolumn{1}{|c}{$\eta y_3'-\lambda
z_3'$} & \multicolumn{1}{c|}{$\delta y_3'-\sigma z_3'$} &
\multicolumn{1}{|c}{$-\eta y_2-\delta z_2+\rho z_3'$} &
\multicolumn{1}{c|}{$\lambda y_2+\sigma z_2-\rho y_3'$} \\
\cline{4-9} \multicolumn{1}{c|}{$y_2$} & \multicolumn{1}{|c}{} &
\multicolumn{1}{c}{} & \multicolumn{1}{c}{} &
\multicolumn{1}{c|}{} & \multicolumn{1}{|c}{$-1$} &
\multicolumn{1}{c|}{$u$} & \multicolumn{1}{|c}{$\beta y_1+\eta
z_1$} &
\multicolumn{1}{c|}{$-\lambda z_1$} \\
\multicolumn{1}{c|}{$z_2$} & \multicolumn{1}{|c}{} &
\multicolumn{1}{c}{} & \multicolumn{1}{c}{} & \multicolumn{1}{c}{}
& \multicolumn{1}{|c}{} & \multicolumn{1}{c|}{$-1$} &
\multicolumn{1}{|c}{$\mu y_1+\delta z_1+\omega z_3'$} &
\multicolumn{1}{c|}{$-\gamma y_1-\sigma z_1-\omega y_3'$}
\\ \cline{6-9} \multicolumn{1}{c|}{$y_3'$} & \multicolumn{1}{|c}{}
& \multicolumn{1}{c}{} & \multicolumn{1}{c}{} &
\multicolumn{1}{c}{} & \multicolumn{1}{c}{} & \multicolumn{1}{c}{}
& \multicolumn{1}{|c}{$-1$} & \multicolumn{1}{c|}{$u+\rho
z_1+\omega z_2$} \\
\multicolumn{1}{c|}{$z_3'$} & \multicolumn{1}{|c}{} &
\multicolumn{1}{c}{} & \multicolumn{1}{c}{} & \multicolumn{1}{c}{}
& \multicolumn{1}{c}{} & \multicolumn{1}{c}{} &
\multicolumn{1}{|c}{} & \multicolumn{1}{c|}{$-1$}
\\ \cline{2-9}
\end{tabular} \]

\vspace{0.5cm} Ainsi $A$ est de division lin\'eaire si et
seulement si \ $\beta\lambda,$ $\beta\gamma,$
$\mu\eta-\beta\delta>0$ \ et \[
(-\beta\sigma-\lambda\mu+\gamma\eta)^2+\beta\lambda\omega^2+\beta\gamma\rho^2<4\beta\lambda(\sigma\mu-\gamma\delta),
\ i.e. \] si et seulement si \ $\beta\lambda,$ $\beta\gamma,$
$\gamma\lambda>0$ \ et \
$(-\beta\sigma-\lambda\mu+\gamma\eta)^2+\beta\lambda\omega^2+\beta\gamma\rho^2<4\beta\lambda(\sigma\mu-\gamma\delta)$
\ car \
$(-\beta\sigma-\lambda\mu+\gamma\eta)^2-4\gamma\lambda(\mu\eta-\beta\delta)=
(-\beta\sigma-\lambda\mu+\gamma\eta)^2-4\beta\lambda(\sigma\mu-\gamma\delta)<0$
\ i.e. \ $\mu\eta-\beta\delta>0.$

\vspace{0.3cm} {\bf 2) Si \ $a, b, c>0,$} on consid\`ere
l'automorphisme $g$ de $W$ dont la matrice, par rapport \`a la
base ${\cal B}$ est:

\vspace{0.4cm}
\[ \left(
\begin{array}{lllllll}
1 &                  &                  &                  &                  &                  & \\
  & a^{-\frac{1}{2}} &                  &                  &                  &                  & \\
  &                  & a^{-\frac{1}{2}} &                  &                  &                  & \\
  &                  &                  & b^{-\frac{1}{2}} &                  &                  & \\
  &                  &                  &                  & b^{-\frac{1}{2}} &                  & \\
  &                  &                  &                  &                  & c^{-\frac{1}{2}} &    \\
  &                  &                  &                  &                  &                  & c^{-\frac{1}{2}}
\end{array}
\right)  \]

\vspace{0.5cm} On obtient la table de multiplication de
l'alg\`ebre $A(g):$

\vspace{0.2cm}
\[ \begin{tabular}{ccccccccc} \\
\multicolumn{1}{c}{} & \multicolumn{1}{c}{$1$} &
\multicolumn{1}{c}{$u$} & \multicolumn{1}{c}{$y_1$} &
\multicolumn{1}{c}{$z_1$} & \multicolumn{1}{c}{$y_2$} &
\multicolumn{1}{c}{$z_2$} & \multicolumn{1}{c}{$y_3$} &
\multicolumn{1}{c}{$z_3$} \\ \cline{2-9} \multicolumn{1}{c|}{$1$}
& \multicolumn{1}{|c}{$1$} & \multicolumn{1}{c|}{$u$} &
\multicolumn{1}{|c}{$y_1$} & \multicolumn{1}{c|}{$z_1$} &
\multicolumn{1}{|c}{$y_2$} & \multicolumn{1}{c|}{$z_2$} &
\multicolumn{1}{|c}{$y_3$} & \multicolumn{1}{c|}{$z_3$} \\
\multicolumn{1}{c|}{$u$} & \multicolumn{1}{|c}{} &
\multicolumn{1}{c|}{$-1$} & \multicolumn{1}{|c}{$z_1$} &
\multicolumn{1}{c|}{$-y_1$} & \multicolumn{1}{|c}{$z_2$} &
\multicolumn{1}{c|}{$-y_2$} & \multicolumn{1}{|c}{$z_3$} &
\multicolumn{1}{c|}{$-y_3$} \\ \cline{2-9}
\multicolumn{1}{c|}{$y_1$} & \multicolumn{1}{|c}{} &
\multicolumn{1}{c|}{} & \multicolumn{1}{|c}{$-1$} &
\multicolumn{1}{c|}{$u$} &
\multicolumn{1}{|c}{$\alpha_0y_3+\beta_0z_3$} &
\multicolumn{1}{c|}{$\gamma_0y_3+\mu_0z_3$} & \multicolumn{1}{|c}{$-\alpha_0y_2-\gamma_0z_2$} &
\multicolumn{1}{c|}{$-\beta_0y_2-\mu_0z_2$} \\
\multicolumn{1}{c|}{$z_1$} & \multicolumn{1}{|c}{} &
\multicolumn{1}{c|}{} & \multicolumn{1}{|c}{} &
\multicolumn{1}{c|}{$-1$} &
\multicolumn{1}{|c}{$\lambda_0y_3+\eta_0z_3$} &
\multicolumn{1}{c|}{$\sigma_0y_3+\delta_0z_3$} &
\multicolumn{1}{|c}{$-\lambda_0y_2-\sigma_0z_2+\rho_0z_3$} &
\multicolumn{1}{c|}{$\eta_0y_2-\delta_0z_2-\rho_0y_3$} \\
\cline{4-9} \multicolumn{1}{c|}{$y_2$} & \multicolumn{1}{|c}{} &
\multicolumn{1}{c}{} & \multicolumn{1}{c}{} &
\multicolumn{1}{c|}{} & \multicolumn{1}{|c}{$-1$} &
\multicolumn{1}{c|}{$u$} &
\multicolumn{1}{|c}{$\alpha_0y_1+\lambda_0z_1$} &
\multicolumn{1}{c|}{$-z_1$} \\
\multicolumn{1}{c|}{$z_2$} & \multicolumn{1}{|c}{} &
\multicolumn{1}{c}{} & \multicolumn{1}{c}{} & \multicolumn{1}{c}{}
& \multicolumn{1}{|c}{} & \multicolumn{1}{c|}{$-1$} &
\multicolumn{1}{|c}{$\gamma_0y_1+\sigma_0z_1+\omega_0z_3$} &
\multicolumn{1}{c|}{$\mu_0y_1+\delta_0z_1-\omega_0y_3$}
\\ \cline{6-9} \multicolumn{1}{c|}{$y_3$} & \multicolumn{1}{|c}{}
& \multicolumn{1}{c}{} & \multicolumn{1}{c}{} &
\multicolumn{1}{c}{} & \multicolumn{1}{c}{} & \multicolumn{1}{c}{}
& \multicolumn{1}{|c}{$-1$} & \multicolumn{1}{c|}{$u+\rho_0z_1+\omega_0z_2$} \\
\multicolumn{1}{c|}{$z_3$} & \multicolumn{1}{|c}{} &
\multicolumn{1}{c}{} & \multicolumn{1}{c}{} & \multicolumn{1}{c}{}
& \multicolumn{1}{c}{} & \multicolumn{1}{c}{} &
\multicolumn{1}{|c}{} & \multicolumn{1}{c|}{$-1$}
\\ \cline{2-9}
\end{tabular} \]

\vspace{0.5cm} o\`u \ $\alpha_0=d\alpha,$ \ $\beta_0=d\beta,$ \
$\gamma_0=d\gamma,$ \ $\mu_0=d\mu,$ \ $\eta_0=d\eta,$ \
$\sigma_0=d\sigma,$ \ $\delta_0=d\delta,$ \
$\omega_0=b^{-\frac{1}{2}}c^{-1}\omega,$ \
$\rho_0=a^{-\frac{1}{2}}c^{-1}\rho$ \ avec
$d=(abc)^{-\frac{1}{2}}.$ Ceci nous ram\`ene \`a la premi\`ere
situation. Ainsi $A(g)$ est de division lin\'eaire si et seulement
si: \ $\beta\lambda-\alpha\eta,$ $\beta\gamma-\alpha\mu,$
$\gamma\lambda-\alpha\sigma>$ \ et \[
c(\alpha\delta-\beta\sigma-\lambda\mu+\gamma\eta)^2+a(\beta\lambda-\alpha\eta)\omega^2+b(\beta\gamma-\alpha\mu)\rho^2
<4c(\beta\lambda-\alpha\eta)(\sigma\mu-\gamma\delta). \] La
d\'emonstration s'ach\`eve vu que $A$ est de division lin\'eaire
si et seulement si $A(g)$ est de division lin\'eaire.$\Box$

\vspace{0.3cm} \begin{corollary} Soient \ $\lambda,$ $\mu,$
$\alpha,$ $\beta,$ $\delta$ \ des nombres r\'eels arbitraires.
Alors l'alg\`ebre r\'eelle \ $\Big(
E_{-1,\alpha,\beta,\delta,-\beta}(\hit^{(\lambda)})\Big)^{(\mu)}$
\ est de Jordan non commutative. Elle est de division lin\'eaire
si et seulement si \ $\lambda,\mu\neq\frac{1}{2},$
$\beta^2+2\alpha-1>0$ \ et \ $\Big(
(2\alpha-1)\delta-\beta\Big)<4(\beta^2+2\alpha-1)(1+\beta\delta).$
\end{corollary}

\vspace{0.1cm} {\bf Preuve.} On peut supposer, en vertu de la
Proposition {\bf 1.47}, que $\mu=1.$ Si \ $1,e_1,e_2,e_3$ \ est la
base canonique de $\hit,$ alors \ $1,$ $e_1,$ $e_2,$ $e_3,$ $e_4,$
$e_5,$ $e_6,$ $e_7$ \ o\`u $e_4=(0,1)$ et $e_{i+4}=e_ie_4$
$i\in\{1,2,3\}$ \ est une base de \
$A=E_{-1,\alpha,\beta,\delta,-\beta}(\hit^{(\lambda)})$ \ pour
laquelle la multiplication de $A$ est donn\'ee par la table:

\vspace{0.5cm}
\[ \begin{tabular}{ccccccccc} \\
\multicolumn{1}{c}{} & \multicolumn{1}{c}{$1$} &
\multicolumn{1}{c}{$e_1$} & \multicolumn{1}{c}{$e_2$} &
\multicolumn{1}{c}{$e_3$} & \multicolumn{1}{c}{$e_4$} &
\multicolumn{1}{c}{$e_5$} & \multicolumn{1}{c}{$e_6$} &
\multicolumn{1}{c}{$e_7$} \\ \cline{2-9} \multicolumn{1}{c|}{$1$}
& \multicolumn{1}{|c}{$1$} & \multicolumn{1}{c|}{$e_1$} &
\multicolumn{1}{|c}{$e_2$} & \multicolumn{1}{c|}{$e_3$} &
\multicolumn{1}{|c}{$e_4$} & \multicolumn{1}{c|}{$e_5$} &
\multicolumn{1}{|c}{$e_6$} & \multicolumn{1}{c|}{$e_7$} \\
\multicolumn{1}{c|}{$e_1$} & \multicolumn{1}{|c}{} &
\multicolumn{1}{c|}{$-1$} &
\multicolumn{1}{|c}{$\alpha'e_3+\beta'e_7$} &
\multicolumn{1}{c|}{$-\alpha'e_2-\beta'e_6$} &
\multicolumn{1}{|c}{$e_5$} & \multicolumn{1}{c|}{$-e_4$} &
\multicolumn{1}{|c}{$\beta'e_3-\lambda'e_7$} & \multicolumn{1}{c|}{$-\beta'e_2+\lambda'e_6$} \\
\cline{2-9} \multicolumn{1}{c|}{$e_2$} & \multicolumn{1}{|c}{} &
\multicolumn{1}{c|}{} & \multicolumn{1}{|c}{$-1$} &
\multicolumn{1}{c|}{$\alpha'e_1+\beta'e_5$} &
\multicolumn{1}{|c}{$e_6$} &
\multicolumn{1}{c|}{$-\beta'e_3+\lambda'e_7$} & \multicolumn{1}{|c}{$-e_4$} &
\multicolumn{1}{c|}{$\beta'e_1-\lambda'e_5$} \\
\multicolumn{1}{c|}{$e_3$} & \multicolumn{1}{|c}{} &
\multicolumn{1}{c|}{} & \multicolumn{1}{|c}{} &
\multicolumn{1}{c|}{$-1$} & \multicolumn{1}{|c}{$e_7$} &
\multicolumn{1}{c|}{$\beta'e_2-\lambda'e_6$} &
\multicolumn{1}{|c}{$-\beta'e_1+\lambda'e_5$} &
\multicolumn{1}{c|}{$-e_4$} \\
\cline{4-9} \multicolumn{1}{c|}{$e_4$} & \multicolumn{1}{|c}{} &
\multicolumn{1}{c}{} & \multicolumn{1}{c}{} &
\multicolumn{1}{c|}{} & \multicolumn{1}{|c}{$-1$} &
\multicolumn{1}{c|}{$e_1$} & \multicolumn{1}{|c}{$e_2$} &
\multicolumn{1}{c|}{$e_3$} \\
\multicolumn{1}{c|}{$e_5$} & \multicolumn{1}{|c}{} &
\multicolumn{1}{c}{} & \multicolumn{1}{c}{} & \multicolumn{1}{c}{}
& \multicolumn{1}{|c}{} & \multicolumn{1}{c|}{$-1$} &
\multicolumn{1}{|c}{$-\lambda'e_3-\delta'e_7$} &
\multicolumn{1}{c|}{$\lambda'e_2+\delta'e_6$}
\\ \cline{6-9} \multicolumn{1}{c|}{$e_6$} & \multicolumn{1}{|c}{}
& \multicolumn{1}{c}{} & \multicolumn{1}{c}{} &
\multicolumn{1}{c}{} & \multicolumn{1}{c}{} & \multicolumn{1}{c}{}
& \multicolumn{1}{|c}{$-1$} &
\multicolumn{1}{c|}{$-\lambda'e_1-\delta'e_5$} \\
\multicolumn{1}{c|}{$e_7$} & \multicolumn{1}{|c}{} &
\multicolumn{1}{c}{} & \multicolumn{1}{c}{} & \multicolumn{1}{c}{}
& \multicolumn{1}{c}{} & \multicolumn{1}{c}{} &
\multicolumn{1}{|c}{} & \multicolumn{1}{c|}{$-1$}
\\ \cline{2-9}
\end{tabular} \]

\vspace{0.6cm} o\`u \ $\lambda'=2\lambda-1,$ \
$\alpha'=(2\lambda-1)(2\alpha-1),$ \ $\beta'=(2\lambda-1)\beta$ \
et \ $\delta'=(2\lambda-1)\delta.$

\vspace{0.3cm} \hspace{0.3cm} La table de multiplication de $A,$
par rapport \`a la nouvelle base \ $u,$ $y_1,$ $z_1,$ $y_2,$
$z_2,$ $y_3,$ $z_3$ \ o\`u $u=e_4,$ $y_i=e_i,$ $z_i=-e_{i+4},$
$i\in\{1,2,3\}$ \ est donn\'ee par:

\vspace{0.5cm}
\[ \begin{tabular}{ccccccccc} \\
\multicolumn{1}{c}{} & \multicolumn{1}{c}{$1$} &
\multicolumn{1}{c}{$u$} & \multicolumn{1}{c}{$y_1$} &
\multicolumn{1}{c}{$z_1$} & \multicolumn{1}{c}{$y_2$} &
\multicolumn{1}{c}{$z_2$} & \multicolumn{1}{c}{$y_3$} &
\multicolumn{1}{c}{$z_3$} \\ \cline{2-9} \multicolumn{1}{c|}{$1$}
& \multicolumn{1}{|c}{$1$} & \multicolumn{1}{c|}{$u$} &
\multicolumn{1}{|c}{$y_1$} & \multicolumn{1}{c|}{$z_1$} &
\multicolumn{1}{|c}{$y_2$} & \multicolumn{1}{c|}{$z_2$} &
\multicolumn{1}{|c}{$y_3$} & \multicolumn{1}{c|}{$z_3$} \\
\multicolumn{1}{c|}{$u$} & \multicolumn{1}{|c}{} &
\multicolumn{1}{c|}{$-1$} & \multicolumn{1}{|c}{$z_1$} &
\multicolumn{1}{c|}{$-y_1$} & \multicolumn{1}{|c}{$z_2$} &
\multicolumn{1}{c|}{$-y_2$} & \multicolumn{1}{|c}{$z_3$} &
\multicolumn{1}{c|}{$-y_3$} \\ \cline{2-9}
\multicolumn{1}{c|}{$y_1$} & \multicolumn{1}{|c}{} &
\multicolumn{1}{c|}{} & \multicolumn{1}{|c}{$-1$} &
\multicolumn{1}{c|}{$u$} &
\multicolumn{1}{|c}{$\alpha'y_3-\beta'z_3$} &
\multicolumn{1}{c|}{$-\beta'y_3-\lambda'z_3$} &
\multicolumn{1}{|c}{$-\alpha'y_2+\beta'z_2$} &
\multicolumn{1}{c|}{$\beta'y_2+\lambda'z_2$} \\
\multicolumn{1}{c|}{$z_1$} & \multicolumn{1}{|c}{} &
\multicolumn{1}{c|}{} & \multicolumn{1}{|c}{} &
\multicolumn{1}{c|}{$-1$} &
\multicolumn{1}{|c}{$-\beta'y_3-\lambda'z_3$} &
\multicolumn{1}{c|}{$-\lambda'y_3+\delta'z_3$} &
\multicolumn{1}{|c}{$\beta'y_2+\lambda'z_2$} &
\multicolumn{1}{c|}{$\lambda'y_2-\delta'z_2$} \\
\cline{4-9} \multicolumn{1}{c|}{$y_2$} & \multicolumn{1}{|c}{} &
\multicolumn{1}{c}{} & \multicolumn{1}{c}{} &
\multicolumn{1}{c|}{} & \multicolumn{1}{|c}{$-1$} &
\multicolumn{1}{c|}{$u$} &
\multicolumn{1}{|c}{$\alpha'y_1-\beta'z_1$} &
\multicolumn{1}{c|}{$-\beta'y_1-\lambda'z_1$} \\
\multicolumn{1}{c|}{$z_2$} & \multicolumn{1}{|c}{} &
\multicolumn{1}{c}{} & \multicolumn{1}{c}{} & \multicolumn{1}{c}{}
& \multicolumn{1}{|c}{} & \multicolumn{1}{c|}{$-1$} &
\multicolumn{1}{|c}{$-\beta'y_1-\lambda'z_1$} &
\multicolumn{1}{c|}{$-\lambda'y_1+\delta'z_1$}
\\ \cline{6-9} \multicolumn{1}{c|}{$y_3$} & \multicolumn{1}{|c}{}
& \multicolumn{1}{c}{} & \multicolumn{1}{c}{} &
\multicolumn{1}{c}{} & \multicolumn{1}{c}{} & \multicolumn{1}{c}{}
& \multicolumn{1}{|c}{$-1$} & \multicolumn{1}{c|}{$u$} \\
\multicolumn{1}{c|}{$z_3$} & \multicolumn{1}{|c}{} &
\multicolumn{1}{c}{} & \multicolumn{1}{c}{} & \multicolumn{1}{c}{}
& \multicolumn{1}{c}{} & \multicolumn{1}{c}{} &
\multicolumn{1}{|c}{} & \multicolumn{1}{c|}{$-1$}
\\ \cline{2-9}
\end{tabular} \]

\vspace{0.5cm} Ainsi, $A$ est de division lin\'eaire si et
seulement si \ $\beta'+\alpha'\lambda'>0,$ \ et \[
(\alpha'\delta'-\beta'\lambda')^2<4(\beta'^2+\alpha'\lambda')(\lambda'^2+\beta'\delta')
\] \ i.e. \ $\lambda\neq\frac{1}{2},$ \
$\beta^2+2\alpha-1>\frac{1}{2}$ \ et \ $\Big(
(2\alpha-1)\delta-\beta\Big)^2<4(\beta^2+2\alpha-1)(1+\beta\delta).\Box$

\vspace{0.3cm} \begin{lemma} Si $\nu=\pi=0$ et $\theta\neq 0,$
alors $A$ est de division lin\'eaire si et seulement si \
$\beta\gamma-\alpha\mu,$ \ $\beta\lambda-\alpha\eta,$ \
$\gamma\lambda-\alpha\sigma>0,$ \
$c(\gamma\lambda-\alpha\sigma)\theta^2+b(\beta\lambda-\alpha\eta)\omega^2>
b\alpha\rho\omega\theta$ \ et \[
bc(\alpha\delta-\beta\sigma-\lambda\mu+\gamma\eta)^2+ab(\beta\lambda-\alpha\eta)\omega^2
+b^2(\beta\gamma-\alpha\mu)\rho^2+ac(\gamma\lambda-\alpha\sigma)\theta^2<
ab\alpha\rho\omega\theta+4bc(\beta\lambda-\alpha\eta)(\sigma\mu-\gamma\delta).
\]
\end{lemma}

\vspace{0.1cm} {\bf Preuve.} On suppose que $A$ est de division
lin\'eaire et on distingue les deux cas suivants:

\vspace{0.3cm} {\bf 1) Si \ $\alpha\neq 0,$} on consid\`ere
l'alg\`ebre $A_1=A(\varphi),$ o\`u $\varphi$ est l'automorphisme
de $W$ d\'efini dans le Lemme {\bf 4.11}. On pose, selon la Note
{\bf 4.12}, \begin{eqnarray*} \omega_0 &=&
(\omega^2+\theta^2)^{\frac{1}{2}}, \\
u' &=& y_1, \\
y_1' &=& \varepsilon\omega_0^{-1}(\omega y_2+\theta y_3), \\
z_1' &=& -z_1, \\
y_2' &=& \varepsilon\omega_0^{-1}(\omega y_2+\theta y_3), \\
z_2' &=& \omega_0^{-1}(-\theta y_2+\omega y_3), \\
y_3' &=& -\varepsilon z_2.
\end{eqnarray*}

On obtient la table:

\vspace{0.3cm}
{\footnotesize \[ \begin{tabular}{ccccccccc} \\
\multicolumn{1}{c}{} & \multicolumn{1}{c}{$1$} &
\multicolumn{1}{c}{$u'$} & \multicolumn{1}{c}{$y_1'$} &
\multicolumn{1}{c}{$z_1'$} & \multicolumn{1}{c}{$y_2'$} &
\multicolumn{1}{c}{$z_2'$} & \multicolumn{1}{c}{$y_3'$} &
\multicolumn{1}{c}{$z_3$}
\\ \cline{2-9}
\multicolumn{1}{c|}{$1$} & \multicolumn{1}{|c}{$1$} &
\multicolumn{1}{c|}{$u'$} & \multicolumn{1}{|c}{$y_1'$} &
\multicolumn{1}{c|}{$z_1'$} & \multicolumn{1}{|c}{$y_2'$} &
\multicolumn{1}{c|}{$z_2'$} & \multicolumn{1}{|c}{$y_3'$} &
\multicolumn{1}{c|}{$z_3$}
\\
\multicolumn{1}{c|}{$u'$} & \multicolumn{1}{|c}{} &
\multicolumn{1}{c|}{$-1$} & \multicolumn{1}{|c}{$az_1'$} &
\multicolumn{1}{c|}{$-ay_1'$} &
\multicolumn{1}{|c}{$|\alpha|z_2'$} &
\multicolumn{1}{c|}{$-|\alpha|y_2'$} &
\multicolumn{1}{|c}{$|\mu'|z_3$} &
\multicolumn{1}{c|}{$-|\mu'|y_3'$}
\\ \cline{2-9}
\multicolumn{1}{c|}{$y_1'$} & \multicolumn{1}{|c}{} &
\multicolumn{1}{c|}{} & \multicolumn{1}{|c}{$-1$} &
\multicolumn{1}{c|}{$au'$} &
\multicolumn{1}{|c}{$\alpha_0y_3'+\beta_0z_3$} &
\multicolumn{1}{c|}{$\gamma_0y_3'+\mu_0z_3$} &
\multicolumn{1}{|c}{$-\alpha_0y_2'-\gamma_0z_2'$} &
\multicolumn{1}{c|}{$-\beta_0y_2'-\mu_0z_2'$}
\\
\multicolumn{1}{c|}{$z_1'$} & \multicolumn{1}{|c}{} &
\multicolumn{1}{c|}{} & \multicolumn{1}{|c}{} &
\multicolumn{1}{c|}{$-1$} &
\multicolumn{1}{|c}{$\lambda_0y_3'+\eta_0z_3$} &
\multicolumn{1}{c|}{$\sigma_0y_3'+\delta_0z_3$} &
\multicolumn{1}{|c}{$-\lambda_0y_2'-\sigma_0z_2'+\varepsilon\delta'z_3$}
&
\multicolumn{1}{c|}{$-\eta_0y_2'-\delta_0z_2'-\varepsilon\delta'z_3$}
\\ \cline{4-9}
\multicolumn{1}{c|}{$y_2'$} & \multicolumn{1}{|c}{} &
\multicolumn{1}{c}{} & \multicolumn{1}{c}{} &
\multicolumn{1}{c|}{} & \multicolumn{1}{|c}{$-1$} &
\multicolumn{1}{c|}{$|\alpha|u'$} &
\multicolumn{1}{|c}{$\alpha_0y_1'+\lambda_0z_1'$} &
\multicolumn{1}{c|}{$\beta_0y_1'+\eta_0z_1'$}
\\ \multicolumn{1}{c|}{$z_2'$} & \multicolumn{1}{|c}{} &
\multicolumn{1}{c}{} & \multicolumn{1}{c}{} & \multicolumn{1}{c}{}
& \multicolumn{1}{|c}{} & \multicolumn{1}{c|}{$-1$} &
\multicolumn{1}{|c}{$\gamma_0y_1'+\sigma_0z_1'+\varepsilon\omega_0z_3$}
&
\multicolumn{1}{c|}{$\mu_0y_1'+\delta_0z_1'-\varepsilon\omega_0y_3'$}
\\ \cline{6-9} \multicolumn{1}{c|}{$y_3'$} & \multicolumn{1}{|c}{}
& \multicolumn{1}{c}{} & \multicolumn{1}{c}{} &
\multicolumn{1}{c}{} & \multicolumn{1}{c}{} & \multicolumn{1}{c}{}
& \multicolumn{1}{|c}{$-1$} &
\multicolumn{1}{c|}{$|\mu'|u+\varepsilon\delta'z_1'+\varepsilon\omega_0z_2'$} \\
\multicolumn{1}{c|}{$z_3$} & \multicolumn{1}{|c}{} &
\multicolumn{1}{c}{} & \multicolumn{1}{c}{} & \multicolumn{1}{c}{}
& \multicolumn{1}{c}{} & \multicolumn{1}{c}{} &
\multicolumn{1}{|c}{} & \multicolumn{1}{c|}{$-1$}
\\ \cline{2-9}
\end{tabular} \]}

\vspace{0.3cm} o\`u \begin{eqnarray*} \alpha_0 &=&
-\omega_0^{-1}\omega b, \\
\beta_0 &=& \omega_0^{-1}\varepsilon\theta c, \\
\gamma_0 &=& \omega_0^{-1}\varepsilon\theta b, \\
\mu_0 &=& \omega_0^{-1}\omega c, \\
\lambda_0 &=& -\omega_0^{-1}\theta\sigma', \\
\eta_0 &=& -\omega_0^{-1}\varepsilon(\omega\eta'+\theta\rho), \\
\sigma_0 &=& -\omega_0^{-1}\varepsilon\omega\sigma', \\
\delta_0 &=& \omega_0^{-1}(\theta\eta'-\omega\rho).
\end{eqnarray*}

Ainsi $A_1$ est de division lin\'eaire si et seulement si \
$\beta_0\lambda_0-\alpha_0\eta_0,$
$\beta_0\gamma_0-\alpha_0\mu_0,$
$\gamma_0\lambda_0-\alpha_0\sigma_0>0$ \ et \
$|\mu'|(\alpha_0\delta_0-\beta_0\sigma_0-\lambda_0\mu_0+\gamma_0\eta_0)^2+a(\beta_0\lambda_0-\alpha_0\eta_0)
(\varepsilon\omega_0)^2+|\alpha|(\beta_0\gamma_0-\alpha_0\mu_0)(\varepsilon\delta')^2<
4|\mu'|(\beta_0\lambda_0-\alpha_0\eta_0)(\sigma_0\mu_0-\gamma_0\delta_0)$
\ i.e. $A$ est de division lin\'eaire si et seulement si \[
\beta\lambda-\alpha\eta, \ \beta\gamma-\alpha\mu, \
\gamma\lambda-\alpha\sigma>0, \
c(\gamma\lambda-\alpha\sigma)\theta^2+b(\beta\lambda-\alpha\eta)\omega^2>b\alpha\rho\omega\theta
\ \mbox{ et }
\] \[ bc(\alpha\delta-\beta\sigma-\lambda\mu+\gamma\eta)^2+ab(\beta\lambda-\alpha\eta)\omega^2+b^2(\beta\gamma-\alpha\mu)
\rho^2+ac(\gamma\lambda-\alpha\sigma)\theta^2< \hspace{5cm} \] \[
\hspace{8cm}
ab\alpha\rho\omega\theta+4bc(\beta\lambda-\alpha\eta)(\sigma\mu-\gamma\delta).
\]

\vspace{0.3cm} {\bf 2) Si \ $\alpha=0,$} la condition $\gamma\neq
0$ est n\'ecessaire pour que $A$ soit de division lin\'eaire. On
consid\`ere alors l'automorphisme $\varphi'$ de $W$ dont la
matrice, par rapport \`a la base ${\cal B}$ est:

\vspace{0.5cm}
\[ \left(
\begin{array}{lllllll}
1 &   &                    &   &     &     & -\theta b^{-1}     \\
  & 1 & -\sigma\gamma^{-1} &   &     &     & -\omega\gamma^{-1} \\
  &   & \hspace{0.5cm}  1  &   &     &     &                    \\
  &   &                    & 1 &     &     &                    \\
  &   &                    &   & 1   &     &                    \\
  &   &                    &   &     &  1  & -\mu\gamma^{-1}    \\
  &   &                    &   &     &     & \hspace{0.5cm} 1
\end{array}
\right)  \]

\vspace{0.5cm} On obtient la multiplication de l'alg\`ebre
$A(\varphi'):$

\vspace{0.3cm}
{\footnotesize \[ \begin{tabular}{ccccccccc} \\
\multicolumn{1}{c}{} & \multicolumn{1}{c}{$1$} &
\multicolumn{1}{c}{$z_2$} & \multicolumn{1}{c}{$u$} &
\multicolumn{1}{c}{$y_2$} & \multicolumn{1}{c}{$y_1$} &
\multicolumn{1}{c}{$y_3$} & \multicolumn{1}{c}{$z_3$} &
\multicolumn{1}{c}{$z_1$}
\\ \cline{2-9}
\multicolumn{1}{c|}{$1$} & \multicolumn{1}{|c}{$1$} &
\multicolumn{1}{c|}{$z_2$} & \multicolumn{1}{|c}{$u$} &
\multicolumn{1}{c|}{$y_2$} & \multicolumn{1}{|c}{$y_1$} &
\multicolumn{1}{c|}{$y_3$} & \multicolumn{1}{|c}{$z_3$} &
\multicolumn{1}{c|}{$z_1$}
\\
\multicolumn{1}{c|}{$z_2$} & \multicolumn{1}{|c}{} &
\multicolumn{1}{c|}{$-1$} & \multicolumn{1}{|c}{$by_2$} &
\multicolumn{1}{c|}{$-bu$} & \multicolumn{1}{|c}{$-\gamma y_3$} &
\multicolumn{1}{c|}{$\gamma y_1$} &
\multicolumn{1}{|c}{$\delta_1z_1$} &
\multicolumn{1}{c|}{$-\delta_1z_3$}
\\ \cline{2-9}
\multicolumn{1}{c|}{$u$} & \multicolumn{1}{|c}{} &
\multicolumn{1}{c|}{} & \multicolumn{1}{|c}{$-1$} &
\multicolumn{1}{c|}{$bz_2$} & \multicolumn{1}{|c}{$az_1$} &
\multicolumn{1}{c|}{$cz_3$} &
\multicolumn{1}{|c}{$-cy_3-a\omega\gamma^{-1}z_1$} &
\multicolumn{1}{c|}{$-ay_1+a\omega\gamma^{-1}z_3$}
\\ \multicolumn{1}{c|}{$y_2$} & \multicolumn{1}{|c}{} &
\multicolumn{1}{c|}{} & \multicolumn{1}{|c}{} &
\multicolumn{1}{c|}{$-1$} & \multicolumn{1}{|c}{$-\beta z_3$} &
\multicolumn{1}{c|}{$\lambda z_1$} & \multicolumn{1}{|c}{$\beta
y_1+\eta_1z_1$} & \multicolumn{1}{c|}{$-\lambda y_3-\eta_1z_3$}
\\ \cline{4-9}
\multicolumn{1}{c|}{$y_1$} & \multicolumn{1}{|c}{} &
\multicolumn{1}{c}{} & \multicolumn{1}{c}{} &
\multicolumn{1}{c|}{} & \multicolumn{1}{|c}{$-1$} &
\multicolumn{1}{c|}{$-\gamma z_2$} & \multicolumn{1}{|c}{$-\beta
y_2+ab^{-1}\theta z_1$} & \multicolumn{1}{c|}{$au-ab^{-1}\theta
z_3$}
\\ \multicolumn{1}{c|}{$y_3$} & \multicolumn{1}{|c}{} &
\multicolumn{1}{c}{} & \multicolumn{1}{c}{} & \multicolumn{1}{c}{}
& \multicolumn{1}{|c}{} & \multicolumn{1}{c|}{$-1$} &
\multicolumn{1}{|c}{$cu+\rho z_1$} & \multicolumn{1}{c|}{$\lambda
y_2-\rho z_3$}
\\ \cline{6-9} \multicolumn{1}{c|}{$z_3$} & \multicolumn{1}{|c}{}
& \multicolumn{1}{c}{} & \multicolumn{1}{c}{} &
\multicolumn{1}{c}{} & \multicolumn{1}{c}{} & \multicolumn{1}{c}{}
& \multicolumn{1}{|c}{$-1$} &
\multicolumn{1}{c|}{$\delta_1z_2-a\omega\gamma^{-1}u+\eta_1y_2+ab^{-1}\theta y_1+\rho y_3$} \\
\multicolumn{1}{c|}{$z_1$} & \multicolumn{1}{|c}{} &
\multicolumn{1}{c}{} & \multicolumn{1}{c}{} & \multicolumn{1}{c}{}
& \multicolumn{1}{c}{} & \multicolumn{1}{c}{} &
\multicolumn{1}{|c}{} & \multicolumn{1}{c|}{$-1$}
\\ \cline{2-9}
\end{tabular} \]}

\vspace{0.3cm} o\`u \ $\delta_1=\delta-\sigma\mu\gamma^{-1}$ \ et
\ $\eta_1=\eta-(\lambda\mu+\beta\sigma)\gamma^{-1}.$ Comme
$A(\varphi')$ est de division lin\'eaire, on a \
$\gamma\delta_1=\gamma\delta-\sigma\mu<0$ \ car l'\'egalit\'e \
$(xz_2+u)(-cy_1+x\gamma z_3)=(\gamma\delta_1 x^2-\omega ax-ac)z_1$
\ a lieu pour tout $x\in\rit$ et $\delta_1\neq 0.$ Il existe alors
$\varepsilon'\in\{1,-1\}$ tel que
$|\delta_1|=\varepsilon'\delta_1,$ \
$|\gamma|=-\varepsilon'\gamma$ \ et on pose \begin{eqnarray*}
\rho_1 &=& (a^2b^{-2}\theta^2+\rho^2)^{\frac{1}{2}}, \\
y_3' &=& \rho_1^{-1}(ab^{-1}\theta y_1+\rho y_3), \\
y_1' &=& \varepsilon'\rho_1^{-1}(\rho y_1-ab^{-1}\theta y_3), \\
z_3' &=& \varepsilon'z_3 \end{eqnarray*}

et si \ $(\omega,\eta_1)\neq (0,0),$

\begin{eqnarray*}
\eta^* &=& (a^2\omega^2\gamma^{-2}+\eta_1^2)^{\frac{1}{2}}, \\
u' &=& \eta^{*^{-1}}(\eta_1u+a\omega\gamma^{-1}y_2), \\
y_2' &=& \eta^{*^{-1}}(-a\omega\gamma^{-1}u+\eta_1y_2)
\end{eqnarray*}

et on obtient la table:

\vspace{0.3cm}
{\footnotesize \[ \begin{tabular}{ccccccccc} \\
\multicolumn{1}{c}{} & \multicolumn{1}{c}{$1$} &
\multicolumn{1}{c}{$z_2$} & \multicolumn{1}{c}{$u'$} &
\multicolumn{1}{c}{$y_2'$} & \multicolumn{1}{c}{$y_1'$} &
\multicolumn{1}{c}{$y_3'$} & \multicolumn{1}{c}{$z_3'$} &
\multicolumn{1}{c}{$z_1$}
\\ \cline{2-9}
\multicolumn{1}{c|}{$1$} & \multicolumn{1}{|c}{$1$} &
\multicolumn{1}{c|}{$z_2$} & \multicolumn{1}{|c}{$u'$} &
\multicolumn{1}{c|}{$y_2'$} & \multicolumn{1}{|c}{$y_1'$} &
\multicolumn{1}{c|}{$y_3'$} & \multicolumn{1}{|c}{$z_3'$} &
\multicolumn{1}{c|}{$z_1$}
\\
\multicolumn{1}{c|}{$z_2$} & \multicolumn{1}{|c}{} &
\multicolumn{1}{c|}{$-1$} & \multicolumn{1}{|c}{$by_2'$} &
\multicolumn{1}{c|}{$-bu'$} & \multicolumn{1}{|c}{$|\gamma|y_3'$}
& \multicolumn{1}{c|}{$-|\gamma|y_1'$} &
\multicolumn{1}{|c}{$|\delta_1|z_1$} &
\multicolumn{1}{c|}{$-|\delta_1|z_3'$}
\\ \cline{2-9}
\multicolumn{1}{c|}{$u'$} & \multicolumn{1}{|c}{} &
\multicolumn{1}{c|}{} & \multicolumn{1}{|c}{$-1$} &
\multicolumn{1}{c|}{$bz_2$} &
\multicolumn{1}{|c}{$\alpha_1z_3'+\beta_1z_1$} &
\multicolumn{1}{c|}{$\gamma_1z_3'+\mu_1z_1$} &
\multicolumn{1}{|c}{$-\alpha_1y_1'-\gamma_1y_3'$} &
\multicolumn{1}{c|}{$-\beta_1y_1'-\mu_1y_3'$}
\\ \multicolumn{1}{c|}{$y_2'$} & \multicolumn{1}{|c}{} &
\multicolumn{1}{c|}{} & \multicolumn{1}{|c}{} &
\multicolumn{1}{c|}{$-1$} &
\multicolumn{1}{|c}{$\lambda_1z_3'+\eta_2z_1$} &
\multicolumn{1}{c|}{$\sigma_1z_3'+\delta_2z_1$} &
\multicolumn{1}{|c}{$-\lambda_1y_1'-\sigma_1y_3'+\varepsilon'\eta^*z_1$}
&
\multicolumn{1}{c|}{$-\eta_2y_1'-\delta_2y_3'-\varepsilon'\eta^*z_3'$}
\\ \cline{4-9}
\multicolumn{1}{c|}{$y_1'$} & \multicolumn{1}{|c}{} &
\multicolumn{1}{c}{} & \multicolumn{1}{c}{} &
\multicolumn{1}{c|}{} & \multicolumn{1}{|c}{$-1$} &
\multicolumn{1}{c|}{$|\gamma|z_2$} &
\multicolumn{1}{|c}{$\alpha_1u'+\lambda_1y_2'$} &
\multicolumn{1}{c|}{$\beta_1u'+\eta_2y_2'$}
\\ \multicolumn{1}{c|}{$y_3'$} & \multicolumn{1}{|c}{} &
\multicolumn{1}{c}{} & \multicolumn{1}{c}{} & \multicolumn{1}{c}{}
& \multicolumn{1}{|c}{} & \multicolumn{1}{c|}{$-1$} &
\multicolumn{1}{|c}{$\gamma_1u'+\sigma_1y_2'+\varepsilon'\rho_1z_1$}
&
\multicolumn{1}{c|}{$\mu_1u'+\delta_2y_2'-\varepsilon'\rho_1z_3'$}
\\ \cline{6-9} \multicolumn{1}{c|}{$z_3'$} & \multicolumn{1}{|c}{}
& \multicolumn{1}{c}{} & \multicolumn{1}{c}{} &
\multicolumn{1}{c}{} & \multicolumn{1}{c}{} & \multicolumn{1}{c}{}
& \multicolumn{1}{|c}{$-1$} &
\multicolumn{1}{c|}{$|\delta_1|z_2+\varepsilon'\eta^*y_2'+\varepsilon'\rho_1y_3'$} \\
\multicolumn{1}{c|}{$z_1$} & \multicolumn{1}{|c}{} &
\multicolumn{1}{c}{} & \multicolumn{1}{c}{} & \multicolumn{1}{c}{}
& \multicolumn{1}{c}{} & \multicolumn{1}{c}{} &
\multicolumn{1}{|c}{} & \multicolumn{1}{c|}{$-1$}
\\ \cline{2-9}
\end{tabular} \]}

\vspace{0.4cm} o\`u \begin{eqnarray*} \alpha_1 &=&
-\eta^{*^{-1}}\rho_1^{-1}(ab^{-1}c\theta\eta_1+a\beta\rho\omega\gamma^{-1}), \\
\beta_1 &=& \eta^{*^{-1}}\rho_1^{-1}\varepsilon'(a\rho\eta_1-a^2b^{-1}\omega\gamma^{-1}\theta\lambda), \\
\gamma_1 &=& \eta^{*^{-1}}\rho_1^{-1}\varepsilon'(c\rho\eta_1-a^2b^{-1}\omega\gamma^{-1}\theta\beta), \\
\mu_1 &=& \eta^{*^{-1}}\rho_1^{-1}(a^2b^{-1}\theta\eta_1+a\omega\gamma^{-1}\rho\lambda), \\
\lambda_1 &=& \eta^{*^{-1}}\rho_1^{-1}(a^2b^{-1}c\omega\gamma^{-1}\theta-\rho\beta\eta_1), \\
\eta_2 &=& -\eta^{*^{-1}}\rho_1^{-1}\varepsilon'(a^2\omega\gamma^{-1}\rho+ab^{-1}\eta_1\theta\lambda), \\
\sigma_1 &=& -\eta^{*^{-1}}\rho_1^{-1}\varepsilon'(ac\omega\gamma^{-1}\rho+ab^{-1}\eta_1\theta\beta), \\
\delta_2 &=&
\eta^{*^{-1}}\rho_1^{-1}(-a^3b^{-1}\omega\gamma^{-1}\theta+\eta_1\rho\lambda).
\end{eqnarray*}

\vspace{0.3cm} Ainsi $A(\varphi')$ est de division lin\'eaire si
et seulement si \ $\beta_1\lambda_1-\alpha_1\eta_2,$
$\beta_1\gamma_1-\alpha_1\mu_1,$
$\gamma_1\lambda_1-\alpha_1\sigma_1>0$ \ et

\[ |\delta_1|(\alpha_1\delta_2-\beta_1\sigma_1-\lambda_1\mu_1+\gamma_1\eta_2)^2
+b(\beta_1\lambda_1-\alpha_1\eta_2)(\varepsilon'\rho_1)^2+|\gamma|(\beta_1\gamma_1-\alpha_1\mu_1)
(\varepsilon'\eta^*)^2< \hspace{1.5cm} \] \[ \hspace{10cm}
4|\delta_1|(\beta_1\lambda_1-\alpha_1\eta_2)(\sigma_1\mu_1-\gamma_1\delta_2)
\]

i.e. $A$ est de division lin\'eaire si et seulement si \
$\beta\lambda,$ $\beta\gamma,$ $\gamma\lambda>0$ \ et

\[ bc(-\beta\sigma-\lambda\mu+\gamma\eta)^2
+ab\beta\lambda\omega^2+b^2\beta\gamma\rho^2+ac\gamma\lambda\theta^2<
4bc\beta\lambda(\sigma\mu-\gamma\delta). \] Ce r\'esultat englobe
le cas $\eta^*=0.\Box$

\vspace{0.3cm} \begin{theorem} Soit $A$ une $\rit$-alg\`ebre de
Jordan non commutative de division lin\'eaire de dimension $8.$
Alors il existe une base \ ${\cal
B}=\{1,u,y_1,z_1,y_1,z_2,y_3,z_3\}$ de $A,$ trois param\`etres \
$a,$ $b,$ $c>0$ \ et treize autres \ $\alpha,$ $\beta,$ $\gamma,$
$\mu,$ $\lambda,$ $\eta,$ $\sigma,$ $\delta,$ $\nu,$ $\pi,$
$\rho,$ $\theta,$ $\omega$ \ pour lesquelles la multiplication de
$A$ est donn\'ee par la Table {\bf 1}. De plus, une
$\rit$-alg\`ebre dont la multiplication est donn\'ee par la Table
{\bf 1} est de Jordan non commutative. Elle est de division
lin\'eaire si et seulement si \ $\beta\lambda-\alpha\eta,$
$\gamma\lambda-\alpha\sigma >0,$

\begin{eqnarray*}
bc(\alpha\delta-\beta\sigma-\lambda\mu+\gamma\eta)^2+ab(\beta\lambda-\alpha\eta)\omega^2
+ac(\gamma\lambda-\alpha\sigma)\theta^2+(\beta\gamma-\alpha\mu)(b\rho-\pi
c)^2+b^2(\sigma\eta-\lambda\delta)\nu^2 \\
+b\nu(\alpha\delta-\beta\sigma+\lambda\mu-\gamma\eta)(b\rho-\pi
c)< a\theta\omega\Big( \alpha(b\rho-\pi
c)-b\lambda\nu\Big)+4bc(\beta\lambda-\alpha\eta)(\sigma\mu-\gamma\delta)
\end{eqnarray*}

et l'une des quatre situations suivantes a lieu:
\begin{enumerate} \item $\nu=\theta=0$ et
$\beta\gamma-\alpha\mu>0.$ \item $\nu=0,$ $\theta\neq 0,$
$\beta\gamma-\alpha\mu>0$ et
$c(\gamma\lambda-\alpha\sigma)\theta^2+b(\beta\lambda-\alpha\eta)\omega^2>\alpha\omega\theta(b\rho-\pi
c).$ \item $\nu\neq 0,$ $\theta=0$ et
$(\beta\gamma-\alpha\mu)(b\rho-\pi
c)^2+b^2(\sigma\eta-\lambda\delta)\nu^2+b\nu(\alpha\delta-\beta\sigma+\lambda\mu-\gamma\eta)(b\rho-\pi
c)>0.$ \item $\nu\theta\neq 0,$
$c(\gamma\lambda-\alpha\sigma)\theta^2+b(\beta\lambda-\alpha\eta)\omega^2>\omega\theta\Big(
\alpha(b\rho-\pi c)-b\lambda\nu\Big)$ et \[
(\beta\gamma-\alpha\mu)(b\rho-\pi
c)^2+b^2(\sigma\eta-\lambda\delta)\nu^2+b\nu(\alpha\delta-\beta\sigma+\lambda\mu-\gamma\eta)(b\rho-\pi
c)>0. \]
\end{enumerate}
\end{theorem}

\vspace{0.1cm} {\bf Preuve.} Il reste seulement \`a \'etablir la
derni\`ere proposition. On distingue les deux cas suivants:

\vspace{0.2cm} {\bf 1) Si $\pi=0.$} On peut supposer $\nu\neq 0$
et, en posant $\rho'=(\nu^2+\rho^2)^{\frac{1}{2}},$
$y_1'=\rho'^{-1}(\rho y_1-\nu z_1)$ et $z_1'=\rho'^{-1}(\nu
y_1+\rho z_1)$ \ on obtient la table:

\vspace{0.2cm}
{\footnotesize \[ \begin{tabular}{ccccccccc} \\
\multicolumn{1}{c}{} & \multicolumn{1}{c}{$1$} &
\multicolumn{1}{c}{$u$} & \multicolumn{1}{c}{$y_1'$} &
\multicolumn{1}{c}{$z_1'$} & \multicolumn{1}{c}{$y_2$} &
\multicolumn{1}{c}{$z_2$} & \multicolumn{1}{c}{$y_3$} &
\multicolumn{1}{c}{$z_3$}
\\ \cline{2-9}
\multicolumn{1}{c|}{$1$} & \multicolumn{1}{|c}{$1$} &
\multicolumn{1}{c|}{$u$} & \multicolumn{1}{|c}{$y_1'$} &
\multicolumn{1}{c|}{$z_1'$} & \multicolumn{1}{|c}{$y_2$} &
\multicolumn{1}{c|}{$z_2$} & \multicolumn{1}{|c}{$y_3$} &
\multicolumn{1}{c|}{$z_3$}
\\
\multicolumn{1}{c|}{$u$} & \multicolumn{1}{|c}{} &
\multicolumn{1}{c|}{$-1$} & \multicolumn{1}{|c}{$az_1'$} &
\multicolumn{1}{c|}{$-ay_1'$} & \multicolumn{1}{|c}{$bz_2$} &
\multicolumn{1}{c|}{$-by_2$} & \multicolumn{1}{|c}{$cz_3$} &
\multicolumn{1}{c|}{$-cy_3$}
\\ \cline{2-9}
\multicolumn{1}{c|}{$y_1$} & \multicolumn{1}{|c}{} &
\multicolumn{1}{c|}{} & \multicolumn{1}{|c}{$-1$} &
\multicolumn{1}{c|}{$au$} &
\multicolumn{1}{|c}{$\alpha_0y_3+\beta_0z_3$} &
\multicolumn{1}{c|}{$\gamma_0y_3+\mu_0z_3$} &
\multicolumn{1}{|c}{$-\alpha_0y_2-\gamma_0z_2$} &
\multicolumn{1}{c|}{$-\beta_0y_2-\mu_0z_2$}
\\
\multicolumn{1}{c|}{$z_1$} & \multicolumn{1}{|c}{} &
\multicolumn{1}{c|}{} & \multicolumn{1}{|c}{} &
\multicolumn{1}{c|}{$-1$} &
\multicolumn{1}{|c}{$\lambda_0y_3+\eta_0z_3$} &
\multicolumn{1}{c|}{$\sigma_0y_3+\delta_0z_3$} &
\multicolumn{1}{|c}{$-\lambda_0y_2-\sigma_0z_2+\rho' z_3$} &
\multicolumn{1}{c|}{$-\eta_0y_2-\delta_0z_2-\rho' z_3$}
\\ \cline{4-9}
\multicolumn{1}{c|}{$y_2$} & \multicolumn{1}{|c}{} &
\multicolumn{1}{c}{} & \multicolumn{1}{c}{} &
\multicolumn{1}{c|}{} & \multicolumn{1}{|c}{$-1$} &
\multicolumn{1}{c|}{$bu+\theta z_3$} &
\multicolumn{1}{|c}{$\alpha_0y_1'+\lambda_0z_1'$} &
\multicolumn{1}{c|}{$\beta_0y_1'+\eta_0z_1'-\theta z_2$}
\\ \multicolumn{1}{c|}{$z_2$} & \multicolumn{1}{|c}{} &
\multicolumn{1}{c}{} & \multicolumn{1}{c}{} & \multicolumn{1}{c}{}
& \multicolumn{1}{|c}{} & \multicolumn{1}{c|}{$-1$} &
\multicolumn{1}{|c}{$\gamma_0y_1'+\sigma_0z_1'+\omega z_3$} &
\multicolumn{1}{c|}{$\mu_0y_1'+\delta_0z_1'+\theta y_2-\omega
y_3$} \\ \cline{6-9} \multicolumn{1}{c|}{$y_3$} &
\multicolumn{1}{|c}{} & \multicolumn{1}{c}{} &
\multicolumn{1}{c}{} & \multicolumn{1}{c}{} & \multicolumn{1}{c}{}
& \multicolumn{1}{c}{} & \multicolumn{1}{|c}{$-1$} &
\multicolumn{1}{c|}{$cu+\rho' z_1'+\omega z_2$} \\
\multicolumn{1}{c|}{$z_3$} & \multicolumn{1}{|c}{} &
\multicolumn{1}{c}{} & \multicolumn{1}{c}{} & \multicolumn{1}{c}{}
& \multicolumn{1}{c}{} & \multicolumn{1}{c}{} &
\multicolumn{1}{|c}{} & \multicolumn{1}{c|}{$-1$}
\\ \cline{2-9}
\end{tabular} \]}

\vspace{0.3cm} o\`u
\[ \alpha_0=\rho'^{-1}(\alpha\rho-\lambda\nu), \hspace{0.3cm} \beta_0=\rho'^{-1}(\beta\rho-\eta\nu),
\hspace{0.3cm} \gamma_0=\rho'^{-1}(\gamma\rho-\sigma\nu),
\hspace{0.3cm} \mu_0=\rho'^{-1}(\mu\rho-\delta\nu), \] \[
\lambda_0=\rho'^{-1}(\alpha\nu+\lambda\rho), \hspace{0.3cm}
\eta_0=\rho'^{-1}(\beta\nu+\eta\rho), \hspace{0.3cm}
\sigma_0=\rho'^{-1}(\gamma\nu+\sigma\rho), \hspace{0.3cm}
\delta_0=\rho'^{-1}(\mu\nu+\delta\rho). \]

\vspace{0.2cm} Ainsi $A$ est de division lin\'eaire si et
seulement si \ $\beta_0\gamma_0-\alpha_0\mu_0,$ \
$\beta_0\lambda_0-\alpha_0\eta_0,$ \
$\gamma_0\lambda_0-\alpha_0\sigma_0>0$ \ et l'une des deux
situations suivantes a lieu \begin{enumerate} \item $\theta=0$ \
et \[
c(\alpha_0\delta_0-\beta_0\sigma_0-\lambda_0\mu_0+\gamma_0\eta_0)^2+a(\beta_0\lambda_0-\alpha_0\eta_0)\omega^2
+b(\beta_0\gamma_0-\alpha_0\mu_0)\rho'^2<4c(\beta_0\lambda_0-\alpha_0\eta_0)(\sigma_0\mu_0-\gamma_0\delta_0).
\] \item $\theta\neq 0,$ \
$c(\gamma_0\lambda_0-\alpha_0\sigma_0)\theta^2+b(\beta_0\lambda_0-\alpha_0\eta_0)\omega^2>\alpha_0b\rho'\omega\theta$
\ et \[
bc(\alpha_0\delta_0-\beta_0\sigma_0-\lambda_0\mu_0+\gamma_0\eta_0)^2+ac(\gamma_0\lambda_0-\alpha_0\sigma_0)\theta^2
+ab(\beta_0\lambda_0-\alpha_0\eta_0)\omega^2+b^2(\beta_0\gamma_0-\alpha_0\mu_0)\rho'^2
\] \[ <ab\alpha_0\rho'\omega\theta+4bc(\beta_0\lambda_0-\alpha_0\eta_0)(\sigma_0\mu_0-\gamma_0\delta_0).
\] \end{enumerate}

\vspace{0.2cm} Ce qui donne \ $\beta\lambda-\alpha\eta,$ \
$\gamma\lambda-\alpha\sigma>0,$ \[
bc(\alpha\delta-\beta\sigma-\lambda\mu+\gamma\eta)^2+ab(\beta\lambda-\alpha\eta)\omega^2
+ac(\gamma\lambda-\alpha\sigma)\theta^2+b^2\Big(
(\beta\gamma-\alpha\mu)\rho^2+(\alpha\delta-\beta\sigma+\lambda\mu-\gamma\eta)\rho\nu+(\sigma\eta-\lambda\delta)\nu^2\Big)
\] \[ <ab(\alpha\rho-\lambda\nu)\theta\omega+4bc(\beta\lambda-\alpha\eta)(\sigma\mu-\gamma\delta) \]
et l'une des deux derni\`eres situations du Th\'eor\`eme.

\vspace{0.2cm} {\bf 2) Si $\pi$ est quelconque}, on consid\`ere
l'automorphisme $h$ de $W$ dont la matrice par rapport \`a la base
${\cal B}$ est

\vspace{0.3cm}
\[ \left(
\begin{array}{lllllll}
1 &   & -\pi b^{-1}      &   &     &     &   \\
  & 1 &                  &   &     &     &   \\
  &   & \hspace{0.5cm} 1 &   &     &     &   \\
  &   &                  & 1 &     &     &   \\
  &   &                  &   & 1   &     &   \\
  &   &                  &   &     &  1  &   \\
  &   &                  &   &     &     & 1
\end{array}
\right)  \]

\vspace{0.3cm} La multiplication de l'alg\`ebre $A(h)$ s'obtient
\`a partir de la Table {\bf 1} en faisant "$\pi=0$" et en
rempla\c{c}ant \ $\rho$ \ par \ $\rho-\pi b^{-1}c.$ Comme $A$ est
de d.l. si et seulement si $A(h)$ est de d.l., la d\'emonstration
s'ach\`eve en vertu du premier cas.$\Box$

\vspace{0.4cm} \hspace{0.3cm} Nous \'enon\c{c}ons maintenant le
Th\'eor\`eme de classification suivant:

\vspace{0.2cm} \begin{theorem} Les alg\`ebres r\'eelles de Jordan
non commutatives de division lin\'eaire de dimension $8$
s'obtiennent, \`a partie de l'alg\`ebre r\'eelle \ $\oit=\Big(
W,(.|.),\wedge\Big)$ de Cayley-Dickson, par isotopie vectorielle
et sont \`a isomorphisme pr\`es \ $\oit(s)$ o\`u $s$ est un
automorphisme sym\'etrique de l'espace euclidien \ $(W,-(.|.)),$
d\'efini positif. De plus, \ $\oit(s')\simeq$ $\oit(s)$ \ ($s$ et
$s'$ \'etant deux automorphismes sym\'etriques de l'espace
euclidien $(W,-(.|.)),$ d\'efinis positifs) si eu seulement si il
existe \ $f\in G_2$ tel que $\tilde{s}'=f^{-1}\tilde{s}f.\Box$
\end{theorem}

\vspace{0.1cm} {\bf Preuve.} En tenant compte du Corollaire {\bf
4.8}, des Lemmes {\bf 4.13}, {\bf 4.15} et du Th\'eor\`eme {\bf
4.16}, nous avons d\'emontr\'e que les $\rit$-alg\`ebres de Jordan
n.c. de d.l. de dimension $8$ s'obtiennent, \`a partir de
l'alg\`ebre r\'eelle \ $\oit=\Big( W,(.|.),\wedge\Big)$ de
Cayley-Dickson, par isotopie vectorielle. Si $A=\Big(
W,(.|.),\wedge\Big)$ est une telle alg\`ebre, il existe un
automorphisme $\varphi$ de l'espace vectoriel r\'eel $W$ tel que
$A=\oit(\varphi).$ D'apr\`es le Th\'eor\`eme de d\'ecomposition
polaire, $\varphi$ s'exprime comme un produit \ $sr$ \ d'un
automorphisme sym\'etrique $s$ de l'espace euclidien $(W,-(.|.)),$
d\'efini positif, et d'une isom\'etrie $r$ de $(W,-(.|.)).$ Ainsi
$A=\oit(sr)=\Big( \oit(s)\Big)(r)\simeq$ $\oit(s)$ \ (Proposition
{\bf 4.6 1)}. Soient maintenant $s$ et $s'$ deux automorphismes
sym\'etriques, d\'efinis positifs, de l'espace euclidien
$(W,-(.|.)).$ Alors \begin{eqnarray*} \oit(s')\simeq\oit(s)
&\Leftrightarrow& \mbox{ Il existe }
\varphi\in O_7(\rit) \mbox{ tel que } \overline{s'\varphi s^{-1}}\in G_2  \mbox{ (Corollaire {\bf 4.10 3)} } \\
&\Leftrightarrow& \mbox{ Il existe }
f\in G_2 \mbox{ tel que } s'^{-1}.f_{/W}.s\in O_7(\rit) \\
&\Leftrightarrow& \mbox{ Il existe }
f\in G_2 \mbox{ tel que } s'^2=f_{/W}.s^2.f^{-1}_{/W} \\
&\Leftrightarrow& \mbox{ Il existe }
f\in G_2 \mbox{ tel que } s'^2=(f_{/W}.s.f^{-1}_{/W})^2 \\
&\Leftrightarrow& \mbox{ Il existe } f\in G_2 \mbox{ tel que }
s'=f_{/W}sf^{-1}_{/W} \hspace{0.3cm} (f_{/W}sf^{-1}_{/W} \mbox{
sym\'etrique d\'efini positif) }.\Box
\end{eqnarray*}

\vspace{0.3cm} \begin{corollary} Les $\rit$-alg\`ebres de Jordan
non commutatives de division lin\'eaire de dimension finie $\geq
2$ s'obtiennent, \`a partir de \ $\cit,$ $\hit$ et \ $\oit$ par
isotopie vectorielle.$\Box$
\end{corollary}

\vspace{10cm}
\section{$\rit$-alg\`ebres de Jordan n.c. de d.l. de dimension 8 ayant un automorphisme non trivial}

\vspace{0.7cm} \subsection{Etude des $\rit$-alg\`ebres de Jordan
n.c. de d.l de dimension 8 qui poss\`edent une d\'erivation non
triviale}

\vspace{0.5cm} \hspace{0.3cm} En dimension finie et moyennant le
Th\'eor\`eme $(1,2,4,8)$ de Hopf-Kervaire-Milnor-Bott ([H 40], [Ke
58], [BM 58]), Benkart et Osborn [BO 81$_2$] d\'etermin\`erent
toutes les possibilit\'es pour l'alg\`ebre de Lie des
d\'erivations, d'une alg\`ebre r\'eelle de division lin\'eaire de
dimension finie, en \'etablissant le Th\'eor\`eme de
classification suivant:

\vspace{0.3cm} \begin{theorem} Soit $A$ une alg\`ebre r\'eelle de
division lin\'eaire de dimension finie. Alors \begin{enumerate}
\item Si $\dim(A)=1$ ou $2,$ alors \ $Der(A)=0.$ \item Si
$\dim(A)=4,$ alors \ $Der(A)$ est isomorphe \`a \ $su(2)$ ou
$\dim(A)\leq 1.$ \item Si $\dim(A)=8,$ alors \ $Der(A)$ est
isomorphe \`a l'une des alg\`ebres de Lie suivantes:
\begin{enumerate} \item $G_2$ compacte,
\item $su(3),$ \item $su(2)\oplus su(2),$ \item $su(2)\oplus N$
o\`u $N$ est une alg\`ebre ab\'elienne de dimension $\leq 1,$
\item $N,$ une alg\`ebre ab\'elienne de dimension $\leq 2.$
\end{enumerate}

De plus, toutes les possibilit\'es pr\'ec\'edentes se
r\'ealisent.$\Box$
\end{enumerate}
\end{theorem}

\vspace{0.5cm} \hspace{0.3cm} Benkart et Osborn [BO 81$_1$] ont
donn\'e ensuite une classification compl\`ete pour les alg\`ebres
r\'eelles de division lin\'eaire de dimension $4,$ dont
l'alg\`ebre de Lie des d\'erivations est $su(2),$ et pour les
alg\`ebres r\'eelles de division lin\'eaire de dimension $8,$ dont
l'alg\`ebre de Lie des d\'erivations est $G_2$ compacte, $su(3)$
ou $su(2)\oplus su(2).$ Ils ont donn\'e \'egalement des exemples
d'alg\`ebres r\'eelles de division lin\'eaire pour chacun des
autres cas de l'alg\`ebre de Lie des d\'erivations, puis ils ont
pos\'e, entre autres, le probl\`eme particulier de l'existence
d'une alg\`ebre r\'eelle de division lin\'eaire de dimension $8,$
dont l'alg\`ebre de Lie des d\'erivations est $su(2)$ et dont la
d\'ecomposition en $su(2)$-modules irr\'eductibles est de la
forme: \[ 1+1+3+3. \]

\vspace{0.1cm} \hspace{0.3cm} Dans ce paragraphe, nous \'etudions
les alg\`ebres r\'eelles de Jordan non commutatives, de division
lin\'eaire de dimension $8,$ dont l'alg\`ebre de Lie des
d\'erivations est non triviale, puis donnons une r\'eponse
affirmative au probl\`eme particulier pr\'ec\'edent.

\vspace{0.3cm} \hspace{0.3cm} Soient maintenant $A$ une
$\rit$-alg\`ebre de division lin\'eaire de dimension finie,
$\partial$ une d\'erivation de $A,$ $\sigma$ un nombre complexe et
$\overline{\sigma}$ son conjugu\'e. Benkart et Osborn
consid\`erent le sous-espace vectoriel

\vspace{0.2cm} \[ B_\sigma=\{x\in A: (\partial-\sigma
I_A)(\partial-\overline{\sigma}I_A)x=0\} \]

\vspace{0.2cm} et \'etablissent l'inclusion \ $B_\sigma
B_\tau\subseteq B_{\sigma+\tau}+B_{\overline{\sigma}+\tau}.$ \ En
particulier, $B_0$ est une sous-alg\`ebre de $A$ [BO 81$_2$]. Nous
r\'esumons d'autres r\'esultats dans le Lemme suivant:

\vspace{0.3cm} \begin{lemma} .
\begin{enumerate} \item Les valeurs propres de \ $\partial$ sont des imaginaires pures.
\item Si \ $\partial\neq 0,$ alors le rang de $\partial$ est $4$
ou $6,$ et on a
\begin{enumerate} \item Si $rg(\partial)=4,$ il existe \ $\alpha'>0$ tel que \ $A=B_0\oplus B_{\alpha'i}$ \ avec
\[ \dim(B_0)=\dim(B_{\alpha'i})=4, \hspace{0.2cm}
B_0B_{\alpha'i}=B_{\alpha'i}B_0=B_{\alpha'i} \ \mbox{ et } \
B_{\alpha'i}B_{\alpha'i}=B_0. \] \item Si $rg(\partial)=6,$ il
existe $\alpha', \beta'$ avec $0<\alpha'\leq\beta'$ tels que \[
A=B_0\oplus (B_{\alpha'i}+B_{\beta'i})\oplus B_{(\alpha'+\beta')i}
\] avec $\dim(B_0)=\dim(B_{(\alpha'+\beta')i}=2$ et on a

\vspace{0.2cm} {\bf i)} \ $B_0B_{\gamma i}=B_{\gamma
i}B_0=B_{\gamma i}$ \ pour \
$\gamma\in\{\alpha',\beta',\alpha'+\beta'\}.$

\vspace{0.2cm} {\bf ii)} \
$B_{\alpha'i}B_{(\alpha'+\beta')i}=B_{(\alpha'+\beta')i}B_{\alpha'i}=B_{\beta'i}.$

\vspace{0.2cm} {\bf iii)} \
$B_{\beta'i}B_{(\alpha'+\beta')i}=B_{(\alpha'+\beta')i}B_{\beta'i}=B_{\alpha'i}.$

\vspace{0.2cm} {\bf iv)} \
$B_{(\alpha'+\beta')i}=B_{(\alpha'+\beta')i}=B_0.\Box$
\end{enumerate}
\end{enumerate}
\end{lemma}

\vspace{0.1cm} {\bf Preuve.} ([BO 81$_2$] Lemmas {\bf 9,
15}).$\Box$

\vspace{3.5cm} \begin{proposition} Soit \ $A=\Big(
V,(.|.),\wedge\Big)$ une $K$-alg\`ebre quadratique, flexible de
division et soit $\partial\in End_K(A).$ Alors les deux
propri\'et\'es suivantes sont \'equivalentes:.
\begin{enumerate} \item $\partial$ est une d\'erivation de $A.$
\item $\partial$ est anti-sym\'etrique par rapport \`a $(.|.)$ et
induit une d\'erivation $\partial_V$ sur l'alg\`ebre
anti-commutative $(V,\wedge).$
\end{enumerate}
\end{proposition}

\vspace{0.1cm} {\bf Preuve.} 1) $\Rightarrow$ 2) Soit $x\in
V-\{0\},$ il existe $\alpha\in K$ et $u\in V$ tels que $\partial
x=\alpha+u$ et on a: \begin{eqnarray*} 0 &=&
\partial(x^2) \\
&=& (\partial x)x+x\partial x \\ &=& 2(u|x)+2\alpha x.
\end{eqnarray*}

On en d\'eduit que $\partial x$ est un vecteur (orthogonal \`a
$x$) i.e. $\partial A\subseteq V.$ Soient maintenant $x,y\in V,$
on a
\begin{eqnarray*} (x|\partial y)+(\partial x|y)+x\wedge (\partial
y)+(\partial x)\wedge y &=&
x\partial y+(\partial x)y \\
&=& \partial (xy) \\
&=& \partial(x\wedge y)\in V. \end{eqnarray*}

Donc $\partial$ induit une d\'erivation $\partial_V$ sur
l'alg\`ebre $(V,\wedge),$ anti-sym\'etrique par rapport \`a
$(.|.).$ Comme $\partial 1=0,$ $\partial$ est, \`a son tour,
anti-sym\'etrique par rapport \`a $(.|.).$

\vspace{0.2cm} 2) $\Rightarrow$ 1) Il suffit de montrer que
$\partial 1=0.$ En effet, on a: $(\partial 1|1)=-(1|\partial 1)$ \
i.e. $\partial 1\in V.$ Donc $\partial^2 1\in V$ car $\partial$
laisse stable $V,$ et on a:
\begin{eqnarray*} (\partial 1)^2 &=&
(\partial 1|\partial 1) \\
&=& -(1|\partial^2 1) \\
&=& 0.\Box. \end{eqnarray*}

\vspace{0.4cm} \hspace{0.3cm} Dans le reste de ce paragraphe
$A=\Big( W,(.|.),\wedge\Big)$ et $\partial$ d\'esigneront une
$\rit$-alg\`ebre de Jordan non commutative de division lin\'eaire
de dimension $8$ et une d\'erivation non triviale de $A.$ Nous
conserverons alors les notations pr\'ec\'edentes.

\vspace{5cm} \begin{remarks} . \begin{enumerate} \item Si $u$ est
un vecteur non nul de $A$ tel que $\partial u=0,$ alors
$\partial_W$ laisse stable le sous-espace $W(u).$ En effet, si
$v\in W(u),$ on a \ $(\partial v|u)=-(v|\partial u)=0.$ On notera
$\partial_u$ l'application induite par $\partial_W$ sur $W(u).$
\item Les valeurs propres de l'op\'erateur sym\'etrique
$\partial^2$ sont de la forme

\[ 0, \ 0, \ -\alpha'^2, \
-\alpha'^2, \ , -\beta'^2, \ -\beta'^2, \ -(\alpha'+\beta')^2, \
-(\alpha'+\beta')^2 \]

avec $0\leq \alpha'\leq\beta'$ et $(\alpha',\beta')\neq(0,0).$
\item Si $rg(\partial)=6,$ \ i.e. $\alpha'>0,$ alors \ $B_0\oplus
B_{(\alpha'+\beta')i}$ est une sous-alg\`ebre de $A,$ de dimension
$4.$ De plus, $B_{(\alpha'+\beta')i}$ est le sous-espace propre de
$\partial^2$ associ\'e \`a la valeur propre \
$-(\alpha'+\beta')^2$ et s'\'ecrit \ $vect\{x,\partial x\}$ o\`u
$x$ est une vecteur non nul quelconque de $B_{(\alpha'+\beta')i}.$
\begin{enumerate} \item Si, de plus, $\alpha'neq\beta',$ alors
$B_{\alpha'i}$ et $B_{\beta'i}$ sont respectivement les
sous-espaces propres de $\partial^2$ associ\'es aux valeurs
propres $-\alpha'^2$ et $-\beta'^2$ et s'\'ecrivent de la m\^eme
fa\c{c}on que celle de $B_{(\alpha'+\beta')i}.$ \item Si
$\alpha'=\beta',$ alors $B_{\alpha'i}$ est un sous-espace propre
de $\partial^2$ de dimension $4$ et l'\'ecriture pr\'ec\'edente
est valable pour tout sous-espace vectoriel de $B_{\alpha'i}$ de
dimension $2$ stable par $\partial.\Box$
\end{enumerate}
\end{enumerate}
\end{remarks}

\vspace{0.4cm} \hspace{0.3cm} On notera par la suite, $B_0$ et
$B_{(\alpha'+\beta')i},$ par $H_0$ et $H_1$ respectivement et la
sous-alg\`ebre $B_0\oplus B_{(\alpha'+\beta')i},$ par $B.$

\vspace{0.3cm} \begin{lemma} .
\begin{enumerate} \item Pour toute valeur propre non nulle,
$\sigma$ de $\partial,$ on a \ $B_{\sigma}\subset W.$ \item Si
$rg(\partial)=6,$ alors la d\'ecomposition \ $A=B_0\oplus
(B_{\alpha'i}+B_{\beta'i})\oplus B_{(\alpha'+\beta')i}$ \ est
orthogonale. Si, de plus, $u$ est un vecteur norm\'e de $H_0,$
alors $H_0$ s'\'ecrit $vect\{1,u\}$ et $B_{\alpha'i}+B_{\beta'i}$
se d\'ecompose en somme directe orthogonale de sous-espaces $H_2$
et $H_3,$ de dimension $2,$ stables par $f_u$ et $\partial_u$ tels
que $H_2^2=H_3^2=H_0.$ \item Si $rg(\partial)=4,$ alors la
d\'ecomposition \ $A=B_0\oplus B_{\alpha'i}$ \ est orthogonale. De
plus, il existe un vecteur norm\'e $u$ de $B_0$ tel que
$B_{\alpha'i}$ se d\'ecompose en somme directe orthogonale de deux
sous-espaces $H_2$ et $H_3,$ de dimension $2,$ stables par $f_u$
et $\partial_u$ avec $H_2^2=H_0.\Box$
\end{enumerate}
\end{lemma}

\vspace{0.1cm} {\bf Preuve.} La proposition 1) d\'ecoule de la
Proposition {\bf 5.3}.

\vspace{0.2cm} \hspace{0.3cm} La premi\`ere proposition de 2) et
3) d\'ecoule alors de 1) et des relations, dans le Lemme {\bf
5.2}, entre les sous-espaces propres de $\partial.$

\vspace{0.2cm} \hspace{0.3cm} Pour le reste de 2), nous remarquons
que les sous-espaces $B_{\alpha'i}$ et $B_{\beta'i}$ sont stables
par $\partial_u$ et $f_u$ car $\partial_u$ et $f_u$ commutent
($[\partial, L_u]=L_{\partial u}\equiv 0$). De plus, \
$B_{\alpha'i}+B_{\beta'i}:=E$ n'est autre que l'orthogonal de $B$
dans $A.$ On distingue alors les deux cas suivants:
\begin{enumerate} \item Si $\alpha'\neq\beta',$ alors $B_{\alpha'i}$ s'\'ecrit \ $vect\{x,\partial x\}$
o\`u $x\in B_{\alpha'i}-\{0\}$ et on a

\begin{eqnarray*} \partial(x\partial x) &=& (\partial
x)^2-x\partial^2x \\
&=& (\partial x)^2+\alpha'^2x^2\in\rit. \end{eqnarray*}

Donc $x\partial x\in H_0$ i.e. $B_{\alpha'i}^2=H_0,$ de m\^eme
$B_{\beta'i}^2=H_0.$ On prend alors $H_2=B_{\alpha'i}$ et
$H_3=B_{\beta'i}$ dans ce premier cas. \item Si $\alpha'=\beta',$
on consid\`ere les valeurs propres \ $\lambda,$ $\lambda,$ $\mu,$
$\mu$ \ de l'op\'erateur sym\'etrique $f,$ restriction de $f_u^2$
\`a $B_{\alpha'i}:=E,$ et on distingue les deux sous-cas suivants:
\begin{enumerate} \item Si $\lambda\neq \mu,$ les sous-espaces
propres $E_\lambda$ et $E_\mu$ correspondants sont orthogonaux,
stables par $\partial_u$ et on montre, comme dans le cas
pr\'ec\'edent, que $E_\lambda^2=E_\mu^2=H_0.$ De plus,
$E=E_\lambda\oplus E_\mu.$ \item Si $\lambda=\mu,$ alors
$f=\lambda I_E=-\lambda\alpha'^{-2}(\partial_E)^2$ \ i.e.

\[ \Big( f_{u_{/E}}-(-\lambda)^{\frac{1}{2}}\alpha'^{-1}\partial_E\Big)
\Big(
f_{u_{/E}}+(-\lambda)^{\frac{1}{2}}\alpha'^{-1}\partial_E\Big)\equiv
0. \]

On distingue les deux situations suivantes:

\vspace{0.2cm} {\bf i)} Si \ $f_{u_{/E}}$ et $\partial_E$ sont
colin\'eaires, on se ram\`ene au sous-cas {\bf 2. (a)} en
consid\'erant des sous-espaces de $E,$ de dimension $2,$
orthogonaux, stables par $f_u$ (donc stables par $\partial_u$
\'egalement).

\vspace{0.2cm} {\bf ii)} Si \ $f_{u_{/E}}$ et $\partial_E$ ne sont
pas colin\'eaires, alors le sous-espace \[ \ker\Big(
f_{u_{/E}}-(-\lambda)^{\frac{1}{2}}\alpha'^{-1}\partial_E\Big):=H
\] est de dimension $2,$ stables par $f_u$ et $\partial_u$ et
s'\'ecrit \ $vect\{x,\partial x\}$ o\`u $x\in H-\{0\}.$ Si $y$ est
un vecteur non nul de l'orthogonal $H^{\perp}:=K,$ de $H$ dans
$E,$ alors \ $x\partial y=\partial(xy)-(\partial x)y$ \ et \
$(\partial x)(\partial y)=\partial(x\partial y)+\alpha'^2xy$ \
sont des vecteurs, car $xy$ et $(\partial x)y$ le sont. Donc
$\partial y\in K$ et, par cons\'equent, \ $K$ coincide avec
$vect\{y,\partial y\}$ et est stable par $\partial_u$ aussi bien
que par $f_u.$ Ce qui nous ram\`ene au sous-cas {\bf 2. (a)}.
\end{enumerate} \end{enumerate}

\vspace{0.2cm} \hspace{0.3cm} Pour le reste de 3), il existe
$(u,y_2)\in S(B_0\cap W)\times S(B_{\alpha'i})$ tel que $u$ et
$y_2$ engendrent une sous-alg\`ebre \ $vect\{1,u,y_2,uy_2\}$, de
$A,$ de dimension $4$ et on a \
$(vect\{y_2,uy_2\})^2=vect\{1,u\}.$ On notera $H_0'$ le
sous-espace $vect\{1,u\}$ et l'on utilise les notations de la
Proposition {\bf 3.14}. On distingue alors les deux cas suivants:
\begin{enumerate} \item Si $h_u^2$ poss\`ede une unique valeur propre $\lambda,$ alors \
$-\lambda=|||h_u|||^2=||u\wedge y||^2$ \ pour tout $y\in
S(B_{\alpha'i}).$ Donc $u$ et $y$ engendrent une sous-alg\`ebre de
$A$ de dimension $4$ pour tout $y\in S(B_{\alpha'i}),$ et on
construit les sous-espaces $H_2$ et $H_3$ de la m\^eme fa\c{c}on
que dans {\bf 2.} De plus, $H_2^2=H_3^2=H_0'.$ \item Si $h_u^2$
poss\`ede deux valeurs propres distinctes, alors les sous-espaces
propres correspondants sont stables par $f_u$ et $\partial_u.$ De
plus, le carr\'e du sous-espace propre associ\'e \`a la plus
grande, en valeur absolue, valeur propre coincide avec $H_0'.\Box$
\end{enumerate}

\vspace{0.3cm} \begin{corollary} . \begin{enumerate} \item Si
$rg(\partial)=6,$ alors pour tout $i\in\{1,2,3\}:$ \ $H_0+H_i$ est
une sous-alg\`ebre de $A,$ de dimension $4,$ et on a \
$H_iH_j=H_k$ \ pour toute permutation $(i\ j\ k)$ de $\{1,2,3\}.$
De plus, tout sous-espace $H_i$ s'\'ecrit $vect\{y_i,z_i\}$ o\`u
$y_i,$ $z_i$ sont deux vecteurs orthonormaux de $A$ tels que \
$\partial y_i=\gamma_i z_i$ et $\partial z_i=-\gamma_iy_i$ avec
$\gamma_1=\alpha'+\beta',$ $\gamma_2=\alpha'$ et
$\gamma_3=\beta'.$ \item Si $rg(\partial)=4,$ alors tout
sous-espace $H_i,$ o\`u \ $i\in\{1,2,3\},$ s'\'ecrit \
$vect\{y_i,z_i\}$ o\`u $y_i$ et $z_i$ sont deux vecteurs
orthonormaux de $A$ tels que \ $\partial y_i=\alpha'z_i$ et
$\partial z_i=-\alpha'y_i.\Box$
\end{enumerate}
\end{corollary}

\vspace{0.1cm} {\bf Preuve.} Si $rg(\partial)=6$ et si, pour
$i\in\{1,2,3\},$ $y_i$ est un vecteur norm\'e de $H_i,$ alors
$y_i,$ $\partial y_i$ est une base orthogonale de $H_i$ et on a

\begin{eqnarray*} (\partial y_i)^2 &=& \partial(y_i\partial
y_i)-y_i\partial^2y_i \\
&=& -y_i\partial^2y_i \\
&=& (\gamma_iy_i)^2. \end{eqnarray*}

Ainsi \ $\gamma_i^{-1}\partial y_i:=z_i$ est un vecteur norm\'e de
$H_i.$ Les relations $H_iH_i=H_k$ sont alors cons\'equences de la
propri\'et\'e trace de $(.|.).\Box$

\vspace{0.5cm} \hspace{0.3cm} Nous noterons ${\cal B}_2$ la base
orthonorm\'ee \ $\{1,u,y_1,z_1,y_2,z_2,y_3,z_3\}$ et utiliserons
les notations du Corollaire {\bf 5.6}.

\vspace{3cm} \begin{remark} Si $rg(\partial)=4,$ \'etant donn\'es
$y_3,z_3$ et sachant que \ $H_3^2\subseteq B,$ il existe un
vecteur norm\'e $z_1$ de $B,$ orthogonal \`a $u,$ tel que $y_3z_3$
soit une combinaison lin\'eaire de $u$ et $z_1.$ On note alors
$y_1$ le vecteur, norm\'e, \ $-||uz_1||^{-1}uz_1$ et $H_1'$ le
sous-espace $vect\{y_1,z_1\}.$ La propri\'et\'e trace de $(.|.),$
le fait que $H_2^2=H_0'$ et l'hypoth\`ese de division lin\'eaire
donnent

\begin{eqnarray} yH_2 &=& H_3 \ \mbox{ pour tout } y\in H_1'. \end{eqnarray}

Il existe alors un vecteur norm\'e $y_2$ de $H_2$ tel que $y_1y_2$
et $y_3$ soient colin\'eaires, avec une constante multiplicative
positive.$\Box$ \end{remark}

\vspace{0.4cm} \hspace{0.3cm} La base orthonorm\'ee \
$\{1,u,y_1,z_1,y_2,z_2,y_3,z_3\}$ \ obtenue sera not\'ee ${\cal
B}_1.$ La table de multiplication de $A,$ par rapport \`a ${\cal
B}_1,$ est alors donn\'ee moyennant {\bf (5.6)}, la propri\'et\'e
trace de $(.|.)$ et le Corollaire {\bf 5.6}, par:

\vspace{0.3cm}
\[ \begin{tabular}{ccccccccc} \\
\multicolumn{1}{c}{} & \multicolumn{1}{c}{$1$} &
\multicolumn{1}{c}{$u$} & \multicolumn{1}{c}{$y_1$} &
\multicolumn{1}{c}{$z_1$} & \multicolumn{1}{c}{$y_2$} &
\multicolumn{1}{c}{$z_2$} & \multicolumn{1}{c}{$y_3$} &
\multicolumn{1}{c}{$z_3$} \\ \cline{2-9} \multicolumn{1}{c|}{$1$}
& \multicolumn{1}{|c}{$1$} & \multicolumn{1}{c|}{$u$} &
\multicolumn{1}{|c}{$y_1$} & \multicolumn{1}{c|}{$z_1$} &
\multicolumn{1}{|c}{$y_2$} & \multicolumn{1}{c|}{$z_2$} &
\multicolumn{1}{|c}{$y_3$} & \multicolumn{1}{c|}{$z_3$} \\
\multicolumn{1}{c|}{$u$} & \multicolumn{1}{|c}{} &
\multicolumn{1}{c|}{$-1$} & \multicolumn{1}{|c}{$az_1$} &
\multicolumn{1}{c|}{$-ay_1$} & \multicolumn{1}{|c}{$bz_2$} &
\multicolumn{1}{c|}{$-by_2$} & \multicolumn{1}{|c}{$cz_3$} &
\multicolumn{1}{c|}{$-cy_3$} \\ \cline{2-9}
\multicolumn{1}{c|}{$y_1$} & \multicolumn{1}{|c}{} &
\multicolumn{1}{c|}{} & \multicolumn{1}{|c}{$-1$} &
\multicolumn{1}{c|}{$au$} & \multicolumn{1}{|c}{$\alpha y_3$} &
\multicolumn{1}{c|}{$\alpha z_3$} & \multicolumn{1}{|c}{$-\alpha
y_2$} &
\multicolumn{1}{c|}{$-\alpha z_2$} \\
\multicolumn{1}{c|}{$z_1$} & \multicolumn{1}{|c}{} &
\multicolumn{1}{c|}{} & \multicolumn{1}{|c}{} &
\multicolumn{1}{c|}{$-1$} & \multicolumn{1}{|c}{$\lambda y_3+\eta
z_3$} & \multicolumn{1}{c|}{$-\eta y_3+\lambda z_3$} &
\multicolumn{1}{|c}{$-\lambda y_2+\eta z_2+\rho z_3$} &
\multicolumn{1}{c|}{$-\eta y_2-\lambda z_2-\rho y_3$} \\
\cline{4-9} \multicolumn{1}{c|}{$y_2$} & \multicolumn{1}{|c}{} &
\multicolumn{1}{c}{} & \multicolumn{1}{c}{} &
\multicolumn{1}{c|}{} & \multicolumn{1}{|c}{$-1$} &
\multicolumn{1}{c|}{$bu$} & \multicolumn{1}{|c}{$\alpha
y_1+\lambda z_1$} &
\multicolumn{1}{c|}{$\eta z_1$} \\
\multicolumn{1}{c|}{$z_2$} & \multicolumn{1}{|c}{} &
\multicolumn{1}{c}{} & \multicolumn{1}{c}{} & \multicolumn{1}{c}{}
& \multicolumn{1}{|c}{} & \multicolumn{1}{c|}{$-1$} &
\multicolumn{1}{|c}{$-\eta z_1$} & \multicolumn{1}{c|}{$\alpha
y_1+\lambda z_1$}
\\ \cline{6-9} \multicolumn{1}{c|}{$y_3$} & \multicolumn{1}{|c}{}
& \multicolumn{1}{c}{} & \multicolumn{1}{c}{} &
\multicolumn{1}{c}{} & \multicolumn{1}{c}{} & \multicolumn{1}{c}{}
& \multicolumn{1}{|c}{$-1$} & \multicolumn{1}{c|}{$cu+\rho z_1$} \\
\multicolumn{1}{c|}{$z_3$} & \multicolumn{1}{|c}{} &
\multicolumn{1}{c}{} & \multicolumn{1}{c}{} & \multicolumn{1}{c}{}
& \multicolumn{1}{c}{} & \multicolumn{1}{c}{} &
\multicolumn{1}{|c}{} & \multicolumn{1}{c|}{$-1$}
\\ \cline{2-9}
\end{tabular} \]

\begin{center} {\bf Table 3} \end{center}

\vspace{0.5cm} o\`u $a,$ $b,$ $c,$ $\alpha,$ $\lambda,$ $\eta,$
$\rho\in\rit$ \ avec $a,$ $\alpha >0.$

\vspace{6cm} \begin{proposition} Soit $A$ une $\rit$-alg\`ebre de
Jordan non commutative de division lin\'eaire de dimension $8,$
poss\'edant une d\'erivation de rang $6.$ Alors il existe une base
orthonorm\'ee \ $1,$ $u,$ $y_1,$ $z_1,$ $y_2,$ $z_2,$ $y_3,$ $z_3$
\ de $A,$ et quatre param\`etres \ $a,$ $b,$ $c,$ $\alpha>0$ \
pour lesquels la multiplication de $A$ est donn\'ee par la Table
{\bf 4} suivante:

\vspace{0.3cm}
\[ \begin{tabular}{ccccccccc} \\
\multicolumn{1}{c}{} & \multicolumn{1}{c}{$1$} &
\multicolumn{1}{c}{$u$} & \multicolumn{1}{c}{$y_1$} &
\multicolumn{1}{c}{$z_1$} & \multicolumn{1}{c}{$y_2$} &
\multicolumn{1}{c}{$z_2$} & \multicolumn{1}{c}{$y_3$} &
\multicolumn{1}{c}{$z_3$} \\ \cline{2-9} \multicolumn{1}{c|}{$1$}
& \multicolumn{1}{|c}{$1$} & \multicolumn{1}{c|}{$u$} &
\multicolumn{1}{|c}{$y_1$} & \multicolumn{1}{c|}{$z_1$} &
\multicolumn{1}{|c}{$y_2$} & \multicolumn{1}{c|}{$z_2$} &
\multicolumn{1}{|c}{$y_3$} & \multicolumn{1}{c|}{$z_3$} \\
\multicolumn{1}{c|}{$u$} & \multicolumn{1}{|c}{} &
\multicolumn{1}{c|}{$-1$} & \multicolumn{1}{|c}{$az_1$} &
\multicolumn{1}{c|}{$-ay_1$} & \multicolumn{1}{|c}{$bz_2$} &
\multicolumn{1}{c|}{$-by_2$} & \multicolumn{1}{|c}{$cz_3$} &
\multicolumn{1}{c|}{$-cy_3$} \\ \cline{2-9}
\multicolumn{1}{c|}{$y_1$} & \multicolumn{1}{|c}{} &
\multicolumn{1}{c|}{} & \multicolumn{1}{|c}{$-1$} &
\multicolumn{1}{c|}{$au$} & \multicolumn{1}{|c}{$\alpha y_3$} &
\multicolumn{1}{c|}{$-\alpha z_3$} & \multicolumn{1}{|c}{$-\alpha
y_2$} &
\multicolumn{1}{c|}{$\alpha z_2$} \\
\multicolumn{1}{c|}{$z_1$} & \multicolumn{1}{|c}{} &
\multicolumn{1}{c|}{} & \multicolumn{1}{|c}{} &
\multicolumn{1}{c|}{$-1$} & \multicolumn{1}{|c}{$-\alpha z_3$} &
\multicolumn{1}{c|}{$-\alpha y_3$} & \multicolumn{1}{|c}{$\alpha
z_2$} &
\multicolumn{1}{c|}{$\alpha y_2$} \\
\cline{4-9} \multicolumn{1}{c|}{$y_2$} & \multicolumn{1}{|c}{} &
\multicolumn{1}{c}{} & \multicolumn{1}{c}{} &
\multicolumn{1}{c|}{} & \multicolumn{1}{|c}{$-1$} &
\multicolumn{1}{c|}{$bu$} & \multicolumn{1}{|c}{$\alpha y_1$} &
\multicolumn{1}{c|}{$-\alpha z_1$} \\
\multicolumn{1}{c|}{$z_2$} & \multicolumn{1}{|c}{} &
\multicolumn{1}{c}{} & \multicolumn{1}{c}{} & \multicolumn{1}{c}{}
& \multicolumn{1}{|c}{} & \multicolumn{1}{c|}{$-1$} &
\multicolumn{1}{|c}{$-\alpha z_1$} & \multicolumn{1}{c|}{$-\alpha
y_1$}
\\ \cline{6-9} \multicolumn{1}{c|}{$y_3$} & \multicolumn{1}{|c}{}
& \multicolumn{1}{c}{} & \multicolumn{1}{c}{} &
\multicolumn{1}{c}{} & \multicolumn{1}{c}{} & \multicolumn{1}{c}{}
& \multicolumn{1}{|c}{$-1$} & \multicolumn{1}{c|}{$cu$} \\
\multicolumn{1}{c|}{$z_3$} & \multicolumn{1}{|c}{} &
\multicolumn{1}{c}{} & \multicolumn{1}{c}{} & \multicolumn{1}{c}{}
& \multicolumn{1}{c}{} & \multicolumn{1}{c}{} &
\multicolumn{1}{|c}{} & \multicolumn{1}{c|}{$-1$}
\\ \cline{2-9}
\end{tabular} \]

\begin{center} {\bf Table 4} \end{center}

\vspace{0.3cm} De plus, une $\rit$-alg\`ebre dont la
multiplication est donn\'ee par la Table {\bf 4} est de Jordan non
commutative de division lin\'eaire et poss\`ede une d\'erivation
de rang $6.$
\end{proposition}

\vspace{0.1cm} {\bf Preuve.} Suivant les notations
pr\'ec\'edentes, il existe des param\`etres \ $\alpha,$ $\beta,$
$\gamma,$ $\mu,$ $\lambda,$ $\eta,$ $\sigma,$ $\delta$ \ tels que

\begin{eqnarray*} y_1y_2 &=& \alpha y_3+\beta z_3, \\
y_1z_2 &=& \gamma y_3+\mu z_3, \\
z_1y_2 &=& \lambda y_3+\eta z_3, \\
z_1z_2 &=& \sigma y_3+\delta z_3. \end{eqnarray*}

Les \'egalit\'es

\begin{eqnarray*} \partial(y_1y_2) &=& (\partial y_1)y_2+y_1\partial y_2, \\
\partial(y_1z_2) &=& (\partial y_1)z_2+y_1\partial z_2 \\
\partial(z_1y_2) &=& (\partial z_1)y_2+z_1\partial y_2, \\
\partial(z_1z_2) &=& (\partial z_1)z_2+z_1\partial z_2 \end{eqnarray*}

donnent respectivement:

\vspace{0.5cm}
\begin{enumerate} \item
$\gamma\alpha'+\beta\beta'=-\lambda(\alpha'+\beta'),$ \item
$-\mu\alpha'+\alpha\beta'=\eta(\alpha'+\beta'),$ \item
$\alpha\alpha'-\mu\beta'=\sigma(\alpha'+\beta'),$ \item
$\beta\alpha'+\gamma\beta'=\delta(\alpha'+\beta'),$ \item
$\sigma\alpha'+\eta\beta'=\alpha(\alpha'+\beta'),$ \item
$\delta\alpha'-\lambda\beta'=\beta(\alpha'+\beta'),$ \item
$\lambda\alpha'-\delta\beta'=-\gamma(\alpha'+\beta'),$ \item
$\eta\alpha'+\sigma\beta'=-\mu(\alpha'+\beta').$
\end{enumerate}

\vspace{0.3cm} {\bf 1.} et {\bf 7.} donnent {\bf 1'.} \
$\beta+\delta=\gamma-\lambda$ \ et \ {\bf 2'.} \
$2(\gamma+\lambda)\alpha'=(\delta-\beta-\gamma-\lambda)\beta'.$

\vspace{0.1cm} {\bf 2.} et {\bf 8.} donnent {\bf 3.} \
$\alpha+\sigma=\eta-\mu$ \ et \ {\bf 4'.} \
$2(\mu+\eta)\alpha'=(\alpha-\eta-\mu-\sigma)\beta'.$

\vspace{0.1cm} {\bf 3.} et {\bf 5.} donnent {\bf 5'.} \
$2(\alpha-\sigma)\alpha'=(\sigma-\alpha\mu+\eta)\beta'.$

\vspace{0.1cm} {\bf 4.} et {\bf 6.} donnent {\bf 6'.} \
$2(\beta-\delta)\alpha'=(\delta-\beta-\gamma-\lambda)\beta'.$

\vspace{0.3cm} Les nouvelles \'equations {\bf 2'.} et {\bf 6'.}
donnent \ {\bf 7'.} \ $\gamma+\lambda=\beta-\delta.$

\vspace{0.1cm} De m\^eme {\bf 4'.} et {\bf 5'.} donnent \ {\bf
8'.} \ $\mu+\eta=\sigma-\alpha.$

\vspace{0.3cm} Les \'equations {\bf 1'.}, {\bf 3'.}, {\bf 7'.} et
{\bf 8'.} donnent alors: \
$(\gamma,\mu,\sigma,\delta)=(\beta,-\alpha,\eta,-\lambda).$ De ce
fait, les \'equations \ {\bf 5'.} et {\bf 6'.} donnent \
$(\alpha-\sigma)(\alpha'+\beta')=(\beta-\delta)(\alpha'+\beta')=0$
\ i.e. \ $\alpha=\sigma$ et $\beta=\delta.$ Ainsi, \
$-\mu=\alpha=\sigma=\eta$ et $\gamma=\beta=\delta=-\lambda.$ On
obtient, en tenant compte de la propri\'et\'e trace de $(.|.),$ la
table de multiplication de $A$ par rapport \`a la base ${\cal
B}_2,$ en se limitant \`a la partie triangulaire sup\'erieure, vu
l'anti-commutativit\'e de $"\wedge":$

\vspace{0.3cm}
\[ \begin{tabular}{ccccccccc} \\
\multicolumn{1}{c}{} & \multicolumn{1}{c}{$1$} &
\multicolumn{1}{c}{$u$} & \multicolumn{1}{c}{$y_1$} &
\multicolumn{1}{c}{$z_1$} & \multicolumn{1}{c}{$y_2$} &
\multicolumn{1}{c}{$z_2$} & \multicolumn{1}{c}{$y_3$} &
\multicolumn{1}{c}{$z_3$} \\ \cline{2-9} \multicolumn{1}{c|}{$1$}
& \multicolumn{1}{|c}{$1$} & \multicolumn{1}{c|}{$u$} &
\multicolumn{1}{|c}{$y_1$} & \multicolumn{1}{c|}{$z_1$} &
\multicolumn{1}{|c}{$y_2$} & \multicolumn{1}{c|}{$z_2$} &
\multicolumn{1}{|c}{$y_3$} & \multicolumn{1}{c|}{$z_3$} \\
\multicolumn{1}{c|}{$u$} & \multicolumn{1}{|c}{} &
\multicolumn{1}{c|}{$-1$} & \multicolumn{1}{|c}{$az_1$} &
\multicolumn{1}{c|}{$-ay_1$} & \multicolumn{1}{|c}{$bz_2$} &
\multicolumn{1}{c|}{$-by_2$} & \multicolumn{1}{|c}{$cz_3$} &
\multicolumn{1}{c|}{$-cy_3$} \\ \cline{2-9}
\multicolumn{1}{c|}{$y_1$} & \multicolumn{1}{|c}{} &
\multicolumn{1}{c|}{} & \multicolumn{1}{|c}{$-1$} &
\multicolumn{1}{c|}{$au$} & \multicolumn{1}{|c}{$\alpha y_3+\beta
z_3$} & \multicolumn{1}{c|}{$\beta y_3-\alpha z_3$} &
\multicolumn{1}{|c}{$-\alpha y_2-\beta z_2$} &
\multicolumn{1}{c|}{$-\beta y_2+\alpha z_2$} \\
\multicolumn{1}{c|}{$z_1$} & \multicolumn{1}{|c}{} &
\multicolumn{1}{c|}{} & \multicolumn{1}{|c}{} &
\multicolumn{1}{c|}{$-1$} & \multicolumn{1}{|c}{$-\beta y_3+\alpha
z_3$} & \multicolumn{1}{c|}{$\alpha y_3+\beta z_3$} &
\multicolumn{1}{|c}{$\beta y_2-\alpha z_2$} &
\multicolumn{1}{c|}{$-\alpha y_2-\beta z_2$} \\
\cline{4-9} \multicolumn{1}{c|}{$y_2$} & \multicolumn{1}{|c}{} &
\multicolumn{1}{c}{} & \multicolumn{1}{c}{} &
\multicolumn{1}{c|}{} & \multicolumn{1}{|c}{$-1$} &
\multicolumn{1}{c|}{$bu$} & \multicolumn{1}{|c}{$\alpha y_1-\beta
z_1$} &
\multicolumn{1}{c|}{$\beta y_1+\alpha z_1$} \\
\multicolumn{1}{c|}{$z_2$} & \multicolumn{1}{|c}{} &
\multicolumn{1}{c}{} & \multicolumn{1}{c}{} & \multicolumn{1}{c}{}
& \multicolumn{1}{|c}{} & \multicolumn{1}{c|}{$-1$} &
\multicolumn{1}{|c}{$\beta y_1+\alpha z_1$} &
\multicolumn{1}{c|}{$-\alpha y_1+\beta z_1$}
\\ \cline{6-9} \multicolumn{1}{c|}{$y_3$} & \multicolumn{1}{|c}{}
& \multicolumn{1}{c}{} & \multicolumn{1}{c}{} &
\multicolumn{1}{c}{} & \multicolumn{1}{c}{} & \multicolumn{1}{c}{}
& \multicolumn{1}{|c}{$-1$} & \multicolumn{1}{c|}{$cu$} \\
\multicolumn{1}{c|}{$z_3$} & \multicolumn{1}{|c}{} &
\multicolumn{1}{c}{} & \multicolumn{1}{c}{} & \multicolumn{1}{c}{}
& \multicolumn{1}{c}{} & \multicolumn{1}{c}{} &
\multicolumn{1}{|c}{} & \multicolumn{1}{c|}{$-1$}
\\ \cline{2-9}
\end{tabular} \]

\vspace{0.3cm} o\`u $a,b,c\in\rit^*.$ De plus, quitte \`a changer
le signe de $u$ si n\'ecessaire, on peut supposer $a>0$ et l'on
r\'eduit, en posant

\begin{eqnarray*} \alpha_1 &=& (\alpha^2+\beta^2)^{\frac{1}{2}}, \\
y_1' &=& \alpha_1^{-1}(\alpha y_1-\beta z_1), \\
z_1' &=& \alpha_1^{-1}(\beta y_1+\alpha z_1), \end{eqnarray*}

la table de multiplication de $A,$ \`a quatre param\`etres:

\vspace{0.3cm}
\[ \begin{tabular}{ccccccccc} \\
\multicolumn{1}{c}{} & \multicolumn{1}{c}{$1$} &
\multicolumn{1}{c}{$u$} & \multicolumn{1}{c}{$y_1'$} &
\multicolumn{1}{c}{$z_1'$} & \multicolumn{1}{c}{$y_2$} &
\multicolumn{1}{c}{$z_2$} & \multicolumn{1}{c}{$y_3$} &
\multicolumn{1}{c}{$z_3$} \\ \cline{2-9} \multicolumn{1}{c|}{$1$}
& \multicolumn{1}{|c}{$1$} & \multicolumn{1}{c|}{$u$} &
\multicolumn{1}{|c}{$y_1'$} & \multicolumn{1}{c|}{$z_1'$} &
\multicolumn{1}{|c}{$y_2$} & \multicolumn{1}{c|}{$z_2$} &
\multicolumn{1}{|c}{$y_3$} & \multicolumn{1}{c|}{$z_3$} \\
\multicolumn{1}{c|}{$u$} & \multicolumn{1}{|c}{} &
\multicolumn{1}{c|}{$-1$} & \multicolumn{1}{|c}{$az_1'$} &
\multicolumn{1}{c|}{$-ay_1'$} & \multicolumn{1}{|c}{$bz_2$} &
\multicolumn{1}{c|}{$-by_2$} & \multicolumn{1}{|c}{$cz_3$} &
\multicolumn{1}{c|}{$-cy_3$} \\ \cline{2-9}
\multicolumn{1}{c|}{$y_1'$} & \multicolumn{1}{|c}{} &
\multicolumn{1}{c|}{} & \multicolumn{1}{|c}{$-1$} &
\multicolumn{1}{c|}{$au$} & \multicolumn{1}{|c}{$\alpha_1y_3$} &
\multicolumn{1}{c|}{$-\alpha_1z_3$} &
\multicolumn{1}{|c}{$-\alpha_1y_2$} &
\multicolumn{1}{c|}{$\alpha_1z_2$} \\
\multicolumn{1}{c|}{$z_1'$} & \multicolumn{1}{|c}{} &
\multicolumn{1}{c|}{} & \multicolumn{1}{|c}{} &
\multicolumn{1}{c|}{$-1$} & \multicolumn{1}{|c}{$\alpha_1z_3$} &
\multicolumn{1}{c|}{$\alpha_1y_3$} &
\multicolumn{1}{|c}{$-\alpha_1z_2$} &
\multicolumn{1}{c|}{$-\alpha_1y_2$} \\
\cline{4-9} \multicolumn{1}{c|}{$y_2$} & \multicolumn{1}{|c}{} &
\multicolumn{1}{c}{} & \multicolumn{1}{c}{} &
\multicolumn{1}{c|}{} & \multicolumn{1}{|c}{$-1$} &
\multicolumn{1}{c|}{$bu$} & \multicolumn{1}{|c}{$\alpha_1y_1'$} &
\multicolumn{1}{c|}{$\alpha_1z_1'$} \\
\multicolumn{1}{c|}{$z_2$} & \multicolumn{1}{|c}{} &
\multicolumn{1}{c}{} & \multicolumn{1}{c}{} & \multicolumn{1}{c}{}
& \multicolumn{1}{|c}{} & \multicolumn{1}{c|}{$-1$} &
\multicolumn{1}{|c}{$\alpha_1z_1'$} &
\multicolumn{1}{c|}{$-\alpha_1y_1'$}
\\ \cline{6-9} \multicolumn{1}{c|}{$y_3$} & \multicolumn{1}{|c}{}
& \multicolumn{1}{c}{} & \multicolumn{1}{c}{} &
\multicolumn{1}{c}{} & \multicolumn{1}{c}{} & \multicolumn{1}{c}{}
& \multicolumn{1}{|c}{$-1$} & \multicolumn{1}{c|}{$cu$} \\
\multicolumn{1}{c|}{$z_3$} & \multicolumn{1}{|c}{} &
\multicolumn{1}{c}{} & \multicolumn{1}{c}{} & \multicolumn{1}{c}{}
& \multicolumn{1}{c}{} & \multicolumn{1}{c}{} &
\multicolumn{1}{|c}{} & \multicolumn{1}{c|}{$-1$}
\\ \cline{2-9}
\end{tabular} \]

\vspace{0.4cm} \hspace{0.3cm} Les param\`etres $b$ et $c$
s'\'ecrivent, respectivement, \ $\varepsilon|b|$ \ et \
$\varepsilon'|c|$ \ o\`u $\varepsilon, \varepsilon'\in\{1,-1\}.$
La table de multiplication de $A,$ par rapport \`a la base \ $1,$
$u,$ $y_1$ $z_1',$ $y_2,$ $\varepsilon z_2:=z_2',$ $y_3,$
$\varepsilon'z_3:=z_3'$ \ est donn\'ee par:

\vspace{0.3cm}
\[ \begin{tabular}{ccccccccc} \\
\multicolumn{1}{c}{} & \multicolumn{1}{c}{$1$} &
\multicolumn{1}{c}{$u$} & \multicolumn{1}{c}{$y_1'$} &
\multicolumn{1}{c}{$z_1'$} & \multicolumn{1}{c}{$y_2$} &
\multicolumn{1}{c}{$z_2'$} & \multicolumn{1}{c}{$y_3$} &
\multicolumn{1}{c}{$z_3'$} \\ \cline{2-9} \multicolumn{1}{c|}{$1$}
& \multicolumn{1}{|c}{$1$} & \multicolumn{1}{c|}{$u$} &
\multicolumn{1}{|c}{$y_1'$} & \multicolumn{1}{c|}{$z_1'$} &
\multicolumn{1}{|c}{$y_2$} & \multicolumn{1}{c|}{$z_2'$} &
\multicolumn{1}{|c}{$y_3$} & \multicolumn{1}{c|}{$z_3'$} \\
\multicolumn{1}{c|}{$u$} & \multicolumn{1}{|c}{} &
\multicolumn{1}{c|}{$-1$} & \multicolumn{1}{|c}{$az_1'$} &
\multicolumn{1}{c|}{$-ay_1'$} & \multicolumn{1}{|c}{$|b|z_2'$} &
\multicolumn{1}{c|}{$-|b|y_2'$} & \multicolumn{1}{|c}{$|c|z_3'$} &
\multicolumn{1}{c|}{$-|c|y_3$} \\ \cline{2-9}
\multicolumn{1}{c|}{$y_1'$} & \multicolumn{1}{|c}{} &
\multicolumn{1}{c|}{} & \multicolumn{1}{|c}{$-1$} &
\multicolumn{1}{c|}{$au$} & \multicolumn{1}{|c}{$\alpha_1y_3$} &
\multicolumn{1}{c|}{$-\varepsilon\varepsilon'\alpha_1z_3'$} &
\multicolumn{1}{|c}{$-\alpha_1y_2$} &
\multicolumn{1}{c|}{$\varepsilon\varepsilon'\alpha_1z_2'$} \\
\multicolumn{1}{c|}{$z_1'$} & \multicolumn{1}{|c}{} &
\multicolumn{1}{c|}{} & \multicolumn{1}{|c}{} &
\multicolumn{1}{c|}{$-1$} &
\multicolumn{1}{|c}{$\varepsilon'\alpha_1z_3'$} &
\multicolumn{1}{c|}{$\varepsilon\alpha_1y_3$} &
\multicolumn{1}{|c}{$-\varepsilon\alpha_1z_2$} &
\multicolumn{1}{c|}{$-\varepsilon'\alpha_1y_2$} \\
\cline{4-9} \multicolumn{1}{c|}{$y_2$} & \multicolumn{1}{|c}{} &
\multicolumn{1}{c}{} & \multicolumn{1}{c}{} &
\multicolumn{1}{c|}{} & \multicolumn{1}{|c}{$-1$} &
\multicolumn{1}{c|}{$|b|u$} & \multicolumn{1}{|c}{$\alpha_1y_1'$}
&
\multicolumn{1}{c|}{$\varepsilon'\alpha_1z_1'$} \\
\multicolumn{1}{c|}{$z_2'$} & \multicolumn{1}{|c}{} &
\multicolumn{1}{c}{} & \multicolumn{1}{c}{} & \multicolumn{1}{c}{}
& \multicolumn{1}{|c}{} & \multicolumn{1}{c|}{$-1$} &
\multicolumn{1}{|c}{$\varepsilon\alpha_1z_1'$} &
\multicolumn{1}{c|}{$-\varepsilon\varepsilon'\alpha_1y_1'$}
\\ \cline{6-9} \multicolumn{1}{c|}{$y_3$} & \multicolumn{1}{|c}{}
& \multicolumn{1}{c}{} & \multicolumn{1}{c}{} &
\multicolumn{1}{c}{} & \multicolumn{1}{c}{} & \multicolumn{1}{c}{}
& \multicolumn{1}{|c}{$-1$} & \multicolumn{1}{c|}{$|c|u$} \\
\multicolumn{1}{c|}{$z_3'$} & \multicolumn{1}{|c}{} &
\multicolumn{1}{c}{} & \multicolumn{1}{c}{} & \multicolumn{1}{c}{}
& \multicolumn{1}{c}{} & \multicolumn{1}{c}{} &
\multicolumn{1}{|c}{} & \multicolumn{1}{c|}{$-1$}
\\ \cline{2-9}
\end{tabular} \]

\vspace{0.5cm} Comme $A$ est de division lin\'eaire, ceci
\'equivaut \`a \ $-\varepsilon'\alpha_1^2,$
$-\varepsilon\alpha_1^2,$ $\varepsilon\varepsilon'\alpha_1^2>0,$
en vertu du Lemme {\bf 4.13} i.e. \
$\varepsilon'=\varepsilon=-1.\Box$

\vspace{4cm} \begin{theorem} Soit $A$ une $\rit$-alg\`ebre de
Jordan non commutative de division lin\'eaire de dimension $8,$
dont l'alg\`ebre de Lie des d\'erivations n'est pas triviale.
Alors il existe une base \ $1,$ $u,$ $y_1,$ $z_1,$ $y_2,$ $z_2,$
$y_3,$ $z_3$ \ de $A,$ quatre param\`etres \ $a,$ $b,$ $c,$
$\alpha>0$ et trois autres \ $\eta,$ $\lambda,$ $\rho$ \ pour
lesquels la multiplication de $A$ est donn\'ee par la Table {\bf
5} suivante:

\vspace{0.3cm}
\[ \begin{tabular}{ccccccccc} \\
\multicolumn{1}{c}{} & \multicolumn{1}{c}{$1$} &
\multicolumn{1}{c}{$u$} & \multicolumn{1}{c}{$y_1$} &
\multicolumn{1}{c}{$z_1$} & \multicolumn{1}{c}{$y_2$} &
\multicolumn{1}{c}{$z_2$} & \multicolumn{1}{c}{$y_3$} &
\multicolumn{1}{c}{$z_3$} \\ \cline{2-9} \multicolumn{1}{c|}{$1$}
& \multicolumn{1}{|c}{$1$} & \multicolumn{1}{c|}{$u$} &
\multicolumn{1}{|c}{$y_1$} & \multicolumn{1}{c|}{$z_1$} &
\multicolumn{1}{|c}{$y_2$} & \multicolumn{1}{c|}{$z_2$} &
\multicolumn{1}{|c}{$y_3$} & \multicolumn{1}{c|}{$z_3$} \\
\multicolumn{1}{c|}{$u$} & \multicolumn{1}{|c}{} &
\multicolumn{1}{c|}{$-1$} & \multicolumn{1}{|c}{$az_1$} &
\multicolumn{1}{c|}{$-ay_1$} & \multicolumn{1}{|c}{$bz_2$} &
\multicolumn{1}{c|}{$-by_2$} & \multicolumn{1}{|c}{$cz_3$} &
\multicolumn{1}{c|}{$-cy_3$} \\ \cline{2-9}
\multicolumn{1}{c|}{$y_1$} & \multicolumn{1}{|c}{} &
\multicolumn{1}{c|}{} & \multicolumn{1}{|c}{$-1$} &
\multicolumn{1}{c|}{$au$} & \multicolumn{1}{|c}{$\alpha y_3$} &
\multicolumn{1}{c|}{$-\alpha z_3$} & \multicolumn{1}{|c}{$-\alpha
y_2$} &
\multicolumn{1}{c|}{$\alpha z_2$} \\
\multicolumn{1}{c|}{$z_1$} & \multicolumn{1}{|c}{} &
\multicolumn{1}{c|}{} & \multicolumn{1}{|c}{} &
\multicolumn{1}{c|}{$-1$} & \multicolumn{1}{|c}{$\lambda y_3+\eta
z_3$} & \multicolumn{1}{c|}{$\eta y_3-\lambda z_3$} &
\multicolumn{1}{|c}{$-\lambda y_2-\eta z_2+\rho z_3$} &
\multicolumn{1}{c|}{$-\eta y_2+\lambda z_2-\rho y_3$} \\
\cline{4-9} \multicolumn{1}{c|}{$y_2$} & \multicolumn{1}{|c}{} &
\multicolumn{1}{c}{} & \multicolumn{1}{c}{} &
\multicolumn{1}{c|}{} & \multicolumn{1}{|c}{$-1$} &
\multicolumn{1}{c|}{$bu$} & \multicolumn{1}{|c}{$\alpha
y_1+\lambda z_1$} &
\multicolumn{1}{c|}{$\eta z_1$} \\
\multicolumn{1}{c|}{$z_2$} & \multicolumn{1}{|c}{} &
\multicolumn{1}{c}{} & \multicolumn{1}{c}{} & \multicolumn{1}{c}{}
& \multicolumn{1}{|c}{} & \multicolumn{1}{c|}{$-1$} &
\multicolumn{1}{|c}{$\eta z_1$} & \multicolumn{1}{c|}{$-\alpha
y_1-\lambda z_1$}
\\ \cline{6-9} \multicolumn{1}{c|}{$y_3$} & \multicolumn{1}{|c}{}
& \multicolumn{1}{c}{} & \multicolumn{1}{c}{} &
\multicolumn{1}{c}{} & \multicolumn{1}{c}{} & \multicolumn{1}{c}{}
& \multicolumn{1}{|c}{$-1$} & \multicolumn{1}{c|}{$cu+\rho z_1$} \\
\multicolumn{1}{c|}{$z_3$} & \multicolumn{1}{|c}{} &
\multicolumn{1}{c}{} & \multicolumn{1}{c}{} & \multicolumn{1}{c}{}
& \multicolumn{1}{c}{} & \multicolumn{1}{c}{} &
\multicolumn{1}{|c}{} & \multicolumn{1}{c|}{$-1$}
\\ \cline{2-9}
\end{tabular} \]

\begin{center} {\bf Table 5} \end{center}

\vspace{0.3cm} De plus, une alg\`ebre r\'eelle dont la
multiplication est donn\'ee par la Table {\bf 5} est de Jordan non
commutative de division lin\'eaire et poss\`ede une d\'erivation
non triviale. Elle est de division lin\'eaire si et seulement si \
$\eta<0$ et $b\rho^2<4c\eta^2.$
\end{theorem}

\vspace{0.1cm} {\bf Preuve.} Suivant les notations de la Table
{\bf 3}, $b$ et $c$ s'\'ecrivent, respectivement, \
$\varepsilon|b|,$ $\varepsilon|c|$ \ o\`u $\varepsilon,
\varepsilon'\in\{1,-1\}.$ La table de multiplication de $A,$ par
rapport \`a la base \ $1,$ $u,$ $y_1,$ $z_1,$ $y_2,$ $\varepsilon
z_2:=z_2',$ $y_3,$ $\varepsilon'z_3:=z_3'$ \ est donn\'ee par:

\vspace{0.3cm}
{\footnotesize \[ \begin{tabular}{ccccccccc} \\
\multicolumn{1}{c}{} & \multicolumn{1}{c}{$1$} &
\multicolumn{1}{c}{$u$} & \multicolumn{1}{c}{$y_1$} &
\multicolumn{1}{c}{$z_1$} & \multicolumn{1}{c}{$y_2$} &
\multicolumn{1}{c}{$z_2'$} & \multicolumn{1}{c}{$y_3$} &
\multicolumn{1}{c}{$z_3'$} \\ \cline{2-9} \multicolumn{1}{c|}{$1$}
& \multicolumn{1}{|c}{$1$} & \multicolumn{1}{c|}{$u$} &
\multicolumn{1}{|c}{$y_1$} & \multicolumn{1}{c|}{$z_1$} &
\multicolumn{1}{|c}{$y_2$} & \multicolumn{1}{c|}{$z_2'$} &
\multicolumn{1}{|c}{$y_3$} & \multicolumn{1}{c|}{$z_3'$} \\
\multicolumn{1}{c|}{$u$} & \multicolumn{1}{|c}{} &
\multicolumn{1}{c|}{$-1$} & \multicolumn{1}{|c}{$az_1$} &
\multicolumn{1}{c|}{$-ay_1$} & \multicolumn{1}{|c}{$|b|z_2'$} &
\multicolumn{1}{c|}{$-|b|y_2$} & \multicolumn{1}{|c}{$|c|z_3'$} &
\multicolumn{1}{c|}{$-|c|y_3$} \\ \cline{2-9}
\multicolumn{1}{c|}{$y_1$} & \multicolumn{1}{|c}{} &
\multicolumn{1}{c|}{} & \multicolumn{1}{|c}{$-1$} &
\multicolumn{1}{c|}{$au$} & \multicolumn{1}{|c}{$\alpha y_3$} &
\multicolumn{1}{c|}{$\varepsilon\varepsilon'\alpha z_3'$} &
\multicolumn{1}{|c}{$-\alpha y_2$} &
\multicolumn{1}{c|}{$-\varepsilon\varepsilon'\alpha z_2'$} \\
\multicolumn{1}{c|}{$z_1$} & \multicolumn{1}{|c}{} &
\multicolumn{1}{c|}{} & \multicolumn{1}{|c}{} &
\multicolumn{1}{c|}{$-1$} & \multicolumn{1}{|c}{$\lambda
y_3+\varepsilon'\eta z_3'$} &
\multicolumn{1}{c|}{$-\varepsilon\eta
y_3+\varepsilon\varepsilon'\lambda z_3'$} &
\multicolumn{1}{|c}{$-\lambda y_2+\varepsilon\eta
z_2'+\varepsilon'\rho z_3'$} &
\multicolumn{1}{c|}{$-\varepsilon'\eta y_2-\varepsilon\varepsilon'\lambda z_2'-\varepsilon'\rho y_3$} \\
\cline{4-9} \multicolumn{1}{c|}{$y_2$} & \multicolumn{1}{|c}{} &
\multicolumn{1}{c}{} & \multicolumn{1}{c}{} &
\multicolumn{1}{c|}{} & \multicolumn{1}{|c}{$-1$} &
\multicolumn{1}{c|}{$|b|u$} & \multicolumn{1}{|c}{$\alpha
y_1+\lambda z_1$} &
\multicolumn{1}{c|}{$\varepsilon'\eta z_1$} \\
\multicolumn{1}{c|}{$z_2'$} & \multicolumn{1}{|c}{} &
\multicolumn{1}{c}{} & \multicolumn{1}{c}{} & \multicolumn{1}{c}{}
& \multicolumn{1}{|c}{} & \multicolumn{1}{c|}{$-1$} &
\multicolumn{1}{|c}{$-\varepsilon\eta z_1$} &
\multicolumn{1}{c|}{$\varepsilon\varepsilon'\alpha
y_1+\varepsilon\varepsilon'\lambda z_1$}
\\ \cline{6-9} \multicolumn{1}{c|}{$y_3$} & \multicolumn{1}{|c}{}
& \multicolumn{1}{c}{} & \multicolumn{1}{c}{} &
\multicolumn{1}{c}{} & \multicolumn{1}{c}{} & \multicolumn{1}{c}{}
& \multicolumn{1}{|c}{$-1$} & \multicolumn{1}{c|}{$|c|u+\varepsilon'\rho z_1$} \\
\multicolumn{1}{c|}{$z_3'$} & \multicolumn{1}{|c}{} &
\multicolumn{1}{c}{} & \multicolumn{1}{c}{} & \multicolumn{1}{c}{}
& \multicolumn{1}{c}{} & \multicolumn{1}{c}{} &
\multicolumn{1}{|c}{} & \multicolumn{1}{c|}{$-1$}
\\ \cline{2-9}
\end{tabular} \]}

\vspace{0.3cm} Comme $A$ est de division lin\'eaire, ceci
\'equivaut \`a  \ $-\varepsilon'\eta,$ $\varepsilon\eta,$
$-\varepsilon\varepsilon'>0$ \ et \
$|b|(-\varepsilon\varepsilon'\alpha^2)(\varepsilon'\rho)^2<4|c|(-\varepsilon'\alpha\eta)^2$
\ i.e. \ $\varepsilon\varepsilon'=-1,$ $\varepsilon\eta>0$ \ et \
$|b|\rho^2<4|c|(\varepsilon\eta)^2.\Box$

\vspace{0.4cm} \begin{remark} Dans {\em [BO 81$_1$]} Benkart et
Osborn ont montr\'e que les $\rit$-alg\`ebres de division
lin\'eaire de dimension $8,$ dont l'alg\`ebre de Lie des
d\'erivations est $su(3),$ sont ou bien \begin{enumerate} \item
des $su(3)$-modules irr\'eductibles et sont dans ce cas des
pseudo-octonions g\'en\'eralis\'es {\em (G-P Algebras [BBO 81],
Theorem {\bf 3.2})}, ou bien \item des $su(3)$-modules qui se
d\'ecomposent en somme de deux $su(3)$-modules irr\'eductibles de
dimension $1:$ \ $vect\{u\},$ $vect\{v\}$ \ et d'un $su(3)$-module
\ $Z=vect\{z_1,\dots,z_6\},$ irr\'eductible, de dimension $6,$ et
il existe dans ce cas des nombres r\'eels \ $\eta_i,$ $\theta_i,$
$\sigma_i,$ $tau_i$ o\`u $i\in\{1,2,3,4\}$ pour lesquels la
multiplication de telles alg\`ebres est donn\'ee par la table {\em
(Theorem {\bf 4.1})}:
\end{enumerate}
\end{remark}

\vspace{0.3cm}
{\tiny \[ \begin{tabular}{ccccccccc} \\
\multicolumn{1}{c}{} & \multicolumn{1}{c}{$u$} &
\multicolumn{1}{c}{$v$} & \multicolumn{1}{c}{$z_1$} &
\multicolumn{1}{c}{$z_2$} & \multicolumn{1}{c}{$z_3$} &
\multicolumn{1}{c}{$z_4$} & \multicolumn{1}{c}{$z_5$} &
\multicolumn{1}{c}{$z_6$} \\ \cline{2-9} \multicolumn{1}{c|}{$u$}
& \multicolumn{1}{|c}{$\eta_1u+\theta_1v$} &
\multicolumn{1}{c|}{$\eta_2u+\theta_2v$} &
\multicolumn{1}{|c}{$\sigma_1z_1+\sigma_2z_3$} &
\multicolumn{1}{c}{$\sigma_1z_2+\sigma_2z_6$} &
\multicolumn{1}{c}{$-\sigma_2z_1+\sigma_1z_3$} &
\multicolumn{1}{c}{$\sigma_1z_4+\sigma_2z_5$} &
\multicolumn{1}{c}{$-\sigma_2z_4+\sigma_1z_5$} & \multicolumn{1}{c|}{$-\sigma_2z_2+\sigma_1z_6$} \\
\multicolumn{1}{c|}{$v$} &
\multicolumn{1}{|c}{$\eta_3u+\theta_3v$} &
\multicolumn{1}{c|}{$\eta_4u+\theta_4v$} &
\multicolumn{1}{|c}{$\sigma_3z_4+\sigma_4z_3$} &
\multicolumn{1}{c}{$\sigma_1z_2+\sigma_4z_6$} &
\multicolumn{1}{c}{$-\sigma_4z_1+\sigma_3z_3$} &
\multicolumn{1}{c}{$\sigma_3z_4+\sigma_4z_5$} &
\multicolumn{1}{c}{$-\sigma_4z_4+\sigma_3z_5$} & \multicolumn{1}{c|}{$-\sigma_4z_2+\sigma_3z_6$} \\
\cline{2-9} \multicolumn{1}{c|}{$z_1$} &
\multicolumn{1}{|c}{$\tau_1z_1+\tau_2z_3$} &
\multicolumn{1}{c|}{$\tau_3z_1+\tau_4z_3$} &
\multicolumn{1}{|c}{$-u$} & \multicolumn{1}{c}{$z_4$} &
\multicolumn{1}{c}{$v$} & \multicolumn{1}{c}{$-z_2$} &
\multicolumn{1}{c}{$z_6$} &
\multicolumn{1}{c|}{$-z_5$} \\
\multicolumn{1}{c|}{$z_2$} &
\multicolumn{1}{|c}{$\tau_1z_2+\tau_2z_6$} &
\multicolumn{1}{c|}{$\tau_3z_2+\tau_4z_6$} &
\multicolumn{1}{|c}{$-z_4$} & \multicolumn{1}{c}{$-u$} &
\multicolumn{1}{c}{$z_5$} & \multicolumn{1}{c}{$z_1$} &
\multicolumn{1}{c}{$-z_3$} &
\multicolumn{1}{c|}{$v$} \\
 \multicolumn{1}{c|}{$z_3$} &
\multicolumn{1}{|c}{$-\tau_2z_1+\tau_1z_3$} &
\multicolumn{1}{c|}{$-\tau_4z_1+\tau_3z_3$} &
\multicolumn{1}{|c}{$-v$} & \multicolumn{1}{c}{$-z_5$} &
\multicolumn{1}{c}{$-u$} & \multicolumn{1}{c}{$z_6$} &
\multicolumn{1}{c}{$z_2$} &
\multicolumn{1}{c|}{$-z_4$} \\
\multicolumn{1}{c|}{$z_4$} &
\multicolumn{1}{|c}{$\tau_1z_4+\tau_2z_5$} &
\multicolumn{1}{c|}{$\tau_3z_4+\tau_4z_5$} &
\multicolumn{1}{|c}{$z_2$} & \multicolumn{1}{c}{$-z_1$} &
\multicolumn{1}{c}{$-z_6$} & \multicolumn{1}{c}{$-u$} &
\multicolumn{1}{c}{$v$} & \multicolumn{1}{c|}{$z_3$}
\\ \multicolumn{1}{c|}{$z_5$} &
\multicolumn{1}{|c}{$-\tau_2z_4+\tau_1z_5$} &
\multicolumn{1}{c|}{$-\tau_4z_4+\tau_3z_5$} &
\multicolumn{1}{|c}{$-z_6$} & \multicolumn{1}{c}{$z_3$} &
\multicolumn{1}{c}{$-z_2$} & \multicolumn{1}{c} {$-v$} &
\multicolumn{1}{c}{$-u$} & \multicolumn{1}{c|}{$z_1$}
\\ \multicolumn{1}{c|}{$z_6$} &
\multicolumn{1}{|c}{$-\tau_2z_2+\tau_1z_6$} &
\multicolumn{1}{c|}{$-\tau_4z_2+\tau_3z_6$} &
\multicolumn{1}{|c}{$z_5$} & \multicolumn{1}{c}{$-v$} &
\multicolumn{1}{c}{$z_4$} & \multicolumn{1}{c}{$-z_3$} &
\multicolumn{1}{c}{$-z_1$} & \multicolumn{1}{c|}{$-u$}
\\ \cline{2-9}
\end{tabular} \]} \begin{center} {\bf Table 6} \end{center}

\vspace{0.5cm} \hspace{0.3cm} Nous constatons ici que ces deux
situations ne peuvent avoir lieu pour une alg\`ebre de Jordan non
commutative. En effet, la premi\`ere est imm\'ediatement
\'elimin\'ee car les pseudo-octonions g\'en\'eralis\'es ne sont
pas de Jordan non commutatives [BBO 81]. L\'elimination de la
seconde est cons\'equence du r\'esultat suivant:

\vspace{0.3cm} \begin{proposition} Soit \ $A=\Big(
W,(.|.),\wedge\Big)$ \ une $\rit$-alg\`ebre de Jordan non
commutative de division lin\'eaire de dimension $8$ dont la
multiplication est donn\'ee par la Table {\bf 6}. Alors \
$A\simeq$ $\oit^{(\lambda)},$ o\`u
$\lambda\in\rit-\{\frac{1}{2}\}.$ \end{proposition}

\vspace{0.1cm} {\bf Preuve.} Si $\partial\in Der(A)=su(3)$ est une
d\'erivation de rang $6,$ elle est nulle sur $vect\{u\}$ et
$vect\{v\},$ et laisse stable $Z.$ Donc \ $Z=\partial Z\subseteq
W$ Proposition {\bf 5.3}. Ainsi les $z_i$ sont des vecteurs
(orthogonaux) aussi bien que $v=z_1z_3$ et $u$ est un multiple
scalaire de l'\'el\'ement unit\'e $1,$ de $A.$ Donc

\[ \theta_1=\eta_2=\sigma_2=\eta_3=\theta_4=\tau_2=0, \hspace{0.2cm}
\eta_1=\theta_2=\sigma_1=\theta_3=\tau_1 \ \mbox{ et } \ \eta_4<0.
\]

De plus, la propri\'et\'e trace de $(.|.)$ et
l'anti-commutativit\'e de $"\wedge"$ donnent \
$\sigma_3=\tau_3=0,$ \ $\sigma_4=-\eta_4>.$

\vspace{0.2cm} Enfin, quitte \`a changer le signe de $u,$ si
n\'ecessaire, on peut supposer $\eta_1>0.$ Ainsi \ $\eta_1^{-1}u$
est l'\'el\'ement unit\'e de $A,$ et en posant

\begin{eqnarray*} e_1 &=& \sigma_4^{-1}z_1, \\
e_2 &=& \sigma_4^{-1}z_2, \\
e_3 &=& \sigma_4^{-1}z_4, \\
e_4 &=& \sigma_4^{-1}v, \\
e_5 &=& -\sigma_4^{-1}z_3, \\
e_6 &=& -\sigma_4^{-1}z_6, \\
e_7 &=& -\sigma_4^{-1}z_5, \end{eqnarray*}

on obtient la table:

\vspace{0.5cm}
\[ \begin{tabular}{ccccccccc} \\
\multicolumn{1}{c}{} & \multicolumn{1}{c}{$1$} &
\multicolumn{1}{c}{$e_1$} & \multicolumn{1}{c}{$e_2$} &
\multicolumn{1}{c}{$e_3$} & \multicolumn{1}{c}{$e_4$} &
\multicolumn{1}{c}{$e_5$} & \multicolumn{1}{c}{$e_6$} &
\multicolumn{1}{c}{$e_7$} \\ \cline{2-9} \multicolumn{1}{c|}{$1$}
& \multicolumn{1}{|c}{$1$} & \multicolumn{1}{c|}{$e_1$} &
\multicolumn{1}{|c}{$e_2$} & \multicolumn{1}{c|}{$e_3$} &
\multicolumn{1}{|c}{$e_4$} & \multicolumn{1}{c|}{$e_5$} &
\multicolumn{1}{|c}{$e_6$} & \multicolumn{1}{c|}{$e_7$} \\
\multicolumn{1}{c|}{$e_1$} & \multicolumn{1}{|c}{} &
\multicolumn{1}{c|}{$-\beta$} & \multicolumn{1}{|c}{$e_3$} &
\multicolumn{1}{c|}{$-e_2$} & \multicolumn{1}{|c}{$e_5$} &
\multicolumn{1}{c|}{$-e_4$} &
\multicolumn{1}{|c}{$-e_7$} & \multicolumn{1}{c|}{$e_6$} \\
\cline{2-9} \multicolumn{1}{c|}{$e_2$} & \multicolumn{1}{|c}{} &
\multicolumn{1}{c|}{} & \multicolumn{1}{|c}{$-\beta$} &
\multicolumn{1}{c|}{$e_1$} & \multicolumn{1}{|c}{$e_6$} &
\multicolumn{1}{c|}{$e_7$} & \multicolumn{1}{|c}{$-e_4$} &
\multicolumn{1}{c|}{$-e_5$} \\
\multicolumn{1}{c|}{$e_3$} & \multicolumn{1}{|c}{} &
\multicolumn{1}{c|}{} & \multicolumn{1}{|c}{} &
\multicolumn{1}{c|}{$-\beta$} & \multicolumn{1}{|c}{$e_7$} &
\multicolumn{1}{c|}{$-e_6$} & \multicolumn{1}{|c}{$e_5$} &
\multicolumn{1}{c|}{$-e_4$} \\
\cline{4-9} \multicolumn{1}{c|}{$e_4$} & \multicolumn{1}{|c}{} &
\multicolumn{1}{c}{} & \multicolumn{1}{c}{} &
\multicolumn{1}{c|}{} & \multicolumn{1}{|c}{$-\beta$} &
\multicolumn{1}{c|}{$e_1$} & \multicolumn{1}{|c}{$e_2$} &
\multicolumn{1}{c|}{$e_3$} \\
\multicolumn{1}{c|}{$e_5$} & \multicolumn{1}{|c}{} &
\multicolumn{1}{c}{} & \multicolumn{1}{c}{} & \multicolumn{1}{c}{}
& \multicolumn{1}{|c}{} & \multicolumn{1}{c|}{$-\beta$} &
\multicolumn{1}{|c}{$-e_3$} & \multicolumn{1}{c|}{$e_2$}
\\ \cline{6-9} \multicolumn{1}{c|}{$e_6$} & \multicolumn{1}{|c}{}
& \multicolumn{1}{c}{} & \multicolumn{1}{c}{} &
\multicolumn{1}{c}{} & \multicolumn{1}{c}{} & \multicolumn{1}{c}{}
& \multicolumn{1}{|c}{$-\beta$} &
\multicolumn{1}{c|}{$-e_1$} \\
\multicolumn{1}{c|}{$e_7$} & \multicolumn{1}{|c}{} &
\multicolumn{1}{c}{} & \multicolumn{1}{c}{} & \multicolumn{1}{c}{}
& \multicolumn{1}{c}{} & \multicolumn{1}{c}{} &
\multicolumn{1}{|c}{} & \multicolumn{1}{c|}{$-\beta$}
\\ \cline{2-9}
\end{tabular} \]

\vspace{0.3cm} o\`u \ $\beta=\eta_1(-\eta_4)^{-2}.$ Donc \
$A\simeq$ $\oit^{(\lambda)}$ \ avec \
$\lambda=\frac{1}{2}(-\eta_4\eta_1^{-\frac{1}{2}}+1).\Box$

\vspace{0.5cm} \hspace{0.3cm} Nous terminons ce paragraphe par les
deux r\'esultats suivants:

\vspace{0.3cm} \begin{theorem} Soit \ $A=\Big(
W,(.|.),\wedge\Big)$ \ une $\rit$-alg\`ebre de Jordan non
commutative de division lin\'eaire de dimension $8.$ Alors \
$Der(A)=G_2$ compacte si et seulement si \ $A\simeq$
$\oit^{(\lambda)}$ o\`u $\lambda\in\rit-\{\frac{1}{2}\}.$
\end{theorem}

\vspace{0.1cm} {\bf Preuve.} $\Leftarrow/$ \ Si $\lambda$ est un
r\'eel distinct de $\frac{1}{2},$ alors \ $Der\Big(
\oit^{(\lambda)}\Big)=Der(\oit)=G_2$ compacte.

\vspace{0.2cm} $\Rightarrow/$ \ Benkart et Osborn [BO 81$_1$] ont
montr\'e que si $A$ est une $\rit$-alg\`ebre de division
lin\'eaire de dimension $8$ dont l'alg\`ebre de Lie des
d\'erivations est $G_2$ compacte, alors il existe une base \
$u,e_1,\dots,e_7$ \ de $A$ et trois param\`etres \ $\eta,$
$\zeta,$ $\beta$ \ avec $\eta\zeta\beta>0,$ pour lesquels la
multiplication de $A$ est donn\'ee par la table (Theorem {\bf
2.2}):

\vspace{0.5cm}
\[ \begin{tabular}{ccccccccc} \\
\multicolumn{1}{c}{} & \multicolumn{1}{c}{$u$} &
\multicolumn{1}{c}{$e_1$} & \multicolumn{1}{c}{$e_2$} &
\multicolumn{1}{c}{$e_3$} & \multicolumn{1}{c}{$e_4$} &
\multicolumn{1}{c}{$e_5$} & \multicolumn{1}{c}{$e_6$} &
\multicolumn{1}{c}{$e_7$} \\ \cline{2-9} \multicolumn{1}{c|}{$u$}
& \multicolumn{1}{|c|}{$u$} & \multicolumn{1}{c}{$\eta e_1$} &
\multicolumn{1}{c}{$\eta e_2$} & \multicolumn{1}{c}{$\eta e_3$} &
\multicolumn{1}{c}{$\eta e_4$} & \multicolumn{1}{c}{$\eta e_5$} &
\multicolumn{1}{c}{$\eta e_6$} & \multicolumn{1}{c|}{$\eta e_7$}
\\ \cline{2-9} \multicolumn{1}{c|}{$e_1$} & \multicolumn{1}{|c|}{$\zeta e_1$} &
\multicolumn{1}{|c}{$-\beta u$} & \multicolumn{1}{c}{$e_4$} &
\multicolumn{1}{c}{$e_7$} & \multicolumn{1}{c}{$-e_2$} &
\multicolumn{1}{c}{$e_6$} & \multicolumn{1}{c}{$-e_5$} &
\multicolumn{1}{c|}{$-e_3$} \\ \multicolumn{1}{c|}{$e_2$} &
\multicolumn{1}{|c|}{$\zeta e_2$} & \multicolumn{1}{c}{$-e_4$} &
\multicolumn{1}{c}{$-\beta u$} & \multicolumn{1}{c}{$e_5$} &
\multicolumn{1}{c}{$e_1$} & \multicolumn{1}{c}{$-e_3$} &
\multicolumn{1}{c}{$e_7$} &
\multicolumn{1}{c|}{$-e_6$} \\
\multicolumn{1}{c|}{$e_3$} & \multicolumn{1}{|c|}{$\zeta e_3$} &
\multicolumn{1}{c}{$-e_7$} & \multicolumn{1}{c}{$-e_5$} &
\multicolumn{1}{c}{$-\beta u$} & \multicolumn{1}{c}{$e_6$} &
\multicolumn{1}{c}{$e_2$} & \multicolumn{1}{c}{$-e_4$} &
\multicolumn{1}{c|}{$e_1$} \\ \multicolumn{1}{c|}{$e_4$} &
\multicolumn{1}{|c|}{$\zeta e_4$} & \multicolumn{1}{c}{$e_2$} &
\multicolumn{1}{c}{$-e_1$} & \multicolumn{1}{c}{$-e_6$} &
\multicolumn{1}{c}{$-\beta u$} & \multicolumn{1}{c}{$e_7$} &
\multicolumn{1}{c}{$e_3$} &
\multicolumn{1}{c|}{$-e_5$} \\
\multicolumn{1}{c|}{$e_5$} & \multicolumn{1}{|c|}{$\zeta e_5$} &
\multicolumn{1}{c}{$-e_6$} & \multicolumn{1}{c}{$e_3$} &
\multicolumn{1}{c}{$-e_2$} & \multicolumn{1}{c}{$-e_7$} &
\multicolumn{1}{c}{$-\beta u$} & \multicolumn{1}{c}{$e_1$} &
\multicolumn{1}{c|}{$e_4$}
\\ \multicolumn{1}{c|}{$e_5$} & \multicolumn{1}{|c|}{$\zeta e_6$}
& \multicolumn{1}{c}{$e_5$} & \multicolumn{1}{c}{$-e_7$} &
\multicolumn{1}{c}{$e_4$} & \multicolumn{1}{c}{$-e_3$} &
\multicolumn{1}{c}{$-e_1$} & \multicolumn{1}{c}{$-\beta u$} &
\multicolumn{1}{c|}{$e_2$} \\
\multicolumn{1}{c|}{$e_7$} & \multicolumn{1}{|c|}{$\zeta e_7$} &
\multicolumn{1}{c}{$e_3$} & \multicolumn{1}{c}{$e_6$} &
\multicolumn{1}{c}{$-e_1$} & \multicolumn{1}{c}{$e_5$} &
\multicolumn{1}{c}{$-e_4$} & \multicolumn{1}{c}{$-e_2$} &
\multicolumn{1}{c|}{$-\beta u$}
\\ \cline{2-9}
\end{tabular} \]

\vspace{0.3cm} Si, de plus, $A$ est unitaire, on a \
$\eta=\zeta=1$ \ i.e. $\beta>0$ et $A\simeq$ $\oit^{(\lambda)}$ \
o\`u \ $\lambda=\frac{1}{2}(\beta^{-\frac{1}{2}}+1).\Box$

\vspace{0.3cm} \begin{theorem} Soit \ $A=\Big(
W,(.|.),\wedge\Big)$ \ une $\rit$-alg\`ebre de Jordan non
commutative de division lin\'eaire de dimension $8.$ Alors \
$Der(A)=su(2)\oplus su(2)$ si et seulement si \ $A\simeq$ $\Big(
E_{-1,\alpha,0}(\hit)\Big)^{(\lambda)}$ avec
$1\neq\alpha>\frac{1}{2}$ et $\lambda\neq\frac{1}{2}.$
\end{theorem}

\vspace{0.1cm} {\bf Preuve.} Benkart et Osborn ([BO 81$_1$]
Theorem {\bf 5.1}) qu'une $\rit$-alg\`ebre $A,$ de division
lin\'eaire de dimension $8,$ dont l'alg\`ebre de Lie des
d\'erivations contient $su(2)\oplus su(2),$ poss\`ede une base \
$u,$ $x_1,$ $x_2,$ $x_3,$ $y_1,$ $y_2,$ $y_3,$ $y_4$ \ pour
laquelle la multiplication est donn\'ee par la table:

\vspace{0.3cm}
\[ \begin{tabular}{ccccccccc} \\
\multicolumn{1}{c}{} & \multicolumn{1}{c}{$u$} &
\multicolumn{1}{c}{$x_1$} & \multicolumn{1}{c}{$x_2$} &
\multicolumn{1}{c}{$x_3$} & \multicolumn{1}{c}{$y_1$} &
\multicolumn{1}{c}{$y_2$} & \multicolumn{1}{c}{$y_3$} &
\multicolumn{1}{c}{$y_4$} \\ \cline{2-9} \multicolumn{1}{c|}{$u$}
& \multicolumn{1}{|c|}{$u$} & \multicolumn{1}{|c}{$\zeta x_1$} &
\multicolumn{1}{c}{$\zeta x_2$} & \multicolumn{1}{c|}{$\zeta x_3$}
& \multicolumn{1}{|c}{$\rho y_1$} & \multicolumn{1}{c}{$\rho y_2$}
& \multicolumn{1}{c}{$\rho y_3$} & \multicolumn{1}{c|}{$\rho y_4$}
\\ \cline{2-9} \multicolumn{1}{c|}{$x_1$} & \multicolumn{1}{|c|}{$\theta x_1$} &
\multicolumn{1}{|c}{$\beta u$} & \multicolumn{1}{c}{$x_3$} &
\multicolumn{1}{c|}{$-x_2$} & \multicolumn{1}{|c}{$\varepsilon
y_4$} & \multicolumn{1}{c}{$\varepsilon y_3$}
& \multicolumn{1}{c}{$-\varepsilon y_2$} & \multicolumn{1}{c|}{$-\varepsilon y_1$} \\
\multicolumn{1}{c|}{$x_2$} & \multicolumn{1}{|c|}{$\theta x_2$} &
\multicolumn{1}{|c}{$-x_3$} & \multicolumn{1}{c}{$\beta u$} &
\multicolumn{1}{c|}{$x_1$} & \multicolumn{1}{|c}{$\varepsilon
y_2$} & \multicolumn{1}{c}{$-\varepsilon y_1$} &
\multicolumn{1}{c}{$\varepsilon y_4$} &
\multicolumn{1}{c|}{$-\varepsilon y_3$} \\
\multicolumn{1}{c|}{$x_3$} & \multicolumn{1}{|c|}{$\theta x_3$} &
\multicolumn{1}{|c}{$x_2$} & \multicolumn{1}{c}{$-x_1$} &
\multicolumn{1}{c|}{$\beta u$} & \multicolumn{1}{|c}{$-\varepsilon
y_3$} & \multicolumn{1}{c}{$\varepsilon y_4$} &
\multicolumn{1}{c}{$\varepsilon y_1$} &
\multicolumn{1}{c|}{$-\varepsilon y_2$} \\
\cline{2-9} \multicolumn{1}{c|}{$y_1$} &
\multicolumn{1}{|c|}{$\sigma y_1$} & \multicolumn{1}{|c}{$-\eta
y_4$} & \multicolumn{1}{c}{$-\eta y_2$} &
\multicolumn{1}{c|}{$\eta y_3$} & \multicolumn{1}{|c}{$\delta u$}
& \multicolumn{1}{c}{$\gamma x_2$} & \multicolumn{1}{c}{$-\gamma
x_3$} &
\multicolumn{1}{c|}{$\gamma x_1$} \\
\multicolumn{1}{c|}{$y_2$} & \multicolumn{1}{|c|}{$\sigma y_2$} &
\multicolumn{1}{|c}{$-\eta y_3$} & \multicolumn{1}{c}{$\eta y_1$}
& \multicolumn{1}{c|}{$-\eta y_4$} & \multicolumn{1}{|c}{$-\gamma
x_2$} & \multicolumn{1}{c}{$\delta u$} &
\multicolumn{1}{c}{$\gamma x_1$} & \multicolumn{1}{c|}{$\gamma
x_3$}
\\ \multicolumn{1}{c|}{$y_3$} & \multicolumn{1}{|c|}{$\sigma y_3$}
& \multicolumn{1}{|c}{$\eta y_2$} & \multicolumn{1}{c}{$-\eta
y_4$} & \multicolumn{1}{c|}{$\eta y_1$} &
\multicolumn{1}{|c}{$\gamma x_3$} & \multicolumn{1}{c}{$-\gamma
x_1$} & \multicolumn{1}{c}{$\delta u$} &
\multicolumn{1}{c|}{$\gamma x_2$} \\
\multicolumn{1}{c|}{$y_4$} & \multicolumn{1}{|c|}{$\sigma y_4$} &
\multicolumn{1}{|c}{$\eta y_1$} & \multicolumn{1}{c}{$\eta y_3$} &
\multicolumn{1}{c|}{$\eta y_2$} & \multicolumn{1}{|c}{$-\gamma
x_1$} & \multicolumn{1}{c}{$-\gamma x_3$} &
\multicolumn{1}{c}{$-\gamma x_2$} & \multicolumn{1}{c|}{$\delta
u$}
\\ \cline{2-9}
\end{tabular} \]

\vspace{0.4cm} o\`u \ $\beta,$ $\gamma,$ $\delta,$ $\varepsilon,$
$\eta,$ $\zeta,$ $\theta,$ $\rho,$ $\sigma$ \ sont des nombres
r\'eels. Si, de plus, $A=\Big( W,(.|.),\wedge\Big)$ est de Jordan
non commutative, alors $u$ (idempotent non nul) est l'\'el\'ement
unit\'e de $A.$ Donc \ $\zeta=\rho=\theta=\sigma=1,$ \ de plus, la
fait que $(.|.)$ soit une forme trace d\'efinie n\'egative, et
l'anti-commutativit\'e de $"\wedge",$ donnent \ $\beta, \delta<0,$
$\varepsilon=\eta$ et $\varepsilon\delta=\gamma\beta.$ On pose
alors

\begin{eqnarray*} x_i' &=& (-\beta)^{-\frac{1}{2}}x_i, \ i=1,2,3; \\
y_j' &=& (-\delta)^{-\frac{1}{2}}y_j, \ j=1,2,4; \\
y_3' &=& -(-\delta)^{-\frac{1}{2}}y_3. \end{eqnarray*}

On obtient la table:

\vspace{0.3cm}
{\footnotesize \[ \begin{tabular}{ccccccccc} \\
\multicolumn{1}{c}{} & \multicolumn{1}{c}{$u$} &
\multicolumn{1}{c}{$x_1'$} & \multicolumn{1}{c}{$x_2'$} &
\multicolumn{1}{c}{$x_3'$} & \multicolumn{1}{c}{$y_1'$} &
\multicolumn{1}{c}{$y_4'$} & \multicolumn{1}{c}{$y_2'$} &
\multicolumn{1}{c}{$y_3'$} \\ \cline{2-9} \multicolumn{1}{c|}{$u$}
& \multicolumn{1}{|c}{$u$} & \multicolumn{1}{c|}{$x_1'$} &
\multicolumn{1}{|c}{$x_2'$} & \multicolumn{1}{c|}{$x_3'$} &
\multicolumn{1}{|c}{$y_1'$} & \multicolumn{1}{c|}{$y_4'$} &
\multicolumn{1}{|c}{$y_2'$} &
\multicolumn{1}{c|}{$y_3'$} \\
\multicolumn{1}{c|}{$x_1'$} & \multicolumn{1}{|c}{} &
\multicolumn{1}{c|}{$-u$} &
\multicolumn{1}{|c}{$(-\beta)^{-\frac{1}{2}}x_3'$} &
\multicolumn{1}{c|}{$-(-\beta)^{-\frac{1}{2}}x_2'$} &
\multicolumn{1}{|c}{$\varepsilon(-\beta)^{-\frac{1}{2}}y_4'$} &
\multicolumn{1}{c|}{$-\varepsilon(-\beta)^{-\frac{1}{2}}y_1'$} &
\multicolumn{1}{|c}{$-\varepsilon(-\beta)^{-\frac{1}{2}}y_3'$} &
\multicolumn{1}{c|}{$\varepsilon(-\beta)^{-\frac{1}{2}}y_2'$} \\
\cline{2-9} \multicolumn{1}{c|}{$x_2'$} & \multicolumn{1}{|c}{} &
\multicolumn{1}{c|}{} & \multicolumn{1}{|c}{$-u$} &
\multicolumn{1}{c|}{$(-\beta)^{-\frac{1}{2}}x_1'$} &
\multicolumn{1}{|c}{$\varepsilon(-\beta)^{-\frac{1}{2}}y_2'$} &
\multicolumn{1}{c|}{$\varepsilon(-\beta)^{-\frac{1}{2}}y_3'$} &
\multicolumn{1}{|c}{$-\varepsilon(-\beta)^{-\frac{1}{2}}y_1'$} &
\multicolumn{1}{c|}{$-\varepsilon(-\beta)^{-\frac{1}{2}}y_4'$} \\
\multicolumn{1}{c|}{$x_3'$} & \multicolumn{1}{|c}{} &
\multicolumn{1}{c|}{} & \multicolumn{1}{|c}{} &
\multicolumn{1}{c|}{$-u$} &
\multicolumn{1}{|c}{$\varepsilon(-\beta)^{-\frac{1}{2}}y_3'$} &
\multicolumn{1}{c|}{$-\varepsilon(-\beta)^{-\frac{1}{2}}y_2'$} &
\multicolumn{1}{|c}{$\varepsilon(-\beta)^{-\frac{1}{2}}y_4'$} &
\multicolumn{1}{c|}{$-\varepsilon(-\beta)^{-\frac{1}{2}}y_1'$} \\
\cline{4-9} \multicolumn{1}{c|}{$y_1'$} & \multicolumn{1}{|c}{} &
\multicolumn{1}{c}{} & \multicolumn{1}{c}{} &
\multicolumn{1}{c|}{} & \multicolumn{1}{|c}{$-u$} &
\multicolumn{1}{c|}{$\varepsilon(-\beta)^{-\frac{1}{2}}x_1'$} &
\multicolumn{1}{|c}{$\varepsilon(-\beta)^{-\frac{1}{2}}x_2'$} &
\multicolumn{1}{c|}{$\varepsilon(-\beta)^{-\frac{1}{2}}x_3'$} \\
\multicolumn{1}{c|}{$y_4'$} & \multicolumn{1}{|c}{} &
\multicolumn{1}{c}{} & \multicolumn{1}{c}{} & \multicolumn{1}{c}{}
& \multicolumn{1}{|c}{} & \multicolumn{1}{c|}{$-u$} &
\multicolumn{1}{|c}{$-\varepsilon(-\beta)^{-\frac{1}{2}}x_3'$} &
\multicolumn{1}{c|}{$\varepsilon(-\beta)^{-\frac{1}{2}}x_2'$}
\\ \cline{6-9} \multicolumn{1}{c|}{$y_2'$} & \multicolumn{1}{|c}{}
& \multicolumn{1}{c}{} & \multicolumn{1}{c}{} &
\multicolumn{1}{c}{} & \multicolumn{1}{c}{} & \multicolumn{1}{c}{}
& \multicolumn{1}{|c}{$-u$} & \multicolumn{1}{c|}{$-\varepsilon(-\beta)^{-\frac{1}{2}}x_1'$} \\
\multicolumn{1}{c|}{$y_3'$} & \multicolumn{1}{|c}{} &
\multicolumn{1}{c}{} & \multicolumn{1}{c}{} & \multicolumn{1}{c}{}
& \multicolumn{1}{c}{} & \multicolumn{1}{c}{} &
\multicolumn{1}{|c}{} & \multicolumn{1}{c|}{$-u$}
\\ \cline{2-9}
\end{tabular} \]}

\vspace{0.5cm} Donc \ $A^{\Big(
\frac{\varepsilon+(-\beta)^{\frac{1}{2}}}{2\varepsilon}\Big)}\simeq$
$E_{-1,\frac{1+\varepsilon}{2\varepsilon},0}(\hit)$ \ i.e.

\begin{eqnarray*} A &=& \Big( A^{\Big(
\frac{\varepsilon+(-\beta)^{\frac{1}{2}}}{2\varepsilon}\Big)}\Big)^{\Big(
\frac{1}{2}(-\beta)^{-\frac{1}{2}}
(\varepsilon+(-\beta)^{\frac{1}{2}})\Big)} \\
&\simeq& \Big(
E_{-1,\frac{1+\varepsilon}{2\varepsilon},0}(\hit)\Big)^{\Big(
\frac{1}{2}(-\beta)^{-\frac{1}{2}}
(\varepsilon+(-\beta)^{\frac{1}{2}})\Big)}. \end{eqnarray*}

De plus, $A$ est de division lin\'eaire si et seulement si \
$\varepsilon>0$ (Corollaire {\bf 4.14}) et on a \ $\varepsilon\neq
1$ car $Der(A)\neq G_2$ compacte i.e. \
$1\neq\alpha=\frac{1+\varepsilon}{2\varepsilon}>\frac{1}{2}.\Box$

\vspace{0.5cm} \begin{remark} Benkart et Osborn {\em [BO 81$_1$]}
posent le probl\`eme de l'existence d'une $\rit$-alg\`ebre de
division lin\'eaire de dimension $8,$ dont l'alg\`ebre de Lie des
d\'erivations est $su(2),$ et dont la d\'ecomposition en
$su(2)$-modules irr\'eductibles est de la forme: \ $1+1+3+3.$ \
Nous donnons ici une r\'eponse affirmative avec une alg\`ebre de
Jordan non commutative. En effet, nous avons vu que si $\delta\neq
0$ et $\alpha\neq\frac{1}{2},$ alors

\[ Der\Big(
E_{-1,\alpha,\delta}(\hit)\Big)=Der(\hit)=su(2) \]

{\em (Proposition {\bf 3.20})}. Nous remarquons de plus, que \
$vect\{1\},$ $vect\{f\},$ $vect\{i,j,k\},$ $vect\{if,jf,kf\},$
o\`u $f=(0,1),$ sont des $su(2)$-modules irr\'eductibles de
$E_{-1,\alpha,\delta}(\hit),$ qui est leur somme directe {\em
(Proposition {\bf 3.23})}. Les conditions suppl\'ementaires \
$\alpha>\frac{1}{2}$ et $(2\alpha-1)\delta^2<4$ \ assurent que
l'alg\`ebre $E_{-1,\alpha,\delta}(\hit)$ est de division
lin\'eaire.$\Box$
\end{remark}

\vspace{1cm}
\subsection{Caract\'erisation des $\rit$-alg\`ebres de Jordan n.c.
de d.l. de dimension $8$ ayant un automorphisme non trivial}

\vspace{0.5cm} \hspace{0.3cm} Dans ce dernier paragraphe nous
donnons une caract\'erisation des alg\`ebres r\'eelles de Jordan
non commutatives de division lin\'eaire de dimension 8 dont le
groupe des automorphismes est non trivial. Nous donnons ensuite un
exemple o\`u le groupe des automorphismes est trivial, ce qui met
en \'evidence l'immensit\'e de la classe des $\rit$-alg\`ebres de
Jordan non commutatives de division lin\'eaire de dimension $8.$

\vspace{0.5cm} \begin{definition} Soit $A$ une $K$-alg\`ebre. On
appelle reflexion de $A$ tout automorphisme involutif de $A,$ non
identique.$\Box$ \end{definition}

\vspace{0.5cm}
\begin{proposition} Soit \ $A=\Big( W,(.|.),\wedge\Big)$ une
$\rit$-alg\`ebre de Jordan non commutative de division lin\'eaire
de dimension $8$ et soit $f$ un automorphisme de l'espace
vectoriel $A.$ Alors les deux affirmations suivantes sont
\'equivalentes: \begin{enumerate} \item $f$ est une reflexion de
$A.$ \item $f$ est sym\'etrique et isom\'etrique par rapport \`a
$(.|.)$ et le sous-espace vectoriel \ $\ker(f-I_A):=B$ \ est une
sous-alg\`ebre de $A,$ de dimension $4$ qui contient le carr\'e de
orthogonal $B^{\perp}.$ \end{enumerate}
\end{proposition}

\vspace{0.1cm} {\bf Preuve.} 1) $\Rightarrow$ 2). Soient $x,y\in
B,$ on a \ $f(xy)=f(x)f(y)=xy\in B.$ Donc $B$ est une
sous-alg\`ebre de $A.$ Comme $f$ est une isom\'etrie, elle est
sym\'etrique: ${^tf}=f^{-1}=f.$ Donc \ $\ker(f-I_A)$ et
$\ker(f+I_A)$ \ sont suppl\'ementaires et on a: \[
\ker(f+I_A)=B^{\perp}. \] Si $x,y\in B^{\perp},$ on a \
$f(xy)=f(x)f(y)=(-x)(-y)=xy,$ \ i.e. $B^{\perp}B^{\perp}\subseteq
B.$ De plus, $BB^{\perp}$ et $B^{\perp}B$ sont inclus dans
$B^{\perp},$ en vertu de la propri\'et\'e trace de $(.|.).$ Comme
$A$ est de division lin\'eaire, on a: \ $\dim(B^{\perp})\leq
\dim(B)\leq\dim(B^{\perp}).$ Donc \ $\dim(B^{\perp})=\dim(B)=4.$

\vspace{0.2cm} 2) $\Rightarrow$ 1). Soit $x\in A,$ alors \
$f(x)=x$ si et seulement si $x\in B.$ Ainsi, pour tout $(y,z)\in
B^{\perp}\times B,$ on a \ $(f(y)|z)=(y|f(z))=(y|z)=0$ \ i.e.
$f(B^{\perp})\subseteq B^{\perp}$ $(\subseteq W).$ Si $x\in
B^{\perp}$ est un vecteur propre de $f,$ associ\'e \`a la valeur
propre $\lambda,$ on a

\begin{eqnarray*} x^2 &=& (x|x) \\
&=& \Big( f(x)|f(x)\Big) \ \mbox{ car } f \mbox{ est isom\'etrique par rapport \`a } (.|.) \\
&=& f(x)^2 \\
&=& (\lambda x)^2 \\
&=& \lambda^2x^2 \end{eqnarray*}

i.e. $\lambda^2=1.$ Comme $x\notin B,$ on a $\lambda=-1.$ Ainsi,
l'application $f:A=B\oplus B^{\perp}\rightarrow A$ est d\'efinie
par \ $x+y\mapsto x-y$ \ et est une reflexion de $A.$ En effet, on
a $f^2=I_A$ et pour tous $x,x'\in B$ et $y,y'\in B^{\perp}:$

\begin{eqnarray*} f\Big( (x+y)(x'+y')\Big) &=& f\Big(
(xx'+yy')+(xy'+yx')\Big) \\
&=& (xx'+yy')-(xy'+yx') \hspace{0.2cm} \mbox{ car }
B^{\perp}B^{\perp}\subseteq
B \ \mbox{ et } \ BB^{\perp}, B^{\perp}B\subseteq B^{\perp} \\
&=& (x-y)(x'-y') \\
&=& f(x+y)f(x'+y'). \end{eqnarray*}

De plus, $f\neq I_A.\Box$

\vspace{0.5cm}
\begin{remark} Soit $B$ une sous-alg\`ebre de $A,$ de dimension $4,$ qui contient le carr\'e de
orthogonal $B^{\perp}.$ Alors l'application lin\'eaire d\'efinie
par \[ B\oplus B^{\perp}=A\rightarrow A\hspace{0.3cm} x+y\mapsto
x-y \] est une reflexion de $A.\Box$ \end{remark}

\vspace{0.5cm}
\begin{lemma} $A=\Big( W,(.|.),\wedge\Big)$ une
$\rit$-alg\`ebre de Jordan non commutative de division lin\'eaire
de dimension $8$ et soit $f$ une isom\'etrie de l'espace euclidien
$(W,-(.|.))$ qui laisse fixe un vecteur non nul $u$ de $W$ et qui
commute avec l'op\'erateur $L_u^*$ {\em (Lemma {\bf 3.1})}. Alors
l'espace vectoriel $W(u)$ se d\'ecompose en une somme directe
orthogonale de sous-espaces vectoriels $H_i,$ o\`u
$i\in\{1,2,3\},$ de la forme $vect\{y_i,u\wedge y_i\}$ stables par
$f$ et par $f_u.$ De plus, la restriction de $f$ \`a chaque $H_i$
est une rotation. \end{lemma}

\vspace{0.1cm} {\bf Preuve.} Puisque $f$ laisse fixe $u,$ elle
induit une isom\'etrie \ $f_0:W(u)\rightarrow W(u),$ et l'on
distingue les trois cas suivants: \begin{enumerate} \item Si
$f_u^2$ poss\`ede au moins deux valeurs propres distinctes, alors
tout sous-espace propre $E,$ ainsi que son orthogonal $E^{\perp}$
dans $W(u),$ sont stables par $f_u,$ et \'egalement par $f_0$
puisque $f_u$ et $f_0$ commutent. L'un de ces sous-espaces est de
dimension $2,$ l'autre de dimension $4.$ \item Si $f_u^2$
poss\`ede une unique valeur propre et $f_0$ n'en poss\`ede aucune,
on consid\`ere le polyn\^ome minimal $P(X)$ de $f_0$ et l'on
distingue les deux sous-cas suivants:
\begin{enumerate} \item Si $P(X)$ est de degr\'e $2,$ de la forme \
$X^2-\alpha X-\beta$ \ avec $\alpha^2+4\beta<0,$ alors il existe
$\gamma>0$ tel que \ $(f_0-\frac{1}{2}\alpha I_{W(u)})^2=\gamma
f_u^2$ \ i.e.

\[ (f_0-\frac{1}{2}\alpha I_{W(u)}-\sqrt{\gamma}f_u)(f_0-\frac{1}{2}\alpha I_{W(u)}+\sqrt{\gamma}f_u)\equiv 0.
\]

\vspace{0.3cm} {\bf i)} Si $f_0-\frac{1}{2}\alpha
I_{W(u)}=\pm\sqrt{\gamma}f_u$ alors tout sous-espace vectoriel de
$W(u)$ stable par $f_u$ est \'egalement stable par $f_0,$ aussi
bien que son suppl\'ementaire orthogonal.

\vspace{0.3cm} {\bf ii)} Si $f_0-\frac{1}{2}\alpha
I_{W(u)}\neq\pm\sqrt{\gamma}f_u$ alors le sous-espace vectoriel \[
\ker(f_0-\frac{1}{2}\alpha I_{W(u)}-\sqrt{\gamma}f_u):=H, \] de
$W(u),$ ainsi que son orthogonal $H^{\perp},$ sont propres et
stables par $f_u$ et $f_0.$ Leur dimension est un nombre pair, \`a
savoir $2$ ou $4.$ \item Si $P(X)$ est de degr\'e $>2,$ on
consid\`ere une composante irr\'eductible $P_1(X)$ de $P(X)$ et
l'on obtient un sous-espace $\ker\Big( P_1(f_0)\Big)$ de $W(u),$
propre et stable par $f_0$ et $f_u,$ aussi bien que son
orthogonal.
\end{enumerate} \item Si $f_u^2$
poss\`ede une unique valeur propre et $f_0$ admet un vecteur
propre $y,$ on a \ $f_0(y)=\pm y$ et $f_0(uy)=f_0\Big(
f_u(y)\Big)=f_u\Big( f_0(y)\Big)=\pm uy.$ Donc $vect\{y,uy\}:=H$
est stable par $f_0$ et \'egalement par $f_u$ car $f_u^2$ est une
homoth\'etie. Son orthogonal $H^{\perp}$ est \'egalement stable
par $f_0$ et $f_u.$
\end{enumerate}

\vspace{0.2cm} \hspace{0.3cm} On se ram\`ene alors \`a un
sous-espace vectoriel de $W(u),$ de dimension $4$ stable par $f_0$
et $f_u,$ et l'on rep\`ete le m\^eme raisonnement. Soit maintenant
$H_i=vect\{y_i,z_i\},$ $y_i, z_i$ orthonormaux, un sous-espace
stable par $f_0$ et $f_u.$ Alors $H_i$ est form\'e de vecteurs
propres de $f_u^2$ associ\'es \`a une valeur propre $\lambda_i<0$
et $z_i=(-\lambda_i)^{-\frac{1}{2}}uy_i.$ D'autre part, il existe
des scalaires $a_i, b_i\in\rit,$ avec $a_i^2+b_i^2=1,$ tels que
$f_0(y_i)=a_iy_i+b_iz_i,$ et on a

\begin{eqnarray*} f_0(z_i) &=& (-\lambda_i)^{-\frac{1}{2}}f_0(uy_i) \\
&=& (-\lambda_i)^{-\frac{1}{2}}uf_0(y_i) \\
&=& (-\lambda_i)^{-\frac{1}{2}}u(a_iy_i+b_iz_i) \\
&=& (-\lambda_i)^{-\frac{1}{2}}\Big( a_i(-\lambda_i)^{\frac{1}{2}}z_i-b_i(-\lambda_i)^{\frac{1}{2}}y_i\Big) \\
&=& -b_iy_i+a_iz_i \\
\end{eqnarray*}

i.e. $f_{0_{/H_i}}$ est une rotation.$\Box$

\vspace{1.5cm}
\begin{proposition} Soit \ $A=\Big( W,(.|.),\wedge\Big)$ une
$\rit$-alg\`ebre de Jordan non commutative de division lin\'eaire
de dimension $8.$ Alors les deux propri\'et\'es suivantes sont
\'equivalentes:
\begin{enumerate} \item Le groupe des automorphismes de $A,$ \ $Aut(A),$ est non trivial. \item $Aut(A)$
contient une reflexion.
\end{enumerate}
\end{proposition}

\vspace{0.1cm} {\bf Preuve.} 1) $\Rightarrow$ 2). Soit $g^*$ un
automorphisme de l'alg\`ebre $A$ qui n'est pas une reflexion,
alors $g^*$ induit une isom\'etrie $g$ de l'espace euclidien
$-(W,(.|.))$ qui admet, \'evidemment, un vecteur propre $u$
associ\'e \`a la valeur propre $\lambda=\pm 1.$ L'automorphisme
$f^*=g^{*^2}$ de $A,$ non identique que l'on peut supposer non
involutif, laisse fixe $u.$ En outren l'isom\'etrie $f$ qu'il
induit sur $(-W,(.|.))$ commute avec $L_u^*.$ En utilisant les
notations du Lemme {\bf 5.18}, on a \begin{eqnarray*} f(y_iz_i) &=& f^*(y_iz_i) \\
&=& f^*(y_)f^*(z_i) \\
&=& f_0(y_i)f_0(z_i) \\
&=& (a_i^2+b_i^2)y_iz_i \\
&=& y_iz_i.
\end{eqnarray*}

On distingue les deux cas suivants: \begin{enumerate} \item Si $f$
n'admet pas de vecteurs propres lin\'eairement ind\'ependants \`a
$u,$ alors $y_iz_i=(-\lambda_i)^{\frac{1}{2}}u.$ Ceci \'etant pour
tout $i\in\{1,2,3\},$ on a: \ $B_i^{\perp}B_i^{\perp}\subseteq
B_i$ o\`u $B_i$ est la sous-alg\`ebre de $A,$ de dimension $4,$
engendr\'ee par $u$ et $y_i$ et qui coincide avec
$vect\{1,u,y_i,z_i\}.$ En effet, en notant $H=vect\{1,u\}$ et,
pour $i\in\{1,2,3\},$ $vect\{y_i,z_i\}$ on a:
\[ B_i^{\perp}=\bigoplus_{j\neq i}H_j. \] De plus, $H_iH_j$ est,
d'apr\`es la propri\'et\'e trace de $(.|.),$ orthogonal \`a \
$H+H_i+H_j$ si $i\neq j$ \ i.e. $H_iH_j=H_k$ pour toute
permutation $(i\ j\ k)$ de $\{1,2,3\}.$ Donc \
$B_i^{\perp}B_i^{\perp}\subseteq B_i.$ Ainsi, $Aut(A)$ contient
une reflexion dans ce premier cas. \item Si $f$ admet un vecteur
propre norm\'e $y_1$ lin\'eairement ind\'ependant \`a $u,$ que
l'on peut supposer orthogonal \`a $u,$ alors $y_1\in W(u)$ est un
vecteur propre de $f_0.$ L'\'el\'ement $uy_1$ est \'egalement un
vecteur propre de $f_0$ associ\'e \`a la m\^eme valeur propre $\pm
1$ que celle de $y_1,$ que l'on peut supposer \'egale \`a $1$ en
consid\'erant l'automorphisme $f^{*^2}.$ Comme $f^{*^2}$ est
suppos\'e non identique, l'espace propre \ $B=\ker(f^{*^2}-I_A)$ \
est une sous-alg\`ebre de $A$ de dimension $4,$ qui coincide avec
$vect\{1,u,y_1,uy_1\},$ et qui contient le carr\'e de son
orthogonal. En effet, en reprenant les notations pr\'ec\'edentes,
on a \ $B^{\perp}=H_2+H_3.$ De plus, pour $i=2,3$ on a \
$H_i^2\subset B$ car $f^{*^2}(y_iz_i)=y_iz_i.$ La propri\'et\'e
trace de $(.|.)$ montre alors que $H_2H_3=H_3H_2\subset B$ pour
$i=2,3.$ Donc $B^{\perp}B^{\perp}\subseteq B,$ et $Aut(A)$
contient une reflexion dans ce cas \'egalement.$\Box$
\end{enumerate}

\vspace{0.1cm}
\begin{theorem} Soit $A$ une $\rit$-alg\`ebre de Jordan non commutative de division lin\'eaire
de dimension $8,$ ayant un automorphisme non trivial. Alors il
existe une base \ $1,$ $u,$ $y_1,$ $z_1,$ $y_2,$ $z_2,$ $y_3,$
$z_3$ \ de $A,$ trois param\`etres \ $a,$ $b,$ $c>0$ \ et neuf
autres \ $\alpha,$ $\beta,$ $\gamma,$ $\mu,$ $\lambda,$ $\eta,$
$\sigma,$ $\delta,$ $\rho$ \ pour lesquels la multiplication de
$A$ est donn\'ee par la Table {\bf 6} suivante:

\vspace{0.1cm}
\[ \begin{tabular}{ccccccccc} \\
\multicolumn{1}{c}{} & \multicolumn{1}{c}{$1$} &
\multicolumn{1}{c}{$u$} & \multicolumn{1}{c}{$y_1$} &
\multicolumn{1}{c}{$z_1$} & \multicolumn{1}{c}{$y_2$} &
\multicolumn{1}{c}{$z_2$} & \multicolumn{1}{c}{$y_3$} &
\multicolumn{1}{c}{$z_3$}
\\ \cline{2-9}
\multicolumn{1}{c|}{$1$} & \multicolumn{1}{|c}{$1$} &
\multicolumn{1}{c|}{$u$} & \multicolumn{1}{|c}{$y_1$} &
\multicolumn{1}{c|}{$z_1$} & \multicolumn{1}{|c}{$y_2$} &
\multicolumn{1}{c|}{$z_2$} & \multicolumn{1}{|c}{$y_3$} &
\multicolumn{1}{c|}{$z_3$}
\\ \multicolumn{1}{c|}{$u$} & \multicolumn{1}{|c}{} &
\multicolumn{1}{c|}{$-1$} & \multicolumn{1}{|c}{$az_1$} &
\multicolumn{1}{c|}{$-ay_1$} & \multicolumn{1}{|c}{$bz_2$} &
\multicolumn{1}{c|}{$-by_2$} & \multicolumn{1}{|c}{$cz_3$} &
\multicolumn{1}{c|}{$-cy_3$}
\\ \cline{2-9}
\multicolumn{1}{c|}{$y_1$} & \multicolumn{1}{|c}{} &
\multicolumn{1}{c|}{} & \multicolumn{1}{|c}{$-1$} &
\multicolumn{1}{c|}{$au$} & \multicolumn{1}{|c}{$\alpha y_3+\beta
z_3$} & \multicolumn{1}{c|}{$\gamma y_3+\mu z_3$} &
\multicolumn{1}{|c}{$-\alpha y_2-\gamma z_2$} &
\multicolumn{1}{c|}{$-\beta y_2-\mu z_2$}
\\ \multicolumn{1}{c|}{$z_1$} & \multicolumn{1}{|c}{} &
\multicolumn{1}{c|}{} & \multicolumn{1}{|c}{} &
\multicolumn{1}{c|}{$-1$} & \multicolumn{1}{|c}{$\lambda y_3+\eta
z_3$} & \multicolumn{1}{c|}{$\sigma y_3+\delta z_3$} &
\multicolumn{1}{|c}{$-\lambda y_2-\sigma z_2+\rho z_3$} &
\multicolumn{1}{c|}{$-\eta y_2-\delta z_2-\rho y_3$}
\\ \cline{4-9}
\multicolumn{1}{c|}{$y_2$} & \multicolumn{1}{|c}{} &
\multicolumn{1}{c}{} & \multicolumn{1}{c}{} &
\multicolumn{1}{c|}{} & \multicolumn{1}{|c}{$-1$} &
\multicolumn{1}{c|}{$bu$} & \multicolumn{1}{|c}{$\alpha
y_1+\lambda z_1$} & \multicolumn{1}{c|}{$\beta y_1+\eta z_1$}
\\ \multicolumn{1}{c|}{$z_2$} & \multicolumn{1}{|c}{} &
\multicolumn{1}{c}{} & \multicolumn{1}{c}{} & \multicolumn{1}{c}{}
& \multicolumn{1}{|c}{} & \multicolumn{1}{c|}{$-1$} &
\multicolumn{1}{|c}{$\gamma y_1+\sigma z_1$} &
\multicolumn{1}{c|}{$\mu y_1+\delta z_1$}
\\ \cline{6-9}
\multicolumn{1}{c|}{$y_3$} & \multicolumn{1}{|c}{} &
\multicolumn{1}{c}{} & \multicolumn{1}{c}{} & \multicolumn{1}{c}{}
& \multicolumn{1}{c}{} & \multicolumn{1}{c}{} &
\multicolumn{1}{|c}{$-1$} & \multicolumn{1}{c|}{$cu+\rho z_1$}
\\ \multicolumn{1}{c|}{$z_3$} & \multicolumn{1}{|c}{} &
\multicolumn{1}{c}{} & \multicolumn{1}{c}{} & \multicolumn{1}{c}{}
& \multicolumn{1}{c}{} & \multicolumn{1}{c}{} &
\multicolumn{1}{|c}{} & \multicolumn{1}{c|}{$-1$}
\\ \cline{2-9}
\end{tabular} \]
\begin{center} {\bf Table 6} \end{center}

\vspace{0.1cm} De plus, une $\rit$-alg\`ebre dont la
multiplication est donn\'ee par la Table {\bf 6} est de Jordan non
commutative et poss\`ede un automorphisme non trivial. Elle est de
division lin\'eaire si et seulement si \ $\beta\gamma-\alpha\mu,$
$\beta\lambda-\alpha\eta,$ $\gamma\lambda-\alpha\sigma>0$ \ et \[
c(\alpha\delta-\beta\sigma-\lambda\mu+\gamma\eta)^2+b(\beta\gamma-\alpha\mu)\rho^2<4c(\beta\lambda-\alpha\eta)
(\sigma\mu-\gamma\delta). \]
\end{theorem}

\vspace{0.1cm} {\bf Preuve.} Il existe une sous-alg\`ebre \
$B=vect\{1,u,y_1,z_1\}$ \ de $A$ de dimension $4$ qui contient le
carr\'e de son orthogonal, les vecteurs $u,$ $y_1,$ $z_1$ \'etant
orthonormaux et tels que \ $uy_1=a z_1,$ $uz_1=-ay_1$ \ et \
$y_1z_1=au$ \ o\`u $a$ est un param\`etre $>0.$ Ainsi, la
Proposition {\bf 3.14} assure l'existence de deuc vecteurs
orthonormaux $y_2,$ $z_2$ orthogonaux aux pr\'ec\'edents, tels que
$y_2z_2=bu$ o\`u $b$ est un param\`etre que l'on peut choisir
$>0,$ quitte \`a changer le signe de $y_2,$ si n\'ecessaire. En
compl\'etant \ $1,$ $u,$ $y_1,$ $z_1,$ $y_2,$ $z_2$ \ en une base
orthonorm\'ee $1,$ $u,$ $y_1,$ $z_1,$ $y_2,$ $z_2,$ $y_3,$ $z_3$ \
de $A$ on obtient, en tenant compte de la propri\'et\'e trace de
$(.|.),$ une premi\`ere table de multiplication:

\vspace{0.2cm}
\[ \begin{tabular}{ccccccccc} \\
\multicolumn{1}{c}{} & \multicolumn{1}{c}{$1$} &
\multicolumn{1}{c}{$u$} & \multicolumn{1}{c}{$y_1$} &
\multicolumn{1}{c}{$z_1$} & \multicolumn{1}{c}{$y_2$} &
\multicolumn{1}{c}{$z_2$} & \multicolumn{1}{c}{$y_3$} &
\multicolumn{1}{c}{$z_3$}
\\ \cline{2-9}
\multicolumn{1}{c|}{$1$} & \multicolumn{1}{|c}{$1$} &
\multicolumn{1}{c|}{$u$} & \multicolumn{1}{|c}{$y_1$} &
\multicolumn{1}{c|}{$z_1$} & \multicolumn{1}{|c}{$y_2$} &
\multicolumn{1}{c|}{$z_2$} & \multicolumn{1}{|c}{$y_3$} &
\multicolumn{1}{c|}{$z_3$}
\\
\multicolumn{1}{c|}{$u$} & \multicolumn{1}{|c}{} &
\multicolumn{1}{c|}{$-1$} & \multicolumn{1}{|c}{$az_1$} &
\multicolumn{1}{c|}{$-ay_1$} & \multicolumn{1}{|c}{$bz_2$} &
\multicolumn{1}{c|}{$-by_2$} & \multicolumn{1}{|c}{$cz_3$} &
\multicolumn{1}{c|}{$-cy_3$}
\\ \cline{2-9}
\multicolumn{1}{c|}{$y_1$} & \multicolumn{1}{|c}{} &
\multicolumn{1}{c|}{} & \multicolumn{1}{|c}{$-1$} &
\multicolumn{1}{c|}{$au$} & \multicolumn{1}{|c}{$\alpha y_3+\beta
z_3$} & \multicolumn{1}{c|}{$\gamma y_3+\mu z_3$} &
\multicolumn{1}{|c}{$-\alpha y_2-\gamma z_2+\nu z_3$} &
\multicolumn{1}{c|}{$-\beta y_2-\mu z_2-\nu y_3$}
\\ \multicolumn{1}{c|}{$z_1$} & \multicolumn{1}{|c}{} &
\multicolumn{1}{c|}{} & \multicolumn{1}{|c}{} &
\multicolumn{1}{c|}{$-1$} & \multicolumn{1}{|c}{$\lambda y_3+\eta
z_3$} & \multicolumn{1}{c|}{$\sigma y_3+\delta z_3$} &
\multicolumn{1}{|c}{$-\lambda y_2-\sigma z_2+\rho z_3$} &
\multicolumn{1}{c|}{$-\eta y_2-\delta z_2-\rho y_3$}
\\ \cline{4-9}
\multicolumn{1}{c|}{$y_2$} & \multicolumn{1}{|c}{} &
\multicolumn{1}{c}{} & \multicolumn{1}{c}{} &
\multicolumn{1}{c|}{} & \multicolumn{1}{|c}{$-1$} &
\multicolumn{1}{c|}{$bu$} & \multicolumn{1}{|c}{$\alpha
y_1+\lambda z_1$} & \multicolumn{1}{c|}{$\beta y_1+\eta z_1$}
\\ \multicolumn{1}{c|}{$z_2$} & \multicolumn{1}{|c}{} &
\multicolumn{1}{c}{} & \multicolumn{1}{c}{} & \multicolumn{1}{c}{}
& \multicolumn{1}{|c}{} & \multicolumn{1}{c|}{$-1$} &
\multicolumn{1}{|c}{$\gamma y_1+\sigma z_1$} &
\multicolumn{1}{c|}{$\mu y_1+\delta z_1$}
\\ \cline{6-9}
\multicolumn{1}{c|}{$y_3$} & \multicolumn{1}{|c}{} &
\multicolumn{1}{c}{} & \multicolumn{1}{c}{} & \multicolumn{1}{c}{}
& \multicolumn{1}{c}{} & \multicolumn{1}{c}{} &
\multicolumn{1}{|c}{$-1$} & \multicolumn{1}{c|}{$cu+\nu y_1+\rho
z_1$} \\ \multicolumn{1}{c|}{$z_3$} & \multicolumn{1}{|c}{} &
\multicolumn{1}{c}{} & \multicolumn{1}{c}{} & \multicolumn{1}{c}{}
& \multicolumn{1}{c}{} & \multicolumn{1}{c}{} &
\multicolumn{1}{|c}{} & \multicolumn{1}{c|}{$-1$}
\\ \cline{2-9}
\end{tabular} \]

\vspace{0.2cm} o\`u $c,$ $\alpha,$ $\beta,$ $\gamma,$ $\mu,$
$\lambda,$ $\eta,$ $\sigma,$ $\delta,$ $\rho$ $\nu\in\rit.$ De
plus, quitte \`a changer le signe de $y_3,$ si n\'ecessaire, on
peut supposer $c>0.$ Si $\nu'=(\nu^2+\rho^2)^{\frac{1}{2}}\neq 0,$
on pose \[ y_1'=\nu'^{-1}(\rho y_1-\nu z_1), \hspace{0.3cm} z_1'=
\nu'^{-1}(\nu+\rho z_1) \] et l'on obtient la table d\'esir\'ee en
consid\'erant la nouvelle base \ $1,$ $u,$ $y_1',$ $z_1',$ $y_2,$
$z_2,$ $y_3,$ $z_3.$ R\'eciproquement, une $\rit$-alg\`ebre dont
la multiplication est donn\'ee par la Table {\bf 6} est de Jordan
non commutative. Le reste de la d\'emonstration est assur\'e par
la Remarque {\bf 5.17} et le Th\'eor\`eme {\bf 4.16}.$\Box$

\vspace{0.3cm}
\begin{proposition} Soit \ $\oit=\Big( W,(.|.),\wedge\Big)$ l'alg\`ebre r\'eelle de Cayley-Dickson et soient
$\varphi$ un automorphisme de l'espace vectoriel r\'eel $W$ et $f$
une isom\'etrie de l'espace euclidien $(W,-(.|.)).$ Alors
\begin{enumerate} \item $\tilde{f}\in Aut\Big( \oit(\varphi)\Big)$ si et seulement si $\overline{\varphi f\varphi^{-1}}
\in G_2.$ Dans ce cas, $f$ commute avec $\varphi^*\varphi.$ \item
$\tilde{f}$ est une reflexion de $\oit(\varphi):=\Big(
W,(.|.),\Delta\Big)$ si et seulement si $\overline{\varphi
f\varphi^{-1}}$ est une reflexion de \ $\oit.$ Dans ce cas, $f$
commute avec $\varphi^*\varphi,$ et $\varphi^*\varphi$ laisse
stable une sous-alg\`ebre de $(W,\Delta)$ de dimension $3$
\end{enumerate}
\end{proposition}

\vspace{0.1cm} {\bf Preuve.} \begin{enumerate} \item On a
\begin{eqnarray*} \tilde{f}\in Aut\Big( \oit(\varphi)\Big) &\Leftrightarrow& \Big( \oit(\varphi)\Big)(f)=\oit(\varphi) \\
&\Leftrightarrow& \oit(\varphi f\varphi^{-1})=\oit \\
&\Leftrightarrow& \overline{\varphi f\varphi^{-1}}\in G_2.
\end{eqnarray*}

Dans ces conditions \ $\varphi f\varphi^{-1}\in O_7(\rit)$ \ i.e.
\ $f(\varphi^*\varphi)=(\varphi^*\varphi)f.$ \item On v\'erifie
facilement que $\tilde{f}$ est involutif non identique, si et
seulement si $\overline{\varphi f\varphi^{-1}}$ est involutif non
identique. De plus, $f$ laisse fixes, uniquement, les \'el\'ements
d'une sous-alg\`ebre $W_0$ de $(W,\Delta),$ de dimension $3,$ et
comme $f(\varphi^*\varphi)=(\varphi^*\varphi)f,$ l'application
$\varphi^*\varphi$ laisse stable $W_0.\Box$
\end{enumerate}

\vspace{0.3cm}
\begin{lemma} Toute sous-alg\`ebre de \ $\oit$ de dimension $4$ contient le carr\'e de son orthogonal.
\end{lemma}

\vspace{0.1cm} {\bf Preuve.} Soit $B$ une telle sous-alg\`ebre
contenant un vecteur non nul $u$ et soient $v\in B^{\perp}-\{0\}$
et $w\in(B+vect\{v,uv\})^{\perp}-\{0\}.$ Alors \ $v,$ $uv,$ $w,$
$uw$ \ est une base orthogonale de $B^{\perp}$ avec laquelle on
v\'erifie facilement que \ $B^{\perp}B^{\perp}\subseteq B.\Box$

\vspace{0.3cm}
\begin{lemma} Soit \ $\oit=\Big( W,(.|.),\wedge\Big)$ l'alg\`ebre r\'eelle de Cayley-Dickson et soient
$\varphi$ un automorphisme de l'espace vectoriel r\'eel $W$ qui
laisse stable une sous-alg\`ebre $W_0$ de $(W,\wedge)$ de
dimension $3,$ ainsi que son orthogonal dans $W.$ Alors \ $\rit
1\oplus W_0:=B$ \ est une sous-alg\`ebre de \
$\oit(\varphi):=\Big( W,(.|.),\Delta\Big),$ \ qui contient la
carr\'e de son orthogonal.
\end{lemma}

\vspace{0.1cm} {\bf Preuve.} $\varphi^*$ laisse stable $W_0$ et
pour tous $x,y\in W_0$ \ on a \ $\varphi^*\Big(
\varphi(x)\wedge\varphi(y)\Big)\in\varphi^*(W_0)\subseteq W_0.$ Si
$x,y\in B^{\perp}=W_0^{\perp},$ on a \ $\varphi(x)\wedge
\varphi(y)\in B\cap W=W_0$ (Lemme {\bf 5.22}). Donc $x\Delta y\in
W_0.\Box$

\vspace{0.3cm}
\begin{lemma} Soit \ $\oit=\Big( W,(.|.),\wedge\Big)$ l'alg\`ebre r\'eelle de Cayley-Dickson et soit
$s$ un automorphisme sym\'etrique de l'espace euclidien
$(W,-(.|.)),$ d\'efini positif, qui laisse stable une
sous-alg\`ebre $(W,\wedge)$ de dimension $3.$ Alors le groupe des
automorphismes de l'alg\`ebre \ $\oit(s)$ n'est pas trivial.
\end{lemma}

\vspace{0.1cm} {\bf Preuve.} Soit $W_0$ une sous-alg\`ebre de
$(W,\wedge)$ de dimension $3,$ stable par $s.$ Alors son
orthogonal dans $W$ est \'egalement stable par $s.$ La Remarque
{\bf 5.17} et le Lemme {\bf 5.23} montrent que $\oit(s)$ contient
une reflexion. $\Box$

\vspace{0.3cm}
\begin{theorem} Soit \ $\oit=\Big( W,(.|.),\wedge\Big)$ l'alg\`ebre r\'eelle de Cayley-Dickson et
soit $s$ un automorphisme sym\'etrique de l'espace euclidien
$(W,-(.|.)),$ d\'efini positif. Alors les deux propri\'et\'es
suivantes sont \'equivalentes: \begin{enumerate} \item $Aut\Big(
\oit(s)\Big)$ n'est pas trivial. \item $\tilde{s}$ laisse stable
une sous-alg\`ebre de $\oit$ de dimension $4.$ \end{enumerate}
\end{theorem}

\vspace{0.1cm} {\bf Preuve.} 1) $\Rightarrow$ 2). L'existence
d'une reflexion dans $\oit(s)=\Big( W,(.|.),\Delta\Big)$ montre
que $s^2$ laisse stable une sous-alg\`ebre $W_0,$ de $(W,\Delta)$
de dimension $3.$ Comme $s$ et $s^2$ ont les m\^emes vecteurs
propres, l'application $s$ laisse stable $W_0.$ Cette derni\`ere
est \'egalement une sous-alg\`ebre de $(W,\wedge).$

\vspace{0.2cm} L'implication 2) $\Rightarrow$ 1) est cons\'equence
du Lemme {\bf 5.24}.$\Box$

\vspace{0.3cm}
\begin{note} Soient \ ${\cal C}$ \ la classe de toutes les \ $\rit$-alg\`ebres de Jordan non commutatives de
division lin\'eaire de dimension $8,$ \ ${\cal G}$ \ la
sous-classe de ${\cal C}$ constitu\'e des alg\`ebres qui
s'obtiennent, \`a partir de $\rit,$ par mutations et par extension
cayleyenne g\'en\'eralis\'ee {\em (Corollaire {\bf 3.4})}, \
${\cal D}$ la sous-classe de ${\cal C}$ constitu\'e des alg\`ebres
ayant une d\'erivation non triviale, et enfin \ ${\cal A}$ la
sous-classe de ${\cal C}$ constitu\'e des alg\`ebres ayant un
automorphisme non trivial. Nous avons les inclusion suivantes: \
${\cal D}\subset {\cal A}\subset {\cal G}\subset {\cal C}.$
\end{note}

\vspace{0.3cm} \hspace{0.3cm} Nous allons montrer, pour terminer
ce dernier paragraphe, que les inclusion pr\'ec\'edentes sont
strictes. Nous consid\'erons l'alg\`ebre r\'eelle $A$ ayant une
base \ ${\cal B}=\{1,u1,y_1,z_1,y_2,z_2,y_3,z_3\}$ \ pour laquelle
la multiplication est donn\'ee par la Table {\bf 7}
suivante:

\vspace{0.1cm}
\[ \begin{tabular}{ccccccccc} \\
\multicolumn{1}{c}{} & \multicolumn{1}{c}{$1$} &
\multicolumn{1}{c}{$u$} & \multicolumn{1}{c}{$y_1$} &
\multicolumn{1}{c}{$z_1$} & \multicolumn{1}{c}{$y_2$} &
\multicolumn{1}{c}{$z_2$} & \multicolumn{1}{c}{$y_3$} &
\multicolumn{1}{c}{$z_3$}
\\ \cline{2-9}
\multicolumn{1}{c|}{$1$} & \multicolumn{1}{|c}{$1$} &
\multicolumn{1}{c|}{$u$} & \multicolumn{1}{|c}{$y_1$} &
\multicolumn{1}{c|}{$z_1$} & \multicolumn{1}{|c}{$y_2$} &
\multicolumn{1}{c|}{$z_2$} & \multicolumn{1}{|c}{$y_3$} &
\multicolumn{1}{c|}{$z_3$}
\\ \multicolumn{1}{c|}{$u$} & \multicolumn{1}{|c}{} &
\multicolumn{1}{c|}{$-1$} & \multicolumn{1}{|c}{$az_1$} &
\multicolumn{1}{c|}{$-ay_1$} & \multicolumn{1}{|c}{$bz_2$} &
\multicolumn{1}{c|}{$-by_2$} & \multicolumn{1}{|c}{$cz_3$} &
\multicolumn{1}{c|}{$-cy_3$}
\\ \cline{2-9}
\multicolumn{1}{c|}{$y_1$} & \multicolumn{1}{|c}{} &
\multicolumn{1}{c|}{} & \multicolumn{1}{|c}{$-1$} &
\multicolumn{1}{c|}{$au$} & \multicolumn{1}{|c}{$\alpha y_3$} &
\multicolumn{1}{c|}{$\mu z_3$} & \multicolumn{1}{|c}{$-\alpha
y_2$} & \multicolumn{1}{c|}{$-\mu z_2$}
\\ \multicolumn{1}{c|}{$z_1$} & \multicolumn{1}{|c}{} &
\multicolumn{1}{c|}{} & \multicolumn{1}{|c}{} &
\multicolumn{1}{c|}{$-1$} & \multicolumn{1}{|c}{$\mu z_3$} &
\multicolumn{1}{c|}{$\mu y_3$} & \multicolumn{1}{|c}{$-\mu z_2$} &
\multicolumn{1}{c|}{$-\mu y_2$}
\\ \cline{4-9}
\multicolumn{1}{c|}{$y_2$} & \multicolumn{1}{|c}{} &
\multicolumn{1}{c}{} & \multicolumn{1}{c}{} &
\multicolumn{1}{c|}{} & \multicolumn{1}{|c}{$-1$} &
\multicolumn{1}{c|}{$bu$} & \multicolumn{1}{|c}{$\alpha y_1$} &
\multicolumn{1}{c|}{$\mu z_1$}
\\ \multicolumn{1}{c|}{$z_2$} & \multicolumn{1}{|c}{} &
\multicolumn{1}{c}{} & \multicolumn{1}{c}{} & \multicolumn{1}{c}{}
& \multicolumn{1}{|c}{} & \multicolumn{1}{c|}{$-1$} &
\multicolumn{1}{|c}{$\mu z_1$} & \multicolumn{1}{c|}{$\mu y_1$}
\\ \cline{6-9}
\multicolumn{1}{c|}{$y_3$} & \multicolumn{1}{|c}{} &
\multicolumn{1}{c}{} & \multicolumn{1}{c}{} & \multicolumn{1}{c}{}
& \multicolumn{1}{c}{} & \multicolumn{1}{c}{} &
\multicolumn{1}{|c}{$-1$} & \multicolumn{1}{c|}{$cu$}
\\ \multicolumn{1}{c|}{$z_3$} & \multicolumn{1}{|c}{} &
\multicolumn{1}{c}{} & \multicolumn{1}{c}{} & \multicolumn{1}{c}{}
& \multicolumn{1}{c}{} & \multicolumn{1}{c}{} &
\multicolumn{1}{|c}{} & \multicolumn{1}{c|}{$-1$}
\\ \cline{2-9}
\end{tabular} \]
\begin{center} {\bf Table 7} \end{center}

\vspace{0.3cm} o\`u \ $a,$ $b,$ $c$ \ sont des param\`etres $>0$ \
et \ $\alpha,$ $\mu\in\rit.$ Alors $A$ est de Jordan non
commutative, elle est de division lin\'eaire si et seulement si \
$\alpha\mu<0.$ Nous utiliserons ces notations dans ce qui suit:

\vspace{0.2cm} \begin{proposition} Si $\alpha=-\mu=1$ et
$c>a=b>1,$ alors \ $A\in{\cal D}$ \ et \ $A\notin{\cal G}.$
\end{proposition}

\vspace{0.1cm} {\bf Preuve.} L'alg\`ebre $A$ admet une
d\'erivation, de rang $6,$ en vertu de la Proposition {\bf 5.8}.
Supposons maintenant que $A$ s'obtient, \`a partir de $\rit,$ par
mutations et par extension cayleyenne g\'en\'eralis\'ee, alors il
existerait une mutation $B=A^{(\lambda)},$
$\lambda\neq\frac{1}{2},$ de $A$ et un vecteur non nul $u$ de $A$
pour lequel l'op\'erateur $f_u^2$ soit une homoth\'etie de $B$
(Remarque {\bf 3.16 4)}). On d\'efinit une nouvelle base \ ${\cal
B}_0=\{1,e_1,\dots,e_7\}$ \ de $A$ en posant

\[ e_i=y_i, \hspace{0.2cm} e_{i+4}=-z_i \ \mbox{ pour } i=1,2,3 \ \mbox{ et } e_4=u. \]

Soit, alors, \ $x=\sum_{1\leq i\leq 7}a_ie_i$ \ un vecteur
quelconque de $B.$ Les cinq premiers termes diagonaux de la
matrice de l'op\'erateur \ $f_x^2$ \ de $B,$ par rapport \`a la
base ${\cal B}_0,$ sont respectivement

\begin{eqnarray*} a_{11} &=& -\sum_{1\leq i\leq 7}a_i^2, \\
a_{22} &=& -a_1^2-\lambda'^2(a_2^2+a_3^2+a_6^2+a_7^2)-\lambda'^2a^2(a_4^2+a_5^2), \\
a_{33} &=& -a_2^2-\lambda'^2(a_1^2+a_3^2+a_5^2+a_7^2)-\lambda'^2a^2(a_4^2+a_6^2), \\
a_{44} &=& -a_3^2-\lambda'^2(a_1^2+a_2^2+a_5^2+a_6^2)-\lambda'^2c^2(a_4^2+a_7^2), \\
a_{55} &=&
-a_4^2-\lambda'^2a^2(a_1^2+a_2^2+a_5^2+a_6^2)-\lambda'^2c^2(a_3^2+a_7^2).
\end{eqnarray*}

o\`u  $\lambda'=2\lambda-1.$ On suppose que $L_x^2$ est une
homoth\'etie, i.e. les $a_{ii}$ sont tous \'egaux, et l'on
distingue les trois cas suivants:

\begin{enumerate}
\item Si $\lambda'^2=a^{-2},$ alors les \'egalit\'es: \
$a_{11}=a_{22}=a_{33}$ donnent \ $a_2=a_3=a_6=a_7=a_1=a_5=0$ car
$a^2\neq 1.$ Ainsi $a_{11}=a_{44}$ donne $a_4=0$ car $a^2\neq
c^2,$ i.e. $x=0.$ \item Si $\lambda'^2< a^{-2},$ alors \
$a_{11}=a_{22}$ donne \ $a_2=\dots=a_7=0$ et $a_{11}=a_{33}$ donne
$a_1=0$ i.e. $x=0.$ \item Si $\lambda'^2> a^{-2},$ alors \
$a_{11}=a_{55}$ donne \ $a_1=a_2=a_3=a_5=a_6=a_7=0$ et
$a_{11}=a_{22}$ donne $a_4=0,$ i.e. $x=0.$
\end{enumerate}

\vspace{0.2cm} \hspace{0.3cm} Ainsi, aucune mutation de $A$ ne
poss\`ede un vecteur non nul $v$ qui satisfait \`a la
propri\'et\'e {\bf 4)} de la Remarque {\bf 3.16} i.e.
$A\notin{\cal G}.\Box$

\vspace{0.4cm} \hspace{0.3cm} La Proposition {\bf 5.27} montre que
le {\bf Proc\'ede de Cayley-Dickson G\'en\'eralis\'e} est
insuffisant pour la d\'etermination de toutes les
$\rit$-alg\`ebres de Jordan non commutatives de division
lin\'eaire de dimension $8.$

\vspace{0.3cm} \begin{proposition} Si $\mu=-a=-1,$ $b, c, \alpha
>0,$ $b<1<c$ et $\alpha\neq 1, b, c$ alors \ $A\in{\cal A}$ \ et \ $A\notin{\cal D}.$
\end{proposition}

\vspace{0.1cm} {\bf Preuve.} En effet, $A\in{\cal A}$ et on a la
table:

\vspace{0.1cm}
\[ \begin{tabular}{ccccccccc} \\
\multicolumn{1}{c}{} & \multicolumn{1}{c}{$1$} &
\multicolumn{1}{c}{$e_1$} & \multicolumn{1}{c}{$e_2$} &
\multicolumn{1}{c}{$e_3$} & \multicolumn{1}{c}{$e_4$} &
\multicolumn{1}{c}{$e_5$} & \multicolumn{1}{c}{$e_6$} &
\multicolumn{1}{c}{$e_7$}
\\ \cline{2-9}
\multicolumn{1}{c|}{$1$} & \multicolumn{1}{|c}{$1$} &
\multicolumn{1}{c|}{$e_1$} & \multicolumn{1}{|c}{$e_2$} &
\multicolumn{1}{c|}{$e_3$} & \multicolumn{1}{|c}{$e_4$} &
\multicolumn{1}{c|}{$e_5$} & \multicolumn{1}{|c}{$e_6$} &
\multicolumn{1}{c|}{$e_7$}
\\ \multicolumn{1}{c|}{$e_1$} & \multicolumn{1}{|c}{} &
\multicolumn{1}{c|}{$-1$} & \multicolumn{1}{|c}{$\alpha e_3$} &
\multicolumn{1}{c|}{$-\alpha e_2$} & \multicolumn{1}{|c}{$e_5$} &
\multicolumn{1}{c|}{$-e_4$} & \multicolumn{1}{|c}{$-e_7$} &
\multicolumn{1}{c|}{$e_6$}
\\ \cline{2-9}
\multicolumn{1}{c|}{$e_2$} & \multicolumn{1}{|c}{} &
\multicolumn{1}{c|}{} & \multicolumn{1}{|c}{$-1$} &
\multicolumn{1}{c|}{$\alpha e_1$} & \multicolumn{1}{|c}{$be_6$} &
\multicolumn{1}{c|}{$e_7$} & \multicolumn{1}{|c}{$-be_4$} &
\multicolumn{1}{c|}{$-e_5$}
\\ \multicolumn{1}{c|}{$e_3$} & \multicolumn{1}{|c}{} &
\multicolumn{1}{c|}{} & \multicolumn{1}{|c}{} &
\multicolumn{1}{c|}{$-1$} & \multicolumn{1}{|c}{$ce_7$} &
\multicolumn{1}{c|}{$-e_6$} & \multicolumn{1}{|c}{$e_5$} &
\multicolumn{1}{c|}{$-ce_4$}
\\ \cline{4-9}
\multicolumn{1}{c|}{$e_4$} & \multicolumn{1}{|c}{} &
\multicolumn{1}{c}{} & \multicolumn{1}{c}{} &
\multicolumn{1}{c|}{} & \multicolumn{1}{|c}{$-1$} &
\multicolumn{1}{c|}{$e_1$} & \multicolumn{1}{|c}{$be_2$} &
\multicolumn{1}{c|}{$ce_3$}
\\ \multicolumn{1}{c|}{$e_5$} & \multicolumn{1}{|c}{} &
\multicolumn{1}{c}{} & \multicolumn{1}{c}{} & \multicolumn{1}{c}{}
& \multicolumn{1}{|c}{} & \multicolumn{1}{c|}{$-1$} &
\multicolumn{1}{|c}{$-e_3$} & \multicolumn{1}{c|}{$e_2$}
\\ \cline{6-9}
\multicolumn{1}{c|}{$e_6$} & \multicolumn{1}{|c}{} &
\multicolumn{1}{c}{} & \multicolumn{1}{c}{} & \multicolumn{1}{c}{}
& \multicolumn{1}{c}{} & \multicolumn{1}{c}{} &
\multicolumn{1}{|c}{$-1$} & \multicolumn{1}{c|}{$-e_1$}
\\ \multicolumn{1}{c|}{$z_3$} & \multicolumn{1}{|c}{} &
\multicolumn{1}{c}{} & \multicolumn{1}{c}{} & \multicolumn{1}{c}{}
& \multicolumn{1}{c}{} & \multicolumn{1}{c}{} &
\multicolumn{1}{|c}{} & \multicolumn{1}{c|}{$-1$}
\\ \cline{2-9}
\end{tabular} \]

\vspace{0.3cm} Si $\partial\in Der(A),$ alors pour tout
$i\in\{1,2,3\},$

\[ \partial e_i=\sum_{j\neq i}a_{ij}e_j \ \mbox{ o\`u } a_{ij}\in\rit. \]

Le fait que $\partial$ soit anti-sym\'etrique par rapport \`a
$(.|.)$ donne \ $a_{ji}=-a_{ij}.$ Ainsi, les
\'egalit\'es:

\begin{eqnarray*}
\partial e_5 &=& \partial(e_1e_4)=(\partial e_1)e_4+e_1\partial e_4, \\
-\partial e_4 &=& \partial(e_1e_5)=(\partial e_1)e_5+e_1\partial e_5, \\
\partial e_1 &=& \partial(e_4e_5)=(\partial e_4)e_5+e_4\partial e_5 \end{eqnarray*}

donnent

\vspace{0.2cm} {\bf (1)} \ $ba_{16}-a_{25}-\alpha a_{34}=0,$
\hspace{0.3cm} {\bf (5)} \ $a_{17}+a_{24}-\alpha a_{35}=0,$
\hspace{0.3cm} {\bf (9)} \ $a_{12}+a_{47}-ba_{56}=0,$

\vspace{0.1cm} {\bf (2)} \ $ca_{17}+\alpha a_{24}-a_{35}=0,$
\hspace{0.3cm} {\bf (6)} \ $a_{16}-\alpha a_{25}-a_{34}=0,$
\hspace{0.3cm} {\bf (10)} \ $a_{13}-a_{46}-ca_{57}=0,$

\vspace{0.1cm} {\bf (3)} \ $ba_{12}+a_{47}-a_{56}=0,$
\hspace{0.6cm} {\bf (7)} \ $a_{13}-a_{46}-a_{57}=0,$
\hspace{0.5cm} {\bf (11)} \ $a_{16}-ba_{25}-a_{34}=0,$

\vspace{0.1cm} {\bf (4)} \ $ca_{13}-a_{46}-a_{57}=0,$
\hspace{0.6cm} {\bf (8)} \ $a_{12}+a_{47}-a_{56}=0,$
\hspace{0.5cm} {\bf (12)} \ $a_{17}+a_{24}-ca_{35}=0.$

\vspace{0.5cm} {\bf (1), (6)} et {\bf (11)} donnent \hspace{0.2cm}
$a_{25}=a_{16}=a_{34}=0,$ \ car $\alpha\neq b.$

\vspace{0.1cm} {\bf (2), (5)} et {\bf (12)} donnent \hspace{0.2cm}
$a_{35}=a_{17}=a_{24}=0,$ \ car $\alpha\neq c.$

\vspace{0.1cm} {\bf (3), (8)} et {\bf (9)} donnent \hspace{0.2cm}
$a_{56}=a_{12}=a_{47}=0,$ \ car $b\neq 1.$

\vspace{0.1cm} {\bf (4), (7)} et {\bf (10)} donnent \hspace{0.2cm}
$a_{57}=a_{13}=a_{46}=0,$ \ car $c\neq 1.$

\vspace{0.5cm} Ainsi

\begin{eqnarray*}
\partial e_1 &=& a_{14}e_4+a_{15}e_5, \\
\partial e_4 &=& -a_{14}e_1+a_{45}e_5, \\
\partial e_5 &=& -a_{15}e_1-a_{45}e_4, \\
\partial e_2 &=& a_{23}e_3+a_{26}e_6+a_{27}e_7, \\
\partial e_3 &=& -a_{23}e_2+a_{36}e_6+a_{37}e_7, \\
\partial e_6 &=& -a_{26}e_2-a_{36}e_3+a_{67}e_7, \\
\partial e_7 &=& -a_{27}e_2-a_{37}e_3-a_{67}e_6. \end{eqnarray*}

Les \'egalit\'es:

\begin{eqnarray*}
-\partial e_7 &=& \partial(e_1e_6)=(\partial e_1)e_6+e_1\partial e_6, \\
\partial e_6 &=& \partial(e_1e_7)=(\partial e_1)e_7+e_1\partial e_7, \\
-\partial e_1 &=& \partial(e_6e_7)=(\partial e_6)e_7+e_6\partial
e_7 \end{eqnarray*}

donnent alors

\vspace{0.3cm} {\bf (1')} \ $ba_{14}-a_{27}+\alpha a_{36}=0,$
\hspace{0.3cm} {\bf (3')} \ $a_{15}+a_{26}+\alpha a_{37}=0,$
\hspace{0.3cm} {\bf (5')} \ $a_{14}-ba_{27}+ca_{36}=0,$

\vspace{0.1cm} {\bf (2')} \ $a_{15}+\alpha a_{26}+a_{37}=0,$
\hspace{0.3cm} {\bf (4')} \ $ca_{14}-\alpha a_{27}+a_{36}=0,$
\hspace{0.3cm} {\bf (6')} \ $a_{15}+a_{26}+a_{37}=0.$

\vspace{0.5cm} {\bf (2'), (3')} et {\bf (6')} donnent
\hspace{0.2cm} $a_{26}=a_{15}=a_{37}=0,$ \ car $\alpha\neq 1.$

\vspace{0.2cm} {\bf (1'), (4')} et {\bf (5')} donnent

\vspace{0.4cm}
\[ \left(
\begin{array}{lll}
b & -1      & \alpha \\
c & -\alpha & 1      \\
1 & -b      & c
\end{array} \right)\left(
\begin{array}{l}
a_{14} \\
a_{27} \\
a_{36}
\end{array} \right)=
\left(
\begin{array}{l}
0 \\
0 \\
0 \\
\end{array}
\right)  \]

\vspace{0.3cm} Le d\'eterminant de ce syst\`eme est: \
$\alpha^2+b^2+c^2-2\alpha bc-1=(\alpha-bc)^2-(1-b^2)(1-c^2)\neq
0.$ Donc \ $a_{14}=a_{27}=a_{36}=0.$ Ainsi

\begin{eqnarray*}
\partial e_1 &=& 0, \\
\partial e_4 &=& a_{45}e_5, \\
\partial e_5 &=& -a_{45}e_4, \\
\partial e_2 &=& a_{23}e_3, \\
\partial e_3 &=& -a_{23}e_2, \\
\partial e_6 &=& a_{67}e_7, \\
\partial e_7 &=& -a_{67}e_6. \end{eqnarray*}

Enfin les \'egalit\'es

\begin{eqnarray*}
\partial e_7 &=& \partial(e_2e_5)=(\partial e_2)e_5+e_2\partial e_5, \\
-\partial e_5 &=& \partial(e_2e_7)=(\partial e_2)e_7+e_2\partial e_7, \\
\partial e_2 &=& \partial(e_5e_7)=(\partial e_5)e_7+e_5\partial e_7 \end{eqnarray*}

donnent

\[ a_{23}+ba_{45}-a_{67}=ca_{23}+a_{45}-ba_{67}=a_{23}+ca_{45}-a_{67}=0. \]

Ce qui entra\^ine: \ $a_{45}=a_{23}=a_{67}$ car $b\neq c$ i.e. \
$\partial\equiv 0.\Box$

\vspace{7cm} \hspace{0.3cm} Finalement, pour \ ${\cal A}\neq{\cal
C},$ \ on fixe des param\`etres \ $\lambda_1,\dots,\lambda_7$ \
positifs et distincts deux \`a deux. Soit alors \ $\oit=\Big(
W,(.|.),\times\Big)$ \ l'alg\`ebre r\'eelle de Cayley-Dickson et
soit \ $e_1,\dots,e_7$ \ la base canonique de l'espace euclidien
$(W,-(.|.)).$ On pose

\begin{eqnarray*} x_1 &=& e_1+e_2, \\
x_2 &=& e_1-e_2, \\
x_3 &=& e_3+e_4, \\
x_4 &=& e_3-e_4, \\
x_5 &=& e_5, \\
x_6 &=& e_6+e_7, \\
x_7 &=& e_6-e_7. \end{eqnarray*}

On obtient une nouvelle base ${\cal B}=\{x_1,\dots,x_7\}$
orthogonale de $(W,-(.|.))$ pour laquelle le produit vectoriel \
$x_i\times x_j$ \ n'est colin\'eaire \`a aucun \'el\'ement de
${\cal B},$ pour tout couple $(i,j)$ d'indices distincts, de
$\{1,\dots,7\}.$ On d\'efinit alors l'automorphisme, sym\'etrique,
$s$ de l'espace euclidien $(W,-(.|.)),$ par \
$s(x_i)=\lambda_ix_i$ pour tout $i\in\{1,\dots,7\}.$ Le
prolongement $\tilde{s}$ de $s$ ne laisse invariante aucune
sous-alg\`ebre de $\oit,$ de dimension $4,$ car sinon $s$
laiserait invariante une sous-alg\`ebre de $(W,\times)$ de
dimension $3.$ Une telle sous-alg\`ebre serait engendr\'ee
lin\'eairement par trois vecteurs propres de $s,$ et l'un d'eux
serait colin\'eaire au produit vectoriel des deux autres. Ceci
contredirait le choix de la base ${\cal B}.$ D'apr\`es le
Th\'eor\`eme {\bf 5.25}, $Aut\Big( \oit(s)\Big)$ est
trivial.$\Box$

\vspace{20cm} {\LARGE\bf Index}

\vspace{1cm} \hspace{2.5cm} {\bf A} \hspace{7.2cm} {\bf H}

\vspace{0.3cm} Albert (th\'eor\`emes d\^us \`a), 14, 29
\hspace{1.7cm} Hopf-Kervaire-Milnor-Bott (th. de), 46

\vspace{0.1cm} Alg\`ebre, 13 \hspace{5.2cm}

\vspace{0.1cm} Alg\`ebre de Lie des d\'erivations, 18
\hspace{4.6cm} {\bf I}

\vspace{0.1cm} Alg\'ebrique (\'el\'ement, alg\`ebre), 16
\hspace{1.6cm}  Isom\'etries (groupe des), 67

\vspace{0.1cm} Alternative (alg\`ebre), 16 \hspace{3.1cm} Isotope
d'une alg\`ebre, 40

\vspace{0.1cm} Anti-commutative (alg\`ebre), 13

\vspace{0.1cm} Artin (th\'eor\`eme de), 16 \hspace{6.3cm} {\bf K}

\vspace{0.1cm} Associateur (de $x,y,z$), 13

\vspace{0.1cm} Automorphismes (groupe des), 18 \hspace{1.5cm}
Kaidi (travaux d\^us \`a), 33, 37-40, 44-46

\vspace{0.8cm} \hspace{2.5cm} {\bf B} \hspace{7.3cm} {\bf M}

\vspace{0.3cm} Base (d'une alg\`ebre), 14 \hspace{3.2cm} Moufang
(identit\'es de), 17

\vspace{0.1cm} \hspace{7.4cm} Multiplication (table de), 14

\vspace{0.3cm} \hspace{2.5cm} {\bf C}

\vspace{0.3cm} Cayley-Dickson (proc\'ed\'e de), 23 \hspace{4.8cm}
{\bf N}

\vspace{0.2cm} Cayleyenne (alg\`ebre), 22 \hspace{3.1cm} Norme
(d'un \'el\'ement), 22

\vspace{0.2cm} Commutateur (de $x,y$), 13

\vspace{0.1cm} Cuenca (sur un travail de), 46 \hspace{5.2cm} {\bf
O}

\vspace{0.3cm} \hspace{7.4cm} Octonions (alg\`ebre des) \ $\oit,$
24

\vspace{0.1cm} \hspace{2.5cm} {\bf D} \hspace{4.3cm} Orthogonaux
(\'el\'ements), 27

\vspace{0.3cm} Dimension (d'une alg\`ebre), 14

\vspace{0.1cm} Division lin\'eaire (\`a droite, \`a gauche), 20
\hspace{3.5cm} {\bf P}

\vspace{0.2cm} \hspace{7.4cm} P-octonions g\'en\'eralis\'ee
(alg\`ebre des), 95

\vspace{0.1cm} \hspace{2.5cm} {\bf F} \hspace{4.3cm} Puissances
\`a droite (d'un \'el\'ement), 14

\vspace{0.3cm} Flexible (alg\`ebre), 16 \hspace{3.6cm} Puissances
\`a gauche (d'un \'el\'ement), 14

\vspace{0.1cm} Frobenius (th\'eor\`eme de), 29  \hspace{2.4cm}
Puissances associatives (alg\`ebre \`a), 14

\vspace{0.7cm} \hspace{2.5cm} {\bf Q} \hspace{7.4cm} {\bf U}

\vspace{0.3cm} Quadratique (alg\`ebre), 24 \hspace{3cm} Unit\'e
\`a gauche (\'el\'ement), 13

\vspace{0.1cm} Quaternions (alg\`ebre des) \ $\hit,$ 13, 28
\hspace{1.2cm} Urbanik-Wright (r\'esultats d\^us \`a), 45

\vspace{0.7cm} \hspace{2.5cm} {\bf R} \hspace{7.4cm} {\bf W}

\vspace{0.3cm} Reflexion, 99 \hspace{5.2cm} Wright (sur les
travaux de), 31

\vspace{0.2cm} Rodr\'iguez (travaux d\^us \`a), 31 \hspace{5.4cm}
{\bf Z}

\vspace{0.2cm} \hspace{7.6cm} Zorn (th\'eor\`eme d\^u \`a), 29

\vspace{0.3cm} \hspace{2.5cm} {\bf S}

\vspace{0.3cm} Sph\`ere unit\'e (d'un e.v. norm\'e), 48

\vspace{0.7cm} \hspace{2.5cm} {\bf T}

\vspace{0.3cm} Trace (forme), 27

\vspace{22cm} \begin{center} {\Large\bf R\'esum\'e}
\end{center}

\vspace{0.5cm} \hspace{0.4cm} Dans ce travail nous nous
int\'eressons au probl\`eme g\'en\'eral de la d\'etermination des
alg\`ebres norm\'ees de division lin\'eaire. Nos r\'esultats
fondamentaux sont obtenus dans la sous-classe particuli\`ere des
alg\`ebres r\'eelles de Jordan non commutatives de division
lin\'eaire de dimension $8.$

\vspace{0.2cm} \hspace{0.4cm} Nous donnons un nouveau proc\'ed\'e
qui g\'en\'eralise celui de Cayley-Dickson et qui permet
l'obtention d'une nouvelle famille d'alg\`ebres de Jordan non
commutatives de division lin\'eaire de dimension $8.$ Nous donnons
des exemples d'alg\`ebres r\'eelles de Jordan non commutatives de
division lin\'eaire de dimension $8$ qui ne peuvent pas s'obtenir
par ce premier proc\'ed\'e de "duplication" et \`a l'aide d'un
second proc\'ed\'e, qui consiste \`a faire une d\'eformation
appropri\'ee du produit de l'alg\`ebre des octonions de
Cayley-Dickson, nous d\'eterminons ces derni\`eres et nous
r\'esolvons le probl\`eme d'isomorphisme. Nous \'etudions ensuite
les alg\`ebres r\'eelles de Jordan non commutatives de division
lin\'eaire de dimension $8$ qui poss\`edent une d\'erivation non
triviale moyennant le proc\'ed\'e de Cayley-Dickson
g\'en\'eralis\'e. Nous donnons \'egalement des exemples
d'alg\`ebres r\'eelles de Jordan non commutatives de division
lin\'eaire de dimension $8$ dont le groupe des automorphismes est
trivial, et caract\'erisons celles dont le groupe des
automorphismes est non trivial. Ceci met en \'evidence
l'immensit\'e de cette sous-classe d'alg\`ebres.

\vspace{0.8cm} \hspace{0.3cm} {\large\bf Mots cl\'es.} Alg\`ebre
norm\'ee, diviseur topologique lin\'eaire de z\'ero, alg\`ebre non
associative (alternative, flexible, de Jordan non commutative,
quadratique, cayleyenne), quaternions, octonions, proc\'ed\'e de
Cayley-Dickson "g\'en\'eralis\'e", isotopie vectorielle,
automorphismes, d\'erivations.

\vspace{16cm}
\begin{thebibliography}{100}

\bibitem{1} [A 47] A. A. Albert, {\em Absolute valued real algebras.} Ann. Math. {\bf 48},
(1947) 495-501.

\bibitem{2} [A 48] A. A. Albert, {\em Power associative rings.} T. AMS {\bf 64},
(1948) 552-593.

\bibitem{3} [A 49] A. A. Albert, {\em Absolute
valued algebraic algebras.} Bull. Amer. Math. Soc. {\bf 55},
(1949) 763-768. Anote of correction. Ibid. {\bf 55}, (1949) 1191.

\bibitem{4} [A] A. A. Albert, {\em A note of correction.} Bull. Amer. Math. Soc. {\bf 55},
(1949) 1191.

\bibitem{5} [AK 83] S. C. Althoen and L. D. Kugler, {\em When
is $\rit^2$ a division algebra?} Amer. Math. Monthly {\bf 90},
(1983) 625-635.

\bibitem{6} [AHK 86] S. C. Althoen, K. D. Hansen, and L. D. Kugler,
{\em $\cit$-associative algebras of dimension $4$ over $\rit.$}
Alg., Gr., Geom. {\bf 3}, (1986) 329-360.

\bibitem{7} [AHK 87] S. C.
Althoen, K. D. Hansen, and L. D. Kugler, {\em Four-dimensional
real algebras satisfying two $\cit$-associative conditions.} Alg.,
Gr., Geom. {\bf 4}, (1987) 395-419.

\bibitem{8} [AHK 92] S. C.
Althoen, K. D. Hansen, and L. D. Kugler, {\em Rotationnel scaled
quaternion division algebras.} J. Algebra {\bf 146}, (1992)
124-143.

\bibitem{9} [BO 81$_1$] G. M. Benkart and J. M. Osborn, {\em An investigation of real division algebras using
derivations.} Pacific Journal of Mathematics {\bf 96}, (1981)
265-300.

\bibitem{10} [BO 81$_2$] G. M. Benkart and J. M. Osborn, {\em The derivation algebra of a real division algebra.}
Amer. J. Math. {\bf 103}, (1981) 1135-1150.

\bibitem{11} [BBO 82] G. M. Benkart, D. J. Britten and J. M. Osborn, {\em Real flexible division algebras.} Can. J. Math.
{\bf XXXIV}, (1982) 550-588.

\bibitem{12} [Ber 73] S. K. Berberian, {\em Lectures in Functional Analysis and Operator Theory.} Springer-Verlag, (1973).

\bibitem{13} [BD 73] F. F. Bonsall and J. Duncan, {\em Complete Normed Algebras.} Springer-Verlag, (1973).

\bibitem{14} [BM 58] R. Bott and J. Milnor, {\em On the parallelizability of the spheres.} Bull. Amer.
Math. Soc. {\bf 64} (1958), 87–89.

\bibitem{15} [Bou 70] N. Bourbaki, {\em El\'ements de
Math\'ematiques-Alg\`ebre- Chapitre I-III.} Hermann (1970).

\bibitem{16} [BKo 66] H. Braun and M. Koecher, {\em Jordan Algebren.} Springer-Verlag, (1966).

\bibitem{17} [BKl 51] R. H. Bruck and E. Kleinfeld, {\em The structure of alternative division rings.} Proc. AMS
{\bf 2}, (1951) 878-890.

\bibitem{18} [CR] M. Cabrera Garc\'ia and A. Rodr\'iguez Palacios,
{\em A new simple proof of the Gelfand-Mazur-Kaplansky Theorem.}
(A para\^itre dans "Proc. AMS")

\bibitem{19} [Cu 92] J. A. Cuenca, {\em On one-sided division infinite-dimensional normed real
algebras.} Publications Matem\'atiques. {\bf 36} (1992), 485-488.

\bibitem{20} [CDKR] J. A. Cuenca, R. De Los Santos Villodres, A. Kaidi and A. Rochdi, {\em Classification of the real
quadratic flexible division algebras of dimension $8.$} (Soumis
\`a la publication).

\bibitem{21} [Cz 76] R. A. Czerwinski, {\em Bonded quadratic division algebras.} Pac. J. Math. {\bf 64}
(1976), 341-351.

\bibitem{22} [E-R 91] H. D. Ebbinghaus, H. Hermes, F. Hirzebruch, M. Koecher, K. Mainzer, J. Neukirch, A. Prestel and
R. Remmert, {\em Numbers.} Springer-Verlag, (1991).

\bibitem{23} [Elm 83] M. L. El-Mallah, {\em Sur les alg\`ebres absolument valu\'ees qui v\'erifient
l'identit\'e $(x,x,x)=0.$} J. Algebra {\bf 80} (1983), 314-322.

\bibitem{24} [Elm 87] M. L. El-Mallah, {\em On finite dimensional absolute valued algebras satisfying
$(x,x,x)=0.$} Arch. Math. {\bf 49} (1987), 16-22.

\bibitem{25} [Elm 88] M. L. El-Mallah, {\em Absolute valued algebras with an involution.} Arch. Math.
{\bf 51} (1988), 39-49.

\bibitem{26} [Elm 90] M. L. El-Mallah, {\em Absolute valued algebras containing a central idempotent.} J.
Algebra {\bf 128} (1990), 180-187.

\bibitem{27} [ElM 80] M. L. El-Mallah and A. Micali, {\em Sur les alg\`ebres norm\'ees sans diviseurs topologiques
de z\'ero.} Bolet\'in de la Sociedad Matem\'atica Mexicana {\bf
25} (1980), 23-28.

\bibitem{28} [ElM 81] M. L. El-Mallah and A. Micali, {\em Sur les dimensions des alg\`ebres absolument valu\'ees.} J.
Algebra {\bf 68} (1981), 237-246.

\bibitem{29} [ErMO 75] T. S. Erickson, W. S. Martindale III and J. M. Osborn, {\em Prime nonassociative algebras.}
Pac. J. Math. {\bf 60} (1975), 49-63.

\bibitem{30} [GG 73] M. G{\" u}naydin and F. G{\" u}rsey, {\em Quark structure and octonions.} J. Math. Phys.
{\bf 14} (1973), 1651-1667.

\bibitem{31} [H 40] H. Hopf, {\em Ein topologischer beitrag zur reellen algebra.} Comment. Math. Helvet. (1940), 219-239.

\bibitem{32} [J 37] N. Jacobson, {\em Abstract derivation and Lie algebras.} trans. AMS {\bf 42} (1937), 206-224.

\bibitem{33} [J 62] N. Jacobson, {\em Lie algebras.} Interscience Publischers (1962).

\bibitem{34} [J 68] N. Jacobson, {\em Structure and representations of Jordan algebras.} Amer. Math. Soc. Coll. Publ.
{\bf 39} (1968).

\bibitem{35} [J 85] N. Jacobson, {\em Basic algebra I.} Freeman and company, New York (1985).

\bibitem{36} [Kai 77] A. M. Kaidi, {\em Bases para una teoria de las algebras no asociativas normadas.} Tesis Doctoral,
Universidad de Granada, Spain (1977).

\bibitem{37} [Kai 91] A. M. Kaidi, {\em Structures des alg\`ebres de Jordan Banach non commutatives r\'eelles de division.}
Ann. Sci. univ. "Blaise Pascal", Cl\'ermont II, S\'er. Math. Fasc.
{\bf 27} (1991), 119-124.

\bibitem{38} [KR 92] A. M. Kaidi and A. Rochdi, {\em Sur les alg\`ebres r\'eelles de Jordan non commutatives de division
lin\'eaire de dimension $8.$} Nonassociative Algebraic Models.
Nova Sciences Publishers, Inc. new York (1992), 183-193.

\bibitem{39} [KS] A. Kaidi and A. S\'anchez S\'anchez, {\em J-diviseurs topologiques de z\'ero dans une alg\`ebre de
Jordan non commutative norm\'ee.} Preprint. Universidad de
M\'alaga.

\bibitem{40} [Kap 49] I. Kaplansky, {\em Normed algebras.} Duke Math. J. {\bf 16} (1949), 399-418.

\bibitem{41} [Ke 58] M. Kervaire, {\em Non-parallelizability of the $n$-sphere for $n>7.$}
Proc. Nat. Acad. Sci. USA {\bf 44} (1958), 280–283.

\bibitem{42} [Ma 77] J. Martinez Moreno, {\em Sobre algebras de jordan normadas completas.} Tesis Doctoral, Universidad
de Granada, Spain (1977).

\bibitem{43} [Mc 65] K. Mc Crimmon, {\em Norms and noncommutative Jordan algebras.} Pacific J. Math. {\bf 15} (1965),
925-956.

\bibitem{44} [Mc 66] K. Mc Crimmon, {\em Structure and representations of noncommutative Jordan algebras.} Trans. AMS
{\bf 121} (1966), 187-199.

\bibitem{45} [Ok 80] S. Okubo, {\em Some new classes of
division algebras.} J. Algebra {\bf 67} (1980), 479-490.

\bibitem{46} [OO 81$_1$] S. Okubo and J. M. Osborn, {\em Algebras with nondegenerate associative symmetric bilinear forms
permitting composition I.} Com. Alg. {\bf 9} (1981), 1233-1261.

\bibitem{47} [OO 81$_2$] S. Okubo and J. M. Osborn, {\em Algebras with nondegenerate associative symmetric bilinear forms
permitting composition II.} Com. Alg. {\bf 9} (1981), 2015-2073.

\bibitem{48} [Os 62] J. M. Osborn, {\em Quadratic division algebras.} Trans. AMS {\bf 105} (1962),
202-221.

\bibitem{49} [Pete 81] H. P. Petersson, {\em On Linear and quadratic Jordan Division Algebras.} Math. Z. {\bf 177}
(1981), 541-548.

\bibitem{50} [Petr 87] J. Petro, {\em Real division algebras of dimension $>1$ contain $\cit.$} Amer. math. Monthly
{\bf 94} (1987), 445-449.

\bibitem{51} [Pi 82] R. S. Pierce, {\em Associative Algebras.} Springer-Verlag
(1982).

\bibitem{52} [Po 85] M. Postnikov, {\em Le\c{c}ons de G\'eom\'etrie. Groupes et alg\`ebres de Lie.} Editions Mir, Moscow,
(1985).

\bibitem{53} [Pr 51] C. M. Price, {\em Jordan division algebras and the algebras $A^{(\lambda)}.$} Trans. AMS {\bf 70}
(1951), 291-300.

\bibitem{54} [Ri 60] C. E. Rickart, {\em General Theory of Banach Algebras.} Princeton (1960)

\bibitem{55} [Roc 87] A. Rochdi, {\em Sur les alg\`ebres non associatives norm\'ees de
division.} Th\`ese de $3^e$ cycle, Facult\'e des Sciences de
Rabat. Maroc (1987).

\bibitem{56} [Roc$_1$] A. Rochdi, {\em Etude des
alg\`ebres r\'eelles de Jordan non commutatives de division de
dimension $8$ dont l'alg\`ebre de Lie des d\'erivations n'est pas
triviale.} (A para\^itre dans "Journal of Algebra").

\bibitem{57} [Roc$_2$] A. Rochdi, {\em Sur les $\rit$-alg\`ebres de Jordan non commutatives de division de dimension $8$
poss\'edant un automorphisme ou une d\'erivation non triviaux.} (A
para\^itre dans "Ed. Kluwer").

\bibitem{58} [Rod 92$_1$] A. Rodr\'iguez Palacios, {\em One-sided division absolute valued algebras.} Pub. Mat., {\bf
36} (1992), 925-954.

\bibitem{59} [Rod 92$_2$] A. Rodr\'iguez Palacios, {\em Jordan Structures in Analysis.} Departamento de An\'alisis
Matem\'atico. Facultad de Ciencias, Universidad de Granada.
18071-Granada (1992).

\bibitem{60} [Rod 93] A. Rodr\'iguez Palacios, {\em N\'umeros hipercomplejos en dimensi\'on infinita.} Universidad de
Granada. Granada (1993).

\bibitem{61} [Sc 55] R. D. Schafer, {\em
Noncommutative Jordan algebras of characteristic $0.$} Proc. AMS
{\bf 6} (1955), 472-475.

\bibitem{62} [Sc 58] R. D. Schafer, {\em On noncommutative Jordan algebras.} Proc. AMS {\bf 9} (1958), 110-117.

\bibitem{63} [Sc 66] R. D. Schafer, {\em An introduction to nonassociative algebras.} Academic Press,
New York (1966).

\bibitem{64} [Sp 54] T. A. Springer, {\em An algebraic proof of a Theorem of H. Hopf.} Indagationes Mathematicae {\bf 16}
(1954), 33-35.

\bibitem{65} [UW 60] K. Urbanik and F. B. Wright,
{\em Absolute valued algebras.} Proc. Amer. Math. Soc. 11 (1960),
861-866.

\bibitem{66} [V 71] D. C. Viola, {\em Jordan algebras with continuous inverse.} Math. Japon {\bf 16} (1971),
115-125.

\bibitem{67} [War 83] F. W. Warner, {\em Foundations of Differentiable Manifolds and Lie Groups.} Springer-Verlag (1983).

\bibitem{68} [Wat 87] W. C. Waterhouse, {\em Nonassociative quaternions algebras.} Alg. Gr. Geo. {\bf 4} (1987), 365-378.

\bibitem{69} [Wr 53] F. B. Wright, {\em Absolute valued
algebras.} Proc. Nat. Acad. Sci. USA. {\bf 39} (1953), 330-332.

\bibitem{70} [Y 81] C. T. Yang, {\em Division algebras and fibrations of spheres by great spheres.} Journal of
Differential Geometry {\bf 16} (1981), 577-593.

\bibitem{71} [ZSSS 82] K. A. Zhevlakov, A. M. Slin'ko, I. P. Shestakov and A. I. Shirshov, {\em Rings that are nearly
associative.} Academic Press (1982).
\end{thebibliography}
\end{document}